\documentclass[11pt, reqno]{amsart}

\usepackage{amsmath, amstext, amssymb, amsthm, bbm,lmodern, enumerate, esint, todonotes}
\usepackage[pdftex]{hyperref}
\usepackage{tikz}
\usepackage{subcaption}
\usepackage{graphicx}
\usepackage[english]{babel}

\usepackage[utf8]{inputenc}     

\usepackage[toc,page]{appendix}
\usepackage{fancyhdr}
\DeclareFontFamily{U}{mathx}{}
\DeclareFontShape{U}{mathx}{m}{n}{<-> mathx10}{}
\DeclareSymbolFont{mathx}{U}{mathx}{m}{n}
\DeclareMathAccent{\widehat}{0}{mathx}{"70}
\DeclareMathAccent{\widecheck}{0}{mathx}{"71}

\usepackage[ right= 1.7cm, left=1.7cm, a4paper,top=2.7cm, bottom= 3cm]{geometry}

\usepackage{mathrsfs}
\usepackage{xpatch}

\makeatletter
\xpatchcmd{\@startsection}{\@m}{\@M}{}{}
\makeatother

\newcommand{\E}{\mathbb {E}}   
\newcommand{\R}{\mathbb {R}}   
\newcommand{\Z}{\mathbb {Z}}   
\newcommand{\N}{\mathbb {N}}   
\newcommand{\ind}{\mathbbm{1}} 

\newcommand{\AAA}{\mathcal {A}}  
\newcommand{\FFF}{\mathcal {F}} 
\newcommand{\EEE}{\mathcal {E}}

\newcommand{\CCC}{\mathcal {C}}
\newcommand{\PPP}{\mathcal {P}}
\newcommand{\MMM}{\mathcal {M}}
\newcommand{\RRR}{\mathcal {R}}
\newcommand{\BBB}{\mathcal {B}} 

\newcommand{\Ent}{\mathrm{Ent}}

\newcommand{\dP}{\mathbf{P}}
\newcommand{\I}{\mathrm{I}}
\newcommand{\dR}{\mathbf{R}}
\newcommand{\dt}{\frac{\mathrm{d}}{\mathrm{dt}}}

\newcommand{\A}{K}

\newcommand{\xil}{\xi^M}
\newcommand{\xl}{x^M}

\newcommand{\Leb}{\mathrm{Leb}}

\newcommand{{\muP}}{{{\mathbf{P}}}}

\newcommand{{\dV}}{{\mathbf{V}}}
\newcommand{{\dPl}}{{{\mathtt{P}}}}
\newcommand{{\dVl}}{{{\mathtt{V}}}}

\def\Poi {{\mathbf{Poi}}}  

\numberwithin{equation}{section}

\theoremstyle{plain}
\newtheorem{satz}{Theorem}[section]

\newtheorem{kor}[satz]{Corollary}
\newtheorem{lem}[satz]{Lemma}
\newtheorem{prop}[satz]{Proposition}
\newtheorem{dfn}[satz]{Definition}

\theoremstyle{remark}
\newtheorem{bem}[satz]{Remark}
\newtheorem{bsp}[satz]{Example}

\usepackage{blindtext}
\usepackage{enumitem} 
\usepackage{todonotes}

\begin{document}
\title{Non-local Transport Geometry for Papangelou Point Processes}
\author{Matthias Erbar}
\address{Matthias Erbar: Faculty of Mathematics, University of Bielefeld,  33615
Bielefeld, Germany}
\email{erbar@math.uni-bielefeld.de }

\author{Martin Huesmann}
\address{Martin Huesmann: Institute for Mathematical Stochastics,
University of M\"unster, 
48149 M\"unster, Germany}
\email{martin.huesmann@uni-muenster.de}

\author{Hanna Stange}
\address{Hanna Stange: Faculty of Mathematics, Ruhr University Bochum,
44780 Bochum, Germany}
\email{hanna.stange@ruhr-uni-bochum.de}
\date{}

\maketitle
\begin{abstract}
We develop new analytical tools for the analysis of interacting point processes and their associated dynamics by introducing a non-local Benamou–Brenier-type transport distance tailored to (stationary) Papangelou point processes. Under a suitable Dobrushin-type condition, quantified by a constant $\zeta$, we show that the natural spatial birth-and-death dynamics associated with a Papangelou point process  is the EVI$(\zeta)$-gradient flow of the (specific) relative entropy. Moreover, we establish $\zeta$-geodesic convexity of the (specific) relative entropy and derive new functional inequalities, including HWI, Log-Sobolev, and Talagrand inequalities. 
\end{abstract}
\tableofcontents

\section{Introduction}\label{introduction}
For interacting point processes even fundamental questions such as existence and uniqueness are far from trivial, see e.g.\ \cite{leble2015uniqueness,Dereudre2019,Houdebert_Zass_2022}. In view of the discrete structure of point processes and their associated birth-and-death dynamics, the theory of non-local optimal transport, initially developed for Markov chains \cite{maas2011gradientflowsentropyfinite,Mielke2011}, provides a promising approach for the development of new analytical tools for their analysis. Building on recent works for Poisson point processes in the non-stationary setting \cite{schiavo2023wasserstein} and the stationary infinite-volume setting \cite{huesmann2025nonlocalwassersteingeometrygradient}, we develop a non-local Wasserstein geometry tailored to the analysis of (stationary) Papangelou point processes. Exploiting the new induced geometry, we establish under a suitable Dobrushin-type condition, quantified by a constant $\zeta>0$, 
\begin{itemize}
    \item that the natural spatial birth-and-death dynamics associated with a Papangelou point process is the EVI$(\zeta)$-gradient flow of the (specific) relative entropy;
    \item that the (specific) relative entropy is $\zeta$-geodesically convex with respect to this geometry; and
    \item functional inequalities, including HWI inequalities, Log-Sobolev inequalities, and Talagrand inequalities.
\end{itemize}
\subsection{Main Results}
\subsubsection{Setting}
Let $(X,d)$ be a Polish metric space. The space of locally finite counting measures on $X$ is denoted by $\Gamma_X$. A point process on $X$ is a random variable taking values in $\Gamma_X$, and its distribution is an element of $\PPP(\Gamma_X)$ the space of probability measures on $\Gamma_X$. We call $\dP\in\PPP(\Gamma_X)$ simple if $\dP({\xi\in\Gamma_X:\xi({x})\leq 1\text{ for all }x\in X})=1$, and locally integrable if $\int\xi(B)\dP(\mathrm{d}\xi)<\infty$ for all bounded Borel sets $B\in\BBB_0(X)$. The set of simple and locally integrable probability measures on $\Gamma_X$ is  $\PPP_1(\Gamma_X)$.

In the Euclidean setting $X=\R^d$, we define the shift operator $\Theta_z:\Gamma_{\R^d}\to\Gamma_{\R^d}$ by $\Theta_z\xi(B)=\xi(B+z)$ for all Borel sets $B\in\BBB(\R^d)$ and 
$z\in\R^d$. A measure $\dP\in\PPP_1(\Gamma_{\R^d})$ is called \emph{stationary} if it is invariant under spatial shifts, i.e.\
$(\Theta_z)_\#\dP=\dP$ for all $z\in\R^d$, where $(\Theta_z)_\#\dP$ denotes the pushforward of $\dP$ under $\Theta_z$. The set of all 
stationary measures is denoted by $\PPP_s(\Gamma_{\R^d})$, and a point process is called \emph{stationary} if its distribution belongs to $\PPP_s(\Gamma_{\R^d})$.

In the following, we introduce a transport distance tailored to the analysis of \textit{Papangelou measures} $\dR\in\PPP_1(\Gamma_{X})$ with Papangelou intensity $\kappa:\Gamma_{X}\times X\to[0,\infty)$ and reference measure $\lambda\in\MMM^+(X)$ satisfying the GNZ equation, named after Georgii \cite{Geogii76}, Nguyen, and Zessin \cite{NguyenZessin1979},
\begin{align}\label{ieq:GNZ}
    \int\int F(\xi-\delta_x,x)\xi(\mathrm{d}x)\dR(\mathrm{d}\xi)
    =
    \int\int F(\xi,x)\kappa(\xi,x)\lambda(\mathrm{d}x)\dR(\mathrm{d}\xi),
\end{align}
for all non-negative measurable functions $F:\Gamma_{X}\times X\to[0,\infty)$.

\subsubsection{Non-Local Transport Geometry}
The construction of our transport distance follows the approach of non-local transport distances for Markov chains \cite{maas2011gradientflowsentropyfinite,Mielke2011}, jump processes \cite{erbar2012gradientflowsentropyjump}, and Poisson point processes \cite{schiavo2023wasserstein, huesmann2025nonlocalwassersteingeometrygradient}, and is inspired by the dynamical Benamou--Brenier formulation of the classical $L^2$-Wasserstein distance \cite{BeBr00}.  For  $\dP_0,\dP_1\in\PPP_1(\Gamma_{X})$, the non-local transport distance takes the form
\begin{align*}
\mathcal{W}_{\kappa,\lambda}^2(\dP_0,\dP_1)
:=
\inf\left\{
\mathcal{A}_{\kappa,\lambda}(\bar{\dP},\bar{\mathbf V})
:
(\bar{\dP},\bar{\mathbf V})\in
\CCC\EEE^X(\dP_0,\dP_1)
\right\},
\end{align*}
where $\CCC\EEE^X(\dP_0,\dP_1)$ denotes the set of  sufficiently regular curves $\bar{\dP}=(\bar{\dP}_t)_{t\in[0,1]}\subset\PPP_1(\Gamma_X)$ connecting $\dP_0$ and $\dP_1$ and velocity fields $\bar{\mathbf V}=(\bar{\mathbf V}_t)_{t\in[0,1]}\subset\MMM(\Gamma_X\times X)$ of signed measures on $\Gamma_X\times X$, satisfying the \emph{non-local continuity equation}
\begin{align}\label{ieq:ce}
\partial_t\bar{\dP}_t+D\bar{\mathbf V}_t=0
\end{align}
in the distributional sense, where the classical gradient is replaced by the \textit{discrete difference operator} $D_xF(\xi):=F(\xi+\delta_x)-F(\xi)$ for functionals $F:\Gamma_X\to\R$.
Further, $\mathcal{A}_{\kappa,\lambda}$ is a non-local \textit{action functional} tailored to the Papangelou intensity $\kappa$ and the reference measure $\lambda$, see Section \ref{sec:transportdistance} for more details about the construction and properties of $\mathcal{W}_{\kappa,\lambda}$. 

However, similar to the relative entropy, the action functional $\AAA_{\kappa,\lambda}$ typically explodes in the stationary Euclidean setting $X=\R^d$. Following the approach of \cite{EHJM25} and \cite{huesmann2025nonlocalwassersteingeometrygradient}, we define for stationary measures $\dP_0,\dP_1\in\PPP_s(\Gamma_{\R^d})$ the \textit{specific} non-local transport distance
\begin{align*}
\mathcal{W}_{\kappa,s}^2(\dP_0,\dP_1)
:=
\inf\left\{\sup_{n\in\N}\frac{1}{\Leb(\Lambda_n)}
\mathcal{A}_{\kappa,\Leb_{|\Lambda_n}}(\bar{\dP},\bar{\mathbf V})
:
(\bar{\dP},\bar{\mathbf V})\in
\CCC\EEE^{s}(\dP_0,\dP_1)
\right\},
\end{align*}
where $\Lambda_n:=[-n/2,n/2]^d$, $\Leb$ is the $d$-dimensional Lebesgue measure, $\bar{\dP}_{|\Lambda_n}$ and $\bar{\mathbf V}_{|\Lambda_n}$ denote the restrictions of $\bar{\dP}$ and $\bar{\mathbf V}$ to $\Gamma_{\Lambda_n}$ and $\Gamma_{\Lambda_n}\times\Lambda_n$ respectively, and $\CCC\EEE^s(\dP_0,\dP_1)$ is the set of stationary solutions to the non-local continuity equation \eqref{ieq:ce}, see Section \ref{sec:stationarydistance} for more details.\\

In the following, we focus on the stationary setting and consider stationary Papangelou measures $\dR\in\PPP_s(\Gamma_{\R^d})$ with reference measure $\Leb$ and Papangelou intensity $\kappa$, which  is locally stable with constant $1$ \hyperlink{iprop:locallystable}{$(\mathbf{S}(1))$}, repulsive \hyperlink{iprop:repulsive}{$(\mathbf{RP})$}, has finite range \hyperlink{iprop:finiterange}{$(\mathbf{FR})$}, and satisfies a certain continuity assumption \hyperlink{iprop:usc+}{$(\mathbf{USC+})$}. Moreover, we impose a technical approximation assumption \hyperlink{iproperty:ALB}{$(\mathbf{ALB})$} and a certain variant of the Dobrushin condition defined by \hyperlink{iprop:dobrushin}{$(\mathbf{D}(\zeta(\kappa)))$}
\[
0<\zeta(\kappa)
:=
\frac12\bigg(1-
\sup_{\substack{\xi\in\Gamma_{\R^d}\\ x\in\R^d}}
\int
\bigl|\kappa(\xi,y)-\kappa(\xi+\delta_x,y)\bigr|
\Leb(\mathrm{d}y)\bigg).
\]
and originally used to establish uniqueness of Gibbs measures, see e.g.\ \cite{Dobruschin1968TheDO,Ruelle,Houdebert_Zass_2022}.  We refer to Section~\ref{sec:proofidea} for a precise formulation and discussion of these assumptions. Examples satisfying all assumptions include the area-interaction model, the hard-core model, and finite-range pair interactions with a hard-core exclusion at the origin, see Example~\ref{exp:Papangelou}. The Poisson case is also covered with $\zeta(\kappa)=\frac12$, thereby recovering the results of \cite{schiavo2023wasserstein} und \cite{huesmann2025nonlocalwassersteingeometrygradient}.  \\

Exploiting the geometry induced by the specific non-local distance $\mathcal{W}_{\kappa,s}$, we establish an evolution variational inequality (EVI) for the specific relative entropy
\[
    \mathcal{E}(\dP|\dR)
    :=
    \liminf_{n\to\infty}
    \frac{1}{\Leb(\Lambda_n)}
    \mathrm{Ent}\bigl(\dP_{|\Lambda_n}|\dR_{|\Lambda_n}\bigr),\quad \dP\in\PPP_s(\Gamma_{\R^d}),
\]
identifying the gradient flow of $\mathcal{E}(\cdot|\dR)$
as the spatial birth-and-death dynamics with birth rate $\kappa$ and death rate $1$. The latter is formally described by the Markov generator
\[
    \mathfrak{L}_{\kappa}F(\xi)
    :=
    \int_{\R^d}
    \kappa(\xi,x)
    \bigl(F(\xi+\delta_x)-F(\xi)\bigr)
    \Leb(\mathrm{d}x)
    +
    \sum_{x\in\xi}
    \bigl(F(\xi-\delta_x)-F(\xi)\bigr),
    \quad F:\Gamma_{\R^d}\to\R, 
\]
and the associated dual semigroup acting on $\PPP_s(\Gamma_{\R^d})$ is denoted by $(\mathrm{S}_t^{\kappa*})_{t\geq0}$.
\begin{satz}[EVI]\label{ithm:evistationary}
The dual spatial birth-and-death semigroup $(\mathrm{S}_t^{\kappa*})_{t\geq0}$ satisfies the evolution variational inequality with parameter $\zeta(\kappa)$ (EVI$(\zeta(\kappa))$) for the specific relative entropy $\mathcal{E}(\cdot|\dR)$, i.e.\  for all $\dP_0,\dP_1\in\PPP_s(\Gamma_{\R^d})$ such that $ \mathcal{E}(\dP_1|\dR)<\infty$ and $\mathcal{W}_{\kappa,s}(\dP_0,\dP_1)<\infty$ it holds
\[
    \frac{1}{2}\frac{\mathrm{d}}{\mathrm{d}t}
    \mathcal{W}_{\kappa,s}^2\bigl(\mathrm{S}_t^{\kappa*}\dP_0,\dP_1\bigr)
    +
    \frac{\zeta(\kappa)}{2}
    \mathcal{W}_{\kappa,s}^2\bigl(\mathrm{S}_t^{\kappa*}\dP_0,\dP_1\bigr)
    \leq
    \mathcal{E}(\dP_1|\dR)
    -
    \mathcal{E}\bigl(\mathrm{S}_t^{\kappa*}\dP_0|\dR\bigr),
    \quad t>0 .
\]
\end{satz}
Among the various formulations of gradient flows, the EVI formulation is one of the strongest, as it entails several fundamental consequences for both the associated semigroup and the underlying transport geometry, see e.g.~\cite{Ambrosio2008GradientFI}. In particular, it implies $\zeta(\kappa)$-geodesic convexity of the specific relative entropy $\Ent(\cdot|\dR)$, which may be regarded as the counterpart of McCann's displacement convexity in our setting.

\begin{satz}[Geodesic convexity]\label{ithm:geodconvexitystationary}
The specific relative entropy $\mathcal{E}(\cdot|\dR)$ is $\zeta(\kappa)$-geodesically convex with respect to $\mathcal{W}_{\kappa,s}$. More precisely, if $\dP_0,\dP_1\in\PPP_s(\Gamma)$ satisfy $\mathcal{W}_{\kappa,s}(\dP_0,\dP_1)<\infty$, then for every $\mathcal{W}_{\kappa,s}$-geodesic $\bar{\dP}$ connecting $\dP_0$ and $\dP_1$, it holds 
\[
    \mathcal{E}(\bar{\dP}_t|\dR)
    \leq
    (1-t)\mathcal{E}(\dP_0|\dR)
    +
    t\mathcal{E}(\dP_1|\dR)
    -
    \frac{\zeta(\kappa)}{2}t(1-t)
    \mathcal{W}_{\kappa,s}^2(\dP_0,\dP_1),
    \quad 0\leq t\leq1 .
\]
\end{satz}
In view of the definition of synthetic Ricci curvature bounds due to Sturm~\cite{St06} and Lott--Villani~\cite{LoVi09}, and their non-local counterpart due to the first author and Maas~\cite{Erbar_2012}, the previous theorem may be viewed as providing a synthetic lower Ricci curvature bound $\zeta(\kappa)$ for the Papangelou space $(\Gamma_{\R^d},\dR)$.

Beyond these geometric consequences, the EVI$(\zeta(\kappa))$ (Theorem~\ref{ithm:evistationary}) also yields several functional inequalities. As a first consequence, we obtain exponential contraction of  $(\mathrm{S}_t^{\kappa*})_{t\geq0}$.
\begin{satz}[Contractivity]\label{ithm:contractivity}
The dual spatial birth-and-death semigroup 
$(\mathrm{S}_t^{\kappa*})_{t\geq0}$ is $\zeta(\kappa)$-exponentially contractive  with respect to $\mathcal{W}_{\kappa,s}$: 
\begin{align*}
    \mathcal{W}_{\kappa,s}
    (
        \mathrm{S}_t^{\kappa*}\dP_0,
        \mathrm{S}_t^{\kappa*}\dP_1
    )
    &\leq
    e^{-\zeta(\kappa)t}
    \mathcal{W}_{\kappa,s}(\dP_0,\dP_1),
\end{align*}
for all $t\geq0$ and $\dP_0,\dP_1\in\PPP_s(\Gamma_{\R^d})$ with \(\mathcal{E}(\dP_0|\dR),\mathcal{E}(\dP_1|\dR)<\infty\). Further, $(\mathrm{S}_t^{\kappa*})_{t\geq0}$ is $2\zeta(\kappa)$-exponentially contractive  with respect to $\mathcal{E}(\cdot|\dR)$:
\begin{align*}
    \mathcal{E}
    (
        \mathrm{S}_t^{\kappa*}\dP_0|\dR
    )
    &\leq
    e^{-2\zeta(\kappa)t}
    \mathcal{E}(\dP_0|\dR)
\end{align*}
for all $t\geq0$ and $\dP_0\in\PPP_s(\Gamma_{\R^d})$.  
\end{satz}
Following the approach of
\cite{ErbarMaas2014GradientFlow,JuengelYue2017}, we will show in Section \ref{sec:phientropies} that it is possible 
to extend the contraction results to more general $\Phi$-entropies.

Another consequence of the EVI$(\zeta(\kappa))$ (Theorem~\ref{ithm:evistationary}) is a non-local analog of the famous $\mathrm{HWI}$ inequality of Otto and Villani~\cite{OtVi00}. Classically, this inequality links the relative entropy, typically denoted by $H$, the $L^2$-Wasserstein distance $W_2$, and the classical Fisher information $I$. In the present stationary Papangelou setting, these quantities are replaced by the specific relative entropy $\EEE(\cdot|\dR)$, the non-local transport distance $\mathcal{W}_{\kappa,s}$, and the specific modified Fisher information 
\begin{align*}
\mathcal{I}(\dP|\dR)
:=\liminf_{n\to\infty}\frac{1}{\Leb(\Lambda_n)}\mathrm{I}(\dP_{|\Lambda_n}|\dR_{|\Lambda_n})
\end{align*}
with \begin{align*}
   \mathrm{I}(\dP_{|\Lambda_n}|\dR_{|\Lambda_n})= \begin{cases}
\displaystyle
\int\int D_x\rho_n(\xi) D_x\log(\rho_n)(\xi)\tilde\kappa_n(\xi,x)\dR_{|\Lambda_n}(\mathrm{d}\xi)\Leb_{|\Lambda_n}(\mathrm{d}x),
&\text{if } \dP_{|\Lambda_n}=\rho_n\dR_{|\Lambda_n},\\
\infty,
&\text{otherwise,}
\end{cases}
\end{align*}
where $\tilde\kappa_n$ is the finite-volume Papangelou intensity of $\mathbf{R}_{|\Lambda_n}$, see Section \ref{sec:entropy} for more details. We only mention here the useful and for us important result that $\mathrm{I}(\dP_{|\Lambda_n}|\dR_{|\Lambda_n}) < \infty$ implies that $\dP_{|\Lambda_n}$ is a Papangelou measure (see Lemma \ref{lem:regularisationFisher}).
\begin{satz}[HWI inequality]\label{ithm:HWI}
 For all
$\dP\in\PPP_s(\Gamma_{\R^d})$  such that
$\EEE(\dP|\dR)<\infty$, the non-local HWI inequality with constant $\zeta(\kappa)$ is satisfied: 
\begin{align*}
    \EEE(\dP|\dR)
    \leq
    \mathcal{W}_{\kappa,s}(\dP,\dR)
    \sqrt{\mathcal{I}(\dP|\dR)}
    -
    \frac{\zeta(\kappa)}{2}
    \mathcal{W}_{\kappa,s}^2(\dP,\dR).
\end{align*}
\end{satz}

 A notable consequence of the HWI inequality is the following modified Log-Sobolev
inequality. In particular, it extends Wu's classical Log-Sobolev
inequality~\cite{Wu00} beyond the Poisson framework to Papangelou measures.

\begin{satz}[Modified Log-Sobolev inequality]\label{ithm:logsobolev}
For all
$\dP\in\PPP_s(\Gamma_{\R^d})$  such that
$\EEE(\dP|\dR)<\infty$, the modified Log-Sobolev inequality with constant $\zeta(\kappa)$ is satisfied:
\begin{align*}
    \EEE(\dP|\dR)
    \leq
    \frac{1}{2\zeta(\kappa)}
    \mathcal{I}(\dP|\dR).
\end{align*}
\end{satz}

Using the classical Herbst argument \cite{GrossRothaus1998}, this modified Log-Sobolev inequality may be used to derive new concentration inequalities for Papangelou functionals, see \cite{steenbeck2026modifiedlogsobolevinequalitiesconcentration}. Moreover, the modified Log-Sobolev inequality also implies the following Talagrand inequality.

\begin{satz}[Talagrand inequality]\label{ithm:talagrand}
For all
$\dP\in\PPP_s(\Gamma_{\R^d})$  such that
$\EEE(\dP|\dR)<\infty$, the non-local Talagrand inequality with constant $\zeta(\kappa)$ is satisfied: 
\begin{align*}
    \mathcal{W}_{\kappa,s}(\dP,\dR)
    \leq
    \sqrt{
        \frac{2}{\zeta(\kappa)}
        \EEE(\dP|\dR)
    }.
\end{align*}
\end{satz}

The previous Talagrand inequality bounds the non-local transport distance $\mathcal{W}_{\kappa,s}$ in terms of $\EEE(\cdot|\dR)$ and, in particular, implies that $\mathcal{W}_{\kappa,s}$ is finite on the proper domain of $\EEE(\cdot|\dR)$.
\subsubsection{Connection to the Literature}

Non-local transport distances were initially introduced for the analysis of finite Markov chains by Maas \cite{maas2011gradientflowsentropyfinite} and independently by Mielke \cite{Mielke2011}, and were subsequently extended to jump processes \cite{erbar2012gradientflowsentropyjump}. The induced transport geometries provide a powerful framework for studying gradient flows and deriving new functional inequalities and convexity properties for various discrete stochastic models, including birth-and-death processes, spin models, and discrete mean-field models, see e.g.\ \cite{Erbar_2012,FM16,EHMT17,ErbarFathiSchlichting2020}. Extending these ideas to the infinite-dimensional setting of configuration spaces, Dello Schiavo, Herry, and Suzuki \cite{schiavo2023wasserstein} introduced a non-local transport distance tailored to Poisson point processes and identified the Ornstein--Uhlenbeck process as the EVI$(1/2)$-gradient flow of the relative entropy. Building on this, in \cite{huesmann2025nonlocalwassersteingeometrygradient} this approach is extended to stationary point processes by introducing a non-local transport distance per unit volume and establishing a corresponding gradient flow  result and functional inequalities in the infinite-volume setting.

\medskip

In contrast to the non-local approach, transport geometries on configuration spaces can also be constructed by lifting the classical Wasserstein distance to the space of point configurations, leading to a local transport geometry in which points move continuously through the underlying space. The investigation of such Wasserstein geometries was initiated  by the first and second authors \cite{ErHu15} and later substantially generalized in \cite{DeSu21,DeSu22}. Building on \cite{HuSt13,Hu16} and motivated by \cite{erbar2018onedimensionalloggasfreeenergy}, a local transport distance per unit volume tailored to stationary point processes was introduced in \cite{EHJM25}. In the resulting geometry, the evolution of infinitely many independent Brownian motions was identified as the EVI$(0)$-gradient flow of the specific relative entropy with respect to the Poisson measure. In dimension one, the second author and Müller \cite{HuMu25} recently introduced a related Wasserstein-type gap distance and identified gradient flows of free energy functionals. Such local transport distances also reveal deep insights into long-range order properties of stationary point processes, including connections to hyperuniformity \cite{dereudre2024nonhyperuniformityperturbedlattices} and number rigidity \cite{DerDig}. In a complementary direction, the second author and Leblé \cite{HuLe24} established an infinite-volume analog of the classical $W_2$-$H^{-1}$ estimate, in which the $H^{-1}$ norm takes the form of a regularized Coulomb energy.
Similar transport distances for point processes were also constructed and studied in \cite{Dec08, DJS10}.

Beyond optimal transport, spatial birth-and-death dynamics have been studied by other means, for instance via Dirichlet form methods in e.g.\ \cite{li2010noteequilibriumglauberkawasaki,KoLy05,Lytvynov-Ohlerich07} and through stochastic equation representations in \cite{garcia2006spatialbirthdeathprocesses}. In the stationary setting, recent works \cite{steenbeck2025reversiblebirthanddeathdynamicscontinuum,jahnel2026reversiblebirthanddeathdynamicscontinuum} establish entropy-dissipation results and de Bruijn identities for reversible birth-and-death dynamics with Gibbs invariant measures having bounded and continuous energy functionals.\\ 

Functional inequalities for Poisson processes have a long history. For instance Poincaré inequalities w.r.t.\ the discrete difference operator $D_x$ have been derived in \cite{Che85, HPA95, LP11}. The modified Log-Sobolev inequality for Poisson processes was proved by Wu \cite{Wu00} and used to derive concentration inequalities, for instance in \cite{BP16}. A new elegant proof of Wu's modified Log-Sobolev inequality via the Poisson transport map was recently given in \cite{ALS24}.  The relation between the modified Log-Sobolev inequality and the Poincaré inequality, both for Poisson processes and in more general discrete settings, was studied in \cite{BL98, BT06}. Direct linearisation yields a Poincaré inequality with a loss of a factor of two in the constant. The sharp constant can be recovered using the intertwining relation $DP_t=e^{-t}P_tD$, see \cite{CJ13}. The same work also treats more general $\Phi$-entropy inequalities, see also \cite{Cha04, Cha06}. Further, using a different transport metric, $T_1$ transport inequalities for Poisson measures and Gibbs measures in the Dobrushin regime were established in \cite{MSWW11, MWW11}. $T_2$ or Marton-type inequalities have been derived in \cite{GHP21} using tensorisation techniques. Moreover,  \cite{Rei13} established a Poisson analogue of Talagrand's convex-distance inequality.

In case of interacting systems, i.e.\ Gibbs point processes the literature is more sparse. Stein's method has been used to obtain estimates in Kantorovich--Rubinstein or $W_1$-type distances \cite{SS14, DV18, Tor17}. Poincaré inequalities have been derived e.g.\ in \cite{BCC02, KoLy05, KKO13} under Dobrushin conditions. Entropy decay along Glauber dynamics was proved in \cite{DPP13} using ideas similar to those employed here, and for Markov chains in \cite{Con22} using a probabilistic coupling approach. Both approaches extend to $\Phi$-entropies, but do not yield Talagrand or HWI inequalities or an EVI. The recent work \cite{steenbeck2026modifiedlogsobolevinequalitiesconcentration} connects modified Log-Sobolev inequalities with uniqueness of Gibbs measures and concentration properties.
\subsection{Proof Idea}\label{sec:proofidea}
We now sketch the proof of the EVI stated in Theorem~\ref{ithm:evistationary} and explain the necessary assumptions in more detail. Once the EVI is established, geodesic convexity (Theorem~\ref{ithm:geodconvexitystationary})  and the functional inequalities (Theorems~ \ref{ithm:contractivity}, \ref{ithm:HWI}, \ref{ithm:logsobolev}, and~\ref{ithm:talagrand}) follow immediately from the classical theory of gradient flows in Wasserstein spaces, see e.g.\ \cite{Ambrosio2008GradientFI}.

We begin with a brief overview of the general proof strategy and then provide more details about the individual steps of the proof.

In view of the analogy between non-local optimal transport for Markov chains \cite{maas2011gradientflowsentropyfinite} and the new non-local transport geometry for Papangelou measures $\mathcal{W}_{\kappa,s}$, one might hope to adapt techniques already developed in the non-local Markov chain setting. However, such a transfer is not directly feasible: the existing non-local transport theory for Markov chains is formulated on finite base spaces, whereas the configuration spaces considered here are infinite-dimensional.

Hence, we first consider discrete base spaces $X=\Lambda_n^M =(\frac{1}{M}\Z^d-(\frac{1}{2M}-\frac{n}{2})\mathbf{1})\cap[-n/2,n/2)^d,$ where $\mathbf{1}=(1,\ldots,1)$. In this setting, the configuration space $\Gamma_{\Lambda_n^M}$ only contains finitely many simple point configurations, making the tools from the non-local theory of optimal transport for Markov chains available. Having established the EVI in the discrete setting, in a second step, we let the mesh size $1/M$ tend to zero and lift the EVI to continuous but bounded base spaces $\Lambda_n$. Finally, we pass to the normalized infinite-volume limit and establish the EVI for stationary measures.

\subsubsection{Discrete Base Spaces}
As noted above, the advantage of restricting the base space to the finite lattice $X=\Lambda_n^M$ is that the number of simple point configurations on $\Lambda_n^M$ is finite. Moreover, spatial birth-and-death dynamics can be viewed as finite-state Markov chains and the non-local transport distance for discrete Papangelou measures coincides with the non-local transport distance for finite Markov chains \cite{maas2011gradientflowsentropyfinite}. This allows us to apply tools from the theory of non-local optimal transport for Markov chains.

Following \cite{FM16,EHMT17}, we use a perturbative approach to establish a discrete Bochner inequality, which is one of several equivalent characterizations of the EVI in the finite Markov chain setting, see \cite[Theorem 4.5]{Erbar_2012}. The resulting EVI holds for discrete Papangelou measures with Papangelou intensities $\kappa^M:\Gamma_{\Lambda_n^M}\times \Lambda_n^M\to\R$ that are locally stable with constant $1$ \hypertarget{iprop:locallystable}{$(\mathbf{S}(1))$}, i.e.\ $\kappa^M(\xi^M,x^M)\leq 1$ for all $\xi^M\in\Gamma_{\Lambda_n^M}$, $x^M\in\Lambda_n^M$, repulsive \hypertarget{iprop:repulsive}{$(\mathbf{RP})$}, i.e.\ $\kappa^M(\xi^M+\delta_{x^M},y^M)\leq\kappa^M(\xi^M,y^M)$  for all $\xi^M\in\Gamma_{\Lambda_n^M}$, $x^M,y^M\in\Lambda_n^M$, and satisfy a discrete Dobrushin-type condition \hypertarget{iprop:dobrushin}{$(\mathbf{D}(\zeta(\kappa^M)))$} of the form
\begin{align*}
    0<\zeta(\kappa^M)
:=\frac{1}{2}\bigg(
1-
\sup_{\substack{\xi^M\in\Gamma_{\Lambda_n^M}\\ x^M\in\Lambda_n^M}}
\int|
\kappa^M(\xi^M,y^M)
-
\kappa^M(\xi^M+\delta_{x^M},y^M)|
\Leb^M_{|\Lambda_n^M}(\mathrm{d}y^M)\bigg),
\end{align*}
with $\Leb^M_{|\Lambda_n^M}:=\frac{1}{M^d}\sum_{z\in\Lambda_n^M}\delta_z$.

Given a finite-volume Papangelou measure $\dR\in\PPP_1(\Gamma_{\Lambda_n})$ with reference measure $\Leb_{|\Lambda_n}$ and Papangelou intensity $\kappa:\Gamma_{\Lambda_n}\to\R$, there are two natural ways to choose a discrete Papangelou measure to which the discrete EVI can be applied:

\begin{itemize}
\item The first approach is to choose the natural discretization $\dR^M=(T^M)_\#\dR$, where $T^M:\Gamma_{\Lambda_n}\to\Gamma_{\Lambda_n^M}$ is defined by
\[
T^M(\xi)=\xi^M
:=
\sum_{z\in\Lambda_n^M}
\delta_z\ind_{\{\xi(z+\Lambda_{\frac{1}{M}})\geq 1\}}.
\]
 The discretization $\dR^M$ is again a Papangelou measure on ${\Gamma}_{\Lambda_n^M}$ with reference measure $\Leb^M_{|\Lambda_n^M}$ and Papangelou intensity
\[
\kappa^{M}_*(\xi^M,{x^M}):=
M^d
\frac{
\dR(\{\eta:\eta^M=\xi^M+\delta_{x^M}\})
}{
\dR(\{\eta:\eta^M=\xi^M\})
}\ind_{\{\xi^M(\{x^M\})=0\}},
\quad \xi^M\in\Gamma_{\Lambda_n^M},\; x^M\in\Lambda_n^M  .
\]
However, repulsiveness of $\kappa$  does not, in general, imply repulsiveness of $\kappa^{M}_*$. Therefore, the discrete EVI is not directly applicable.
\item The second approach is to consider the finite-volume Papangelou measure $\tilde{\dR}^M$ with reference measure $\Leb^M_{|\Lambda_n^M}$ and Papangelou intensity
\[
\kappa^M(\xi^M,x^M):=\kappa(\xi^M,x^M)\ind_{\{\xi^M(\{x^M\})=0\}}.
\]
If $\kappa$ satisfies \hyperlink{iprop:locallystable}{$(\mathbf{S}(1))$} and \hyperlink{iprop:repulsive}{$(\mathbf{RP})$},
then $\kappa^M$ satisfies these assumptions as well.
\end{itemize}

Hence, we follow the second approach and obtain the following discrete EVI$(\zeta(\kappa^M))$:
\begin{align*}
\frac{1}{2}\frac{\mathrm{d}}{\mathrm{d}t}
\mathcal{W}_{\kappa^M,\Leb_{|\Lambda_n^M}^M}^2\!(\mathrm{S}_t^{\kappa^M*}\dP_0^M,\dP_1^M)
\!+\!
\frac{\zeta(\kappa^M)}{2}
\mathcal{W}_{\kappa^M,\Leb_{|\Lambda_n^M}^M}^2\!(\mathrm{S}_t^{\kappa^M*}\dP_0^M,\dP_1^M)
\leq
\Ent(\dP_1^M|\tilde{\dR}^M)
\!-\!
\Ent(\mathrm{S}_t^{\kappa^M*}\dP_0^M|\tilde{\dR}^M),
\end{align*}
whenever all terms are finite. This finiteness assumption is non-trivial, roughly speaking, because we combine two different discretization procedures. More precisely, if $\Ent(\dP_0|\dR)<\infty$, then the contraction property of the relative entropy yields $\Ent(\dP_0^M|\dR^M)\leq \Ent(\dP_0|\dR)<\infty$, see \cite[Lemma 9.4.5]{Ambrosio2008GradientFI}, but  nevertheless it might happen that $\Ent(\dP_0^M|\tilde{\dR}^M)=\infty$. This discrepancy arises from the fact that under $\tilde{\dR}^M$ points are sampled according to the Papangelou intensity $\kappa$ while, in contrast, $\dR^M$ is obtained by first sampling point configurations according to $\dR$ and subsequently projecting them onto the lattice. If the Papangelou intensity $\kappa$ is allowed to vanish, as for instance in the hard-core model, the latter projection procedure may assign positive probability to configurations that have zero probability under $\tilde{\dR}^M$.\\

These degeneracies could be avoided by assuming that the Papangelou intensity $\kappa$ is uniformly bounded away from zero, ensuring that all configurations have positive probability. However, this assumption is rather restrictive and excludes, for instance, pair-interaction models. We therefore impose the weaker assumption \hypertarget{iproperty:ALB}{$(\mathbf{ALB})$}, which requires that $\kappa$ admits an approximation from above by Papangelou intensities $(\kappa_a)_{a\in\mathbb N}$ such that each $\kappa_a$ is uniformly bounded from below by a positive constant and satisfies \hyperlink{iprop:locallystable}{$(\mathbf{S}(1))$} and \hyperlink{iprop:repulsive}{$(\mathbf{RP})$}.

Instead of applying the discrete EVI to $\tilde \dR^M$, we apply it  to the Papangelou measure $\tilde{\dR}_a^M$ with Papangelou intensity $\kappa_a^M=(\kappa_a)^M$ and reference measure $\Leb^M_{|\Lambda_n^M}$. We obtain that  \begin{align*}
    \frac{1}{2}\frac{\mathrm{d}}{\mathrm{d}t}
\mathcal{W}_{\kappa^M_a}^2(\mathrm{S}_t^{\kappa^M_a*}\dP_0^M,\dP_1^M)
+
\frac{\zeta(\kappa^M_a)}{2}
\mathcal{W}_{\kappa^M_a}^2(\mathrm{S}_t^{\kappa^M_a*}\dP_0^M,\dP_1^M)
\leq
\Ent(\dP_1^M|\tilde{\dR}^M_a)
-
\Ent(\mathrm{S}_t^{\kappa^M_a*}\dP_0^M|\tilde{\dR}^M_a),
\end{align*}
for all $t>0$. 
\subsubsection{Convergence of the Discretization}
Next, we first take the continuum limit by letting $M$ go to infinity while keeping $a$ fixed, and then let $a$ go to infinity. The main technical ingredient of this step is the
following $\Gamma$-convergence result. Due to the particular discretization procedure the following result requires that  $\kappa$ is jointly upper semicontinuous and the maps $(x_1,\ldots,x_{j+1})\longmapsto
\kappa\left(\sum_{i=1}^j\delta_{x_i},x_{j+1}\right)$ are $\Leb\otimes\Leb^{\otimes j}$ a.e.\ continuous \hypertarget{iprop:usc+}{$(\mathbf{USC+})$}.

\begin{satz}
The following $\Gamma$-convergence results hold:
\begin{enumerate}[label=$(\roman*)$]
\item
$\Gamma\text{-}\lim_{M\to\infty}\mathcal{W}_{\kappa_a^M,\Leb_{|\Lambda_n^M}^M}(\cdot,\cdot)
=\mathcal{W}_{\kappa_a,\Leb_{|\Lambda_n}}(\cdot,\cdot)$  for every $a\in\N$; 

\item 
$\Gamma\text{-}\!\lim_{M\to\infty}
\mathrm{Ent}(\cdot|\tilde{\dR}_a^M)\!
=\!\mathrm{Ent}(\cdot|{\dR}_a)$ for every $a\in\N$;

\item $\Gamma\text{-}\lim_{a\to\infty}\mathcal{W}_{\kappa_a,\Leb_{|\Lambda_n}}(\cdot,\cdot)
=\mathcal{W}_{\kappa,\Leb_{|\Lambda_n}}(\cdot,\cdot)$;

\item $\Gamma\text{-}\lim_{a\to\infty}
\mathrm{Ent}(\cdot|{\dR}_a)
=\mathrm{Ent}(\cdot|\dR)$.

\end{enumerate}
\end{satz}
\subsubsection{Infinite-Volume Limit}
In the last step of the proof, we lift the EVI for bounded boxes $\Lambda_n$ to stationary measures.
 The guiding idea is to start from the EVI for  $\Lambda_n$, divide the inequality by the volume of the box, and then pass to the large-volume limit.
 More specifically, let $\dR\in\PPP_s(\Gamma_{\R^d})$ now be a stationary Papangelou measure with Papangelou intensity $\kappa:\Gamma_{\R^d}\to\R$ and reference measure $\Leb$. By the previous steps, we obtain 
\begin{align}
&\liminf_{n\to\infty}
\frac{1}{\Leb(\Lambda_n)}
\bigg(
\frac{1}{2}\frac{\mathrm{d}}{\mathrm{d}t}
\mathcal{W}_{\kappa_n,\Leb_{|\Lambda_n}}^2
\bigl(\mathrm{S}_t^{\kappa_n*}\dP_{0|\Lambda_n},
\dP_{1|\Lambda_n}\bigr)
+
\frac{\zeta(\kappa_n)}{2}
\mathcal{W}_{\kappa_n,\Leb_{|\Lambda_n}}^2
\bigl(\mathrm{S}_t^{\kappa_n*}\dP_{0|\Lambda_n},
\dP_{1|\Lambda_n}\bigr)
\bigg)\nonumber
\\
&\qquad\leq
\liminf_{n\to\infty}
\frac{1}{\Leb(\Lambda_n)}
\bigg(
\Ent(\dP_{1|\Lambda_n}|\dR_{n})
-
\Ent\bigl(
\mathrm{S}_t^{\kappa_n*}\dP_{0|\Lambda_n}
|\dR_{n}
\bigr)
\bigg)\label{eq:1111}
\end{align}
 for all $t>0$ and $\dP_0,\dP_1\in\PPP_s(\Gamma_{\R^d})$, where $\dR_n$ is the finite-volume Papangelou measure with Papangelou intensity $\kappa_n:=\kappa\ind_{\Lambda_n}$ and reference measure $\Leb_{|\Lambda_n}$.

The main difficulty lies in proving the following $\Gamma$-convergence result, which requires that $\kappa$ has finite range \hypertarget{iprop:finiterange}{$(\mathbf{FR})$}, i.e.\ there exists some ball $B_r(0)$ with radius $r$ centered at the origin such that $\kappa(\xi,x)=\kappa(\xi_{|B_r(x)},x)$ for all $\xi\in\Gamma_{\R^d}$, $x\in\R^d$. 
\begin{satz}\label{ithm:gammaconvergence}
The following $\Gamma$-convergence results hold:
\begin{itemize}
\item[$(i)$] 
$\Gamma\text{-}\lim_{n\to\infty}
\frac{1}{\Leb(\Lambda_n)}
\mathcal{W}_{\kappa_n,\Leb_{|\Lambda_n}}^2
(\cdot,\cdot)
=
\mathcal{W}_{\kappa,s}^2(\cdot,\cdot);$
\item[$(ii)$] 
$\Gamma\text{-}\lim_{n\to\infty}
\frac{1}{\Leb(\Lambda_n)}
\Ent(\cdot|\dR_n)
=
\mathcal{E}(\cdot|\dR)$;
\end{itemize}
\end{satz}
Applying Theorem \ref{ithm:gammaconvergence} to \eqref{eq:1111}, we obtain the EVI of Theorem \ref{ithm:evistationary}.

We finish this section with the important observation that Theorem \ref{ithm:gammaconvergence} implies
\begin{align*}
     \mathcal{W}_{\kappa,s}^2(\dP_0,\dP_1)&=\inf\left\{\sup_{n\in\N}\frac{1}{\Leb(\Lambda_n)}\mathcal{A}_{\kappa,\Leb_{|\Lambda_n}}({\bar{\dP}_{|\Lambda_n}},{\bar{\dV}_{|\Lambda_n}}):({\bar{\dP}},{\bar{\dV}})\in\CCC\EEE^s(\dP_0,\dP_1)\right\}\\
            &=\lim_{n\to\infty}\frac{1}{\Leb(\Lambda_n)}\inf\left\{\mathcal{A}_{\kappa,\Leb_{|\Lambda_n}}({\bar{\dP}^n},{\bar{\dV}^n}):({\bar{\dP}^n},{\bar{\dV}^n})\in\CCC\EEE^{\Lambda_n}(\dP_{0|\Lambda_n},\dP_{1|\Lambda_n})\right\}\nonumber \\
            &=\lim_{n\to\infty}\frac{1}{\Leb(\Lambda_n)}\mathcal W_{\kappa,\Leb_{|\Lambda_n}}^2(\dP_{0|\Lambda_n},\dP_{1|\Lambda_n}), \quad \dP_0,\dP_1\in\PPP_s(\Gamma_{\R^d}),\nonumber
\end{align*}
justifying the name \textit{specific} non-local transport distance.
A subtle point in this representation is that the finite-volume minimization on the right-hand side does not impose any stationarity constraint on the admissible curves. In particular, it is not clear why finite-volume minimizers should arise as restrictions of stationary curves.
If it were, the right-hand side would be a perfect definition of $\mathcal W_{\kappa,s}$ allowing us to lift the  EVI for  $\Lambda_n$ directly to the setting of stationary measures.

\section{Preliminaries}
In this section, we introduce and recall notation, definitions, and  several standard results that will be used throughout this paper. We begin with an overview of the notation used for functions and measures on a Polish metric space $X$.

 We denote by $\BBB(X)$ the Borel sets of $X$ and by $\BBB_0(X)$ those that are additionally bounded.
 A set $B\in\BBB(X)$ is a Polish subspace of $X$ if and only if $B$ is a $G_\delta$ set, i.e.\ a countable intersection of open sets, see e.g.\ \cite[Theorem 24.12]{Willard}.

 We denote by $\FFF(X)$ the space of measurable functions mapping $X$ to $\R$. For continuous functions, we write $\CCC(X)$. The space of continuous and compactly supported functions is $\CCC_c(X)$ and the space of continuous and bounded functions is $\CCC_b(X)$. The space of locally finite Borel signed measures on $X$ is $\mathcal{M}(X)$ and we endow it with the vague topology generated by the maps $\lambda\mapsto\int f(x)\lambda(\mathrm{d}x)$, $f\in\CCC_c(X)$. 
For $p\in[1,\infty)$ and $\lambda\in\MMM(X)$, $L^p(\lambda)$ denotes the space of all $p$-integrable functions, i.e.\  functions $f\in\FFF(X)$ s.t.\   $\|f\|_{L^p(\lambda)}:=(\int |f(x)|^p\lambda(\mathrm{d}x))^{1/p}<\infty$. If a function or a measure space is assumed to contain only non-negative functions or measures, this will always be indicated by an additional superscript $+$. 
\subsection{Point Processes}
We continue with a brief introduction to the theory of point processes. For a detailed overview, we refer e.g.\ to \cite{daley2003introduction, last_penrose_2017}.
 We denote by $\Gamma_X$ the space of all \textit{locally finite  counting measures} on $X$. If $X=\R^d$, we briefly write $\Gamma$ instead of $\Gamma_{\R^d}$. We endow $\Gamma_X$ with the {vague topology} generated by the maps $\xi\mapsto\int f(x)\xi(\mathrm{d}x),\;f\in\CCC_c(X)$. 
 A \textit{point process} on the base space $X$ is a random variable taking values in $\Gamma_X$. The distribution $\dP$ of a point process is an element of $\PPP(\Gamma_X)$, the space of probability measures on $\Gamma_X$. 
 The \textit{intensity measure} of  $\dP\in\PPP(\Gamma_X)$ is defined by
  \[
\I_\dP(B):=\int\int\ind_B(x)\xi(\mathrm{d}x){\dP}(\mathrm{d}\xi),\quad \text{$B\in\BBB(X)$}.
\]
A measure $\dP$ is \textit{locally integrable} if $\I_\dP(B)<\infty$ for all $B\in\BBB_0(X)$. Further, a point configuration $\xi\in\Gamma_X$  is called \textit{simple} if $\xi(\{x\})\leq1$ for all $x\in X$. A measure $\dP\in\PPP(\Gamma_X)$ is called simple if ${\dP(\{\xi\in\Gamma_X:\xi(\{x\})\leq1 \text{ for all }x\in X\})=1}$. 
 The space of all locally integrable and simple probability measures on $\Gamma_X$ is $\PPP_1(\Gamma_X)$.

 An $\R$-valued measurable function on the configuration space $F\in\FFF(\Gamma_X)$ has sublinear growth if there exists some  constant $c>0$ and some $f\in\CCC_c^+(X)$ s.t.\   
\begin{align*}
    |F(\xi)|\leq c\left(1+\int f(x)\xi(\mathrm{d}x)\right),\quad \xi\in\Gamma_X.
\end{align*}
   The set of all continuous functions having sublinear growth is denoted by $\CCC_1(\Gamma_X)$. Unless stated otherwise, we endow $\PPP_1(\Gamma_X)$ with the $\CCC_1(\Gamma_X)$-weak topology, i.e.\ the topology generated by the maps $\dP\mapsto\int F(\xi)\dP(\mathrm{d}\xi)$, $F\in\CCC_1(\Gamma_X)$.  This topology turns
    $\PPP_1(\Gamma_X)$ into a Polish space, see \cite[Theorem 2.11]{schiavo2023wasserstein}. Closely connected to the intensity measure is the \textit{reduced Campbell measure} 
\begin{align*}
    C_\dP(\A\times B):=\int\int \ind_B(x)\ind_\A(\xi-\delta_x)\xi(\mathrm{d}x)\dP(\mathrm{d}\xi), \quad \A\in\BBB(\Gamma_X), \; B\in\BBB(X),
\end{align*}
where $\delta_x$ is the \textit{Dirac measure} at $x$.  The reduced Campbell measure is a measure on the product space $\Gamma_X\times X$. We endow $\Gamma_X\times X$ with the product topology and 
define $\MMM_{b,0}(\Gamma_X\times X)$ to be the set of all signed measures $\dV\in\MMM(\Gamma_X\times X)$  s.t.\   $|\dV|(\Gamma_X\times B)<\infty$ for all $B\in\BBB_0(X)$. Further, $\FFF_{b,0}(\Gamma_X\times X)$ is the set of all measurable $\R$-valued functions on $\Gamma_X\times X$ that are bounded and vanish outside a set of the form $\Gamma_X\times B$ for some $B\in\BBB_0(X)$.
Moreover,  $\mathcal{C}_{b,0}(\Gamma_X\times X)$ is the set of functions that satisfy the previous condition and are additionally continuous.
Unless stated otherwise, we always equip $\MMM_{b,0}(\Gamma_X\times X)$ with the $\CCC_{b,0}(\Gamma_X\times X)$-weak topology, i.e.\ the topology induced by the maps 
    $\dV\mapsto|\int F(\xi,x)\dV(\mathrm{d}\xi,\mathrm{d}x)|, \ F\in \mathcal{C}_{b,0}(\Gamma_X\times X)$. Since $C_\dP(\Gamma_X\times B)=\I_\dP(B)$, $B\in\BBB(X)$, it follows that  $C_\dP\in\MMM_{b,0}^+(\Gamma_X\times X)$ whenever $\dP\in\PPP_1(\Gamma_X)$.\\

    We call $B\in\mathcal{B}(X)$ a \textit{continuity set} w.r.t.\ $\dP\in\PPP_1(\Gamma_X)$ if $\dP(\{\xi\in\Gamma_X:\xi(\partial B)=0\})=1$, where $\partial B$ denotes the boundary of $B$. A sequence of measures $(\dP_n)_{n\in\N}\subset\PPP_1(\Gamma_X)$ converges in finite-dimensional distributions to a measure $\dP\in\PPP_1(\Gamma_X)$, if for all $m \in \N$, $a_1, \ldots, a_m \in \N$, and   $B_1, \ldots, B_m \in \BBB_0(X)$ continuity sets of $\dP$, it holds that
        \begin{align*}
            \lim_{n \to \infty} \int \ind_{\{\xi(B_1) = a_1, \ldots, \xi(B_m) = a_m\}}  \dP_n(\mathrm{d}\xi)
            = \int \ind_{\{\xi(B_1) = a_1, \ldots, \xi(B_m) = a_m\}}  \dP(\mathrm{d}\xi).
            \end{align*}
    The convergence in $\PPP_1(\Gamma_X)$ can be characterized either by convergence in finite-dimensional distributions, by convergence of the reduced Campbell measure, or by using the \textit{Laplace functional} $$L_\dP(f):=\int e^{-\int f(x) \xi(\mathrm{d}x)}\dP(\mathrm{d}\xi),\quad f\in\CCC_c^+(X).$$

\begin{prop}[{\cite[Proposition 2.12]{schiavo2023wasserstein}, \cite[Proposition 11.1.VIII]{daley2007introduction}}]\label{prop:convergencep1}
    Consider a sequence $(\dP_n)_{n\in\N}\subset\PPP_1(\Gamma_X)$ and let $\dP\in\PPP_1(\Gamma_X)$. The following are equivalent:
    \begin{enumerate}
       \item[$(i)$] $\lim_{n\to\infty}\dP_n=\dP$ in $\PPP_1(\Gamma_X)$;
       \item[$(ii)$] $\lim_{n\to\infty}\I_{\dP_n}=\I_{\dP}$ vaguely in $\MMM^+(X)$ and $\lim_{n\to\infty}L_{\dP_n}(f)=L_{\dP}(f)$ for all $f\in\CCC_c^+(X)$;
       \item[$(iii)$] $\lim_{n\to\infty}\I_{\dP_n}=\I_{\dP}$ vaguely in $\MMM^+(X)$ and $\lim_{n\to\infty}\dP_n= \dP$ in the sense of finite-dimensional distributions;
       \item[$(iv)$] $\lim_{n\to\infty} C_{\dP_n}=C_\dP$ in $\MMM_{b,0}(\Gamma_X\times X)$. 
    \end{enumerate}
\end{prop}
 Prominent examples of point processes are \textit{Poisson point processes} $\pi^\lambda$   with intensity measure $\lambda\in\MMM^+(X)$. They are characterized by the following two properties: 
\begin{enumerate}
    \item [$(i)$]  The point statistics $\pi^\lambda(B):=\int\ind_B(x)\pi^\lambda(\mathrm{d}x)$ is   Poisson distributed with parameter $\lambda(B)$ for any $B\in\BBB(X)$.
    \item [$(ii)$] The point statistics  $\pi^\lambda(B_1),\dots,\pi^\lambda(B_n)$ are independent for all pairwise disjoint sets   $B_1,\dots , B_n \in\BBB(X)$.
\end{enumerate}
We denote the distribution of $\pi^\lambda$ by $\Poi^\lambda$. The Poisson point process $\pi^\lambda$ is called a homogeneous Poisson point process if $\pi^\lambda$  is a Poisson point process on $\R^d$ with $\lambda = c\Leb$, where $c > 0$ and $\Leb$ is the $d$-dimensional Lebesgue measure. In this case, we briefly denote the distribution by $\Poi^c$. If $c=1$, we simply write $\Poi$.
\subsubsection{Stationary Point Processes}For $B\in\BBB(\R^d)$ and $z\in\R^d$, let $B+z:=\{x+z:x\in B\}$ and $B-z:=B+(-z)$. We define the shift operator $\Theta_z:\Gamma\to\Gamma,\; z\in\R^d$, by $$\Theta_z\xi(B):=\xi(B+z), \quad B\in\BBB(\R^d).$$  A measure $\dP\in\PPP_1(\Gamma)$ is  called \textit{{stationary}} if  ${(\Theta_z)}_\#\dP=\dP$ for all $z\in\R^d$, where $(\Theta_z)_\#\dP$ denotes the pushforward measure of $\dP$ under $\Theta_z$.  We denote the space of all stationary measures by $\PPP_s(\Gamma)$ and call a point process stationary if its distribution is a stationary measure. Analogously to $\PPP_1(\Gamma)$, we endow $\PPP_s(\Gamma)$ with the $\CCC_1(\Gamma)$-weak topology  induced by  the maps $\dP\mapsto\int F(\xi)\dP(\mathrm{d}\xi)$, $F\in\CCC_1(\Gamma)$. 
The intensity measure of a stationary measure $\dP\in\PPP_s(\Gamma)$ is always a multiple of the Lebesgue measure, i.e.\  $\I_\dP=c\Leb$ for $c\geq0$. 
\subsubsection{Transformation of Point Processes}\label{subsec:transformation}
We continue by introducing transformations of point configurations on $\R^d$. Most of them can be generalized to point configurations on more general Polish spaces, but for the sake of presentation we focus on the Euclidean setting. 

 Let $B\in\BBB(\R^d)$. We define the \textit{projection} map $\mathrm{pr}_B:\Gamma\to\Gamma_B$ by  $\mathrm{pr}_B(\xi):=\xi(B\cap\cdot)=\xi_{|B}$.  The \textit{restriction} of 
$\dP\in\PPP_1(\Gamma)$ to $\Gamma_B$ is defined by  $\dP_{|B}:=(\mathrm{pr}_B)_{\#} \dP$.

Let $B_1,B_2\in\BBB(\R^d)$. We define the \textit{superposition} of point configurations $\mathrm{sp}:\Gamma_{B_1}\times\Gamma_{B_2}\to\Gamma_{{B_1} \cup B_2}$ by $\mathrm{sp}(\xi^{1},\xi^2):=\xi^1+\xi^2$. For  measures $\dP^1\in\PPP_1(\Gamma_{B_1}),\;\dP ^2\in\PPP_1(\Gamma_{B_2})$, we define the \textit{independent superposition} by $\dP^1\oplus\dP^2:=\mathrm{sp}_\#(\dP^1\otimes\dP^2)$. In particular, it holds that 
\begin{align*}
    \int F(\xi) \dP^1\oplus\dP^2(\mathrm{d}\xi)=\int\int F(\xi^1+\xi^2) \dP^1(\mathrm{d}\xi^1)\dP^2(\mathrm{d}\xi^2), \quad F\in\FFF^+(\Gamma_{B_1\cup B_2}).
\end{align*} 

Let $M,n\in\N$. Consider  $\Lambda_n:=[-n/2,n/2)^d$ and define the shifted lattice with mesh size $\frac{1}{M}$ by $$\Lambda_n^M:=\left(\frac{1}{M}\Z^d-\left(\frac{1}{2M}-\frac{n}{2}\right)\mathbf{1}\right)\cap\Lambda_n,$$ where $\mathbf{1}=(1,\ldots,1)$.
The translated cells $(z+\Lambda_{\frac1M})_{z\in\Lambda_n^M}:=(z+[-\frac{1}{2M},\frac{1}{2M})^d)_{z\in\Lambda_n^M}$ form a partition of $\Lambda_n$, i.e.\
\(
\Lambda_n=\dot\bigcup_{z\in\Lambda_n^M}(z+\Lambda_{\frac1M}).
\) The discretization map $T^M:\Lambda_n\to\Lambda_n^M$ is defined by
\[
T^M(x)
:=x^M=
\sum_{z^M\in\Lambda_n^M}
z^M\ind_{z^M+\Lambda_{\frac{1}{M}}}(x).
\]
The discretization of a measure 
$\lambda\in\MMM^+(\Lambda_n)$ is $\lambda^M:=(T^M)_\#\lambda$. Similarly, we define the discretization of point configurations $T_\Gamma^M:\Gamma_{\Lambda_n}\to\Gamma_{\Lambda_n^M}$
by
\[
T_\Gamma^M(\xi)
:=\xi^M=
\sum_{z^M\in\Lambda_n^M}
\delta_{z^M}\ind_{\{\xi(z^M+\Lambda_{\frac{1}{M}})\ge 1\}}.
\]
The discretization of $\dP\in\PPP_1(\Gamma_{\Lambda_n})$ is defined by $\dP^M:=(T_\Gamma^M)_\#\dP$. The following lemma shows that the discretized measures converge to the original measure as the mesh size  tends to zero.

    \begin{lem}\label{lem:convergenceQM}
        Let   $\dP\in\PPP_1(\Gamma_{\Lambda_n})$. It holds that 
            $\lim_{M\to\infty}\dP^M=\dP\quad \text{in } \PPP_1(\Gamma_{\Lambda_n})$.

    \end{lem}
    \begin{proof}
  By Proposition~\ref{prop:convergencep1}, it suffices to show that  $\lim_{M\to\infty}L_{\dP^M}(f)=L_{\dP}(f)$ for all functions $f\in\CCC^+_c({\Lambda_n})$ and $\lim_{M\to\infty}\I_{\dP^M}=\I_{\dP}$ vaguely. Let $f\in\CCC^+_c({\Lambda_n})$. By the definition of $\dP^M$ and the dominated convergence theorem, we obtain  that 
   \begin{align*}
       \lim_{M\to\infty}L_{\dP^M}(f)&=\lim_{M\to\infty}\int e^{-\int f(y)\xi^M(\mathrm{d}y)}\dP(\mathrm{d}\xi)=\int\lim_{M\to\infty} e^{-\int f(y)\xi^M(\mathrm{d}y)}\dP(\mathrm{d}\xi)=L_{\dP}(f).
       \end{align*}
       For $B\in\BBB_0(\R^d)$, by definition $ \sum_{z^M\in B\cap\Lambda_n^M}
\ind_{\{\xi(z^M+\Lambda_{\frac{1}{M}})\ge1\}}\leq \xi(\cup_{z^M\in B\cap\Lambda_n^M}(z^M+\Lambda_{\frac{1}{M}}))$, which is $\dP$-integrable as $\dP\in\PPP_1(\Gamma_{\Lambda_n})$. Thus, as before, we obtain by the dominated convergence theorem that
\[
\lim_{M\to\infty}\I_{\dP^M}(B)
=
\lim_{M\to\infty}
\int
\sum_{z^M\in B\cap\Lambda_n^M}
\ind_{\{\xi(z^M+\Lambda_{\frac{1}{M}})\ge1\}}
\dP(\mathrm d\xi)
=
\int \xi(B)\dP(\mathrm d\xi)
=
\I_{\dP}(B),
\]
proving that $\lim_{M\to\infty}\I_{\dP^M}=\I_{\dP}$ vaguely and finishing the proof. 
    \end{proof}
\subsubsection{Papangelou Measures}\label{sec:preliminariespapangelou}
In this subsection, we define  \textit{Papangelou measures}, first introduced in \cite{papangelou1974conditional}, and discuss the most relevant properties for this manuscript.
We follow mainly \cite{Betsch23}, which provides a rigorous presentation of Papangelou measures on general Polish spaces $X$. For further references, see e.g.\ \cite{Papangelou, daley2007introduction}. 
\begin{dfn}[Papangelou measures]\label{dfn:papangelou}
     A probability measure $\dR\in\PPP_1(\Gamma_X)$ is called \textit{Papangelou measure} w.r.t.\ some measurable function $\kappa:\Gamma_X\times X\to[0,\infty)$, called \textit{Papangelou intensity}, and a reference measure $\lambda\in\MMM^+(X)$ if  
    \begin{align}\label{eq:papangelou}
        \int\int F(\xi-\delta_x,x)\xi(\mathrm{d}x)\dR(\mathrm{d}\xi)=\int\int F(\xi,x)\kappa(\xi,x)\lambda(\mathrm{d}x)\dR(\mathrm{d}\xi),\quad F\in\FFF^+(\Gamma_X\times X).
    \end{align}
    The set of all Papangelou measures w.r.t.\ the Papangelou intensity $\kappa$ and the reference measure $\lambda$ is denoted by $\RRR_{\kappa,\lambda}(\Gamma_X)$.
\end{dfn}
\begin{bem}
 Equation \eqref{eq:papangelou} is known as the GNZ equation, named after Georgii \cite{Geogii76}, Nguyen, and Zessin \cite{NguyenZessin1979}, and can be interpreted as a detailed balance equation for the spatial birth-and-death process introduced in Section~\ref{sec:spatialbirthanddeath}.  In this process, points are removed at unit rate and added with rate given by $\kappa(\xi,x)\lambda(\mathrm{d}x)$. 
    \end{bem}
Note that Definition~\ref{dfn:papangelou} covers both finite and infinite base spaces $X$. However, the existence and uniqueness of measures satisfying the GNZ equation~\eqref{eq:papangelou} are, in particular for infinite base spaces, non-trivial and will be discussed below. If a Papangelou measure $\dR\in\RRR_{\kappa,\lambda}(\Gamma_X)$ exists, then $\kappa$ satisfies the following {cocycle property}, see e.g.\ \cite[Lemma 2.3]{Betsch23}.
\begin{enumerate}[label=$(\mathbf{CC})$,leftmargin=*, align=left]
    \item \label{property:cocycle}  The Papangelou intensity $\kappa:\Gamma_X\times X\to[0,\infty)$ satisfies the \textit{cocycle property} \begin{align*}
        \kappa(\xi+\delta_x,y)\kappa(\xi,x)=\kappa(\xi+\delta_y,x)\kappa(\xi,y), \quad \text{for all $\xi\in\Gamma_X,\; x,y\in X$}.
    \end{align*}
\end{enumerate}
\begin{kor}\label{kor:CampbellPapangelou}
Let $\lambda \in \MMM^+(X)$ and assume that $\kappa : \Gamma_X \times X \to [0,\infty)$ satisfies \ref{property:cocycle}. If $\dR \in \RRR_{\kappa,\lambda}(\Gamma_X)$ exists, then  the intensity measure of $\dR$ is 
\begin{align*}
\I_\dR(\mathrm{d}x)
=
\int \kappa(\xi,x)\dR(\mathrm{d}\xi)\lambda(\mathrm{d}x),
\end{align*}
and the reduced Campbell measure of $\dR$ is 
\begin{align*}
C_{\dR}(\mathrm{d}\xi,\mathrm{d}x)
=
\kappa(\xi,x)\dR(\mathrm{d}\xi)\lambda(\mathrm{d}x).
\end{align*}
Further, let $\dP \in \PPP_1(\Gamma_X)$ be absolutely continuous w.r.t.\ $\dR$, i.e.\  $\dP\ll\dR$, and  $\dP = \rho \dR$. The reduced Campbell measure of $\dP$ is
\begin{align*}
C_{\dP}(\mathrm{d}\xi,\mathrm{d}x)
=
\rho(\xi+\delta_x)\kappa(\xi,x)\dR(\mathrm{d}\xi)\lambda(\mathrm{d}x).
\end{align*}
\end{kor}
\begin{proof}
    The claim follows directly from the GNZ equation~\eqref{eq:papangelou}.
\end{proof}
The \textit{partition function} plays a central role in the analysis of Papangelou measures. Let $B \in \mathcal{B}(X)$ and define the partition function  $Z_{\kappa,\lambda}^B : \Gamma_{B^C} \to [0,\infty]$ by
\begin{align*}
Z_{\kappa,\lambda}^{B}(\eta)
:=
\sum_{j\geq0}\frac{1}{j!}
\int_B \cdots \int_B
\varphi_\kappa^{\eta}\left(\sum_{i=1}^j\delta_{x_i}\right)
\lambda(\mathrm{d}x_1)\cdots\lambda(\mathrm{d}x_j),
\end{align*}
with
\begin{align*}
\varphi_\kappa^\eta\left(\sum_{i=1}^j\delta_{x_i}\right)
:=
\prod_{i=1}^j
\kappa\left(\sum_{l=1}^{i-1}\delta_{x_l}+\eta,x_i\right),
\quad j\geq0.
\end{align*}
For the zero measure, we set $\varphi_\kappa^\eta\left(\sum_{i=1}^0\delta_{x_i}\right)=\varphi_\kappa^\eta(0):=1$. 
The cocycle property \ref{property:cocycle} implies that $\varphi_\kappa^\eta$ is well-defined and does not depend on the ordering of the points in $\sum_{i=1}^j\delta_{x_i}$.
For empty boundary conditions, i.e.\  $\eta=0$, we briefly write $\varphi_\kappa$ and $Z_{\kappa,\lambda}^{B}$ instead of $\varphi_\kappa^0$ and $Z_{\kappa,\lambda}^{B}(0)$, respectively.
\begin{lem}[{\cite[Lemma 4.2]{Betsch23}}]\label{lem:finitevolumePapangelou}
Let $\lambda\in\MMM^+(X)$ and assume that $\kappa:\Gamma_X\times X\to[0,\infty)$ satisfies \ref{property:cocycle}. Further, let $B \in \mathcal{B}(X)$  s.t.\   $\lambda(B)<\infty$. A measure  $\dR^B(\cdot|\eta)\in\PPP_1(\Gamma_B)$ is called finite-volume Papangelou measure with Papangelou intensity $\kappa$, reference measure $\lambda$, and boundary condition $\eta\in\Gamma_{B^C}$ if it satisfies the GNZ equation~\eqref{eq:papangelou} w.r.t.\  the Papangelou intensity $\kappa^\eta(\xi,x):=\kappa(\xi+\eta,x)$ and the reference measure $\lambda_{|B}$, i.e.\ $\dR^B(\cdot|\eta)\in\RRR_{\kappa^\eta,\lambda_{|B}}(\Gamma_B)$. Such a measure exists  if and only if $Z_{\kappa,\lambda}^B(\eta)<\infty$. Moreover, it is uniquely defined and it holds that
\begin{align*}
\int F(\xi)\dR^B(\mathrm{d}\xi | \eta)
=
\frac{1}{Z_{\kappa,\lambda}^B(\eta)}
\sum_{j\geq0}\frac{1}{j!}
\int_B \cdots \int_B
F\left(\sum_{i=1}^j\delta_{x_i}\right)
\varphi_\kappa^\eta\left(\sum_{i=1}^j\delta_{x_i}\right)
\lambda(\mathrm{d}x_1)\cdots \lambda(\mathrm{d}x_j),
\end{align*}
for all functions $F\in\FFF^+(\Gamma_B)$. In particular,
\begin{align*}
\dR^B(\mathrm{d}\xi | \eta)
=
\frac{e^{\lambda(B)}}{Z_{\kappa,\lambda}^B(\eta)}
\varphi_\kappa^\eta(\xi)
\Poi^{\lambda_{|B}}(\mathrm{d}\xi).
\end{align*}
\end{lem}
Unless stated otherwise, we always assume empty boundary conditions for finite-volume Papangelou measures and briefly write $\dR^B(\cdot)$ instead of $\dR^B(\cdot | 0 )$. 
Besides the GNZ equation \eqref{eq:papangelou}, another characterization of Papangelou measures is given by the \textit{DLR equation}, named after Dobrushin, Lanford, and Ruelle.

\begin{satz}[{\cite[Corollary 6.3]{Betsch23}}]\label{thm:DLR}
Let $\lambda\in\MMM^+(X)$ and assume that $\kappa:\Gamma_X\times X\to[0,\infty)$ satisfies \ref{property:cocycle}.  Then  $\dR \in \RRR_{\kappa,\lambda}(\Gamma_X)$ if and only if
\begin{align*}
\int F(\xi)\dR(\mathrm{d}\xi)
=
\int \int
F(\xi^{B}+\xi_{|B^C})
\dR^B(\mathrm{d}\xi^B | \xi_{|B^C})
\dR(\mathrm{d}\xi),
\end{align*}
for all $B \in \BBB_0(X),\;F \in \FFF^+(\Gamma_X)$.
\end{satz}

We now introduce several assumptions on the Papangelou intensities $\kappa$, the reference measure $\lambda$, and, consequently, on the associated Papangelou measures $\dR\in\RRR_{\kappa,\lambda}(\Gamma_X)$. Unless stated otherwise, we always assume that all Papangelou intensities are measurable and satisfy the following non-degeneracy assumption.

\begin{enumerate}[label=$(\mathbf{ND})$, leftmargin=*,leftmargin=1.7cm,
  labelwidth=1.5cm,
  labelsep=0.1cm,
  align=left]
\item\label{property:nondegenerate}
The Papangelou intensity $\kappa:\Gamma_X\times X\to[0,\infty)$ is \textit{non-degenerate}, i.e.\  for all $x\in X$, it holds that $\kappa(0,x)>0$.
\end{enumerate}

We continue with a classical stability assumption.

\begin{enumerate}[label=$(\mathbf{S}(A))$, leftmargin=*,leftmargin=1.7cm,
  labelwidth=1.5cm,
  labelsep=0.1cm,
  align=left]
\item\label{property:locallystable}
The Papangelou intensity $\kappa:\Gamma_X\times X\to[0,\infty)$ is \textit{locally stable} with parameter $A<\infty$, i.e.\ it holds that
\[
\kappa(\xi,x) \leq A, \quad \text{for all } \xi \in \Gamma_X,\; x \in X.
\]
\end{enumerate}
If $\kappa$ satisfies \ref{property:cocycle} and is locally stable \ref{property:locallystable}, then the finite-volume Papangelou measure  with empty boundary condition $\dR^B$ exists and is unique by Lemma~\ref{lem:finitevolumePapangelou}. To emphasize this uniqueness we write $\dR^B=\mathcal{R}_{\kappa,\lambda}(\Gamma_B)$ instead of $\dR^B\in\mathcal{R}_{\kappa,\lambda}(\Gamma_B)$.  Further, for any $c>0$, it holds that
\begin{align}
\int c^{\xi(B)}\dR^B(\mathrm{d}\xi)
\leq
\int \frac{e^{\lambda(B)}}{Z_{\kappa,\lambda}^B}
(cA)^{\xi(B)}\Poi^{\lambda_{|B}}(\mathrm{d}\xi)
< \infty, \label{property:integrability}
\end{align}
since $\int (cA)^{\xi(B)} \Poi^{\lambda_{|B}}(\mathrm{d}\xi) < \infty$ by the definition of the Poisson measure. Moreover, \ref{property:locallystable} ensures that, for the birth-and-death dynamics introduced in Section~\ref{sec:spatialbirthanddeath}, the rate of adding a point to a configuration is uniformly bounded and, in particular, does not explode. 

Next, we introduce assumptions regarding the continuity of the Papangelou intensity.
\begin{enumerate}[label=$(\mathbf{USC})$, leftmargin=*,leftmargin=1.7cm,
  labelwidth=1.5cm,
  labelsep=0.1cm,
  align=left]
\item\label{property:continuity}
The Papangelou intensity $\kappa:\Gamma_X\times X\to[0,\infty)$ is jointly upper semicontinuous.
\end{enumerate}
\begin{enumerate}[label=$(\mathbf{USC+})$, leftmargin=*,leftmargin=1.7cm,
  labelwidth=1.5cm,
  labelsep=0.1cm,
  align=left] \item\label{property:continuity2} Let $\lambda\in\MMM^+(X)$. The Papangelou intensity $\kappa:\Gamma_X\times X\to[0,\infty)$ is jointly upper semicontinuous and, for all $j\in\N$, the function $$(x_1,\ldots,x_j,x_{j+1})\mapsto\kappa\left(\sum_{i=1}^j\delta_{x_i},x_{j+1}\right)$$ is $\lambda^{\otimes j}\otimes\lambda$-a.e.\ continuous. \end{enumerate}
If \ref{property:continuity2} is satisfied and $X$ is a bounded Polish metric space, then Lemma~\ref{lem:finitevolumePapangelou} directly implies that $\kappa:\Gamma_X\times X\to[0,\infty)$ is $\dR\otimes\lambda$-a.e.\ continuous.  We note that not every upper semicontinuous $\kappa$ satisfies the additional a.e.-continuity assumption in \ref{property:continuity2}. However, the additional assumption still includes a large class of Papangelou measures, see Example~\ref{exp:Papangelou}.
\begin{enumerate}[label=$(\mathbf{RP})$, leftmargin=*,leftmargin=1.7cm,
  labelwidth=1.5cm,
  labelsep=0.1cm,
  align=left]
\item\label{property:repulsive}
The Papangelou intensity $\kappa:\Gamma_X\times X\to[0,\infty)$ is \textit{repulsive}, i.e.\ 
\[
\kappa(\xi_1+\xi_2,x) \leq \kappa(\xi_1,x),\quad \text{for all } \xi_1,\xi_2 \in \Gamma_X,\; x \in X.
\]
\end{enumerate}
In view of the spatial birth-and-death process introduced in Section~\ref{sec:spatialbirthanddeath}, repulsiveness \ref{property:repulsive} implies that the birth rate described by $\kappa(\xi,x)\lambda(\mathrm{d}x)$ decreases as the number of existing points increases. Furthermore, \ref{property:repulsive} implies the so-called \textit{heredity} assumption: if $\kappa(\xi,x)=0$ for $\xi\in\Gamma_X,\;x\in X$, then $\kappa(\xi+\delta_y,x)=0$ for all $y\in X$.

We continue with the discussion of assumptions  dealing with the range of the interactions between points.
\begin{enumerate}[label=$(\mathbf{D}(\zeta(\kappa)))$, leftmargin=*,leftmargin=1.85cm,
  labelwidth=1.7cm,
  labelsep=0.15cm,
  align=left]
\item\label{property:dobrushin}
The Papangelou intensity $\kappa:\Gamma_X\times X\to[0,\infty)$ and the reference measure $\lambda\in\MMM^+(X)$ satisfy the \textit{Dobrushin condition}
\begin{align*}
0<\zeta(\kappa)
:=
\frac12\bigg(1-\sup_{\substack{\xi \in \Gamma_X \\ x \in X}}
\int |\kappa(\xi,y)-\kappa(\xi+\delta_x,y)|\lambda(\mathrm{d}y)\bigg).
\end{align*}
\end{enumerate}

The Dobrushin  condition originates from statistical physics, see e.g.\ \cite{Dobruschin1968TheDO}, and is classically used to ensure existence and uniqueness of Gibbs measures, see e.g.\ \cite{Ruelle, Houdebert_Zass_2022}. 
Every Papangelou intensity satisfying the following finite-range assumption and~\ref{property:locallystable} satisfies the Dobrushin condition \ref{property:dobrushin} after a suitable rescaling.
 \begin{enumerate}[label=$(\mathbf{FR})$, leftmargin=1.75cm,
  labelwidth=1.6cm,
  labelsep=0.15cm,
  align=left]
         \item \label{property:finiterange} The Papangelou intensity $\kappa:\Gamma_X\times X\to[0,\infty)$ has finite range, i.e.\  there exists some $0\leq r<\infty$  s.t.\   $$\kappa(\xi,x)=\kappa(\xi_{|B_r(x)},x),\quad \text{for all } \xi\in\Gamma_X,\; x\in X,$$  where $B_r(x)$ denotes the ball of  radius $r$  centered at $x$.
        \end{enumerate}
In the context of stationary measures, the following shift-invariance is a natural assumption.
\begin{enumerate}
        [label=$(\mathbf{SI})$, leftmargin=1.75cm,
  labelwidth=1.6cm,
  labelsep=0.15cm,
  align=left]
            \item \label{property:shiftinvariance} The Papangelou intensity $\kappa:\Gamma\times\R^d\to[0,\infty)$ is shift-invariant, i.e.\  $$\kappa(\Theta_z\xi,x-z)=\kappa(\xi,x),\quad \text{for all } \xi\in\Gamma,\; x,z\in \R^d.$$ 
       \end{enumerate}
       Assume that $\kappa:\Gamma\times\R^d\to[0,\infty)$ satisfies \ref{property:cocycle}, \ref{property:locallystable}, \ref{property:finiterange}, and~\ref{property:shiftinvariance}. By \cite[Theorem 8.1]{Betsch23}, there exists an infinite-volume stationary Papangelou measure w.r.t.\  the Papangelou intensity $\kappa$ and the reference measure $\Leb$. The set of all these measures will be denoted by $\RRR_{\kappa,s}(\Gamma)$.
       In the stationary setting, we will always assume that \ref{property:cocycle}, \ref{property:locallystable}, \ref{property:finiterange}, and~\ref{property:shiftinvariance} are satisfied. Therefore, we will not mention the existence of a stationary Papangelou measure explicitly each time.
       \\

      Finally, we introduce assumptions that are  specifically required for a discretization procedure, see Section~\ref{sec:gradientflows} for more details. Let $M,n\in\N$ and recall the definition of discretization of point configurations and measures introduced in Subsection~\ref{subsec:transformation}. For $\kappa:\Gamma_{\Lambda_n}\times\Lambda_n\to[0,\infty)$, we define $\kappa^M:\Gamma_{\Lambda_n^M}\times\Lambda_n^M\to[0,\infty)$ by $$\kappa^M(\xi^M,x^M):=\kappa(\xi^M,x^M)\ind_{\{\xi^M(\{x^M\})=0\}}.$$ 
       We introduce a special form of the Dobrushin condition ensuring that $\zeta(\kappa^M)$ converges to $\zeta(\kappa)$ as the mesh size $\frac{1}{M}$ tends to zero. \begin{enumerate}[label=$(\mathbf{D}^\prime(\zeta(\kappa)))$, leftmargin=1.75cm,
  labelwidth=1.6cm,
  labelsep=0.15cm,
  align=left] \item \label{property:dobrushin2} The Papangelou intensity  $\kappa:\Gamma_{\Lambda_n}\times\Lambda_n\to[0,\infty)$ satisfies the Dobrushin condition \ref{property:dobrushin} and it holds that \begin{align*} 0&<\lim_{M\to\infty}\zeta(\kappa^M)=\frac12\bigg(1-\lim_{M\to\infty} \sup_{\substack{\xi^M\in\Gamma_{\Lambda_n^M}\\x^M\in\Lambda_n^M}}\int_{\Lambda_n^M} |\kappa^M(\xi^M,y^M)-\kappa^M(\xi^M+\delta_{x^M},y^M)|\Leb^M(\mathrm{d}y^M)\bigg)=\zeta(\kappa). \end{align*} 
\end{enumerate} 
When discretizing Papangelou measures, technical difficulties arise whenever $\kappa(\xi,x)=0$, see Section~\ref{sec:convergencediscretization} for more details. This can be avoided by assuming that $\kappa$ admits a strictly positive lower bound.
\begin{enumerate}[label=$(\mathbf{LB}(a ))$, leftmargin=1.75cm,
  labelwidth=1.6cm,
  labelsep=0.15cm,
  align=left]
\item\label{property:lowerbound}
The Papangelou intensity $\kappa:\Gamma_{\Lambda_n}\times {\Lambda_n}\to[0,\infty)$ is $\dR$-a.s.\ lower bounded, i.e.\ there exists an $a\in\N$ and a set $E\in\BBB(\Gamma_{\Lambda_n})$ s.t.\   $ \dR(E)=1$, $\xi+\delta_x\in E$ implies $\xi\in E$, $\xi\in E$ implies $\xi^M\in E$, and
\[
\frac{1}{a }\leq\kappa(\xi,x), \quad \text{for all } \xi+\delta_x \in E.
\]
\end{enumerate}
\begin{bem}\label{bem:lowerboundcampbell}
    If \ref{property:lowerbound} is satisfied, then for any $B\in\BBB(\Lambda_n)$ it holds that
    \begin{align*}
        \int \ind_B(x)\ind_{\{\kappa(\xi,x)<\frac{1}{a}\}}C_\dR(\mathrm{d}\xi,\mathrm{d}x)&\leq\int\ind_B(x)\ind_{\{\xi+\delta_x\in E^C\}}C_\dR(\mathrm{d}\xi,\mathrm{d}x)=\int\int\ind_B(x)\ind_{\{\xi\in E^C\}}\xi(\mathrm{d}x
        ) \dR(\mathrm{d}\xi)=0.
    \end{align*}
    This implies that also $\frac 1 a\leq\kappa(\xi,x)$ $C_\dR$-a.e.
\end{bem}

 The condition \ref{property:lowerbound} is rather restrictive as it excludes many classical examples such as pair-interaction models, see Example~\ref{exp:Papangelou}. It can be relaxed by the following approximation condition.
\begin{enumerate}[label=$(\mathbf{ALB})$, leftmargin=1.75cm,
  labelwidth=1.6cm,
  labelsep=0.15cm,
  align=left]
\item\label{property:ALB}
The Papangelou intensity $\kappa:\Gamma_{\Lambda_n}\times {\Lambda_n}\to[0,\infty)$ can be approximated from above, i.e.\  there exists a  sequence of Papangelou intensities $(\kappa_a)_{a\in\N}$ and a monotone decreasing function $g:\N\to[0,\infty)$ with $\lim_{a\to\infty}g(a)=0$ s.t.\  for every $a\in\N$ it holds that: 
\begin{itemize}
\item[$(i)$] $0\leq\kappa_a(\xi,x)-\kappa(\xi,x)\leq g(a)$;
\item[$(ii)$] $\kappa_a$ satisfies  \ref{property:cocycle};
\item[$(iii)$] $\kappa_a$ satisfies \ref{property:lowerbound};

\item[$(iv)$] if $\kappa$ satisfies \ref{property:locallystable}, \ref{property:continuity}, \ref{property:continuity2}, or~\ref{property:repulsive}, then $\kappa_a$ satisfies the same conditions.
\end{itemize}

\end{enumerate}

In the following, we consider several examples of Papangelou measures and discuss under which conditions they satisfy the assumptions introduced above.
\begin{bsp}[Poisson measures]
      A measure $\dR\in\PPP_1(\Gamma_X)$ is a Poisson measure with intensity measure $\lambda\in\MMM^+(X)$ if and only if it satisfies the Mecke equation  
\begin{align}\label{eq:poisson}
        \int\int F(\xi-\delta_x,x)\xi(\mathrm{d}x)\dR(\mathrm{d}\xi)=\int\int F(\xi,x)\lambda(\mathrm{d}x)\dR(\mathrm{d}\xi),\quad F\in\FFF^+(\Gamma_X\times X),
    \end{align}
 see \cite{Mecke1967}. Hence, a Poisson measure $\Poi^\lambda$ is a Papangelou measure w.r.t.\  the Papangelou intensity $\kappa\equiv1$ and reference measure $\lambda$. In particular, the Poisson measure satisfies \ref{property:cocycle}, \hyperref[property:locallystable]{$(\mathbf{S}(1))$}, \hyperref[property:lowerbound]{$(\mathbf{LB}(1))$},   \ref{property:continuity2}, \ref{property:repulsive},  \hyperref[property:dobrushin]{$(\mathbf{D}(1/2))$}, \hyperref[property:dobrushin2]{$(\mathbf{D}^\prime(1/2))$}, and~\ref{property:finiterange}.
 Moreover, in the Euclidean setting $X=\R^d$,  the homogeneous Poisson measure $\Poi^c$  satisfies \ref{property:shiftinvariance}.
 \end{bsp}
 \begin{bsp}[Gibbs measures]\label{exp:Papangelou}
      Let \(\dR\in\PPP_1(\Gamma)\) and \(H\in\FFF(\Gamma)\). Assume that the limit
\[
h(\xi,x)
:=
\lim_{n\to\infty}
(
H(\xi_{|\Lambda_n}+\delta_x)-H(\xi_{|\Lambda_n}))
\]
exists for \(\dR\otimes\Leb\)-a.e.\
\((\xi,x)\in\Gamma\times\mathbb{R}^d\). Then $\dR$
      is an infinite-volume Gibbs measure with activity $z>0$, inverse temperature $\beta>0$, and Hamiltonian $H$ if 
    \begin{align*}
       \int\int
F(\xi-\delta_x,x)\xi(\mathrm{d}x)\dR(\mathrm{d}\xi)
=
\int\int
F(\xi,x)
z e^{-\beta h(\xi,x)}
\dR(\mathrm{d}\xi)\Leb(\mathrm{d}x),\quad F\in\FFF^+(\Gamma\times\R^d),
    \end{align*}
      see e.g.\ \cite{NguyenZessin77}. Hence, an infinite-volume Gibbs measure is a Papangelou measure w.r.t.\  the Papangelou intensity $\kappa(\xi,x)=ze^{-\beta h(\xi,x) }$ and the reference measure $\Leb$.
    \begin{enumerate}
    \item[$(i)$] \textbf{Pair interactions.}
Let $\phi:\mathbb{R}^d\to\mathbb{R}\cup\{\infty\}$ be symmetric, i.e.\  $\phi(x)=\phi(-x)$, and define
\[
H(\xi_{|\Lambda_n})
:=\frac{1}{2} \int\int \ind_{\{x,y\in\Lambda_n, x\not=y\}}\phi(x-y)\xi(\mathrm{d}x)\xi(\mathrm{d}y)=\frac{1}{2}
\sum_{\substack{x,y\in\xi\cap\Lambda_n\\x\not=y}}\phi(x-y).
\]
 The corresponding Papangelou intensity is given by
\[
\kappa(\xi,x)
=
z\exp\left(
-\beta\lim_{n\to\infty}\big(H(\xi_{|\Lambda_n}+\delta_x)-H(\xi_{|\Lambda_n})\big)
\right)
=
z\exp\left(-\beta\sum_{y\in\xi}\phi(x-y)\right),
\]
whenever it is well-defined. A direct computation shows that $\kappa$ satisfies the cocycle property \ref{property:cocycle}.   
If $\phi\geq 0$, then $\kappa$ is repulsive \ref{property:repulsive}. In this case, \ref{property:locallystable} is always satisfied. If additionally \(z\leq 1\), this improves to
\hyperref[property:locallystable]{\((\mathbf{S}(1))\)}. \\
 If $\phi$ is lower semicontinuous, then $\kappa$ satisfies \ref{property:continuity}. If $\phi$ is additionally $\Leb$-a.e.\ continuous, then also \ref{property:continuity2} holds. Note that in contrast to \ref{property:continuity}, \ref{property:continuity2} mainly excludes pair interactions $\phi$ which are too irregular. For instance, in the $1$-dimensional setting $\phi(x)=1-\ind_C(x)$ with $C$ being the fat Cantor set is not $\Leb$-a.e.\ continuous but lower semicontinuous.\\ 
If $\phi\geq0$, we obtain for the Dobrushin condition \ref{property:dobrushin} the estimate
\begin{align*}
\zeta(\kappa)
&=
\frac12\bigg(1-\sup_{\substack{\xi\in\Gamma\\ x\in\mathbb{R}^d}}
z\int e^{-\beta\sum_{v\in\xi}\phi(v-y)}\big(1-e^{-\beta\phi(x-y)}\big)\Leb(\mathrm{d}y)\bigg)\\
&\geq
\frac12\bigg(1-\sup_{x\in\mathbb{R}^d}
z\int \big(1-e^{-\beta\phi(x-y)}\big)\Leb(\mathrm{d}y)\bigg)\\
&
=:\tilde{\zeta}(\kappa).
\end{align*}
Note that if $0<\tilde{\zeta}(\kappa)$, then \ref{property:dobrushin} is satisfied. In particular, in Theorems~\ref{thm:evi}, \ref{thm:geodconvexity}, 
\ref{thm:FIbounded}, \ref{thm:evistationary}, 
\ref{thm:geodconvexitystationary}, and~\ref{thm:FIstationary}, 
the assumption \ref{property:dobrushin2} can be replaced by the assumption that
\begin{align}\label{eq:dobrushinmodified}
    0<\lim_{M\to\infty}\tilde{\zeta}(\kappa^M)
=
\tilde{\zeta}(\kappa),\end{align}
at the cost of replacing $\zeta(\kappa)$ by $\tilde\zeta(\kappa)$. The assumption in \eqref{eq:dobrushinmodified} has the advantage of being easier to
verify. For instance, it holds for $z$ sufficiently small and sufficiently integrable pair potentials \(\phi\) that are \(\Leb\)-a.e.\
continuous.\\
The finite range condition \ref{property:finiterange} is satisfied if $\phi$ has bounded support. Moreover, $\kappa$ is shift-invariant \ref{property:shiftinvariance} by definition.\\
We conclude with a discussion of the approximation condition \ref{property:ALB}.
Examples of pair potentials $\phi$ s.t.\  $\kappa$ satisfies \ref{property:ALB} are potentials having a hard-core exclusion at the origin. More precisely, assume that there exists some $0<s<\infty$ s.t.\  $\phi(x)=\infty$ for all $x\in B_s(0)$. Under this assumption, there exists some constant $C_s$ s.t.\  $\kappa(\xi,x)=0$ if $\xi(\Lambda_n)>C_s\Leb(\Lambda_n)$. Define 
 \begin{align*}
     E:=\{\xi\in\Gamma_{\Lambda_n}:\xi(\Lambda_n)\leq C_s\Leb(\Lambda_n)\}
 \end{align*}
 and
 \begin{align*}
     \kappa_a(\xi,x):=\ind_{E}(\xi+\delta_x)ze^{-\beta\sum_{y\in\xi}(\phi(x-y)\wedge\frac{1}{\beta C_s\Leb(\Lambda_n)}\log(az))}, \quad a\in\N,\; \xi\in\Gamma_{\Lambda_n},\; x\in\Lambda_n.
\end{align*} 
Then $\kappa_a$ satisfies the conditions of \ref{property:ALB}.
\begin{itemize}
    \item[$(i)$] By definition, we directly obtain that
    $$0\leq \kappa_a(\xi,x)- \kappa(\xi,x)\ind_{E}(\xi+\delta_x)=\kappa_a(\xi,x)-\kappa(\xi,x)\leq z\left(\frac{1}{za}\right)^{\frac{1}{C_s\Leb(\Lambda_n)}}=:g(a), $$
    with $\lim_{a\to\infty}g(a)=0$.
    \item[$(ii)$] Using the pair interaction structure, a direct computation shows that $\kappa_a$ satisfies the cocycle property \ref{property:cocycle}.
    \item[$(iii)$] For $\xi+\delta_x\in E$ it holds that
 \begin{align*}
 \kappa_a(\xi,x)= ze^{-\beta\sum_{y\in\xi}(\phi(x-y)\wedge\frac{1}{\beta C_s\Leb(\Lambda_n)}\log(az))}\geq \frac{1}{a}.
\end{align*}
Hence, $\kappa_a$ satisfies \ref{property:lowerbound}.
\item[$(iv)$] If $\phi\geq0$, then both $\kappa$ and $\kappa_a$ are repulsive \ref{property:repulsive}. Further, if $\kappa$ satisfies \ref{property:locallystable}, \ref{property:continuity}, and~\ref{property:continuity2}, then $\kappa_a$ satisfies these conditions by definition.
\end{itemize}
A prominent example satisfying all assumptions discussed above is the hard-core model, defined by
\[
\phi(x)=
\begin{cases}
\infty, & \|x\|< r,\\
0, & \|x\|\geq r,
\end{cases}
\]
where $\|\cdot\|$ denotes the Euclidean norm.

\item[$(ii)$] \textbf{Area interactions.} We continue with the discussion of the area-interaction model. For $\xi \in \Gamma$ and $\tilde r>0$  define $L_{\tilde r}(\xi) := \bigcup_{x \in \xi} B_{\tilde r}(x)$. The energy functional with parameter $c>0$ is defined by
\[
H(\xi_{|\Lambda_n}) := - \log(c)\mathrm{Leb}(L_{\tilde r}(\xi_{|\Lambda_n})),\quad \xi\in\Gamma,
\]
and the associated Papangelou intensity is
\[
\kappa(\xi,x)
=
z c^{\beta\mathrm{Leb}(B_{\tilde r}(x)\setminus L_{\tilde r}(\xi))},\quad\xi\in\Gamma,\;x\in\R^d.
\]
By definition, $\kappa$ satisfies the cocycle property \ref{property:cocycle} and is shift-invariant \ref{property:shiftinvariance}. Further,
 both \ref{property:locallystable} and~\ref{property:lowerbound} are satisfied since
\[
z \min\{1,c^{\beta\mathrm{Leb}(B_{\tilde r}(0))}\}
\leq
\kappa(\xi,x)
\leq
z \max\{1,c^{\beta\mathrm{Leb}(B_{\tilde r}(0))}\}.
\]
Moreover,  $\kappa$ is jointly continuous \ref{property:continuity2} and~\ref{property:finiterange} holds with range $r=2\tilde r$. If $c\geq1$, a direct computation shows that $\kappa$ is also  repulsive \ref{property:repulsive}. Further, in this setting, the Dobrushin condition \ref{property:dobrushin} can be bounded by 
\begin{align*}
\zeta(\kappa)&=\frac12\bigg(1-\sup_{\substack{\xi\in\Gamma\\ x\in\mathbb{R}^d}}
\int_{\mathbb{R}^d}
|\kappa(\xi,y)-\kappa(\xi+\delta_x,y)|\mathrm{d}y\bigg)
\\
&\geq\frac12\bigg(1-
\sup_{\substack{\xi\in\Gamma\\ x\in\mathbb{R}^d}}
\int_{B_{2\tilde r}(x)}
z c^{\beta\mathrm{Leb}(B_{\tilde r}(0))}
\Big(1 - c^{-\beta\mathrm{Leb}(B_{\tilde r}(x)\cap B_{\tilde r}(y))}\Big)
\mathrm{Leb}(\mathrm{d}y)\bigg)\\
&=\frac12\bigg(1- \int_{B_{2{\tilde r}}(0)}
z c^{\beta\mathrm{Leb}(B_{\tilde r}(0))}
\Big(1 - c^{-\beta\mathrm{Leb}(B_{\tilde r}(0)\cap B_{\tilde r}(y))}\Big)
\mathrm{Leb}(\mathrm{d}y)\bigg)\\
&=:\tilde\zeta(\kappa).
\end{align*}
In particular, \ref{property:dobrushin} is satisfied if $z$ and ${\tilde r}$ are sufficiently small and $c\geq1$ is sufficiently close to $1$. Regarding \ref{property:dobrushin2}, we note that  $\lim_{M\to\infty}\tilde\zeta(\kappa^M)=\tilde\zeta(\kappa)$.
 \end{enumerate}
    
\end{bsp} 
\begin{bsp}[Determinantal measures]
    Let $\lambda\in\MMM^+(X)$ and let $K$ be a linear Hermitian  integral operator on $L^2(\lambda)$ with kernel $\tilde{K}$, i.e.\  
    \begin{align*}
        Kf(x)=\int \tilde K(x,y)f(y)\lambda(\mathrm{d}y),\quad f\in L^2(\lambda),\; x\in X.
    \end{align*}
    For $B\in\BBB(X)$, we define the restriction operator by $(\chi_Bf)(x):=\ind_B(x)f(x)$, $f\in L^2(\lambda)$, $x\in X$, and denote $K_B:=\chi_B K\chi_B$. In particular, $K_B$ is an integral operator with kernel $$\tilde{K}_B(x,y)=\ind_{B\times B}(x,y)\tilde{K}(x,y).$$
    In the following, we always assume that $K$ is of locally finite trace class, i.e.\  
    \[
    \operatorname{Tr}(K_B):=\sum_{i\in\mathbb N}\lambda_i^{(B)}<\infty,
\]
for all $B\in\mathcal B_0(X)$, where $(\lambda_i^{(B)})_{i\in\mathbb N}$ are the eigenvalues of $K_B$. Moreover, we assume that $0\leq K\leq1$, i.e.\ 
    \begin{align*}
        0\leq\int\int \bar{f}(x)\tilde{K}(x,y)f(y)\lambda(\mathrm{d}x)\lambda(\mathrm{d}y)\leq \|f\|_{L^2(\lambda)}^2.
    \end{align*}
    A measure $\dP\in\PPP_1(\Gamma_X)$ is called \textit{determinantal measure} w.r.t.\ $K$ if the $n$-point correlation function of $\dP$ is given by $h_\dP^{(n)}(x_1,\ldots,x_n)=\det(\tilde{K}(x_i,x_j)_{i,j=1}^n)$, i.e.\  
    \begin{align*}
        \int \sum_{\{x_1,\ldots,x_n\}\subseteq\xi}f_n(x_1,\ldots,x_n)\dP(\mathrm{d}\xi)=\frac{1}{n!}\int f_n(x_1,\ldots,x_n)\det(\tilde{K}(x_i,x_j)_{i,j=1}^n)\lambda(\mathrm{d}x_1)\ldots\lambda(\mathrm{d}x_n),
    \end{align*}
    for all symmetric functions $f_n:X^n\to[0,\infty)$. Since $0\leq K\leq1$ and $K$ is of locally finite trace class, the determinantal measure w.r.t.\ $K$ exists uniquely, see e.g.\ \cite{Macchi_1975}, \cite[Theorem 3]{Soshnikov2000Determinantal}. If additionally $\|K\|<1$ in operator norm,  \cite[Theorem 3.5]{GeorgiiYoo2005} establishes that the determinantal measure is a Papangelou measure. The associated Papangelou intensity $\kappa$ satisfies \ref{property:cocycle} by the construction of \cite{GeorgiiYoo2005}, \ref{property:locallystable} by \cite[Corollary 3.4]{GeorgiiYoo2005},   and~\ref{property:repulsive} by \cite[Theorem 3.6]{GeorgiiYoo2005}. A rigorous verification of the remaining assumptions introduced in this subsection requires further investigations.   
\end{bsp}
    \subsection{Continuity Equation}\label{subsec:ce}
In this subsection, we recall the definition of a non-local continuity equation for probability measures in $\PPP_1(\Gamma_X)$, introduced in \cite{schiavo2023wasserstein}, and summarize its main properties developed in \cite{schiavo2023wasserstein} and \cite{huesmann2025nonlocalwassersteingeometrygradient}. The role of the gradient will be taken over by the \textit{discrete difference operator}, also known as the \textit{add-one-cost operator},
$DF:\Gamma_X\times X\to\R$ of  $F\in\FFF(\Gamma_X)$, which  is defined by 
\begin{equation*}
        D_xF(\xi):=F(\xi+\delta_x)-F(\xi).
    \end{equation*}

 \begin{dfn}[Non-local continuity equation]\label{dfn:ce}
    Let $\dP_0,\dP_T\in\PPP_1(\Gamma_X)$. The set of all continuous solutions to the \textit{non-local continuity equation} on $X$ connecting $\dP_0$ and $\dP_T$ in time $T$ is denoted by  $\CCC\EEE_T^X(\dP_0,\dP_T)$ and contains  all pairs $({\bar{\dP}_t},{\bar{\dV}_t})_{t\in[0,T]}\subset \PPP_1(\Gamma_X)\times \MMM_{b,0}(\Gamma_X\times X)$ satisfying the following conditions:
  \begin{enumerate}
      \item[${(i)}$] $\bar{\dP}_0=\dP_0$, $\bar{\dP}_T=\dP_T$;
      \item[${(ii)}$]${\bar{\dP}}:[0,T]\to\PPP_1(\Gamma_X)$ is $\PPP_1(\Gamma_X)$-continuous;
      \item[${(iii)}$] for all $f\in\CCC_c^+(X)$ and compactly supported smooth functions $\phi\in \CCC_c^\infty((0,T))$, it holds that 
     \begin{equation}\label{eq:continuityequation}
        \int_0^T\left(\dt \phi(t)\right)\int e^{-\int f(x)\xi(\mathrm{d}x)}\bar{\dP}_t(\mathrm{d}\xi)\mathrm{d}t+ \int_0^T \phi(t)\int D_x e^{-\int f(y)\xi(\mathrm{d}y)}\bar{\dV}_t(\mathrm{d}\xi, \mathrm{d}x)\mathrm{d}t=0;
    \end{equation}
      \item [${(iv)}$] $\int_0^T|\bar{\dV}_t|(\Gamma_X\times B)\mathrm{d}t<\infty$ for all $B\in\BBB_0(X)$.
  \end{enumerate}
  Further, $\CCC\EEE^X_T$ denotes the set of all pairs $(\bar\dP,\bar\dV)$ satisfying properties $(ii)$--$(iv)$.
\end{dfn}
 We will usually drop the brackets around $\dt \phi(t)$ when working with the continuity equation and if $X=\R^d$, we briefly write $\CCC\EEE_T(\dP_0,\dP_T)$ instead of $\CCC\EEE_T^{\R^d}(\dP_0,\dP_T)$. By \cite[Theorem 3.2]{huesmann2025nonlocalwassersteingeometrygradient} the set of solutions $\CCC\EEE_T^X(\dP_0,\dP_T)$ is always non-empty.
\begin{bem}\label{bem:topologyMb0}
The velocity fields $(\Bar{\dV}_t)_{t\in[0,T]}$ can also be considered  as the measure $$\Bar\dV(\mathrm{d}\xi,\mathrm{d}x,\mathrm{d}t)=\int \Bar\dV_t(\mathrm{d}\xi,\mathrm{d}x)\mathrm{d}t\in\MMM(\Gamma_X\times X\times [0,T]).$$
    In analogy with $\mathcal{M}_{b,0}(\Gamma_X \times X)$, we define $\mathcal{M}_{b,0}(\Gamma_X\times (X\times [0,T]))=\mathcal{M}_{b,0}(\Gamma_X\times X\times [0,T])$  as the space of all signed measures
    $\Bar\dV \in \mathcal{M}(\Gamma_X \times X \times [0,T])$  s.t.\   $|\Bar\dV|(\Gamma_X\times B)<\infty$ for all 
    $B\in\BBB_0(X\times[0,T])$. Further, we endow $\MMM_{b,0}(\Gamma_X\times X\times[0,T])$ with the topology induced by the seminorms $\Bar\dV\mapsto|\int G(\xi,x,t)\Bar\dV(\mathrm{d}\xi,\mathrm{d}x,\mathrm{d}t)|$, $G\in\mathcal{C}_{b,0}(\Gamma_X\times (X\times[0,T]))$.  Property $(iv)$ of Definition \ref{dfn:ce} directly implies that $\Bar\dV\in\MMM_{b,0}(\Gamma_X\times X\times [0,T])$. 
\end{bem}
  We introduce the set of bounded \textit{local functions}
\begin{align}\label{eq:H}
    \mathcal{H} := \{ F \in \FFF_b(\Gamma_X):  DF \in \FFF_{b,0}(\Gamma_X \times X)\}.
\end{align} 
\begin{satz}\label{thm:propertiesce}
     Let  $ \dP_0,\dP_T\in\PPP_1(\Gamma_X)$.
     \begin{enumerate}
        \item  [$(i)$]  \emph{(\cite[Lemma 4.5]{schiavo2023wasserstein})} Let  $\tau:[0,\tilde T]\to [0,{T}]$ be an increasing and absolutely continuous function s.t.\   its inverse is also absolutely continuous. Then $({\bar{\dP}},{\bar{\dV}})\in\CCC\EEE^X_{{T}}(\dP_0,\dP_{{T}})$ if and only if $({\bar{\dP}_{\tau(t)}},\dt\tau(t)\cdot\bar{\dV}_{\tau(t)})_{t\in[0,\tilde T]}\in\CCC\EEE^X_{\tilde T}(\dP_0,\dP_T)$. 
         
         \item [$(ii)$] \emph{(\cite[Proposition 4.8]{schiavo2023wasserstein})} Let $({\bar{\dP}},{\bar{\dV}})\in\CCC\EEE^X_T(\dP_0,\dP_T)$. For any $\phi\in \CCC_c^\infty((0,T))$ and $F\in\mathcal{H}$, it holds that 
    \begin{align*}
        \int_0^T\dt \phi(t)\int F(\xi)\bar{\dP}_t(\mathrm{d}\xi)\mathrm{d}t+ \int_0^T \phi(t)\int D_x F(\xi)\bar{\dV}_t(\mathrm{d}\xi, \mathrm{d}x)\mathrm{d}t=0.
    \end{align*}
    \item[$(iii)$] \emph{(\cite[Lemma 4.2]{schiavo2023wasserstein})} Let $({\bar{\dP}}^i,{\bar{\dV}}^i)_{i\in\N}\subset\CCC\EEE^X_T(\dP_0,\dP_T),\;(\bar{\dP}_t)_{t\in[0,T]}\subset\PPP_1(\Gamma_X)$ a $\PPP_1(\Gamma_X)$-continuous curve, and $\bar{\dV}\in\MMM_{b,0}(\Gamma_X\times X\times[0,T])$. If $\lim_{i\to\infty}{\bar{\dP}}_{t}^i={\bar{\dP}}_t$  in $\PPP_1(\Gamma_X)$ for a.e.\ $t\in[0,T]$ and $\lim_{i\to\infty}{\bar{\dV}}^i={\bar{\dV}}$ in $\MMM_{b,0}(\Gamma_X\times X\times [0,T])$, then $({\bar{\dP}},{\bar{\dV}})\in\CCC\EEE^X_T(\dP_0,\dP_T)$.
    \item[$(iv)$] \emph{(\cite[Theorem 4.10, Theorem 4.12]{schiavo2023wasserstein})}
    Let $(\bar\dP,\bar\dV)\in\CCC\EEE^X_T(\dP_0,\dP_T)$. For a.e.\ $t\in[0,T]$ and $F\in\mathcal{H}$, it holds that
    \begin{align*}
        \int F(\xi)\bar{\dP}_t(\mathrm{d}\xi)=\int F(\xi)\bar{\dP}_0(\mathrm{d}\xi)+\int_0^t\int D_xF(\xi)\bar{\dV}_s(\mathrm{d}\xi,\mathrm{d}x)\mathrm{d}s.
    \end{align*}
    In particular, this implies
    \begin{align*}
        \I_{\bar{\dP}_t}(B)=\I_{\bar{\dP}_0}(B)+\int_0^t \bar{\dV}_s(\Gamma_X\times B)\mathrm{d}s,
    \end{align*}
    for all $B\in\BBB_0(X)$ and a.e.\ $t\in[0,T]$.
    \item[$(v)$] \emph{(\cite[Theorem 4.10]{schiavo2023wasserstein})} Let $(\bar{\dP}_t)_{t\in[0,T]}\subset\PPP_1(\Gamma_X)$ and $\bar{\dV}\in\MMM_{b,0}(\Gamma_X\times X\times[0,T])$. If $(\bar\dP,\bar\dV)$ solves the continuity equation~\eqref{eq:continuityequation}, then  there exists a unique continuous representative $(\bar{\dP}^\prime,\bar{\dV})$, i.e.\  $\bar{\dP}_t=\bar{\dP}_t^\prime$ for a.e.\ $t\in[0,T]$ and $(\bar{\dP}_t^\prime)_{t\in[0,T]}$ is $\PPP_1(\Gamma_X)$-continuous.
     \end{enumerate}
     
\end{satz}
\begin{bem}\label{bem:ce}
 Theorem~\ref{thm:propertiesce} $(i)$ shows that solutions to the continuity equation are invariant under time reparameterization. This justifies restricting to solutions to the continuity equation in time 1. To simplify the notation, we briefly write $\CCC\EEE^X(\dP_0,\dP_1)$ instead of $\CCC\EEE^{X}_1(\dP_0,\dP_1)$.
\end{bem}
The following theorem summarizes basic transformations of solutions to the continuity equation on $\R^d$ such as superposition or restriction. Further transformations such as discretization, tiling, and stationarization are discussed throughout the paper, see e.g.\ Lemma~\ref{lem:discretizedCE} and Lemma~\ref{lem:stationarysolutions}. 
\begin{satz}\label{thm:transformationofce}
    Let $B_1,B_2\in\BBB(\R^d)$ be  $G_\delta$ sets, $z\in\R^d$, $\dP_0^1,\dP_1^1\in\PPP_1(\Gamma_{B_1})$, $\dP_0^2,\dP_1^2\in\PPP_1(\Gamma_{B_2})$ and $(\Bar{\dP}^1,\Bar{\dV}^1)\in\CCC\EEE^{B_1}(\dP_0^1,\dP_1^1)$, $(\Bar{\dP}^2,\Bar{\dV}^2)\in\CCC\EEE^{B_2}(\dP_0^2,\dP_1^2)$.
    \begin{itemize}
        \item[$(i)$] Define the shift operator on the product space $\Theta_z^{\mathrm{prod}}:\Gamma\times\R^d\to\Gamma\times\R^d$ by $$\Theta^{\mathrm{prod}}_z(\xi,x):=(\Theta_z\xi,x-z),\quad \xi\in\Gamma,\; x,z\in\R^d.$$
    Then $((\Theta_z)_\#\bar{\dP},(\Theta_z^\mathrm{prod})_\#\bar{\dV})
    \in
    \CCC\EEE^{B_1-z}((\Theta_z)_\#\dP_0,(\Theta_z)_\#\dP_1)$.
         \item[$(ii)$]\emph{(\cite[Lemma 3.4]{huesmann2025nonlocalwassersteingeometrygradient})} Define the independent superposition of $(\bar{\dP}_t^1,\bar{\dV}_t^1)$ and $(\bar{\dP}_t^2,\bar{\dV}_t^2)$ by 
\begin{align}\label{eq:superpositionvelocity}
    (\bar{\dP}_t^1, \bar{\dV}_t^1)\oplus^\mathrm{prod}(\bar{\dP}_t^2, \bar{\dV}_t^2):= (\mathrm{sp}^{\mathrm{prod}})_\#(\bar{\dV}_t^1\otimes \bar{\dP}_t^2)+ (\mathrm{sp}^{\mathrm{prod}})_\#(\bar{\dV}^2_t\otimes \bar{\dP}_t^1),
\end{align}
where  $\mathrm{sp}^{\mathrm{prod}}(\xi^1,x,\xi^2):=(\xi^1+\xi^2,x)$.
         Then, $(\Bar\dP^1,\Bar\dV^1)\oplus^\mathrm{prod}(\Bar\dP^2,\Bar\dV^2)\in\CCC\EEE^{B_1\cup B_2}(\dP_0^{1}\oplus\dP_0^2,\dP_1^1\oplus\dP_1^2)$.
          \item[$(iii)$] \emph{(\cite[Lemma 3.5]{huesmann2025nonlocalwassersteingeometrygradient})} Assume that $B_1\subseteq B_2$.  We define the restriction of $\bar\dV_t^2$  to $\Gamma_{B_1}\times B_1$ by $$\bar{\dV}_{t|B_1}^2:=\ind_{B_1}(x){(\mathrm{pr}_{B_1},\mathrm{id})}_\#\bar{\dV}_t^2.$$ Then $(\Bar{\dP}^2_{|B_1},\Bar{\dV}^2_{|B_1})\in\CCC\EEE^{B_1}(\dP^2_{0|B_1},\dP^2_{1|B_1})$.
    \end{itemize}
    \begin{proof}
        The proof of $(i)$ follows directly by the definition of the respective objects.
    \end{proof}
\end{satz}
\subsubsection{Stationary Solutions to the Continuity Equation}
 For stationary measures $\dP_0,\dP_1 \in \PPP_s(\Gamma)$, we are particularly interested in solutions to the  non-local continuity equation $(\bar{\dP},\bar{\dV}) \in \CCC\EEE(\dP_0,\dP_1)$ (see Definition~\ref{dfn:ce}) for which the curve $\bar{\dP}$ is a curve of stationary measures. To this end, we introduce the set of stationary measures on the product space $\Gamma\times\R^d$ by
\begin{align*}
    \MMM_s(\Gamma\times \R^d):=\{\dV\in\MMM_{b,0}(\Gamma\times \R^d):(\Theta_z^{\mathrm{prod}})_\#\dV=\dV \text{ for all } z\in\R^d\}.
\end{align*}
Similar to $\MMM_{b,0}(\Gamma\times \R^d)$, we endow $\MMM_s(\Gamma\times\R^d)$ with the $\CCC_{b,0}(\Gamma\times\R^d)$-weak topology induced by the maps  $\bar\dV\mapsto|\int F(\xi,x)\bar\dV(\mathrm{d}\xi,\mathrm{d}x)|$, $F\in\CCC_{b,0}(\Gamma\times \R^d)$.
\begin{dfn}[Stationary solution]\label{dfn:stationarysolutions}
    Let $\dP_0,\dP_1\in\PPP_s(\Gamma)$. The set of \textit{stationary solutions} $\CCC\EEE^s(\dP_0,\dP_1)$ is the set of all solutions $(\bar{\dP},\bar{\dV})\in\CCC\EEE(\dP_0,\dP_1)$ s.t.\  $\bar{\dV}_t\in\MMM_s(\Gamma\times\R^d)$ for a.e.\ $t\in[0,1]$. As before, we denote by $\CCC\EEE^s$ the set of all stationary solutions independent of the starting and the final measures.
\end{dfn}
\cite[Corollary 3.9]{huesmann2025nonlocalwassersteingeometrygradient} ensures the existence of stationary solutions, i.e.\ if $\dP_0,\dP_1\in\PPP_s(\Gamma)$ then $\CCC\EEE^s(\dP_0,\dP_1)$ is always  non-empty.
\subsection{Entropy Functionals and Fisher Information}\label{sec:entropy}
In this section, we recall the definitions of the relative and specific relative entropy as well as the modified and specific modified Fisher Information and discuss their most important properties.
\begin{dfn}[Relative entropy]\label{dfn:fentropy}
Let $\dP,\dR\in\PPP_1(\Gamma_X)$. The \textit{relative entropy} is defined by
\begin{align*}
\Ent(\dP|\dR):=\begin{cases}
\int \rho(\xi)\log(\rho(\xi))\dR(\mathrm{d}\xi), & \text{if } \dP\ll\dR \text{ with } \dP=\rho\dR,\\
\infty, & \text{otherwise. }
\end{cases}
\end{align*}
The proper domain of the relative entropy is denoted by $$\mathcal{D}(\Ent(\cdot|\dR)):=\{\dP\in\PPP_1(\Gamma_{X}):\Ent(\dP|\dR)<\infty\}.$$
\end{dfn}
\begin{lem}\label{lem:propertiesent}
The  relative entropy satisfies the following properties:
\begin{itemize}
\item [$(i)$]  \emph{(\cite[Lemma 9.4.3]{Ambrosio2008GradientFI})}  $\Ent(\cdot|\cdot)$ is jointly lower semicontinuous on $\PPP_1(\Gamma_X)\times\PPP_1(\Gamma_X)$.
\item [$(ii)$] \emph{(\cite[Lemma 9.4.5]{Ambrosio2008GradientFI})} Let $Y$ be another Polish metric space and  $T:\Gamma_X\to\Gamma_Y$ measurable. It holds that $$\Ent (T_\#\dP|T_\#\dR)\leq \Ent(\dP|\dR).$$
\end{itemize}
\end{lem}

The relative entropy is typically infinite for stationary measures as the absolute continuity assumption on the full configuration space is rather restrictive. For instance, the only stationary measure that is absolutely continuous with respect to a homogeneous Poisson measure $\Poi^{c}$ is $\Poi^{c}$ itself. Hence, for these measures one considers the specific relative entropy.

\begin{dfn}[Specific relative entropy]
Let $\dP,\dR\in\PPP_s(\Gamma)$. The \textit{specific relative entropy} is defined by
\begin{align}\label{eq:specificrelativeentropy}
\mathcal{E}(\dP|\dR):=\liminf_{n\to\infty}\frac{1}{\Leb(\Lambda_n)}\Ent (\dP_{|\Lambda_n}|\dR_{|\Lambda_n}).
\end{align}
The proper domain of the specific relative entropy is denoted by $$\mathcal{D}(\mathcal{E}(\cdot|\dR)):=\{\dP\in\PPP_s(\Gamma):\mathcal{E}(\dP|\dR)<\infty\}.$$
\end{dfn}
\begin{lem}\label{lem:entropy}
Assume that $\kappa:\Gamma\times\R^d\to[0,\infty)$ satisfies \ref{property:cocycle}, \ref{property:locallystable}, \ref{property:finiterange}, and~\ref{property:shiftinvariance}. The specific relative entropy w.r.t.\  $\dR\in\RRR_{\kappa,s}(\Gamma)$ satisfies the following properties:
\begin{enumerate}
\item[${(i)}$] The limit inferior in \eqref{eq:specificrelativeentropy} is a limit, i.e.\ it holds that
\begin{align*}
\mathcal{E}(\dP|\dR)=\lim_{n\to\infty}\frac{1}{\Leb(\Lambda_n)}\Ent (\dP_{|\Lambda_n}|\dR_{|\Lambda_n}),\quad \dP\in\PPP_s(\Gamma).
\end{align*}
\item[${(ii)}$] The specific relative entropy $\mathcal{E}(\cdot|\dR)$ is lower semicontinuous on $\PPP_s(\Gamma)$.
\end{enumerate}
\end{lem}
\begin{proof}
We first prove $(i)$. Since $\kappa$ satisfies \ref{property:locallystable},   Lemma~\ref{lem:finitevolumePapangelou} implies
\[
    Z_{\kappa,\Leb_{|B}}^B(\eta)
    =
    \sum_{j\geq0}\frac{1}{j!}
    \int_B\ldots\int_B
    \varphi_\kappa^\eta\left(\sum_{l=1}^j\delta_{x_l}\right)
    \Leb(\mathrm{d}x_1)\ldots\Leb(\mathrm{d}x_j)
    \leq
    \sum_{j\geq0}\frac{1}{j!}A^j\Leb(B)^j
    =
    e^{A\Leb(B)},
\]
for all $B\in\BBB_0(\R^d)$, $\eta\in\Gamma_{B^C}$.
Hence,
\begin{align}\label{eq:voidprobabilityentropy}
    \dR^B(\{0\}|\eta)
    =
    \frac{1}{Z_{\kappa,\Leb_{|B}}^B(\eta)}
    \geq
    e^{-A\Leb(B)}.
\end{align}
Recall that $\kappa$ satisfies  \ref{property:finiterange} and let $0\leq r<\infty$ be the range of $\kappa$. For $n>2r$, we define
\[
    \Lambda_n^{r-}:=[-n/2+r,n/2-r)^d.
\]
Further, let $m\in\N$, $I=\{-m,\ldots,m\}^d$, and define
\[
    B_{m,n}:=
    \bigcup_{i\in I}(\Lambda_n^{r-}-ni)
    \quad \text{ and }\quad
    C_{m,n}:=
    \Lambda_{(2m+1)n}
    \setminus
    B_{m,n}.
\]
For every $K\in\BBB(\Gamma_{B_{m,n}})$, the DLR equation (Theorem~\ref{thm:DLR}) and
\eqref{eq:voidprobabilityentropy} yield
\begin{align*}
\dR(\xi(C_{m,n})=0)
\left(
    \bigoplus_{i\in I}
    \dR^{\Lambda_n^{r-}-ni}
\right)(K)&=
\int
\ind_K\left(\xi_{|B_{m,n}}\right)
\ind_{\{\xi(C_{m,n})=0\}}(\xi)
\dR(\mathrm{d}\xi)
\\
&=
\int
\ind_K\left(\xi_{|B_{m,n}}\right)
\dR^{C_{m,n}}
\left(
    \{0\}|
    \xi_{|C_{m,n}^C}
\right)
\dR(\mathrm{d}\xi)
\\
&\geq
e^{-A\Leb(C_{m,n})}
\dR_{|B_{m,n}}(K).
\end{align*}
Since $\dR(\xi(C_{m,n})=0)\leq 1$, it follows that
\begin{align}
\dR_{|B_{m,n}}(K)
\leq
e^{A\Leb(C_{m,n})}
\left(
    \bigoplus_{i\in I}
    \dR^{\Lambda_n^{r-}-ni}
\right)(K).
\label{eq:firstcomparisonentropy}
\end{align}

Moreover, for every $i\in I$ and every
$K_i\in\BBB(\Gamma_{\Lambda_n^{r-}-ni})$, the DLR equation,
the finite-range property \ref{property:finiterange}, and \eqref{eq:voidprobabilityentropy} give
\begin{align*}
\dR_{|\Lambda_n^{r-}-ni}(K_i)
&\geq
\int
\ind_{K_i}\left(
    \xi_{|\Lambda_n^{r-}-ni}
\right)
\ind_{\left\{
    \xi\left(
        (\Lambda_n-ni)\setminus
        (\Lambda_n^{r-}-ni)
    \right)=0
\right\}}(\xi)
\dR(\mathrm{d}\xi)
\\
&=
\dR\left(
    \xi\left(
        (\Lambda_n-ni)\setminus
        (\Lambda_n^{r-}-ni)
    \right)=0
\right)
\dR^{\Lambda_n^{r-}-ni}(K_i)
\\
&\geq
e^{-A\Leb(\Lambda_n\setminus\Lambda_n^{r-})}
\dR^{\Lambda_n^{r-}-ni}(K_i).
\end{align*}
Consequently,
\begin{align*}
\dR^{\Lambda_n^{r-}-ni}(K_i)
\leq
e^{A\Leb(\Lambda_n\setminus\Lambda_n^{r-})}
\dR_{|\Lambda_n^{r-}-ni}(K_i).
\end{align*}
Taking the independent superposition over all $i\in I$ and using
\(
\Leb(C_{m,n})
=
(2m+1)^d
\Leb(\Lambda_n\setminus\Lambda_n^{r-}),
\)
we obtain
\begin{align}
\left(
    \bigoplus_{i\in I}
    \dR^{\Lambda_n^{r-}-ni}
\right)(K)
\leq
e^{A\Leb(C_{m,n})}
\left(
    \bigoplus_{i\in I}
    \dR_{|\Lambda_n^{r-}-ni}
\right)(K)
\label{eq:secondcomparisonentropy}
\end{align}
for every $K\in\BBB(\Gamma_{B_{m,n}})$.
Combining \eqref{eq:firstcomparisonentropy} and
\eqref{eq:secondcomparisonentropy}, we conclude that
\begin{align*}
\dR_{|B_{m,n}}(K)
\leq
e^{2A\Leb(C_{m,n})}
\left(
    \bigoplus_{i\in I}
    \dR_{|\Lambda_n^{r-}-ni}
\right)(K).
\end{align*}
Hence,
\begin{align}\label{eq:productcomparisonentropy}
\frac{
    \mathrm{d}\dR_{|B_{m,n}}
}{
    \mathrm{d}
    \left(
        \bigoplus_{i\in I}
        \dR_{|\Lambda_n^{r-}-ni}
    \right)
}
\leq
e^{2A\Leb(C_{m,n})}.
\end{align}
By the definition of the relative entropy and \eqref{eq:productcomparisonentropy}, we obtain
\begin{align*}
    \Ent \left(
    \dP_{|B_{m,n}}
    |
    \dR_{|B_{m,n}}
    \right)&=
    \Ent \left(
    \dP_{|B_{m,n}}
    |
    \bigoplus_{i\in I}\dR_{|\Lambda_n^{r-}-ni}
    \right)-
    \int
    \log\left(
    \frac{
   \mathrm d\dR_{|B_{m,n}} }{\mathrm d\bigoplus_{i\in I}\dR_{|\Lambda_n^{r-}-ni}
    }(\xi)
    \right)
    \dP_{|B_{m,n}}(\mathrm{d}\xi)\\
    &\geq
    \Ent \left(
    \dP_{|B_{m,n}}
    |
    \bigoplus_{i\in I}\dR_{|\Lambda_n^{r-}-ni}
    \right)
    -
    2A\Leb(C_{m,n}).
\end{align*}
By the super-additivity of the relative entropy for product reference measures, see \cite[Proposition~15.10]{georgii2011gibbs}, and since the sets $\Lambda_n^{r-}-ni$, $i\in I$, are pairwise disjoint, we obtain
\begin{align}\label{eq:quasisuperadditivityentropy}
    \Ent \left(
    \dP_{|B_{m,n}}
    |
    \dR_{|B_{m,n}}
    \right)&\geq
    \sum_{i\in I}
    \Ent (\dP_{|\Lambda_n^{r-}-ni}|\dR_{|\Lambda_n^{r-}-ni})
    -
    2A\Leb(C_{m,n})\nonumber\\
    &=
    (2m+1)^d
    \Ent (\dP_{|\Lambda_n^{r-}}|\dR_{|\Lambda_n^{r-}})
    -
    2A(2m+1)^d
    \Leb(\Lambda_n\setminus\Lambda_n^{r-}).
\end{align}
Since
\(
    B_{m,n}\subset\Lambda_{(2m+1)n},
\)
the contraction property of the relative entropy (Lemma~\ref{lem:propertiesent}~$(ii)$) gives
\[
    \Ent (\dP_{|\Lambda_{(2m+1)n}}|\dR_{|\Lambda_{(2m+1)n}})
    \geq
    \Ent \left(
    \dP_{|B_{m,n}}
    |
    \dR_{|B_{m,n}}
    \right).
\]
Together with \eqref{eq:quasisuperadditivityentropy}, this implies
\begin{align*}
    \frac{1}{\Leb(\Lambda_{(2m+1)n})}
    \Ent (\dP_{|\Lambda_{(2m+1)n}}|\dR_{|\Lambda_{(2m+1)n}})
    \geq
    \frac{\Leb(\Lambda_n^{r-})}{\Leb(\Lambda_n)}
    \frac{1}{\Leb(\Lambda_n^{r-})}
    \Ent (\dP_{|\Lambda_n^{r-}}|\dR_{|\Lambda_n^{r-}})
    -
    2A
    \frac{\Leb(\Lambda_n\setminus\Lambda_n^{r-})}{\Leb(\Lambda_n)}.
\end{align*}

The same argument with all boxes $\Lambda_n^{r-}-ni$ s.t.\  $\Lambda_n-ni\subset\Lambda_m$ gives, for arbitrarily large boxes,
\begin{align*}
    \liminf_{m\to\infty}
    \frac{1}{\Leb(\Lambda_m)}
    \Ent (\dP_{|\Lambda_m}|\dR_{|\Lambda_m})
    &\geq
    \frac{\Leb(\Lambda_n^{r-})}{\Leb(\Lambda_n)}
    \frac{1}{\Leb(\Lambda_n^{r-})}
    \Ent (\dP_{|\Lambda_n^{r-}}|\dR_{|\Lambda_n^{r-}})-
    2A
    \frac{\Leb(\Lambda_n\setminus\Lambda_n^{r-})}{\Leb(\Lambda_n)}.
\end{align*}
Define
\[
a_n:=
\frac{1}{\Leb(\Lambda_n)}
\Ent (\dP_{|\Lambda_n}|\dR_{|\Lambda_n}).
\]
For $n$ large, let $k_n:=\lfloor n-2r\rfloor$. Since
$\Lambda_{k_n}\subset \Lambda_n^{r-}$, the contraction property of the
relative entropy (Lemma~\ref{lem:propertiesent}~$(ii)$) implies
\(
    \Ent (\dP_{|\Lambda_n^{r-}}|\dR_{|\Lambda_n^{r-}})
    \geq
    \Ent (\dP_{|\Lambda_{k_n}}|\dR_{|\Lambda_{k_n}}).
\)
Hence, the previous estimate yields
\begin{align}\label{eq:specificentropylowerbound}
\liminf_{m\to\infty} a_m
\geq
\frac{\Leb(\Lambda_{k_n})}{\Leb(\Lambda_n)}a_{k_n}
-
2A
\frac{\Leb(\Lambda_n\setminus\Lambda_n^{r-})}{\Leb(\Lambda_n)}.
\end{align}
Letting $n\to\infty$, we obtain
\(
\liminf_{m\to\infty} a_m
\geq
\limsup_{n\to\infty} a_n,
\)
proving that the limit in the definition of the specific relative entropy \eqref{eq:specificrelativeentropy} exists.

We proceed with the proof of $(ii)$. Let $(\dP_i)_{i\in\N}\subset\PPP_s(\Gamma)$ be s.t.\  $\lim_{i\to\infty}\dP_i=\dP$ in $\PPP_s(\Gamma)$. For $\tilde{\dP}\in\PPP_s(\Gamma)$, define
\[
a_n(\tilde{\dP}):=
\frac{1}{\Leb(\Lambda_n)}
\Ent (\tilde{\dP}_{|\Lambda_n}|\dR_{|\Lambda_n}).
\]
 We obtain by \eqref{eq:specificentropylowerbound}, 
\[
\mathcal{E}(\tilde{\dP}|\dR)
\geq
\frac{\Leb(\Lambda_{k_n})}{\Leb(\Lambda_n)}a_{k_n}(\tilde{\dP})
-
2A
\frac{\Leb(\Lambda_n\setminus\Lambda_n^{r-})}{\Leb(\Lambda_n)}.
\]
for all sufficiently large $n$.
Applying this estimate to $\dP_i$ and taking the liminf, we get
\[
\liminf_{i\to\infty}\mathcal{E}(\dP_i|\dR)
\geq
\frac{\Leb(\Lambda_{k_n})}{\Leb(\Lambda_n)}
\liminf_{i\to\infty}a_{k_n}(\dP_i)
-
2A
\frac{\Leb(\Lambda_n\setminus\Lambda_n^{r-})}{\Leb(\Lambda_n)}.
\]
By the lower semicontinuity of $\Ent(\cdot|\dR_{|\Lambda_n})$
(Lemma~\ref{lem:propertiesent}~$(i)$), we obtain
\(
\liminf_{i\to\infty}a_{k_n}(\dP_i)
\geq
a_{k_n}(\dP).
\)
Consequently,
\[
\liminf_{i\to\infty}\mathcal{E}(\dP_i|\dR)
\geq
\frac{\Leb(\Lambda_{k_n})}{\Leb(\Lambda_n)}
a_{k_n}(\dP)
-
2A
\frac{\Leb(\Lambda_n\setminus\Lambda_n^{r-})}{\Leb(\Lambda_n)}.
\]
Letting $n\to\infty$ and using $(i)$,
we conclude that
\(
\liminf_{i\to\infty}\mathcal{E}(\dP_i|\dR)
\geq
\mathcal{E}(\dP|\dR),
\)
finishing the proof.
\end{proof}
We finish this section with the definition of a non-local counterpart of the Fisher information.
\begin{dfn}[Modified Fisher information]
        Let $\lambda\in\MMM^+(X)$ and let $\kappa:\Gamma_X\times X\to[0,\infty)$ satisfy \ref{property:cocycle}, \ref{property:locallystable}, and~\ref{property:finiterange}. The \textit{modified Fisher information} w.r.t.\ a Papangelou measure $\dR\in\RRR_{\kappa,\lambda}(\Gamma_X)$ is defined by 
             \begin{align*}
             \mathrm{I}(\dP|\dR):=\begin{cases}
        \int\int D_x\rho(\xi) D_x\log(\rho)(\xi)\kappa(\xi,x)\dR(\mathrm{d}\xi)\lambda(\mathrm{d}x), \quad &\text{if } \dP\ll\dR \text{ with } \dP=\rho \dR,\\
        \infty, & \text{otherwise,}
    \end{cases}
    \end{align*}
    with the convention that $(0-0)(\log(0)-\log(0))=0$ and $(u-0)(\log(u)-\log(0))=+\infty$ for $u>0$.
    \end{dfn}
    Similar to the specific relative entropy, we also introduce the \textit{specific modified Fisher information}.
\begin{dfn}[Specific modified Fisher information]
Assume that $\kappa:\Gamma\times \R^d\to[0,\infty)$ satisfies  \ref{property:cocycle}, \ref{property:locallystable},  \ref{property:finiterange}, and~\ref{property:shiftinvariance}. The \textit{specific modified Fisher information} of $\dP\in\PPP_s(\Gamma)$ w.r.t.\ a Papangelou measure $\dR\in\RRR_{\kappa,s}(\Gamma)$ is defined by 
    \begin{align*}
    \mathcal{I}(\dP|\dR)
    :=
    \lim_{n\to\infty}
    \frac{1}{\Leb(\Lambda_n)}
    \I(\dP_{|\Lambda_n}|\dR_{|\Lambda_n}).
\end{align*}
\end{dfn}
\begin{bem}
    The limit in the definition of $\mathcal{I}$ exists by \cite[Lemma 9.1]{jahnel2026reversiblebirthanddeathdynamicscontinuum}.
\end{bem}
 \section{Transport Geometry for Non-Stationary Measures } \label{sec:gradientflows}
     The main goal of this section is to identify the natural spatial birth-and-death process associated with a Papangelou measure $\dR=\RRR_{\kappa,\lambda}(\Gamma_{\Lambda_n})$  as the gradient flow of the relative entropy $\Ent(\cdot | \dR)$ in a geometry induced by a suitable transport distance  tailored to $\dR$. We prove the gradient flow result by establishing an Evolution Variational Inequality (EVI) which implies several central consequences, including geodesic convexity of $\Ent(\cdot | \dR)$, exponential convergence of the gradient flow, and functional inequalities such as an HWI inequality and a Log-Sobolev inequality.

    \subsection{General Framework and Definitions}\label{subsec:generalsetup}
    In this subsection, we briefly introduce the general framework needed to state the gradient flow result. We refer to  Sections~\ref{sec:transportdistance} and~\ref{subsec:semigroup} for a more detailed discussion. 
    
    Let $\lambda\in\MMM^+(\Lambda_n)$, and assume that  $\kappa:\Gamma_{\Lambda_n}\times \Lambda_n\to[0,\infty) $  satisfies \ref{property:cocycle} and \ref{property:locallystable}. The spatial birth-and-death process with birth rate $\kappa(\xi,x)\lambda(\mathrm{d}x)$, death rate $1$, and initial configuration 
$\xi \in \Gamma_{\Lambda_n}$ is the $\Gamma_{\Lambda_n}$-valued Markov process 
$\bigl(Y_t^{\kappa,\lambda,(\xi)}\bigr)_{t \ge 0}$ whose generator is formally given by
\begin{align*}
\mathfrak{L}_{\kappa,\lambda}F(\xi)
:=
\int_{\Lambda_n} \kappa(\xi,x)\bigl(F(\xi+\delta_x)-F(\xi)\bigr)\,\lambda(\mathrm{d}x)
+ \sum_{x\in\xi}\bigl(F(\xi-\delta_x)-F(\xi)\bigr),
\end{align*}
for a class of test functions $F : \Gamma_{\Lambda_n} \to \mathbb{R}$ specified in Theorem~\ref{thm:SDESemigroup}. The Markov semigroup associated with the process $(Y_t^{\kappa,\lambda,(\xi)})_{t\geq0}$ is defined by
\[
\mathrm{S}_t^{\kappa,\lambda} F(\xi) := \mathbb{E}\bigl[F(Y_t^{\kappa,\lambda,(\xi)})\bigr], 
\quad\ t \geq 0,
\]
for all functions $F\in\FFF(\Gamma_{\Lambda_n})$ s.t.\   the right hand side is well-defined. The associated dual semigroup $(\mathrm{S}_t^{\kappa,\lambda*})_{t\geq0}$ on $\PPP_1(\Gamma_{\Lambda_n})$ is defined by
\begin{align*}
\mathrm{S}_t^{\kappa,\lambda*}\dP(K) 
:= \int \mathrm{S}_t^{\kappa,\lambda} \ind_K(\xi)\, \dP(\mathrm{d}\xi), 
\quad t \geq 0,\ K \in \BBB(\Gamma_{\Lambda_n}).
\end{align*}
   By Theorem~\ref{thm:dualsemigroup}, $(\mathrm{S}_t^{\kappa,\lambda*})_{t\geq0}$ is a $C_0$-semigroup on $\PPP_1(\Gamma_{\Lambda_n})$, and the associated Papangelou measure $\dR=\RRR_{\kappa,\lambda}(\Gamma_{\Lambda_n})$ is invariant, i.e., $\mathrm{S}_t^{\kappa,\lambda*}\dR=\dR$ for all $t\geq0,$ see \eqref{eq:Rinvariant}.

    Further, we follow the approach of \cite{maas2011gradientflowsentropyfinite}, \cite{erbar2012gradientflowsentropyjump}, and \cite{schiavo2023wasserstein}, and define a non-local transport distance tailored to a Papangelou measure $\dR=\RRR_{\kappa,\lambda}(\Gamma_{\Lambda_n})$ in the spirit of the dynamical Benamou-Brenier formulation of the Wasserstein distance $\mathcal{W}_{\kappa,\lambda}:\PPP_1(\Gamma_{\Lambda_n})\times \PPP_1(\Gamma_{\Lambda_n})\to [0,\infty]$ by
\begin{align*}
\mathcal{W}_{\kappa,\lambda}^2(\dP_0, \dP_1) := \inf \left\{ \mathcal{A}_{\kappa,\lambda}(\Bar{\dP}, \Bar{\dV}) : (\Bar{\dP}, \Bar{\dV}) \in \CCC\EEE^{\Lambda_n}(\dP_0, \dP_1) \right\},
\end{align*}
     where $\CCC\EEE^{\Lambda_n}(\dP_0,\dP_1)$ denotes the set of solutions to the non-local continuity equation introduced in Subsection~\ref{subsec:ce}. The action functional $\mathcal{A}_{\kappa,\lambda}$ is defined as follows. Let $\theta:[0,\infty)\times[0,\infty)\to[0,\infty)$ be the logarithmic mean 
     \begin{align}\label{eq:logmean}
           \theta(x,y)=\int_0^1x^ry^{1-r}\mathrm{d}r=\frac{y-x}{\log(y)-\log(x)}\ind_{\{x,y>0, x\not=y\}}+x\ind_{\{x=y\}}
       \end{align}
    (see Section \ref{sec:transportdistance} for a more general choice of $\theta$). We define the \textit{action density} $\alpha: [0,\infty)\times [0,\infty)\times\R\to [0,\infty]$  by
       \begin{align*}
           \alpha(x,y,w):=\frac{|w|^2}{\theta(x,y)},
       \end{align*}
      where  $0/0=0$ by convention.
      The \textit{Lagrange functional}  w.r.t.\ $\kappa$ and $\lambda$  is for $\dP\in\PPP_1(\Gamma_{\Lambda_n})$ and $\dV\in \MMM_{b,0}(\Gamma_{\Lambda_n}\times {\Lambda_n})$ defined by
    \begin{align}\label{eq:dfnlagrange}
           \mathcal{L}_{\kappa,\lambda}({\dP},{\dV}):=\int\alpha\left(\frac{\mathrm{d}\kappa(\dP\otimes \lambda)}{\mathrm{d}\sigma}(\xi,x),\frac{\mathrm{d}C_{\dP}}{\mathrm{d}\sigma}(\xi,x),\frac{\mathrm{d}{\dV}}{\mathrm{d}\sigma}(\xi,x)\right)\sigma(\mathrm{d}\xi,\mathrm{d}x),
       \end{align}
       where $\sigma\in\MMM_{b,0}(\Gamma_{\Lambda_n}\times {\Lambda_n})$ is non-negative s.t.\   $\kappa(\dP\otimes {\lambda}),$ $ C_{\dP},$ and $\dV$ are absolutely continuous w.r.t.\ $\sigma,$ i.e., $\kappa(\dP\otimes {\lambda}),$ $C_{\dP},$ ${\dV}\ll\sigma.$ Since the action density $\alpha$ is a positive homogeneous function, the definition of the Lagrange functional is independent of the choice of $\sigma$, see Section \ref{sec:transportdistance} for more details. The \textit{action functional}  w.r.t.\  $\kappa$ and $\lambda$  is defined by
     \begin{align*}
         \mathcal{A}_{\kappa,\lambda}({\Bar{\dP}},{\Bar{\dV}}):=\int_0^1\mathcal{L}_{\kappa,\lambda}(\Bar{\dP}_t,\Bar{\dV}_t)\mathrm{d}t,\quad (\Bar{\dP},\Bar{\dV})\in\CCC\EEE^{\Lambda_n}(\dP_0,\dP_1).
     \end{align*}
     \subsection{Gradient Flows, Convexity, and Functional Inequalities}
     The main result of this section is the following EVI identifying $(\mathrm{S}_t^{\kappa,\lambda*})_{t\geq0}$ as the gradient flow of $\Ent(\cdot|\dR)$ in the geometry induced by $\mathcal{W}_{\kappa,\lambda}.$
     \begin{satz}\label{thm:evi}
Let $\lambda=\Leb_{|\Lambda_n}$ and assume that $\kappa:\Gamma_{\Lambda_n}\times\Lambda_n\to[0,\infty)$ satisfies \ref{property:cocycle}, \hyperref[property:locallystable]{$(\mathbf{S}(1))$}, \ref{property:ALB}, \ref{property:continuity2}, \ref{property:repulsive}, and~\ref{property:dobrushin2}.  Further, let $\dR = \RRR_{\kappa,\lambda}(\Gamma_{\Lambda_n})$. 
Then $(\mathrm{S}_t^{\kappa,\lambda*})_{t\geq0}$ satisfies the EVI$(\zeta(\kappa))$ for $\Ent(\cdot|\dR)$, i.e.\ for $\dP_0,\dP_1\in\PPP_1(\Gamma_{\Lambda_n})$ s.t.\   $\Ent(\dP_1|\dR)<\infty$ and $\mathcal{W}_{\kappa,\lambda}(\dP_0,\dP_1)<\infty$ it holds  
 \begin{align*}
            \frac{1}{2}\frac{\mathrm{d}}{\mathrm{d}t}
\mathcal{W}_{\kappa,\lambda}^2(\mathrm{S}_t^{\kappa,\lambda*}\dP_0,\dP_1)
+
\frac{\zeta(\kappa)}{2}
\mathcal{W}_{\kappa,\lambda}^2(\mathrm{S}_t^{\kappa,\lambda*}\dP_0,\dP_1)
\leq
\Ent(\dP_1 | \dR)
-
\Ent(\mathrm{S}_t^{\kappa,\lambda*}\dP_0 | \dR),
\quad
t>0,
        \end{align*}
        or equivalently in integrated form
    \begin{align*}
        \frac{e^{\zeta(\kappa)t}}{2}\mathcal{W}_{\kappa,\lambda}^2(\mathrm{S}_t^{\kappa,\lambda*}\dP_0,\dP_1)-\frac{1}{2}\mathcal{W}_{\kappa,\lambda}^2(\dP_0,\dP_1)\leq E_{\zeta(\kappa)}(t)(\Ent(\dP_1 | {\dR})
-
\Ent(\mathrm{S}_t^{\kappa,\lambda*}\dP_0 | {\dR})),\quad t>0,
    \end{align*}
    with $E_{\zeta(\kappa)}(t)=\int_0^te^{\zeta(\kappa)s}\mathrm{d}s$.
    \end{satz}
    Here, and in the sequel, the time derivative is intended as the right upper derivative. The equivalence of the differential and integrated formulation of the EVI is discussed e.g.\ in  \cite[Proposition 3.1]{Daneri_2008}.
    
    Before proving Theorem~\ref{thm:evi}, we first discuss some of its consequences. The following theorem 
    may be regarded as the counterpart of McCann's displacement convexity in our setting.

\begin{satz}\label{thm:geodconvexity}
Let $\lambda=\Leb_{|\Lambda_n}$ and assume that $\kappa:\Gamma_{\Lambda_n}\times\Lambda_n\to[0,\infty)$ satisfies \ref{property:cocycle}, \hyperref[property:locallystable]{$(\mathbf{S}(1))$}, \ref{property:ALB}, \ref{property:continuity2}, \ref{property:repulsive}, and~\ref{property:dobrushin2}. Further, let $\dR=\RRR_{\kappa,\lambda}(\Gamma_{\Lambda_n})$. The relative entropy $\Ent(\cdot|\dR)$ is $\zeta(\kappa)$-geodesically convex w.r.t.\ $\mathcal{W}_{\kappa,\lambda}$. More precisely, if $\dP_0,\dP_1\in\PPP_1(\Gamma_{\Lambda_n})$ are such that  $\mathcal{W}_{\kappa,\lambda}(\dP_0,\dP_1)<\infty$, then for every $\mathcal{W}_{\kappa,\lambda}$-geodesic $\bar{\dP}$ connecting $\dP_0$ and $\dP_1$, it holds that 
   \begin{align*}
        \Ent(\bar{\dP}_t|\dR)\leq (1-t) \Ent({\dP}_0|\dR)+t\Ent({\dP}_1|\dR)-\frac{\zeta(\kappa)t(1-t)}{2}\mathcal{W}_{\kappa,\lambda}^2(\dP_0,\dP_1), \quad 0\leq t\leq 1.
    \end{align*}
    \end{satz}
     \begin{proof}
    The proof follows directly from \cite[Proposition 3.2]{Daneri_2008} and Theorem~\ref{thm:evi}.
\end{proof}
Further consequences of the EVI of Theorem~\ref{thm:evi} are several functional inequalities such as a non-local HWI inequality, a modified Log-Sobolev inequality, and a non-local Talagrand inequality. For a discussion of the relevance of the functional inequalities and possible applications, we refer to Subsection~\ref{introduction}.

\begin{satz}\label{thm:FIbounded}
 Let $\lambda=\Leb_{|\Lambda_n}$ and assume that $\kappa:\Gamma_{\Lambda_n}\times\Lambda_n\to[0,\infty)$ satisfies \ref{property:cocycle}, \hyperref[property:locallystable]{$(\mathbf{S}(1))$}, \ref{property:ALB}, \ref{property:continuity2}, \ref{property:repulsive}, and~\ref{property:dobrushin2}. For $\dR= \RRR_{\kappa,\lambda}(\Gamma_{\Lambda_n})$ and  $\dP_0,\dP_1\in\PPP_1(\Gamma_{\Lambda_n})$ s.t.\   $\Ent(\dP_0|\dR),$ $\Ent(\dP_1|\dR)<\infty$,  the following functional inequalities are satisfied:
    \begin{itemize}
        \item[$(i)$] The dual semigroup $(\mathrm{S}_t^{\kappa,\lambda*})_{t\geq0}$ is $\zeta(\kappa)$-contractive w.r.t.\ $\mathcal{W}_{\kappa,\lambda}$: 
    \begin{align*}
        \mathcal{W}_{\kappa,\lambda}(\mathrm{S}_t^{\kappa,\lambda*}\dP_0,\mathrm{S}_t^{\kappa,\lambda*}\dP_1)\leq e^{-\zeta(\kappa)t} \mathcal{W}_{\kappa,\lambda}(\dP_0,\dP_1),\quad t\geq0.
    \end{align*}
    Further, $(\mathrm{S}_t^{\kappa*})_{t\geq0}$ is $2\zeta(\kappa)$-exponentially contractive  w.r.t.\ $\mathrm{Ent}(\cdot|\dR)$:
\begin{align*}
    \mathrm{Ent}
    (
        \mathrm{S}_t^{\kappa,\lambda*}\dP_0|\dR
    )
    &\leq
    e^{-2\zeta(\kappa)t}
    \mathrm{Ent}(\dP_0|\dR),\quad t\geq0.
\end{align*}
    \end{itemize}
    \begin{itemize}
        \item[$(ii)$] The non-local HWI inequality with parameter $\zeta(\kappa)$ is satisfied: \begin{align*}
           \Ent(\dP_0|\dR)\leq \mathcal{W}_{\kappa,\lambda}(\dP_0,\dR)\sqrt{\mathrm{I}(\dP_0|\dR)}-\frac{\zeta(\kappa)}{2}\mathcal{W}_{\kappa,\lambda}^2(\dP_0,\dR).
\end{align*}
         \item[$(iii)$] The modified Log-Sobolev inequality with parameter $\zeta(\kappa)$ is satisfied:  \begin{align*}
           \Ent(\dP_0|\dR)\leq \frac{1}{2\zeta(\kappa)}\mathrm{I}(\dP_0|\dR).
\end{align*}  
         \item[$(iv)$] The non-local Talagrand inequality with parameter $\zeta(\kappa)$ is satisfied:  \begin{align*}
          \mathcal{W}_{\kappa,\lambda}(\dP_0,\dR)\leq \sqrt{\frac{2}{\zeta(\kappa)}\Ent(\dP_0|\dR)}.
\end{align*}

    \end{itemize}
\end{satz}
\begin{proof}
    $(i)$ The proof of the contractivity w.r.t. $\mathcal{W}_{\kappa,\lambda}$ follows directly from \cite[Proposition 3.1]{Daneri_2008} and Theorem~\ref{thm:evi}. The proof of the contractivity w.r.t.\ $\mathrm{Ent}(\cdot|\dR)$ follows from \cite[Theorem 3.5]{MuSa20} and Theorem~\ref{thm:evi}.\\
    $(ii)$ The proof is identical to the proof of  \cite[Theorem 7.3]{Erbar_2012}   by replacing \cite[Proposition 4.2 (2)]{Erbar_2012} by Lemma~\ref{lem:derivativeFisher}.\\
    $(iii)$ The proof is identical to the proof of \cite[Theorem 7.4]{Erbar_2012} and follows directly from $(ii)$.\\
    $(iv)$ For technical reasons, we postpone the proof until after the proof of Theorem \ref{thm:geodconvexity}
\end{proof}

    \subsection{Proof of Theorem~\ref{thm:evi}}
    \label{subsec:proofevi}
    The proof of Theorem~\ref{thm:evi} proceeds in two steps. First, we focus on discrete base spaces $X=\Lambda_n^M =(\frac{1}{M}\Z^d-(\frac{1}{2M}-\frac{n}{2})\mathbf{1})\cap[-n/2,n/2)^d,$ where $\mathbf{1}=(1,\ldots,1)$. In this setting, the configuration space $\Gamma_{\Lambda_n^M}$ only contains finitely many simple point configurations, allowing us to establish the EVI by proving a discrete Bochner inequality (Proposition~\ref{prop:Bochner}). In a second step, we let the mesh size $\tfrac{1}{M}$ tend to zero and lift the EVI from the discrete to the continuous setting.
    \subsubsection{Discrete Base Spaces}\label{sec:curvaturediscretized}
    Let $M,n\in\N$, $\lambda=\Leb^M_{|\Lambda_n^M}$, and assume that $\kappa:\Gamma_{\Lambda_n^M}\times\Lambda_n^M\to[0,\infty)$ satisfies \ref{property:cocycle} and~\ref{property:locallystable}. In the discrete setting, we always assume that $\kappa(\xi^M,x^M)=0$ if $\xi^M(\{x^M\})=1$, ensuring that the associated Papangelou measure $\dR=\RRR_{\kappa,\lambda}(\Gamma_{\Lambda_n^M})$ is simple. 
    
    The main result of this subsection is the  following EVI for discrete base spaces.
    \begin{satz}\label{thm:evidiscretized}
    Let $\lambda=\Leb^M_{|\Lambda_n^M}$ and assume that  $\kappa:\Gamma_{\Lambda_n^M}\times\Lambda_n^M\to[0,\infty)$ satisfies \ref{property:cocycle}, \hyperref[property:locallystable]{$(\mathbf{S}(1))$}, \ref{property:repulsive}, and~\ref{property:dobrushin}. Further, let $\dR=\RRR_{\kappa,\lambda}(\Gamma_{\Lambda_n^M})$. 
    Then $(\mathrm{S}_t^{\kappa,\lambda*})_{t\geq0}$ satisfies the EVI$(\zeta(\kappa))$ for $\Ent(\cdot|\dR)$, i.e.\ if $\dP_0,\dP_1\in\PPP_1(\Gamma_{\Lambda_n^M})$ are s.t.\   $\Ent(\dP_0|\dR),\Ent(\dP_1|\dR)<\infty$, then
    \begin{align*}
            \frac{1}{2}\frac{\mathrm{d}}{\mathrm{d}t}
\mathcal{W}_{\kappa,\lambda}^2(\mathrm{S}_t^{\kappa,\lambda*}\dP_0,\dP_1)
+
\frac{\zeta({\kappa})}{2}
\mathcal{W}_{\kappa,\lambda}^2(\mathrm{S}_t^{\kappa,\lambda*}\dP_0,\dP_1)
\leq
\Ent(\dP_1 | {\dR})
-
\Ent(\mathrm{S}_t^{\kappa,\lambda*}\dP_0 | {\dR}),
\quad
t>0,
        \end{align*}
or equivalently in integrated form
    \begin{align*}
        \frac{e^{\zeta(\kappa)t}}{2}\mathcal{W}_{\kappa,\lambda}^2(\mathrm{S}_t^{\kappa,\lambda*}\dP_0,\dP_1)-\frac{1}{2}\mathcal{W}_{\kappa,\lambda}^2(\dP_0,\dP_1)\leq E_{\zeta(\kappa)}(t)(\Ent(\dP_1 | {\dR})
-
\Ent(\mathrm{S}_t^{\kappa,\lambda*}\dP_0 | {\dR})),\quad t>0,
    \end{align*}
    with $E_{\zeta(\kappa)}(t)=\int_0^te^{\zeta(\kappa)s}\mathrm{d}s$.
    \end{satz}
   We note that the dual semigroup $(\mathrm{S}^{\kappa,\lambda*}_t)_{t\geq0}$ stems from a continuous-time Markov chain  on the finite set $\Gamma_{\Lambda_n^M}$. The distance $\mathcal W_{\kappa,\lambda}$ coincides with the transport distance constructed by Maas \cite{maas2011gradientflowsentropyfinite}. In the framework developed by the first author and Maas \cite{Erbar_2012}, the EVI in Theorem~\ref{thm:evidiscretized} is one of several equivalent formulations of the property that the Markov chain has curvature bounded below by $\zeta(\kappa)$.

In the following, we will specify the underlying Markov chain and introduce further notation in order to state an equivalent formulation of the curvature bound in terms of a discrete Bochner inequality. Then we will use techniques developed in \cite{CDPP09, DPP13, FM16, EHMT17} in order to establish the discrete Bochner inequality and thus also the EVI.\\

The generator $\mathfrak L_{\kappa,\lambda}$ of $(\mathrm{S}_t^{\kappa,\lambda*})_{t\geq0}$ introduced in Section~\ref{sec:spatialbirthanddeath}  can be expressed as follows. For $x^M\in \Lambda_n^M$  we introduce the annihilation and creation maps $u^\pm_{x^M}:\Gamma_{\Lambda_n^M}\to \Gamma_{\Lambda_n^M}$ as
\[
u^+_{x^M}(\xi^M) := \xi^M +\delta_{x^M}\ind_{\{x^M\not\in\xi^M\}},\qquad u^-_{{x^M}}(\xi^M) := \xi^M-\delta_{x^M}\ind_{\{x^M\in\xi^M\}}.
\]
We set $U=\{u^\pm_{{x^M}} : {x^M}\in \Lambda_n^M\}$ and define the collection of transition rates $c:\Gamma_{\Lambda_n^M}\times U\to [0,\infty)$ by
\[
c(\xi^M,u^+_{{x^M}}) := \kappa(\xi^M,x^M)\lambda(\{x^M\}),\qquad
c(\xi^M,u^-_{{x^M}}) := \begin{cases}
    1, &  \text{ if }x^M\in\xi^M,\\
    0, & \text{ if }x^M\not\in\xi^M.
\end{cases}
\]
Note that for every $u \in U$ there is a unique $u^{-1}\in U$ s.t.\   $u^{-1}(u( \xi^M))=\xi^M$ for all $\xi^M\in\Gamma_{\Lambda_n^M}$ with $c(\xi^M,u)>0$. Namely, $(u_{x^M}^\pm)^{-1}=u_{x^M}^\mp$ and $u\mapsto u^{-1}$ is an involution. Moreover, as a consequence of the Papangelou property, the transition rates $c$ are reversible w.r.t.~$\dR$, i.e.~for any function $\Phi:\Gamma_{\Lambda_n^M}\times U\to \R$ we have
\begin{equation}\label{eq:reversible}
\sum_{\substack{\xi^M\in \Gamma_{\Lambda_n^M}\\u\in U}} \Phi(\xi^M, u)c(\xi^M,u)\dR(\{\xi^M\})=\sum_{\substack{\xi^M\in \Gamma_{\Lambda_n^M}\\u\in U}}\Phi(u(\xi^M), u^{-1})c(\xi^M,u)\dR(\{\xi^M\}).
\end{equation}
For a function $f:\Gamma_{\Lambda_n^M}\to \R$, we define the extension of the discrete difference operator $Df:\Gamma_{\Lambda_n^M}\times U\to \R$ via $D_u f(\xi^M) = f(u(\xi^M))-f(\xi^M)$ allowing also the subtraction of points of a point configuration. We will also write $D^\pm_{{x^M}}$ instead of $D_{u^\pm_{{x^M}}}$. Note that $D^-_{{x^M}}f(\xi^M)=0$ if $x^M\not\in\xi^M$ and $D^+_{{x^M}}f(\xi^M)= 0$ if $x^M\in\xi^M$.

Now, the generator writes as
\[
\mathfrak{L}_{\kappa,\lambda} f (\xi^M) = \sum_{{x^M}\in \xi^M}D^-_{{x^M}}f(\xi^M) + \sum_{{x^M}\in \Lambda_n^M}D^+_{{x^M}}f(\xi^M)\kappa(\xi^M,{x^M})\lambda(\{x^M\}) = \sum_{u\in U}D_u f(\xi^M)c(\xi^M,u).
\]

To formulate the Bochner inequality, we introduce the following quantities.

For functions $\Phi,\Psi:\Gamma_{\Lambda_n^M}\times U\to \R$, we introduce the inner product
\begin{align*}
    \langle \Phi,\Psi\rangle:= \frac12\sum_{\substack{\xi^M
    \in \Gamma_{\Lambda_n^M}\\u\in U}} \Phi(\xi^M,u)\Psi(\xi^M,u)c(\xi^M,u)\dR(\{\xi^M\}).
\end{align*}
For $\dP=\rho\dR\in\PPP_1(\Gamma_{\Lambda_n^M})$, we define $\theta(\rho):\Gamma_{\Lambda_n^M}\times U\to [0,\infty)$ and  $\hat{\mathfrak L}_{\kappa,\lambda} \rho: \Gamma_{\Lambda_n^M}\times U\to \R$ by
\begin{align*}
\theta(\rho)(\xi^M,u) &:= \theta\big(\rho(\xi^M),\rho(u(\xi^M))\big),\\
\hat{\mathfrak L}_{\kappa,\lambda}\rho(\xi^M,u)   &:=\partial_1\theta\big(\rho(\xi^M),\rho(u(\xi^M))\big)\mathfrak L_{\kappa,\lambda}\rho(\xi^M)+\partial_2\theta\big(\rho(\xi^M),\rho(u(\xi^M))\big)\mathfrak L_{\kappa,\lambda}\rho(u(\xi^M)).
\end{align*}
Now, for a probability density $\rho$ and a function $\psi:\Gamma_{\Lambda_n^M}\to \R$, we set
\begin{align*}
\mathcal A(\rho,\psi) &:= \langle\theta(\rho) D\psi, D \psi\rangle,\qquad
\mathcal B(\rho,\psi) :=\frac12\langle \hat{\mathfrak L}_{\kappa,\lambda}\rho D\psi,  D \psi\rangle -\langle\theta(\rho)D\mathfrak{L}_{\kappa,\lambda}\psi,D \psi\rangle. 
\end{align*}

By \cite[Theorem 4.5]{Erbar_2012}, the statement of Theorem~\ref{thm:evidiscretized} is equivalent to the following Bochner inequality.

\begin{prop}\label{prop:Bochner}
  For any probability density $\rho:\Gamma_{\Lambda_n^M}\to [0,\infty)$ w.r.t.~$\dR$ and any function $\psi:\Gamma_{\Lambda_n^M}\to \R$, we have
  \begin{equation}\label{eq:Bochner}
  \mathcal B(\rho,\psi) \ge \zeta(\kappa)\mathcal A(\rho,\psi).
  \end{equation}
\end{prop}

In the rest of this section we provide the proof of the Bochner inequality. To this end we will use techniques developed in \cite{FM16, EHMT17} and start by introducing some further notation.

Given $\rho:\Gamma_{\Lambda_n^M}\to[0,\infty)$, $\psi:\Gamma_{\Lambda_n^M}\to\R$, we set 
\begin{align*}
\mathcal B(\rho,\psi)[\xi^M;u,v] 
:= \frac12 |D_u\psi(\xi^M)|^2\partial_1\theta\big(\rho(\xi^M),\rho(u(\xi^M))\big)D_v\rho(\xi^M) \!+\! D_u\psi(\xi^M)D_v\psi(\xi^M)\theta\big(\rho(\xi^M),\rho(u(\xi^M))\big),
\end{align*}
for $\xi^M\in \Gamma_{\Lambda_n^M},\; u,v\in U$, so that
\[
\mathcal B(\rho,\psi) = \sum_{\xi^M\in \Gamma_{\Lambda_n^M}}\sum_{u,v\in U} \mathcal B(\rho,\psi)[\xi^M;u,v]c(\xi^M,u)c(\xi^M,v)\dR(\{\xi^M\}).
\]

The strategy will be to split the term $\mathcal B$ into a non-negative contribution and a positive part that will be estimated against $\mathcal A$. This will be achieved by suitably splitting the weight $c(\xi^M,u)c(\xi^M,v)$.

We call a function $E:\Gamma_{\Lambda_n^M}\times U\times U\to [0,\infty)$ \emph{admissible}, if for all $\xi^M\in \Gamma_{\Lambda_n^M}$ and $u,v\in U$ the following hold:
\begin{itemize}
    \item[$(i)$] $E(\xi^M,u,v) = E(\xi^M,v,u)$;
    \item[$(ii)$] for all $h:\Gamma_{\Lambda_n^M}\times U\times U\to \R$, we have
    \[
    \sum_{\xi^M\in \Gamma_{\Lambda_n^M}} \sum_{u,v\in U}h(\xi^M,u,v)E(\xi^M,u,v) = 
    \sum_{\xi^M\in \Gamma_{\Lambda_n^M}} \sum_{u,v\in U}h(u(\xi^M),u^{-1},v)E(\xi^M,u,v);
    \]
    \item[$(iii)$] $u (v( \xi^M)) = v(u(\xi^M))$ provided $E(\xi^M,u,v)>0$.
\end{itemize}

The following result provides the first key estimate on the $\mathcal B$-term, see \cite[Theorem 3.5]{FM16} or \cite[Proposition~3.4]{EHMT17}. Note that $E(\xi^M,u,v)= \tilde E(\xi^M,u,v)\dR(\{\xi^M\})$, where $\tilde E$ corresponds to the function $R$ used in the notation of \cite{FM16}.

\begin{lem}\label{lem:splitting}
Assume that there are $E,F:\Gamma_{\Lambda_n^M}\times U\times U\to[0,\infty)$ with $E$ admissible, s.t.\   $$c(\xi^M,u)c(\xi^M,v)\dR(\{\xi^M\})=E(\xi^M,u,v)+F(\xi^M,u,v).$$ Then, for any $\rho$ and $\psi$, we have
\begin{equation*}
 \sum_{\xi^M\in \Gamma_{\Lambda_n^M}}\sum_{u,v\in U}\mathcal B(\rho,\psi)[\xi^M;u,v]E(\xi^M,u,v)\ge 0,
\end{equation*}
and hence
\begin{equation*}
 \mathcal B(\rho,\psi)\geq\sum_{\xi^M\in \Gamma_{\Lambda_n^M}}\sum_{u,v\in U}\mathcal B(\rho,\psi)[\xi^M;u,v]F(\xi^M,u,v).
\end{equation*}
\end{lem}

To compare the remaining part of $\mathcal B$ involving the weight $F$ against $\mathcal A$, it is important to control the negative part of the off-diagonal terms in $\mathcal B$ by the diagonal terms.

\begin{lem}[{\cite[Lemma 3.7]{EHMT17}}]\label{lem:off-diag}
For all $\rho,\psi$ and $\xi^M\in \Gamma_{\Lambda_n^M},\;u,v\in U$, we have
\[
\mathcal B(\rho,\psi)[\xi^M;u,v] + \mathcal B(\rho,\psi)[\xi^M;v,u]
\geq
-\mathcal B(\rho,\psi)[u(\xi^M);u^{-1},u^{-1}]
-\mathcal B(\rho,\psi)[v(\xi^M);v^{-1},v^{-1}].
\]
\end{lem}

For the present model, we choose the following function $E:\Gamma_{\Lambda_n^M}\times U\times U\to [0,\infty)$ given by
\begin{align}\label{eq:R}
E(\xi^M,u^+_{x^M},u^+_{y^M})&=c(u^+_{x^M}(\xi^M), u^+_{y^M})c(\xi^M, u^+_{x^M})\dR(\{\xi^M\})\nonumber\\
    &= \kappa(\xi^M+\delta_{x^M},{y^M})\kappa(\xi^M,{x^M})\lambda(\{y^M\})\lambda(\{x^M\})\dR(\{\xi^M\}),\nonumber\\
E(\xi^M,u^-_{x^M},u^-_{y^M}) &= c(u^-_{x^M}(\xi^M),u^-_{y^M})c(\xi^M,u^-_{x^M})\dR(\{\xi^M\})\nonumber\\
&= \begin{cases}
    \dR(\{\xi^M\}), & \text{ if } x^M,y^M\in\xi^M, x^M\neq y^M,\\
    0, & \text{otherwise,}
\end{cases}\nonumber\\
E(\xi^M, u^-_{x^M},u^+_{y^M}) &= E(\xi^M, u^+_{y^M},u^-_{x^M}) \nonumber\\&= \begin{cases}
   c(\xi^M,u^-_{x^M})c(\xi^M,u^+_{y^M})\dR(\{\xi^M\}), & \text{ if } x^M\neq y^M,\\
   0, & \text{otherwise}.
\end{cases} 
\end{align}
One checks that this function is indeed admissible. In this case, we have
\begin{align}\label{eq:S}
F(\xi^M,u^+_{x^M},u^+_{y^M})&= \big(\kappa(\xi^M,{y^M})-\kappa(\xi^M+\delta_{x^M},{y^M})\big)\kappa(\xi^M,{x^M})\lambda(\{y^M\})\lambda(\{x^M\})\dR(\{\xi^M\}),\nonumber\\
F(\xi^M,u^-_{x^M},u^-_{y^M}) &= \begin{cases}
   \dR(\{\xi^M\}), & \text{ if } {x^M}\in \xi^M, {x^M}= {y^M},\nonumber\\
    0 ,& \text{otherwise},
\end{cases}\\
F(\xi^M, u^-_{x^M},u^+_{y^M}) &= F(\xi^M, u^+_{y^M},u^-_{x^M}) = 0.
\end{align}

We can now prove the Bochner inequality.

\begin{proof}[Proof of Proposition~\ref{prop:Bochner}]
First note that by the reversibility \eqref{eq:reversible} we have
\begin{align}\nonumber
    \mathcal A(\rho,\psi) 
    &=
    \sum_{\substack{\xi^M\in \Gamma_{\Lambda_n^M}\\u\in U}} \frac12 |D_u \psi(\xi^M)|^2 \theta(\rho)(\xi^M,u)c(\xi^M,u)\dR(\{\xi^M\}) \\\label{eq:A-}
    &= 
    \sum_{\substack{\xi^M\in \Gamma_{\Lambda_n^M}\\ {x^M}\in \xi^M}} |D^-_{x^M}\psi(\xi^M)|^2\theta(\rho)(\xi^M, u^-_{x^M})\dR(\{\xi^M\}).
\end{align}

Now, let $E$ be the admissible function in \eqref{eq:R}. Using Lemma~\ref{lem:splitting}, Lemma~\ref{lem:off-diag}, and repeatedly the symmetry $F(\xi^M,u,v)=F(\xi^M,v,u)$,
we can estimate (dropping mostly the dependence on $\rho,\psi$ from the notation)
\begin{align*}
    \mathcal B(\rho,\psi)
    &\ge
    \sum_{\substack{\xi^M\in \Gamma_{\Lambda_n^M}\\ u,v\in U}} \mathcal B(\rho,\psi)[\xi^M;u,v] F(\xi^M,u,v)\\
    &=
    \sum_{\substack{\xi^M\in \Gamma_{\Lambda_n^M}\\ u,v\in U}} \Big[\boldsymbol{1}_{\{u=v\}}\mathcal B[\xi^M;u,v] + \boldsymbol{1}_{\{u\neq v\}}\mathcal B[\xi^M;u,v]\Big]F(\xi^M,u,v)\\
    &=
    \sum_{\substack{\xi^M\in \Gamma_{\Lambda_n^M}\\ u,v\in U}} \Big[\boldsymbol{1}_{\{u=v\}}\mathcal B[\xi^M;u,v] + \frac12 \boldsymbol{1}_{\{u\neq v\}}\Big(\mathcal B[\xi^M;u,v]+\mathcal B[\xi^M;v,u]\Big)\Big]F(\xi^M,u,v)\\
    &\ge
    \sum_{\substack{\xi^M\in \Gamma_{\Lambda_n^M}\\ u,v\in U}} \Big[\boldsymbol{1}_{\{u=v\}}\mathcal B[\xi^M;u,v] - \frac12 \boldsymbol{1}_{\{u\neq v\}}\Big(\mathcal B[u(\xi^M);u^{-1},u^{-1}]+\mathcal B[v(\xi^M);v^{-1},v^{-1}]\Big)\Big]F(\xi^M,u,v)\\
    &=
    \sum_{\substack{\xi^M\in \Gamma_{\Lambda_n^M}\\ u,v\in U}} \Big[\boldsymbol{1}_{\{u=v\}}\mathcal B[\xi^M;u,v] - \boldsymbol{1}_{\{u\neq v\}}\mathcal B[u(\xi^M);u^{-1},u^{-1}]\Big]F(\xi^M,u,v).
\end{align*}
Using the explicit form of $F$ in \eqref{eq:S} and again reversibility, and dropping the non-negative diagonal birth terms, we obtain
\begin{align}\nonumber
    \mathcal B(\rho,\psi) &\geq \sum_{\substack{\xi^M\in \Gamma_{\Lambda_n^M}\\{x^M}\in \xi^M}} \mathcal B[\xi^M;u^-_{x^M},u^-_{x^M}]\dR(\{\xi^M\})\nonumber\\
    &\quad-\sum_{\substack{\xi^M\in \Gamma_{\Lambda_n^M}\\{x^M}\in \xi^M}}\sum_{y^M\in\Lambda_n^M\setminus\xi^M} \big(\kappa(\xi^M-\delta_{x^M},y^M)-\kappa(\xi^M,y^M)\big)\lambda(\{y^M\})\mathcal B[\xi^M;u^-_{x^M},u^-_{x^M}]\dR(\{\xi^M\})\nonumber\\
    &\geq 2\zeta(\kappa)\sum_{\substack{\xi^M\in \Gamma_{\Lambda_n^M}\\ {x^M}\in \xi^M}} \mathcal B[\xi^M;u^-_{x^M},u^-_{x^M}]\dR(\{\xi^M\}) \label{eq:B-final}.
\end{align}

Finally, concavity and non-negativity of $\theta$ imply that
\begin{equation}\label{eq:Theta-props}
p\partial_1\theta(p,q)+q\partial_2\theta(p,q)
\leq\theta(p,q),\qquad p,q>0.
\end{equation}
Moreover, monotonicity of $\theta$ implies that $\partial_1\theta,\partial_2\theta\geq0$. This yields, for any $\xi^M\in \Gamma_{\Lambda_n^M}$ and $u\in U$, that
\begin{align*}
    \mathcal B[\xi^M;u,u]\nonumber
    &=
    \!\frac12 |D_u \psi(\xi^M)|^2 \!\big[\partial_1\theta(\rho(\xi^M),\rho(u(\xi^M)))D_u\rho(\xi^M) +2 \theta(\rho(\xi^M),\rho(u(\xi^M)))\big]\\
    &\ge\!
    \frac12 |D_u \psi(\xi^M)|^2 \!\big[(\partial_1\theta(\rho(\xi^M),\rho(u(\xi^M)))\!+\!\partial_2\theta(\rho(\xi^M),\rho(u(\xi^M))))\rho(u(\xi^M)) \!+\! \theta(\rho(\xi^M),\rho(u(\xi^M)))\big]\\
    &\ge
    \!\frac12 |D_u \psi(\xi^M)|^2\theta(\rho(\xi^M),\rho(u(\xi^M))).
\end{align*}
With this and \eqref{eq:A-}, we finally obtain
\begin{align*}
\sum_{\substack{\xi^M\in \Gamma_{\Lambda_n^M}\\{x^M}\in \xi^M}} \!\!\!\mathcal B[\xi^M;u^-_{x^M},u^-_{x^M}] \dR(\{\xi^M\}) \ge \frac 12 \sum_{\substack{\xi^M\in \Gamma_{\Lambda_n^M}\\{x^M}\in \xi^M}}\!\!\!|D_{x^M}^-\psi(\xi^M)|^2\theta(\rho)(\xi^M,u_{x^M}^-)\dR(\{\xi^M\}) = \frac12 \mathcal A(\rho,\psi).
\end{align*}
Together with \eqref{eq:B-final}, this yields $\mathcal B(\rho,\psi)\geq\zeta(\kappa)\mathcal A(\rho,\psi)$. 
\end{proof}
\subsubsection{Convergence to Continuous Base Spaces} \label{sec:curvaturebounded}
Having established the EVI in the discrete setting  (Theorem~\ref{thm:evidiscretized}), we now let the mesh size tend to zero to obtain the EVI in the continuous setting. The main challenge is to establish the convergence of the quantities appearing in the discrete EVI, such as the transport distance, the dual spatial birth-and-death semigroup, and the relative entropy, to their continuous counterparts. To state the required convergence results precisely, we recall the definition of $\Gamma$-convergence.
\begin{dfn}[$\Gamma$-convergence]
Let $\Omega$ be a topological space and 
$F_i,F:\Omega\to(-\infty,\infty],$ $i\in\N$. The sequence $(F_i)_{i\in\N}$ $\Gamma$-converges to $F$, denoted by
\(
\Gamma\text{-}\lim_{i\to\infty} F_i=F,
\)
if the following two conditions are satisfied:
\begin{itemize}
    \item[$(i)$] For all $(\omega_i)_{i\in\N}\subset\Omega$ and $\omega\in\Omega$ s.t.\  $\lim_{i\to\infty}\omega_i=\omega$, it holds
    \[
    F(\omega)\leq \liminf_{i\to\infty} F_i(\omega_i).
    \]
    \item[$(ii)$] For all $\omega\in \Omega$ there exists a sequence $(\omega_i)_{i\in\N}\subset\Omega$ with $\lim_{i\to\infty}\omega_i=\omega$ s.t.\ 
    \[
    F(\omega)\geq \limsup_{i\to\infty} F_i(\omega_i).
    \]
    A sequence satisfying this condition is called recovery sequence.
\end{itemize}
For \(\Omega_\infty\subset\Omega\), we say that \((F_i)_{i\in\N}\) \(\Gamma\)-converges to
	\(F\) on \(\Omega_\infty\) if conditions $(i)$ and $(ii)$
are satisfied for all \(\omega\in\Omega_\infty\) with sequences $(\omega_i)_{i\in\N}\subset\Omega$.
\end{dfn}

\begin{satz}\label{thm:convergenceresultsdiscretization}
Let $\lambda=\Leb_{|\Lambda_n}$ and assume that
$\kappa:\Gamma_{\Lambda_n}\times\Lambda_n\to[0,\infty)$ satisfies
\hyperref[property:locallystable]{$(\mathbf{S}(1))$}, \ref{property:lowerbound},
\ref{property:continuity2}, and~\ref{property:repulsive}. Further, let
$\dR=\RRR_{\kappa,\lambda}(\Gamma_{\Lambda_n})$. Then 
\begin{enumerate}[label=$(\roman*)$]
    \item 
    $\lim_{M\to\infty}\mathrm{S}_t^{\kappa^M,\lambda^M*}\dP^M
    =\mathrm{S}_t^{\kappa,\lambda*}\dP$ in $\PPP_1(\Gamma_{\Lambda_n})$ for all $\dP\in\PPP_1(\Gamma_{\Lambda_n})$ s.t.\  $\dP\ll\dR$ and $t\geq0$;

    \item $\Gamma\text{-}\lim_{M\to\infty}\Ent(\cdot|\tilde{\dR}^M)
    =\Ent(\cdot|\dR)$ on $\PPP_1(\Gamma_{\Lambda_n})$.  A recovery sequence for $\dP\in\PPP_1(\Gamma_{\Lambda_n})$ is  $(\dP^M)_{M\in\N}$;

\item[$(iii)$] $\Gamma\text{-}\lim_{M\to\infty}\mathcal{W}_{\kappa^M,\lambda^M}(\cdot,\cdot)=\mathcal{W}_{\kappa,\lambda}(\cdot,\cdot)$ on $\PPP_1(\Gamma_{\Lambda_n})\times\{\dP\in\PPP_1(\Gamma_{\Lambda_n}):\dP\ll\dR\}$ (or with the two spaces interchanged). A recovery sequence for $(\dP_0,\dP_1)\in\PPP_1(\Gamma_{\Lambda_n})\times\{\dP\in\PPP_1(\Gamma_{\Lambda_n}):\dP\ll\dR\}$ is $(\dP_0^M,\dP_1^M)_{M\in\N}$.
\end{enumerate}
\end{satz}
For the sake of  a clearer structure, we postpone the proof to Section~\ref{sec:convergencediscretization}. The assumptions of Theorem~\ref{thm:convergenceresultsdiscretization} are stated in a unified form. Some of the individual statements hold under weaker assumptions specified in the subsections of the corresponding proofs.

Theorem~\ref{thm:convergenceresultsdiscretization} is stated under the additional assumption that the Papangelou intensity $\kappa$ is bounded below by a strictly positive constant \ref{property:lowerbound}. In order to weaken the additional assumption \ref{property:lowerbound}, we use \ref{property:ALB} and approximate the Papangelou intensity $\kappa$ from above by Papangelou intensities that are strictly positively bounded from below. The following theorem establishes the necessary convergence results for this approximation.  

 \begin{satz} \label{thm:convergenceresultstruncation}
Let $\lambda=\Leb_{|\Lambda_n}$ and assume that $\kappa:\Gamma_{\Lambda_n}\times\Lambda_n\to[0,\infty)$ satisfies \ref{property:cocycle}, \ref{property:locallystable}, and~\ref{property:continuity2}. Further, assume that $\kappa$ satisfies \ref{property:ALB}, and let $(\kappa_a)_{a\in\N}$ be an associated approximating sequence of $\kappa$. Let $\dR=\RRR_{\kappa,\lambda}(\Gamma_{\Lambda_n})$ and $\dR_a=\RRR_{\kappa_a,\lambda}(\Gamma_{\Lambda_n}),\; a\in\N.$  Then
    \begin{itemize}
        \item[$(i)$]  $\lim_{a \to\infty}\mathrm{S}_t^{\kappa_a ,\lambda*}\dP= \mathrm{S}_t^{\kappa,\lambda*}\dP$  in $\PPP_1(\Gamma_{\Lambda_n})$ for all $\dP\in\PPP_1(\Gamma_{\Lambda_n})$ s.t.\  $\dP\ll\dR$ and $t\geq0$;
        \item [$(ii)$]$\Gamma\text{-}\lim_{a\to\infty}\Ent(\cdot|{\dR}_a)
    =\Ent(\cdot|\dR)$ on $\PPP_1(\Gamma_{\Lambda_n})$. A recovery sequence for $\dP\in\PPP_1(\Gamma_{\Lambda_n})$ is given by the constant sequence $(\dP)_{a\in\N}$;

    \item[$(iii)$] $\Gamma\text{-}\lim_{a\to\infty}\mathcal{W}_{\kappa_a,\lambda}(\cdot,\cdot)
    =\mathcal{W}_{\kappa,\lambda}(\cdot,\cdot)$  on $\PPP_1(\Gamma_{\Lambda_n})\times \PPP_1(\Gamma_{\Lambda_n})$.  A recovery sequence for  $(\dP_0,\dP_1)\in$ $\PPP_1(\Gamma_{\Lambda_n})\times \PPP_1(\Gamma_{\Lambda_n})$ is given by the constant sequence $(\dP_0,\dP_1)_{a\in\N}$.
    \end{itemize}
    \end{satz}
     Again, the proof of Theorem~\ref{thm:convergenceresultstruncation} is not discussed here but given in Section~\ref{sec:truncation}. Likewise, the assumptions of Theorem~\ref{thm:convergenceresultstruncation} are stated in a unified form, while some of the individual statements hold under weaker assumptions specified in the subsections of the corresponding proofs.
     
    Having all necessary approximation results at hand, we finish this subsection with the proof of Theorem~\ref{thm:evi} and Theorem \ref{thm:geodconvexity}.
    \begin{proof}[Proof of Theorem~\ref{thm:evi}]
    Since $\kappa$ satisfies \ref{property:ALB}, there exists a sequence of Papangelou intensities $(\kappa_a)_{a\in\N}$ s.t.\      $\kappa_a:\Gamma_{\Lambda_n}\times\Lambda_n\to[0,\infty),\;a\in\N,$  satisfies \ref{property:cocycle}, \ref{property:lowerbound},  \hyperref[property:locallystable]{$(\mathbf{S}(1))$}, \ref{property:continuity2},  \ref{property:repulsive}, and
      $$0\leq\kappa_a(\xi,x)-\kappa(\xi,x)\leq g(a),$$ for all $\xi\in\Gamma_{\Lambda_n}, x\in {\Lambda_n}$, with $\lim_{a\to\infty}g(a)=0$. Let $\kappa_a ^M(\xi^M,x^M):=\kappa_a(\xi^M,x^M)\ind_{\{\xi^M(\{x^M\})=0\}}$ and $\tilde{\dR}_a ^M=\RRR_{\kappa_a ^M,\lambda^M}(\Gamma_{\Lambda_n^M})$. Since $\kappa$ satisfies \hyperref[property:locallystable]{$(\mathbf{S}(1))$} and~\ref{property:repulsive}, also $\kappa_a ^M$ satisfies these properties by definition.  It holds that
    \begin{align*}
     &|\zeta(\kappa_a^M)-\zeta(\kappa)|\\
     &\quad\leq|\zeta(\kappa_a^M)-\zeta(\kappa^M)|+|\zeta(\kappa^M)-\zeta(\kappa)|\\ &\quad\leq\frac12\sup_{\substack{\xi^M\in\Gamma_{\Lambda_n^M}\\y^M\in\Lambda_n^M}}\int|\kappa_a^M(\xi^M,x^M)- \kappa^M(\xi^M,x^M)|+|\kappa_a^M(\xi^M+\delta_{y^M},x^M)- \kappa^M(\xi^M+\delta_{y^M},x^M)|\lambda^M(\mathrm{d}x^M)\\
      &\qquad +|\zeta(\kappa^M)-\zeta(\kappa)|\\
      &\quad\leq\lambda(\Lambda_n)g(a) +|\zeta(\kappa^M)-\zeta(\kappa)|.
  \end{align*}
    Hence, by \ref{property:dobrushin2}, we obtain that \begin{align} 0<\lim_{M\to\infty}\lim_{a\to\infty}\zeta(\kappa_a^M)=\lim_{a\to\infty}\liminf_{M\to\infty}\zeta(\kappa_a^M)=\zeta(\kappa)\leq\frac12.\label{eq:1705}\end{align}
    In particular, $\kappa_a ^M$  satisfies \hyperref[property:dobrushin]{$(\mathbf{D}(\zeta(\kappa_a ^M)))$} for $a,M\in\N$ large enough. Since both $\Ent(\dP_1|\dR)<\infty$  and $\mathcal{W}_{\kappa,\lambda}(\dP_0,\dP_1)<\infty,$ we obtain by Lemma \ref{lem:preservationabsolutecontinuity} that $\dP_0\ll\dR.$ 
    Further, $\kappa\leq\kappa_a$ implies that $\dR\ll\dR_a$. Consequently, $\dP_i\ll\dR_a$ and hence $\dP_i^M\ll\dR_a^M$. The lower-bound assumption on $\kappa_a$ implies
$\dR_a^M\ll\tilde{\dR}_a^M$ and since the set of simple configurations on $\Lambda_n^M$ is finite,
we deduce $\Ent(\dP_0^M|\tilde{\dR}_a^M),\Ent(\dP_1^M|\tilde{\dR}_a^M)<\infty.$
    Hence, the integrated form of the discrete $\mathrm{EVI}(\zeta(\kappa_a ^M))$ (Theorem~\ref{thm:evidiscretized}) applies, and for $t>0$, it holds that
    \begin{align}\nonumber
    &\liminf_{a\to\infty}\liminf_{M\to\infty}\left(\frac{e^{\zeta(\kappa_a ^M)t}}{2}\mathcal{W}_{\kappa_a ^M,\lambda^M}^2(\mathrm{S}_t^{\kappa_a ^M,\lambda^M*}\dP_0^M,\dP_1^M)-\frac{1}{2}\mathcal{W}_{\kappa_a ^M,\lambda^M}^2(\dP_0^M,\dP_1^M)\right)\\
    &\quad\leq  \liminf_{a\to\infty}\liminf_{M\to\infty}\left(E_{\zeta(\kappa_a ^M)}(t)(\Ent(\dP_1^M | \tilde\dR_a ^M)-
\Ent(\mathrm{S}_t^{\kappa_a ^M,\lambda^M*}\dP_0^M | \tilde\dR_a ^M))\right).\label{eq:3456}
    \end{align}
 
    We note that $\kappa_a $ satisfies by assumption  \hyperref[property:locallystable]{$(\mathbf{S}(1))$}, \ref{property:lowerbound}, \ref{property:continuity2}, and~\ref{property:repulsive}. Thus, the claim follows from \eqref{eq:3456} by applying the convergence results of \eqref{eq:1705},  Theorem~\ref{thm:convergenceresultsdiscretization}, and Theorem~\ref{thm:convergenceresultstruncation}. \end{proof}

\begin{proof}[Proof of Theorem~\ref{thm:FIbounded}~$(iv)$]
We use the same notation as in the proof of Theorem~\ref{thm:evi}. By the EVI for discrete spaces (Theorem~\ref{thm:evidiscretized}) and \cite[Theorems 7.4 and 7.5]{Erbar_2012}, we obtain the discrete Talagrand inequality
\begin{align*}
\mathcal W_{\kappa_a^M,\lambda^M}^2(\dP_0^M,\tilde\dR_a^M) &\leq \frac{2}{\zeta(\kappa_a^M)}\Ent(\dP_0^M|\tilde\dR_a^M),
\end{align*}
for $a,M\in\N$ sufficiently large. Applying the convergence results of Theorems~\ref{thm:convergenceresultsdiscretization} and \ref{thm:convergenceresultstruncation}, together with \eqref{eq:1705}, analogously to the proof of Theorem~\ref{thm:evi}, finishes the proof.
\end{proof}

\subsection{$\Phi$-Entropies}\label{sec:phientropies}
In this section, we extend the exponential contraction of the
relative entropy along the birth-and-death semigroup
$(\mathrm S_t^{\kappa,\lambda*})_{t\geq0}$ to more general
$\Phi$-entropies
\[
\Ent_\Phi(\dP|\dR):=
\begin{cases}
\displaystyle\int\Phi(\rho)\,\mathrm d\dR,
&\text{if }\dP=\rho\dR,\\
\infty,&\text{otherwise},
\end{cases}
\]
where $\Phi\in\CCC([0,\infty))\cap\CCC^\infty((0,\infty))$
is superlinear, satisfies $\Phi^{\prime \prime}>0$, and is normalized by
$\Phi(1)=\Phi^{\prime}(1)=0$.
The relative entropy corresponds to
$\Phi(r)=r\log r-r+1$. Further examples are discussed in Example \ref{ex:renyi}.

In the discrete setting, the Bochner-Bakry-\'Emery method that we used to obtain the discrete EVI and MLSI for $\Ent(\cdot|\tilde\R^M_a)$ has been adapted in \cite{JuengelYue2017} to obtain MLSI for $\Phi$-entropies. 
Here, one replaces the logarithmic mean by the weight function
\[
\theta_\Phi(r,s):=
\begin{cases}
\displaystyle\frac{r-s}{\Phi^{\prime}(r)-\Phi^{\prime}(s)},&r\neq s,\\[1ex]
\displaystyle\frac{1}{\Phi^{\prime \prime}(r)},&r=s,
\end{cases}
\qquad r,s>0.
\]
We assume that $\Phi$ is such that $\theta_\Phi$ extends continuously to a
non-negative function on $[0,\infty)^2$, is non-decreasing
in each variable, and jointly concave. We highlight, that we do not require that $\theta_\Phi$ is
$1$-homogeneous.
Starting from the results in \cite{JuengelYue2017} and following our discretisation procedure, we obtain the following exponential
entropy contraction.

\begin{satz}[Exponential entropy contraction]
\label{thm:Phicontraction}
Let $\lambda=\Leb_{|\Lambda_n}$ and assume that
$\kappa:\Gamma_{\Lambda_n}\times\Lambda_n\to[0,\infty)$
satisfies \ref{property:cocycle},
\hyperref[property:locallystable]{$(\mathbf{S}(1))$},
\ref{property:ALB}, \ref{property:continuity2},
\ref{property:repulsive}, and~\ref{property:dobrushin2}.
Further, let
$\dR=\RRR_{\kappa,\lambda}(\Gamma_{\Lambda_n})$.
Then, for every $\dP\in\PPP_1(\Gamma_{\Lambda_n})$, it holds that
\[
\Ent_\Phi\bigl(
\mathrm S_t^{\kappa,\lambda*}\dP\mid\dR
\bigr)
\leq
e^{-2\zeta(\kappa)t}\Ent_\Phi(\dP|\dR),
\qquad t\geq0.
\]
\end{satz}

\begin{proof}
Under the above assumption on $\Phi$, we can apply \cite[Cor.~4]{JuengelYue2017} to the discretised dynamics to obtain
\begin{align*}
\Ent_\Phi(
\mathrm S_t^{\kappa_a^M,\lambda^M*}\dP^M
\mid\tilde\dR_a^M)
\leq
e^{-2\zeta(\kappa_a^M)t}
\Ent_\Phi(\dP^M|\tilde\dR_a^M).
\end{align*}
We then follow the same line of argumentation as in the proof of Theorem \ref{thm:evi} and pass to the limit first for $a$ fixed and $M$ going to infinity and then for $a$ going to infinity. The necessary entropy $\Gamma$-convergence results in
Theorems~\ref{thm:convergenceresultsdiscretization}~$(ii)$
and~\ref{thm:convergenceresultstruncation}~$(ii)$
extend to $\Ent_\Phi$ with only minor modifications, finishing the proof.
\end{proof}
\begin{bsp}\label{ex:renyi}
The power entropy densities
\[
    \Phi_p(r):=\frac{r^p-r}{p-1}-r+1,
    \qquad 1<p\leq2,
\]
generate the power entropies 
\begin{align*}
   \Ent_{\Phi_p}(\dP|\dR):=
    \begin{cases}
        \displaystyle\frac{1}{p-1}
       \Big(\int\rho(\xi)^p\dR(\mathrm{d}\xi)-1\Big),
        &\text{if }\dP=\rho\dR,\\
        \infty,&\text{otherwise.}
    \end{cases}
\end{align*}
associated with R\'enyi divergences. Their corresponding mean function is given by
\[
    \theta_{\Phi_p}(r,s):=
    \begin{cases}
        \displaystyle\frac{p-1}{p}
        \frac{r-s}{r^{p-1}-s^{p-1}},&r\neq s,\\[1ex]
        \displaystyle\frac1p r^{2-p},&r=s,
    \end{cases}
    \qquad r,s>0.
\]
and satisfy the assumptions introduced above, see \cite[Lemma~16]{JuengelYue2017}).
As $p\downarrow1$, $\Phi_p$ converges to the logarithmic
entropy density, while $p=2$ yields
$\Phi_2(r)=(r-1)^2$, corresponding to the so-called
$\chi^2$-divergence.
\end{bsp}
\section{The Transport Distance}\label{sec:transportdistance}
    Let $X$ be a Polish space, $\lambda\in\MMM^+(X)$, and $\kappa:\Gamma_X\times X\to[0,\infty) $. In this section, we slightly generalize the  transport distance $\mathcal{W}_{\kappa,\lambda}$ introduced in Subsection \ref{subsec:generalsetup} and discuss its central properties. Recall that  $\mathcal{W}_{\kappa,\lambda}:\PPP_1(\Gamma_X)\times\PPP_1(\Gamma_X)\to[0,\infty]$ is defined in
    analogy with the classical Benamou--Brenier formulation of the Wasserstein distance, by \begin{align*}
\mathcal{W}_{\kappa,\lambda}^2(\dP_0, \dP_1) &= \inf \{ \mathcal{A}_{\kappa,\lambda}(\Bar{\dP}, \Bar{\dV}) : (\Bar{\dP}, \Bar{\dV}) \in \CCC\EEE^X(\dP_0, \dP_1)\},
    \end{align*} 
    where $\CCC\EEE^X(\dP_0,\dP_1)$ is the set of solutions satisfying the continuity equation (see Definition \ref{dfn:ce}) and $\mathcal{A}_{\kappa,\lambda}$ is an action functional defined as follows. Recall the definition of  the action density $$\alpha(x,y,w)=\frac{|w|^2}{\theta(x,y)},\quad x,y\geq0,\ w\in\R,$$ where  $0/0=0$  by convention. While  $\theta:[0,\infty)\times[0,\infty)\to[0,\infty)$ was assumed to be the logarithmic mean in Subsection \ref{subsec:generalsetup}, here it is a more general function satisfying the following properties.
   \begin{itemize}
       \item [(A1)] (Regularity) $\theta$ is continuous on $[0,\infty)\times [0,\infty)$ and $C^1$ on $(0,\infty)\times(0,\infty);$
       \item [(A2)] (Symmetry) $\theta(x,y)=\theta(y,x)$ for $x,y\geq0$;
       \item[(A3)] (Positivity and normalisation) $\theta(x,y)>0$ for $x,y>0$ and $\theta(1,1)=1$;
       \item[(A4)] (Zero at the boundary) $\theta(0,x)=0$ for all $x\geq 0$;
       \item[(A5)] (Monotonicity) $\theta(x,y)\leq \theta(z,y)$ for all $0\leq x\leq z$ and $y\geq 0$;
       \item[(A6)] (Positive homogeneity) $\theta(rx,ry)=r\theta(x,y)$ for all $r\geq0$;
       \item[(A7)]  (Concavity) The function $\theta:[0,\infty)\times [0,\infty)\to[0,\infty)$ is concave.
   \end{itemize}
   Assumptions (A1) -- (A7) imply that $\alpha$  is lower semi-continuous, convex, and positively
homogeneous, see e.g.\ \cite[Lemma 2.2]{erbar2012gradientflowsentropyjump}. 

        The \textit{Lagrange functional}  w.r.t.\ the Papangelou intensity $\kappa$ and reference measure $\lambda$  is for $\dP\in\PPP_1(\Gamma_X)$ and $\dV\in \MMM_{b,0}(\Gamma_X\times X)$ defined by
        \begin{align}
           \mathcal{L}_{\kappa,\lambda}({\dP},{\dV})=\int\alpha\left(\frac{\mathrm{d}\kappa(\dP\otimes \lambda)}{\mathrm{d}\sigma}(\xi,x),\frac{\mathrm{d}C_{\dP}}{\mathrm{d}\sigma}(\xi,x),\frac{\mathrm{d}{\dV}}{\mathrm{d}\sigma}(\xi,x)\right)\sigma(\mathrm{d}\xi,\mathrm{d}x),
       \end{align}
       where $\sigma\in\MMM_{b,0}(\Gamma_X\times X)$ is non-negative s.t.\   $\kappa(\dP\otimes {\lambda}),$ $C_{\dP},$ ${\dV}\ll\sigma.$ Due to the positive homogeneity of $\alpha$, the definition of the Lagrange functional $\mathcal{L}_{\kappa,\lambda}$ is independent of the choice of $\sigma$. The \textit{action functional}  w.r.t.\ the Papangelou intensity $\kappa$ and reference measure $\lambda$  is defined by
     \begin{align*}
         \mathcal{A}_{\kappa,\lambda}({\Bar{\dP}},{\Bar{\dV}})=\int_0^1\mathcal{L}_{\kappa,\lambda}(\Bar{\dP}_t,\Bar{\dV}_t)\mathrm{d}t,\quad (\Bar{\dP},\Bar{\dV})\in\CCC\EEE^X(\dP_0,\dP_1).
     \end{align*}
    \subsection{Properties of the Lagrange and the Action Functionals}\label{subsec:Lagrangeandaction}
    Assume that $\dR\in\RRR_{\kappa,\lambda}(\Gamma_X)$ exists.
     If $\dP=\rho\dR\in\PPP_1(\Gamma_X)$ and  $\dV= w(\kappa(\dP\otimes\lambda))\in\MMM_{b,0}(\Gamma_X\times X)$, the Lagrange functional $\mathcal{L}_{\kappa,\lambda}$ can be expressed in terms of the densities $\rho$ and $w$. By choosing $\sigma=\kappa(\dP\otimes\lambda)$ and using Corollary~\ref{kor:CampbellPapangelou}, we obtain
     \begin{align}\label{eq:representationLagrangeGibbs}
             \mathcal{L}_{\kappa,\lambda}({\dP},{\dV})&=\int\int\alpha(\rho(\xi),\rho(\xi+\delta_x),w(\xi,x))\kappa(\xi,x)\dR(\mathrm{d}\xi)\lambda(\mathrm{d}x).
         \end{align}

    The following lemma provides a condition for the velocity field to be absolutely continuous w.r.t.\ $\kappa(\dP\otimes\lambda)$.
    \begin{lem}\label{lem:absolutecontinuityvelocity}
          Let $\lambda\in\MMM^+(X)$ and $\kappa:\Gamma_X\times X\to[0,\infty)$. Assume that $\dR\in\RRR_{\kappa,\lambda}(\Gamma_X)$ exists and let $\dP\in\PPP_1(\Gamma_X)$ s.t.\  $\dP\ll\dR$. If  $\dV\in\MMM_{b,0}(\Gamma_X\times X)$ s.t.\   $\mathcal{L}_{\kappa,\lambda}(\dP,\dV)<\infty$, then  $\dV\ll \kappa(\dP\otimes \lambda)$. 
         \end{lem}
         \begin{proof}
             The proof is identical to that of {\cite[Lemma 5.3]{schiavo2023wasserstein}}.
         \end{proof}

         \begin{lem}\label{lem:BoundviaLagrange}
      Let $\lambda\in\MMM^+(X)$ and assume that $\kappa:\Gamma_X\times X\to[0,\infty)$ satisfies \ref{property:locallystable}. For all $\dP\in\PPP_1(\Gamma_X)$, $\dV\in\MMM_{b,0}(\Gamma_X\times X)$, $K\in\BBB(\Gamma_X)$, and $B\in\BBB_0(X)$, it holds that 
         \begin{align}\label{eq:BoundviaLagrange}
             |\dV|(K\times B)\leq \left(\frac{1}{2}(A\lambda(B)+\I_\dP(B))\right)^{1/2}\mathcal{L}_{\kappa,\lambda}(\dP,\dV)^{1/2}.
         \end{align}
     \end{lem}
     \begin{proof}
         The proof follows the same line of argumentation as the proof of \cite[Lemma 5.5]{schiavo2023wasserstein} and \cite[Lemma 2.6]{erbar2012gradientflowsentropyjump}. Due to the additional Papangelou intensity and for the sake of completeness, we provide a full proof here. Assume that $\mathcal{L}_{\kappa,\lambda}(\dP,\dV)<\infty$, since otherwise the claim is immediate. Let
\[
\sigma=\kappa(\dP\otimes \lambda)+C_\dP+|\dV|
\]
and write $\kappa(\dP\otimes \lambda)=\psi\sigma$, $C_\dP=\phi\sigma$, and $\dV=w\sigma.$ Further, let $K\in\BBB(\Gamma_X)$ and $B\in\BBB_0(X)$. By the definition of $\sigma$, we have
\[
|\dV|(K\times B)=\int_{K\times B}|w(\xi,x)|\sigma(\mathrm{d}\xi,\mathrm{d}x).
\]
Moreover, since $\mathcal{L}_{\kappa,\lambda}(\dP,\dV)<\infty$, it follows that $w=0$ $\sigma$-a.e.\ on $\{\theta(\psi,\phi)=0\}$, and hence
\[
|w(\xi,x)|
=
\mathbf{1}_{\{\theta(\psi(\xi,x),\phi(\xi,x))>0\}}
\theta(\psi(\xi,x),\phi(\xi,x))^{1/2}
\alpha(\psi(\xi,x),\phi(\xi,x),w(\xi,x))^{1/2}
\quad \sigma\text{-a.e.}
\]
Therefore, by the Cauchy--Schwarz inequality,
\begin{align*}
|\dV|(K\times B)
&= \int_{K\times B}\theta(\psi,\phi)^{1/2}\alpha(\psi,\phi,w)^{1/2}\sigma(\mathrm{d}\xi,\mathrm{d}x)\\
&\leq \left(\ \int_{K\times B}\theta(\psi,\phi)\sigma(\mathrm{d}\xi,\mathrm{d}x)\right)^{1/2}
\left(\ \int_{K\times B}\alpha(\psi,\phi,w)\sigma(\mathrm{d}\xi,\mathrm{d}x)\right)^{1/2}\\
&\leq \left(\ \int_{K\times B}\theta(\psi,\phi)\sigma(\mathrm{d}\xi,\mathrm{d}x)\right)^{1/2}
\mathcal{L}_{\kappa,\lambda}(\dP,\dV)^{1/2}.
\end{align*}
              Properties (A1)--(A7) of $\theta$ directly imply that  $\theta(x,y)\leq (x+y)/2$, see  \cite[p. 7]{erbar2012gradientflowsentropyjump}. Together with the local stability of $\kappa$ \ref{property:locallystable}, this yields 
         \begin{align*}
             \int_{K\times B}\theta(\psi(\xi,x),\phi(\xi,x))\sigma(\mathrm{d}\xi,\mathrm{d}x)&\leq \frac{1}{2} \int_{K\times B}\psi(\xi,x)+\phi(\xi,x)\sigma(\mathrm{d}\xi,\mathrm{d}x)\\
             &=\frac{1}{2}\int_B\int_{K}\kappa(\xi,x)\dP(\mathrm{d}\xi)\lambda(\mathrm{d}x)+\frac{1}{2}C_\dP(K\times B)\\
             &\leq \frac{1}{2}\int_B\int_{K}A\dP(\mathrm{d}\xi)\lambda(\mathrm{d}x)+\frac{1}{2}C_\dP(\Gamma_X\times B)\\
             &\leq \frac{1}{2}(A\lambda(B)+\I_{\dP}(B)).
         \end{align*}
         Combining both estimates finishes the proof.
     \end{proof}
     In order to establish lower semicontinuity of the Lagrange functional, we first recall a general result for integral functionals.
        \begin{satz}[{\cite[Theorem 3.4.3]{buttazzo1989semicontinuity}}]\label{athm:lowersemicontinuity}
        Let $\Omega$ be a locally compact Polish space and $f:\Omega\times \R^d\to[0,\infty]$ be a lower semicontinuous function s.t.\    $f(\omega,\cdot)$ is convex and positively homogeneous for every $\omega\in\Omega$. Then the functional \begin{align*}
            F(\lambda)=\int_{\Omega}f\left(\omega,\frac{\mathrm{d}\lambda}{\mathrm{d}|\lambda|}(\omega)\right)|\lambda|(\mathrm{d}\omega)
       \end{align*}
        is vaguely lower semicontinuous on the space of vector-valued signed measures $\mathcal{M}(\Omega,\R^d)$.
    \end{satz}
     \begin{lem}\label{lem:propertiesLagrange}
       Let $\lambda\in\MMM^+(X)$ and assume that $\kappa:\Gamma_X\times X\to[0,\infty)$ satisfies \ref{property:locallystable} and~\ref{property:continuity}. 
      \begin{itemize}
      \item[$(i)$] 
              The Lagrange functional $\mathcal{L}_{\kappa,\lambda}$ is jointly convex, i.e.\  for   $({\bar{\dP}}^1,{\bar{\dV}}^1),({\bar{\dP}}^2,{\bar{\dV}}^2)\!\in\!\PPP_1(\Gamma_X)\times \MMM_{b,0}(\Gamma_X\times X)$ it holds that 
             \begin{align*}
                 \mathcal{L}_{\kappa,\lambda}(c{\bar{\dP}}^1+(1-c){\bar{\dP}}^2,c{\bar{\dV}}^1+(1-c){\bar{\dV}}^2)\leq c\mathcal{L}_{\kappa,\lambda} ({\bar{\dP}}^1,{\bar{\dV}}^1)+(1-c) \mathcal{L}_{\kappa,\lambda}({\bar{\dP}}^2,{\bar{\dV}}^2),\quad c\in[0,1].
            \end{align*}
          \item[$(ii)$] The Lagrange functional $\mathcal{L}_{\kappa,\lambda}$ is lower semicontinuous, i.e.\  if $(\dP^i)_{i\in\N}\subset\PPP_1(\Gamma_X)$ and $\dP\in\PPP_1(\Gamma_X)$ s.t.\   $\lim_{i\to\infty}\dP^i=\dP$ in $\PPP_1(\Gamma_X)$ and $(\dV^i)_{i\in\N}\subset\MMM_{b,0}(\Gamma_X\times X)$ and $\dV\in \MMM_{b,0}(\Gamma_X\times X)$ s.t.\   $\lim_{i\to\infty}\dV^i=\dV$ in $\MMM_{b,0}(\Gamma_X\times X)$, then     \begin{align*}
               \mathcal{L}_{\kappa,\lambda}(\dP,\dV)\leq \liminf_{i\to\infty}\mathcal{L}_{\kappa,\lambda}(\dP^i,\dV^i).
           \end{align*}
           \end{itemize}
       \end{lem}
       \begin{proof}
       $(i)$ The convexity of $\mathcal{L}_{\kappa,\lambda}$ follows directly from the convexity of the action density $\alpha$.
       
       $(ii)$ Assume that $\kappa$ is not only upper semicontinuous \ref{property:continuity} but continuous.
       Let $F\in\CCC_{b,0}(\Gamma_{X}\times X)$. As $\kappa$ is uniformly bounded \ref{property:locallystable}, $\kappa F\in\CCC_{b,0}(\Gamma_X\times X)$. Since $\lim_{i\to\infty}\dP^i=\dP$ in $\PPP_1(\Gamma_X)$, we obtain  
       \begin{align*}
           \lim_{i\to\infty}\int F(\xi,x)\kappa(\xi,x)\dP^i\otimes\lambda(\mathrm{d}\xi,\mathrm{d}x)= \int F(\xi,x)\kappa(\xi,x)\dP\otimes\lambda(\mathrm{d}\xi,\mathrm{d}x).
       \end{align*}
       This implies that $\lim_{i\to\infty}\kappa(\dP^i\otimes\lambda)=\kappa(\dP\otimes\lambda)$ in $\MMM_{b,0}(\Gamma_X\times X)$. Further, by Proposition~\ref{prop:convergencep1}, we obtain that $\lim_{i\to\infty}C_{\dP^i}=C_\dP$ in $\MMM_{b,0}(\Gamma_X\times X)$. We note that, for $\Omega=\Gamma_X\times X$, the vague topology of $\MMM(\Omega)$ coincides with the usual $\CCC_{b,0}(\Gamma_X\times X)$-weak topology on $\MMM_{b,0}(\Gamma_X\times X)$, see \cite[p.\ 968]{schiavo2023wasserstein}. Thus, Theorem~\ref{athm:lowersemicontinuity} implies that 
       \begin{align*}
           \mathcal{L}_{\kappa,\lambda}(\dP,\dV)&=\int \alpha\left(\frac{\mathrm{d}\kappa(\dP\otimes\lambda)}{\mathrm{d}\sigma}(\xi,x),\frac{\mathrm{d}C_\dP}{\mathrm{d}\sigma}(\xi,x), \frac{\mathrm{d}\dV}{\mathrm{d}\sigma}(\xi,x)\right)\sigma(\mathrm{d}\xi,\mathrm{d}x)\\
           &\leq \liminf_{i\to\infty}\int \alpha\left(\frac{\mathrm{d}\kappa(\dP^i\otimes\lambda)}{\mathrm{d}\sigma}(\xi,x),\frac{\mathrm{d}C_{\dP^i}}{\mathrm{d}\sigma}(\xi,x), \frac{\mathrm{d}\dV^i}{\mathrm{d}\sigma}(\xi,x)\right)\sigma(\mathrm{d}\xi,\mathrm{d}x)\\
           &= \liminf_{i\to\infty}\mathcal{L}_{\kappa,\lambda}(\dP^i,\dV^i),
       \end{align*}
       finishing the proof for $\kappa$ being continuous. Now, assume that $\kappa$ is upper semicontinuous. There exists a sequence of functions $(\kappa_j)_{j\in\N}\subset\CCC_b(\Gamma_X\times X)$ s.t.\   $\kappa_j\downarrow\kappa$ for $j\to\infty.$ By the monotone convergence theorem, the monotonicity of $\alpha$, and the previous result for the continuous functions $\kappa_j$, we obtain 
       \begin{align*}
           \mathcal{L}_{\kappa,\lambda}(\dP,\dV)&=\int \alpha\left(\frac{\mathrm{d}\kappa(\dP\otimes\lambda)}{\mathrm{d}\sigma}(\xi,x),\frac{\mathrm{d}C_\dP}{\mathrm{d}\sigma}(\xi,x), \frac{\mathrm{d}\dV}{\mathrm{d}\sigma}(\xi,x)\right)\sigma(\mathrm{d}\xi,\mathrm{d}x)\\
           &=\sup_{j\in\N}\int \alpha\left(\frac{\mathrm{d}\kappa_j(\dP\otimes\lambda)}{\mathrm{d}\sigma}(\xi,x),\frac{\mathrm{d}C_\dP}{\mathrm{d}\sigma}(\xi,x), \frac{\mathrm{d}\dV}{\mathrm{d}\sigma}(\xi,x)\right)\sigma(\mathrm{d}\xi,\mathrm{d}x)\\
           &\leq \sup_{j\in\N}\liminf_{i\to\infty}\int \alpha\left(\frac{\mathrm{d}\kappa_j(\dP^i\otimes\lambda)}{\mathrm{d}\sigma}(\xi,x),\frac{\mathrm{d}C_{\dP^i}}{\mathrm{d}\sigma}(\xi,x), \frac{\mathrm{d}\dV^i}{\mathrm{d}\sigma}(\xi,x)\right)\sigma(\mathrm{d}\xi,\mathrm{d}x)\\
           &\leq \liminf_{i\to\infty}\sup_{j\in\N}\int \alpha\left(\frac{\mathrm{d}\kappa_j(\dP^i\otimes\lambda)}{\mathrm{d}\sigma}(\xi,x),\frac{\mathrm{d}C_{\dP^i}}{\mathrm{d}\sigma}(\xi,x), \frac{\mathrm{d}\dV^i}{\mathrm{d}\sigma}(\xi,x)\right)\sigma(\mathrm{d}\xi,\mathrm{d}x)\\
           &= \liminf_{i\to\infty} \mathcal{L}_{\kappa,\lambda}(\dP^i,\dV^i),
       \end{align*}
       finishing the proof.
       \end{proof}

        In the next lemma, we lift the convexity and the lower semicontinuity result for the Lagrange functional to the action functional. In addition, we also establish a compactness result. Together, these results will be used in Theorem~\ref{thm:minimizer} to prove the existence of minimizers of the action functional.

    \begin{lem}\label{lem:compactnessofsolutions}
     Let $\lambda\in\MMM^+(X)$ and assume that $\kappa:\Gamma_X\times X\to[0,\infty)$ satisfies \ref{property:locallystable} and~\ref{property:continuity}.
    \begin{itemize}
    \item[$(i)$] The action functional $\mathcal{A}_{\kappa,\lambda}$ is jointly  convex, i.e.\  for $\dP_0,\dP_1\in\PPP_1(\Gamma_X)$ and all pairs $({\bar{\dP}}^1,{\bar{\dV}}^1),$ $({\bar{\dP}}^2,{\bar{\dV}}^2)\in\CCC\EEE^X(\dP_0,\dP_1)$, it holds that \begin{align*}
                 \mathcal{A}_{\kappa,\lambda}(c{\bar{\dP}}^1+(1-c){\bar{\dP}}^2,c{\bar{\dV}}^1+(1-c){\bar{\dV}}^2)\leq c\mathcal{A}_{\kappa,\lambda} ({\bar{\dP}}^1,{\bar{\dV}}^1)+(1-c) \mathcal{A}_{\kappa,\lambda}({\bar{\dP}}^2,{\bar{\dV}}^2),\quad c\in[0,1].
             \end{align*} 
    \item[$(ii)$] The action functional $\mathcal{A}_{\kappa,\lambda}$ is lower semicontinuous, i.e.\    let $(\dP_0^i)_{i\in\N},(\dP_1^i)_{i\in\N}\subset\PPP_1(\Gamma_{X})$ and $\dP_0,\dP_1\in\PPP_1(\Gamma_X)$ s.t.\   $\lim_{i\to\infty}\dP_0^i=\dP_0$ and $\lim_{i\to\infty}\dP_1^i=\dP_1$ in $\PPP_1(\Gamma_X)$. Further, assume that $(\bar{\dP}^i,\bar{\dV}^i)\in\CCC\EEE^X(\dP_0^i,\dP_1^i),(\bar{\dP},\bar{\dV})\in\CCC\EEE^X(\dP_0,\dP_1),$ s.t.\   $\lim_{i\to\infty}\bar\dP_t^i=\bar\dP_t$ in $\PPP_1(\Gamma_X)$ for all $t\in[0,1]$ and $\lim_{i\to\infty}\bar\dV^i=\bar\dV$ in $\MMM_{b,0}(\Gamma_X\times X\times[0,1])$, then
    \begin{align*}
        \mathcal{A}_{\kappa,\lambda}(\bar{\dP},\bar{\dV})\leq \liminf_{i\to\infty}\mathcal{A}_{\kappa,\lambda}(\bar{\dP}^i,\bar{\dV}^i).
    \end{align*}
        \item[$(iii)$] 
        Let $(\dP_0^i)_{i\in\N},(\dP_1^i)_{i\in\N}\subset\PPP_1(\Gamma_X)$ and $\dP_0,\dP_1\in\PPP_1(\Gamma_X)$ s.t.\   $\lim_{i\to\infty}\dP_0^i=\dP_0$ and $\lim_{i\to\infty}\dP_1^i=\dP_1$ in $\PPP_1(\Gamma_X)$. Further, let $(\bar{\dP}^i,\bar{\dV}^i)_{i\in\N}\subset\CCC\EEE^X(\dP_0^i,\dP_1^i)$ s.t.\   
        \begin{align*}
            \sup_{i\in\N}\mathcal{A}_{\kappa,\lambda}(\bar{\dP}^i,\bar{\dV}^i)\leq C<\infty.
        \end{align*}
        There exists a subsequence $(\bar{\dP}^{i_j},\bar{\dV}^{i_j})_{j\in\N}$ and a solution $(\bar{\dP},\bar{\dV})\in\CCC\EEE^X(\dP_0,\dP_1)$ s.t.\   
        \begin{align*}
            \lim_{j\to\infty}\bar{\dP}_t^{i_j}=\bar{\dP}_t\quad &\text{ in }\PPP_1(\Gamma_X)\text{ for all } t\in [0,1],\\
            \lim_{j\to\infty}\bar{\dV}^{i_j}=\bar{\dV}\quad &\text{ in }\MMM_{b,0}(\Gamma_X\times X\times[0,1]).
             \end{align*}
        \end{itemize}
    \end{lem}
     \begin{proof}
     $(i)$ The convexity of the action functional $\mathcal{A}_{\kappa,\lambda}$ follows directly from the convexity of the Lagrange functional $\mathcal{L}_{\kappa,\lambda}$ (Lemma~\ref{lem:propertiesLagrange} $(i)$).
     
     $(ii)$ The proof of the lower semicontinuity follows analogously to Lemma~\ref{lem:propertiesLagrange} by using Theorem~\ref{athm:lowersemicontinuity} and considering $\int_0^1\alpha(\bar{\dP}_t,\bar{\dV}_t)\mathrm{d}t$ as an integral functional on $\MMM_{b,0}(\Gamma_X\times X\times [0,1])$.
     
    $(iii)$ The proof is analogous to that of \cite[Lemma 5.8, Lemma 5.9]{schiavo2023wasserstein}. However, due to the additional Papangelou intensity $\kappa$ and since we will adapt this proof to slightly different settings throughout the paper, we provide a full proof here. We first construct a subsequence of $(\bar{\dV}^i)_{i\in\N}$ and a limiting measure $\bar{\dV}\in\MMM_{b,0}(\Gamma_X\times X\times[0,1])$. Let $K\in\BBB_c(\Gamma_X),\;B\in\BBB_c(X),$ and $I\in\BBB_c([0,1])$ be compact sets. Since $\lim_{i\to\infty}\I_{\dP_0^i}=\I_{\dP_0}$ vaguely by Proposition~\ref{prop:convergencep1},  for $\vartheta>0$ there exists some $i_0\in\N$ s.t.\   $\I_{\dP_0^i}(B)\leq \I_{\dP_0}(B)+\vartheta$ for all $i\geq i_0$.  Hence, using Lemma~\ref{lem:BoundviaLagrange}, Theorem~\ref{thm:propertiesce} $(iv)$, and the Cauchy--Schwarz inequality, we have  \begin{align}
                    |\bar{\dV}^{i}|(K\times B\times I)&\leq \int_I \sqrt{\frac{1}{2}(A\lambda(B)+\I_{\dP_0^i}(B)+|\bar{\dV}^{{i}}|(\Gamma_X\times B\times [0,1]))}\mathcal{L}_{\kappa,\lambda}(\bar{\dP}^{{i}}_t,\bar{\dV}^{{i}}_t)^{1/2}\mathrm{d}t\nonumber\\
                    &\leq \int_I \sqrt{\frac{1}{2}(A\lambda(B)+\I_{\dP_0}(B)+\vartheta+|\bar{\dV}^{{i}}|(\Gamma_X\times B\times [0,1]))}\mathcal{L}_{\kappa,\lambda}(\bar{\dP}^{{i}}_t,\bar{\dV}^{{i}}_t)^{1/2}\mathrm{d}t\nonumber\\
                    &\leq \sqrt{\frac{1}{2}(A\lambda(B)+\I_{\dP_0}(B)+\vartheta+|\bar{\dV}^{{i}}|(\Gamma_X\times B\times [0,1]))}\sqrt{\Leb(I)}\sqrt{C}.\label{eq:10011}
                \end{align}
                In particular, this yields that 
                \begin{align*}
                    |\bar{\dV}^{i}|(\Gamma_X\times B\times [0,1])\leq \sqrt{\frac{1}{2}(A\lambda(B)+\I_{\dP_0}(B)+\vartheta+|\bar{\dV}^{{i}}|(\Gamma_X\times B\times [0,1]))}\sqrt{C}.
                \end{align*}
                By rearranging this equation, we obtain 
                \begin{align*}
                    |\bar{\dV}^{i}|(\Gamma_X\times B\times [0,1])^2- \frac{C}{2}(A\lambda(B)+\I_{\dP_0}(B)+\vartheta)-\frac{C}{2}|\bar{\dV}^{{i}}|(\Gamma_X\times B\times [0,1])\leq 0,
                \end{align*}
                which is a quadratic inequality and directly implies
                \begin{align*}
                    |\bar{\dV}^{i}|(\Gamma_X\times B\times [0,1])\leq \frac{C+\sqrt{C^2+8C(A\lambda(B)+\I_{\dP_0}(B)+\vartheta)}}{4}.
                \end{align*}
               Hence, choosing suitable constants, we get 
            \begin{align}\label{eq:10010}
             |\bar{\dV}^{i}|(K\times B\times I)\leq    |\bar{\dV}^{i}|(\Gamma_X\times B\times [0,1])\leq a_C+b_{C,A}\sqrt{\lambda(B)+\I_{\dP_0}(B)+\vartheta},
            \end{align}
            i.e.\ $ |\bar{\dV}^{i}|$ has total variation uniformly bounded on every compact subset of $\Gamma_X\times X\times I$. Hence, there exists a measure $\bar{\dV}\in\MMM_{b,0}(\Gamma_{X}\times X\times [0,1])$ and a sequence of natural numbers $(i_j)_{j\in\N}$  s.t.\   $\lim_{j\to\infty}\bar{\dV}^{i_j}=\bar{\dV}$  in $\MMM_{b,0}(\Gamma_{X}\times X\times [0,1])$. By \eqref{eq:10011}, if $I\in\BBB([0,1])$ s.t.\   $\Leb(I)=0$, then $|\bar{\dV}^{i}|(K\times B\times I)=0$. Hence, by disintegration  we can find a Borel family  $(\bar{\dV}_t)_{t\in[0,1]}$ satisfying $\bar{\dV}=\int\bar{\dV}_t\mathrm{d}t$. Further, let 
            \begin{align*}
                \bar{\dP}_t(K):={\dP_0}(K)+ \int_0^t \int D_x\ind_K(\xi)\bar{\dV}_s(\mathrm{d}\xi, \mathrm{d}x)\mathrm{d}s,
            \end{align*}
            for sets in the local algebra $K\in\mathfrak{A}=\bigcup_{B\in\BBB_0(X)}\mathfrak{A}_B$, where $$\mathfrak{A}_B=\{K\in\BBB(\Gamma_{X}):  \xi\in K\iff \xi+\delta_x\in K\quad \forall x\in X\setminus B\},$$ see \cite[Section 4.2.1]{schiavo2023wasserstein}. Note that if $K\in\mathfrak{A},$ then $\ind_K\in\FFF_b(\Gamma_{X})$ and $D\ind_K\in\FFF_{b,0}(\Gamma_{X}\times X).$ Hence, by definition $\ind_K\in\mathcal{H},$ see \eqref{eq:H}.  Since $\lim_{i\to\infty}\dP_0^i=\dP_0$ in $\PPP_1(\Gamma_X)$, Theorem~\ref{thm:propertiesce} $(iv)$ implies that $\lim_{j\to\infty}\bar{\dP}^{i_j}_t(K)= \bar{\dP}_t(K)$. Consequently, $\bar\dP_t$ is a non-negative set function on the algebra $\mathfrak{A}$ with total mass $1$. By Hahn's extension theorem, the measure can be extended to a probability measure $\bar{\dP}_t$ on $\sigma(\mathfrak{A})=\BBB(\Gamma_{X})$. Let $f\in\CCC_c^+(X).$ It holds that $e^{-\int f(y)\xi(\mathrm{d}y)}\in\mathcal{H}.$ Again, by Theorem~\ref{thm:propertiesce} $(iv)$, we conclude that 
            \begin{align*}
                \lim_{j\to\infty}\int e^{-\int f(y)\xi(\mathrm{d}y)} \bar\dP^{i_j}_t(\mathrm{d}\xi)&=\lim_{j\to\infty}\int e^{-\int f(y)\xi(\mathrm{d}y)} \bar\dP^{i_j}_0(\mathrm{d}\xi)+\int_0^t\int D_xe^{-\int f(y)\xi(\mathrm{d}y)} \bar{\dV}_s^{i_j}(\mathrm{d}\xi,\mathrm{d}x)\mathrm{d}s\\
                &= \int e^{-\int f(y)\xi(\mathrm{d}y)} \bar\dP_0(\mathrm{d}\xi)+\int_0^t\int D_xe^{-\int f(y)\xi(\mathrm{d}y)} \bar{\dV}_s(\mathrm{d}\xi,\mathrm{d}x)\mathrm{d}s\\
                &= \int e^{-\int f(y)\xi(\mathrm{d}y)} \bar\dP_t(\mathrm{d}\xi).
            \end{align*}
            A similar argument proves that $\lim_{j\to\infty}\I_{\bar{\dP}^{i_j}_t}=\I_{\bar{\dP}_t}$ vaguely. Thus,  $\lim_{j\to\infty}\bar\dP_t^{i_j}=\bar{\dP}_t$ in $\PPP_1(\Gamma_{X})$ by Proposition~\ref{prop:convergencep1}. Finally, Theorem~\ref{thm:propertiesce} $(iii)$ implies that $(\bar{\dP},\bar{\dV})\in\CCC\EEE^{X}(\dP_0,\dP_1)$. 
     \end{proof}
The following lemma shows that absolute continuity with respect to $\dR$ is preserved along solutions to the continuity equation with finite action. 
     \begin{lem}\label{lem:preservationabsolutecontinuity}
             Let $X$ be a bounded Polish space, $\lambda\in\MMM^+(X)$, and
             $\kappa:\Gamma_X\times X\to[0,\infty)$ such that
             $\dR\in\RRR_{\kappa,\lambda}(\Gamma_X)$ exists.
             Let $\dP_0,\dP_1\in\PPP_1(\Gamma_X)$ with $\dP_0\ll\dR$ or $\dP_1\ll\dR$ and
             $(\bar\dP,\bar\dV)\in\CCC\EEE^X(\dP_0,\dP_1)$ such that
             \begin{align*}
                 \int_0^1\mathcal{L}_{\kappa,\lambda}(\bar\dP_t,\bar\dV_t)\mathrm{d}t<\infty.
             \end{align*}
             Then $\bar\dP_t\ll\dR$ for all $t\in[0,1]$.
         \end{lem}
         \begin{proof}
             For $t\in[0,1]$, let
$\bar{\dP}_t=\bar{\dP}_t^{\mathrm{ac}}+\bar{\dP}_t^{\mathrm{s}}$
be the Lebesgue decomposition w.r.t.\ $\dR$, i.e.\
$\bar{\dP}_t^{\mathrm{ac}}\ll\dR$ and $\bar{\dP}_t^{\mathrm{s}}$ and $\dR$ are singular which will be denoted by  $\bar{\dP}_t^{\mathrm{s}}\perp\dR$. 

Choose a countable dense set $D\subset[0,1]$ and, for each $r\in D$, a set
$N_r\in\BBB(\Gamma_X)$ such that $\dR(N_r)=0$ and $\bar{\dP}_r^{\mathrm{s}}(N_r^c)=0$.
Set $N:=\bigcup_{r\in D}N_r$.
Then $\dR(N)=0$ and $\bar{\dP}_r(N)=\bar{\dP}_r^{\mathrm{ac}}(N)+\bar{\dP}_r^{\mathrm{s}}(N)=\bar{\dP}_r^{\mathrm{s}}(\Gamma_X)$ for all $r\in D$. Using an integrated version of the  continuity equation  \cite[Corollary 4.11]{schiavo2023wasserstein} yields
\begin{align}\label{eq:1895}
   |\bar{\dP}_t^{\mathrm{s}}(\Gamma_X)-\bar{\dP}_s^{\mathrm{s}}(\Gamma_X)|
    \leq \|\bar{\dP}_t-\bar{\dP}_s\|_{\mathrm{TV}}
    \leq 2\int_s^t |\bar{\dV}_r|(\Gamma_X\times X)\mathrm{d}r,
\end{align}
for $0\leq s\leq t\leq1$, with $\|\cdot\|_{\mathrm{TV}}$ denoting the classical total variation norm. Since $\int_0^1|\bar{\dV}_r|(\Gamma_X\times X)\mathrm{d}r<\infty$ by property $(iv)$ of Definition~\ref{dfn:ce}, \eqref{eq:1895} implies that both $t\mapsto\bar{\dP}_t(N)$ and $t\mapsto\bar{\dP}_t^{\mathrm{s}}(\Gamma_X)$ are continuous. Hence, their equality on $D$ extends to $[0,1]$ and we have that
\begin{align}\label{eq:commonnullsetac}
    \dR(N)=0\quad \text{ and }\quad
    \bar{\dP}_t^{\mathrm{s}}(N^c)=0\quad\text{ for all  }t\in[0,1].
\end{align}
Let
\begin{align*}
    A&:=\{(\xi,x):\xi\in N,\ \xi+\delta_x\notin N\},\\
    B&:=\{(\xi,x):\xi\notin N,\ \xi+\delta_x\in N\}.
\end{align*}
Since $D_x\ind_N(\xi)$ vanishes outside $A\cup B$, the integrated version of the continuity equation \cite[Corollary 4.11]{schiavo2023wasserstein} yields
\begin{align}\label{eq:massbalancecommonnullset}
    \bar\dP_t(N)-\dP_0(N)
    =\int_0^t\int_{\Gamma_X\times X}D_x\ind_N(\xi)\bar\dV_s(\mathrm{d}\xi,\mathrm{d}x)\mathrm{d}s=\int_0^t\int_{A\cup B}D_x\ind_N(\xi)\bar\dV_s(\mathrm{d}\xi,\mathrm{d}x)\mathrm{d}s.
\end{align}

If we can prove that $|\bar\dV_s|(A\cup B)=0$ for a.e.\ $s\in[0,1]$, then the claim follows from \eqref{eq:massbalancecommonnullset}, since
\begin{align*}
    \bar\dP_t^{\mathrm{s}}(\Gamma_X)=\bar\dP_t(N)=\dP_0(N)=0
\end{align*}
for all $t\in[0,1]$, where the last equality uses $\dP_0\ll\dR$ and $\dR(N)=0$.

We note that the same argument as in the proof of Lemma~\ref{lem:absolutecontinuityvelocity} yields that
\begin{align}\label{eq:flowacbothmeasures}
    \bar\dV_t\ll\kappa(\bar\dP_t\otimes\lambda),\qquad
    \bar\dV_t\ll C_{\bar\dP_t}
\end{align}
for a.e.\ $t\in[0,1]$. 
By the definition of the reduced Campbell measure and \eqref{eq:commonnullsetac},
$C_{\bar\dP_t^{\mathrm{s}}}(A)=0$, while $C_{\bar\dP_t^{\mathrm{ac}}}\ll C_{\dR}$.
The GNZ equation \eqref{eq:papangelou} gives
\begin{align*}
    C_{\dR}(A)\leq C_{\dR}(N\times X)
    =\int_N\int_X\kappa(\xi,x)\lambda(\mathrm{d}x)\dR(\mathrm{d}\xi)=0.
\end{align*}
Thus, $C_{\bar\dP_t}(A)=0$, and \eqref{eq:flowacbothmeasures} yields $|\bar\dV_t|(A)=0$ for a.e.\ $t$. Next, $\kappa(\bar\dP_t^{\mathrm{s}}\otimes\lambda)(B)=0$ by \eqref{eq:commonnullsetac}, whereas
$\kappa(\bar\dP_t^{\mathrm{ac}}\otimes\lambda)\ll\kappa(\dR\otimes\lambda)=C_{\dR}$.
By the definition of the reduced Campbell measure,
\begin{align*}
    C_{\dR}(B)
    \leq\int\int_X\ind_N(\xi)\xi(\mathrm{d}x)\dR(\mathrm{d}\xi)
    =\int_N\xi(X)\dR(\mathrm{d}\xi)=0.
\end{align*}
Hence, $\kappa(\bar\dP_t\otimes\lambda)(B)=0$, and \eqref{eq:flowacbothmeasures} yields $|\bar\dV_t|(B)=0$ for a.e.\ $t$. All in all, we conclude $|\bar\dV_t|(A\cup B)=0$, finishing the proof.
         \end{proof}
         \subsection{Properties of the Transport Distance}\label{subsec:propertiesdistance}
      In this subsection, we study the transport distance \(\mathcal{W}_{\kappa,\lambda}\) in more detail and establish its central properties. We first prove existence of minimizers for the optimization problem defining \(\mathcal{W}_{\kappa,\lambda}\). \begin{satz} \label{thm:minimizer}
 Let $\lambda\in\MMM^+(X)$ and assume that $\kappa:\Gamma_X\times X\to[0,\infty)$ satisfies \ref{property:locallystable} and~\ref{property:continuity}.     For all $\dP_0,\dP_1\in\PPP_1(\Gamma_X)$   there exists a solution $({\bar{\dP}},{\bar{\dV}})\in\CCC\EEE^X(\dP_0,\dP_1)$ s.t.\   
    \begin{align*}
      \mathcal{W}_{\kappa,\lambda}^2(\dP_0,\dP_1)=\mathcal{A}_{\kappa,\lambda}({\bar{\dP}},{\bar{\dV}}).
    \end{align*}
    \end{satz}
    \begin{proof}
        If $\mathcal W_{\kappa,\lambda}(\dP_0,\dP_1)=\infty$, the claim is immediate by choosing any solution in
$\CCC\EEE^X(\dP_0,\dP_1)$. If $\mathcal W_{\kappa,\lambda}(\dP_0,\dP_1)<\infty$, the claim follows from the direct method of calculus of variations and Lemma~\ref{lem:compactnessofsolutions}.
    \end{proof}
 By a standard reparametrization argument, see e.g.\ \cite[Theorem 5.4]{Dolbeault_2008}, \cite[Lemma 1.1.4]{Ambrosio2008GradientFI}, we obtain the following representation for $\mathcal{W}_{\kappa,\lambda}$.
      \begin{lem}
       Let $\lambda\in\MMM^+(X)$ and assume that $\kappa:\Gamma_X\times X\to[0,\infty)$ satisfies \ref{property:locallystable} and~\ref{property:continuity}. Let $\dP_0,\dP_1\in\PPP_1(\Gamma_X).$ It holds that 
        \begin{align}\label{eq:31}
            \mathcal{W}_{\kappa,\lambda}(\dP_0,\dP_1)=\inf\left\{\int_0^T\mathcal{L}_{\kappa,\lambda}^{1/2}(\bar{\dP}_{t},{\bar{\dV}_{t}})\mathrm{d}t:({\bar{\dP}},{\bar{\dV}})\in\CCC\EEE^X_T(\dP_0,\dP_1)\right\},
        \end{align}
        for all $T>0.$
        \end{lem}
        Using the previous representation, the following theorem establishes the central properties of $\mathcal{W}_{\kappa,\lambda}$. \begin{satz}\label{thm:propertiesdistance}
       Let $\lambda\in\MMM^+(X)$ and assume that $\kappa:\Gamma_X\times X\to[0,\infty)$ satisfies \ref{property:locallystable} and~\ref{property:continuity}.
             \begin{enumerate}
                \item[$(i)$] The map $\mathcal{W}_{\kappa,\lambda}:\PPP_1(\Gamma_X)\times\PPP_1(\Gamma_X)\to[0,\infty]$ is an extended distance on $\PPP_1(\Gamma_X).$
             
                 \item[${(ii)}$]  The topology induced by $\mathcal{W}_{\kappa,\lambda}$ on $\PPP_1(\Gamma_X)$ is stronger than the usual $\CCC_1$-weak topology  on $\PPP_1(\Gamma_X)$. 
                 \item[${(iii)}$] $\mathcal{W}_{\kappa,\lambda}$ is lower semicontinuous on $\PPP_1(\Gamma_X)\times \PPP_1(\Gamma_X)$.
                 \item[${(iv)}$] 
                 Bounded sets w.r.t.\ $\mathcal{W}_{\kappa,\lambda}$ are $\PPP_1(\Gamma_X)$-relatively compact.
                 \item[${(v)}$]  The accessible component $(\{\mathcal{W}_{\kappa,\lambda}(\dP, \cdot) < \infty\}, \mathcal{W}_{\kappa,\lambda}) $ of $\dP\in\PPP_1(\Gamma_X)$ is a complete geodesic space.
                 \item[${(vi)}$]  If $X$ is bounded and $\mathcal{W}_{\kappa,\lambda}(\dP,\dR)<\infty$, then $\dP\ll\dR$. 
            \end{enumerate}
         \end{satz}
         \begin{proof}
             The proofs of $(i)$ and $(ii)$ are identical to those of \cite[Theorem 5.15 $(i)$, $(ii)$]{schiavo2023wasserstein}.

             The statements in $(iii)$ and  $(iv)$ directly follow from Lemma~\ref{lem:compactnessofsolutions}.
             
               The completeness stated in $(v)$ follows from $(iii)$ and $(iv).$ The geodesics are given by the minimizers of Theorem~\ref{thm:minimizer}.

               $(vi)$ follows directly from Lemma \ref{lem:preservationabsolutecontinuity}.
         \end{proof}
          Having established the topological properties of $\mathcal W_{\kappa,\lambda}$, we next study the corresponding notion of absolutely continuous curves and their infinitesimal descriptions.
    \begin{dfn}[Absolutely continuous curves]\label{dfn:metricderivative}
        A curve $(\bar{\dP}_t)_{t\in[0,T]}\subset\PPP_1(\Gamma_X)$ is called \textit{absolutely continuous} w.r.t.\ $\mathcal{W}_{\kappa,\lambda}$ if there exists a function $g\in L^1(\Leb_{|(0,T)})$ s.t.\   \begin{align}\label{eq:abscts}
            \mathcal{W}_{\kappa,\lambda}(\bar{\dP}_s,\bar{\dP}_t)\leq \int_s^t g(u)\mathrm{d}u,\quad 0\leq s\leq t\leq T.
        \end{align}
        The minimal function $g$ satisfying \eqref{eq:abscts} is called metric derivative and is denoted by $|\dot{\bar{\dP}}_t|$.
    \end{dfn}
    \begin{bem}
        Assume that $\kappa:\Gamma_X\times X\to[0,\infty)$ satisfies \ref{property:locallystable} and~\ref{property:continuity}. By \cite[Theorem 1.1.2]{Ambrosio2008GradientFI}, it holds that 
        \begin{align*}
            |\dot{\bar{\dP}}_t|=\lim_{\varepsilon\to0}\frac{\mathcal{W}_{\kappa,\lambda}(\bar{\dP}_{t+\varepsilon},\bar{\dP}_t)}{\varepsilon}\quad \text{ for a.e.\ } t\in(0,T).
        \end{align*}
        
    \end{bem}
    \begin{prop}\label{prop:metricderivative}
      Let $\lambda\in\MMM^+(X)$ and assume that $\kappa:\Gamma_X\times X\to[0,\infty)$ satisfies \ref{property:locallystable} and~\ref{property:continuity}. A $\PPP_1(\Gamma_X)$-continuous curve $(\bar{\dP}_t)_{t\in[0,T]}\subset\PPP_1(\Gamma_X)$ is absolutely continuous w.r.t.\ $\mathcal{W}_{\kappa,\lambda}$ if and only if there exists a $(\bar{\dV}_t)_{t\in[0,T]}\subset\mathcal{M}_{b,0}(\Gamma_X\times X)$ s.t.\   $(\bar{\dP},\bar{\dV})\in\CCC\EEE_T^X(\bar{\dP}_0,\bar{\dP}_T)$ and \begin{align*}
            \int_0^T\sqrt{\mathcal{L}_{\kappa,\lambda}(\bar{\dP}_t,\bar{\dV}_t)}\mathrm{d}t<\infty.
        \end{align*}
        If this condition is satisfied, then $|\dot{\bar{\dP}}_t|^2\leq \mathcal{L}_{\kappa,\lambda}(\bar{\dP}_t,\bar{\dV}_t)$ for a.e.\ $t\in[0,T].$ Moreover, there exists an a.e.\ unique velocity field $(\bar{\dV}^\prime_t)_{t\in[0,T]}\subset\mathcal{M}_{b,0}(\Gamma_X\times X)$ s.t.\   $(\bar{\dP},\bar{\dV}^\prime)\in\CCC\EEE_T^X(\bar{\dP}_0,\bar{\dP}_T)$ and \begin{align}\label{eq:uniquevelocityfield}
            |\dot{\bar{\dP}}_t|^2= \mathcal{L}_{\kappa,\lambda}(\bar{\dP}_t,\bar{\dV}^\prime_t)\  \text{ for a.e.\ }t\in[0,T].\end{align}
    \end{prop}
    \begin{proof}
        The proof is analogous to the proof of \cite[Theorem 5.17]{Dolbeault_2008}. In our setting, \cite[Corollary 4.10]{Dolbeault_2008} corresponds to Lemma~\ref{lem:compactnessofsolutions} $(iii)$.
    \end{proof}
    The previous proposition shows that for every absolutely continuous curve there exists an a.e.\ unique velocity field s.t.\  the Lagrange functional realizes the metric derivative. We now describe these particular velocity fields in more detail.
   To this end, we call a velocity field $\dV\in\MMM_{b,0}(\Gamma_X\times X)$ \textit{divergence-free} if 
      $$ \int D_xF(\xi)\dV(\mathrm{d}\xi,\mathrm{d}x)=0, \quad \text{for all } F\in\mathcal{H}.$$
    The \textit{tangent space} of $\PPP_1(\Gamma_X)$ at the measure $\dP\in\PPP_1(\Gamma_X)$ is defined by
   \begin{align*}
       T_\dP^{\kappa,\lambda}\PPP_1(\Gamma_X):=\{\dV\in\MMM_{b,0}(\Gamma_X\times X):\mathcal{L}_{\kappa,\lambda}(\dP,\dV)\leq& \mathcal{L}_{\kappa,\lambda}(\dP,\dV+\dV^\prime)\\
       &\forall \text{ divergence-free } \dV^\prime\in\MMM_{b,0}(\Gamma_X\times X)\}.
       \end{align*}
       \begin{kor}
         Let $\lambda\in\MMM^+(X)$ and assume that $\kappa:\Gamma_X\times X\to[0,\infty)$ satisfies \ref{property:locallystable} and~\ref{property:continuity}. Let $(\bar\dP,\bar\dV)\in\CCC\EEE^X(\dP_0,\dP_1)$ s.t.\   $\bar\dP$ is absolutely continuous w.r.t.\ $\mathcal{W}_{\kappa,\lambda}$. Then $\bar\dV$ is the unique solution of \eqref{eq:uniquevelocityfield} if and only if $\bar{\dV}_t\in T_{\bar\dP_t}^{\kappa,\lambda}\PPP_1(\Gamma_X)$ for a.e.\ $t\in[0,1]$.
       \end{kor}
       Analogously to the Euclidean setting, see e.g.\ \cite[Section 8]{Ambrosio2008GradientFI}, and the Poisson setting $\kappa\equiv1$ \cite[Proposition 5.25]{schiavo2023wasserstein}, we derive an explicit representation of the tangent space.
       \begin{lem}
         Let $\lambda\in\MMM^+(X)$ and assume that $\kappa:\Gamma_X\times X\to[0,\infty)$ satisfies \ref{property:locallystable} and~\ref{property:continuity}, and let $\dR\in\RRR_{\kappa,\lambda}(\Gamma_X)$.  Let $\dP=\rho\dR\in\PPP_1(\Gamma_X)$. Then $\dV\in T_{\dP}^{\kappa,\lambda}\PPP_1(\Gamma_X)$ if and only if $\dV=w(\kappa(\dP\otimes\lambda))$ for some
$w$ belonging to the $L^2(\kappa(\dP\otimes\lambda))$-closure  of $\{DF:F\in\mathcal{H}\}.$ 
       \end{lem}
       \begin{proof}
           The proof is identical to the proof of \cite[Proposition 5.25]{schiavo2023wasserstein} and \cite[Proposition 4.8]{erbar2012gradientflowsentropyjump}.
       \end{proof}
 \section{Spatial Birth-and-Death Processes}\label{sec:spatialbirthanddeath}

Let  \(\lambda\in\mathcal M^+(X)\) and assume that
\(\kappa:\Gamma_X\times X\to[0,\infty)\) satisfies \ref{property:cocycle},
\ref{property:locallystable}, and~\ref{property:finiterange}. In this section, we study the \textit{spatial birth-and-death
process} with birth rate \(\kappa\) and death rate \(1\), whose formal generator is given by
\[
\mathfrak L_{\kappa,\lambda}F(\xi)
:=\int_X \kappa(\xi,x)\bigl(F(\xi+\delta_x)-F(\xi)\bigr)\lambda(\mathrm dx)
+\sum_{x\in\xi}\bigl(F(\xi-\delta_x)-F(\xi)\bigr),
\quad F:\Gamma_X\to\R .
\]

The following result establishes the existence of a Markov process solving the martingale problem of the generator $\mathfrak{L}_{\kappa,\lambda}$ acting on the class of \textit{cylinder functions}  \[
\mathcal{G}:=\left\{ F(\xi)=\phi\left({\int f_1(x)\xi(\mathrm{d}x)},\ldots,\int f_n(x)\xi(\mathrm{d}x)\right): f_1,\ldots, f_n\in\CCC_c(X), \phi\in\CCC_b(\R^n) \right\}.
\] 

     \begin{satz}[{\cite[Lemma 2.5, Theorem 2.13]{garcia2006spatialbirthdeathprocesses}}]\label{thm:SDESemigroup}
       Let $\lambda\in\MMM^+(X)$ and let $\kappa:\Gamma_X\times X\to[0,\infty)$ satisfy \ref{property:cocycle}, \ref{property:locallystable} and~\ref{property:finiterange}.  Let $\pi$ be a Poisson point process on $X\times [0,\infty)^3$ with intensity measure $\lambda\otimes \Leb_{|[0,\infty)}\otimes e^{-r}\Leb_{|[0,\infty)}\otimes \Leb_{|[0,\infty)}$, $\xi=\sum_{i\geq 1}\delta_{x_i}\in\Gamma_X$, and $\hat{\xi}=\sum_{i\geq1}\delta_{(x_i,\tau_i)}$ be the point process on $X\times [0,\infty)$ obtained by associating to each $x_i$ an independent unit-exponential random variable $\tau_i$. Suppose $\FFF=(\FFF_{t})_{t\geq 0}$ is a filtration s.t.\  $\pi$ is $\FFF$-compatible, i.e.\  s.t.\  for each $A\in\BBB_0(X\times [0,\infty)^2)$, $\pi(A,\cdot)$ is $\FFF$-adapted and $\pi(A,t+s)-\pi(A,t)$ is independent of $\FFF_t$ for any $s,t\geq 0$. Then there exists a unique \(\Gamma_X\)-valued càdlàg process
\(
    (Y_t^{\kappa,\lambda,(\xi)})_{t\geq0}
\)
adapted to $\FFF$ s.t.\  $Y_0^{\kappa,\lambda,(\xi)}=\xi$ and
\begin{align}
            Y_t^{\kappa,\lambda,(\xi)}(B)=\int_{B\times[0,\infty)^2\times [0,t]}\ind_{[0,\kappa(Y_{s-}^{\kappa,\lambda,(\xi)},x)]}(u)\ind_{(t-s,\infty)}(r)\pi(\mathrm{d}x,\mathrm{d}u,\mathrm{d}r,\mathrm{d}s)&+\int_{B\times [0,\infty)}\ind_{(t,\infty)}(r)\hat{\xi}(\mathrm{d}x,\mathrm{d}r),\label{eq:solutionsde}
        \end{align}
         $B\in\BBB(X)$, $t\geq0$. Moreover, $(Y_t^{\kappa,\lambda,(\xi)})_{t\geq0}$ is the unique solution to the martingale problem of the generator $\mathfrak{L}_{\kappa,\lambda}$ acting on $\mathcal{G}$, i.e.\  for all $F\in\mathcal{G}$, \begin{align}\label{eq:martingaleproblem}
           M_t^{\kappa,\lambda,(\xi)}(F):=F(Y_t^{\kappa,\lambda,(\xi)})-F(Y_0^{\kappa,\lambda,(\xi)})-\int_0^t\mathfrak{L}_{\kappa,\lambda}F(Y_s^{\kappa,\lambda,(\xi)})\mathrm{d}s,\quad t\geq0,
       \end{align}
       is a martingale w.r.t.\ the filtration $\FFF$.
    \end{satz}
    \subsection{Semigroup}\label{subsec:semigroup}
    In the following, we assume that the assumptions of Theorem~\ref{thm:SDESemigroup} are satisfied.
      The Markov semigroup associated with the unique solution  of the stochastic equation \eqref{eq:solutionsde} $(Y_t^{\kappa,\lambda,(\xi)})_{t\geq0}$ is defined by
         $$\mathrm{S}_t^{\kappa,\lambda} F(\xi):=\E[F(Y_t^{\kappa,\lambda,(\xi)})],\quad t\geq0,$$
        for all functions $F\in\FFF(\Gamma_X)$ s.t.\   the right-hand side is well-defined.
   The proof of \cite[Lemma 5.1]{jahnel2026reversiblebirthanddeathdynamicscontinuum}
carries over to the present setting and yields
\begin{align}\label{eq:derivativegenerator}
    \mathrm{S}_{t}^{\kappa,\lambda}F
    =
    F+\int_0^t \mathrm{S}_{s}^{\kappa,\lambda}
    \mathfrak{L}_{\kappa,\lambda}F\mathrm{d}s,
    \quad t\geq0,
\end{align}
for all $F\in\mathcal{G}$. The same applies to the reversibility argument in
\cite[Proposition 8.6]{jahnel2026reversiblebirthanddeathdynamicscontinuum},
which shows that every Papangelou measure
$\dR\in\RRR_{\kappa,\lambda}(\Gamma_X)$ is reversible with respect to
$(\mathrm{S}_t^{\kappa,\lambda})_{t\geq0}$, i.e.
\begin{align}\label{eq:Rinvariant}
    \int F(\xi)\mathrm{S}_t^{\kappa,\lambda}G(\xi)\dR(\mathrm{d}\xi)
    =
    \int \mathrm{S}_t^{\kappa,\lambda}F(\xi)G(\xi)\dR(\mathrm{d}\xi),
    \quad t\geq0,
\end{align}
for all $F,G\in\mathcal{G}$. By Jensen's inequality and
\eqref{eq:Rinvariant}, we obtain that
\begin{align}
   \|\mathrm{S}_t^{\kappa,\lambda}F\|_{L^2(\dR)}^2
   &\leq
   \int \mathrm{S}_t^{\kappa,\lambda}|F|^2(\xi)\dR(\mathrm{d}\xi)=
   \int |F(\xi)|^2\dR(\mathrm{d}\xi)
   =
   \|F\|_{L^2(\dR)}^2,
   \label{eq:1010}
\end{align}
for all $F\in\mathcal{G}$. Further, $\mathcal{G}$ is dense in
$L^2(\dR)$, see e.g.\ \cite[p. 690]{KoLy05}. Hence, we can uniquely extend $(\mathrm{S}_t^{\kappa,\lambda})_{t\geq0}$ to a semigroup on $L^2(\dR)$. In abuse of notation, we will also denote the extension by  
       $(\mathrm{S}_t^{\kappa,\lambda})_{t\geq0}$. 
    We begin our discussion of the semigroup by establishing several central properties, including contractivity, $C_0$-continuity (Lemma~\ref{lem:strongcontinuity}), and its long-term behavior (Corollary~\ref{kor:convergencesemigroup}).
    \begin{lem} \label{lem:strongcontinuity}
       Let $\lambda\in\MMM^+(X)$ and assume that $\kappa:\Gamma_X\times X\to[0,\infty)$ satisfies \ref{property:cocycle}, \ref{property:locallystable} and~\ref{property:finiterange}. Further, let $\dR \in \RRR_{\kappa,\lambda}(\Gamma_{X})$. The semigroup $(\mathrm{S}_t^{\kappa,\lambda})_{t\geq0}$ is contractive, i.e.\ 
        \begin{align*}
            \|\mathrm{S}_t^{\kappa,\lambda}F\|_{L^2(\dR)}\leq \|F\|_{L^2(\dR)}, \quad t\geq0,
        \end{align*}
       for all $F\in L^2(\dR)$, and  a $C_0$-semigroup, i.e.\  
        \begin{align*}
           \lim_{t\to t_0 } \|\mathrm{S}_t^{\kappa,\lambda}F-\mathrm{S}_{t_0}^{\kappa,\lambda}F\|_{L^2(\dR)}=0,
        \end{align*}
         for all $F\in L^2(\dR)$. Moreover, the measure $\dR$ is reversible for the semigroup, i.e.\  
         \begin{align*}
              \int F(\xi)\mathrm{S}_t^{\kappa,\lambda} G(\xi)\dR(\mathrm{d}\xi)=\int \mathrm{S}_t^{\kappa,\lambda} F(\xi) G(\xi)\dR(\mathrm{d}\xi),\quad t\geq0,
         \end{align*}
         for all $F,G\in L^2(\dR)$.
    \end{lem}
   \begin{proof}
  The contractivity follows from equation \eqref{eq:1010} and the construction of the semigroup.  For the proof of the $C_0$-continuity, we recall that $(Y_t^{\kappa,\lambda,(\xi)})_{t\geq0}$ is càdlàg. Hence, using the dominated convergence theorem, we obtain 
       \begin{align*}
           \lim_{t\to 0} \|\mathrm{S}_t^{\kappa,\lambda}F-F\|_{L^2(\dR)}^2=\lim_{t\to 0}\int \left(\mathrm{S}_t^{\kappa,\lambda}F(\xi)-F(\xi)\right)^2\dR(\mathrm{d}\xi)=0,
       \end{align*}
       for all $F\in\mathcal{G}$. Now, assume that $F\in L^2(\dR)$. Since $\mathcal{G}$ is dense in $L^2(\dR)$, we can find a sequence  of  functions $(F_i)_{i\in\N}\subset\mathcal{G}$ s.t.\   $\lim_{i\to\infty}F_i= F$ in $L^2(\dR)$.  Using the contractivity of $(\mathrm{S}_t^{\kappa,\lambda})_{t\geq0}$ and the previous result for $F_i\in\mathcal{G}$, we conclude that
      \begin{align*}
          \lim_{t\to 0}\|\mathrm{S}_{t}^{\kappa,\lambda} F- F\|_{L^2(\dR)}&\leq  \lim_{t\to 0} \|\mathrm{S}_{t}^{\kappa,\lambda} F- \mathrm{S}_{t}^{\kappa,\lambda} F_i\|_{L^2(\dR)}+ \| \mathrm{S}_{t}^{\kappa,\lambda} F_i-F_i\|_{L^2(\dR)}+\|  F_i- F\|_{L^2(\dR)}\\
          &\leq  \lim_{t\to 0}2\|F-F_i\|_{L^2(\dR)}+ \|\mathrm{S}_t^{\kappa,\lambda} F_i-F_i\|_{L^2(\dR)}\\
          &=2\|F-F_i\|_{L^2(\dR)},
      \end{align*}
      for all $i\in\N$. Finally,  since $\lim_{i\to\infty}F_i= F$ in $L^2(\dR),$ we obtain that 
      \begin{align*}
          \lim_{t\to 0}\|\mathrm{S}_{t}^{\kappa,\lambda} F- F\|_{L^2(\dR)}=0.
      \end{align*}
      Using the right-continuity in $t=0$  and  the contractivity, we obtain that 
        \begin{align*}
            \lim_{t\downarrow t_0} \|\mathrm{S}_t^{\kappa,\lambda}F-\mathrm{S}_{t_0}^{\kappa,\lambda}F\|_{L^2(\dR)}^2= \lim_{t\downarrow t_0} \|\mathrm{S}_{t_0}^{\kappa,\lambda}\mathrm{S}_{t-t_0}^{\kappa,\lambda}F-\mathrm{S}_{t_0}^{\kappa,\lambda}F\|_{L^2(\dR)}^2\leq \lim_{t\downarrow t_0} \|\mathrm{S}_{t-t_0}^{\kappa,\lambda}F-F\|_{L^2(\dR)}^2=0,
        \end{align*}
         for all $t_0>0$. Similarly, it holds that 
        \begin{align*}
            \lim_{t\uparrow t_0} \|\mathrm{S}_t^{\kappa,\lambda}F-\mathrm{S}_{t_0}^{\kappa,\lambda}F\|_{L^2(\dR)}^2= \lim_{t\uparrow t_0} \|\mathrm{S}_{t}^{\kappa,\lambda}F-\mathrm{S}_{t}^{\kappa,\lambda}\mathrm{S}_{t_0-t}^{\kappa,\lambda}F\|_{L^2(\dR)}^2\leq \lim_{t\uparrow t_0} \|F-\mathrm{S}_{t_0-t}^{\kappa,\lambda}F\|_{L^2(\dR)}^2=0,
        \end{align*}
        for all $t_0>0$, finishing the proof of the $C_0$-continuity. The reversibility of the Papangelou measures $\dR\in\RRR_{\kappa,\lambda}(\Gamma_X)$ for functions $F,G\in L^2(\dR)$ follows directly from equation \eqref{eq:Rinvariant}, the density of $\mathcal{G}$ in $L^2(\dR)$, and the contractivity of the semigroup.
   \end{proof}
    The following lemma establishes upper bounds for the first and second moments of the number of points of the process $(Y_t^{\kappa,\lambda,(\xi)})_{t\geq0}$.
    \begin{lem}\label{lem:intensitymeasuresemigroup}
       Let $\lambda\in\MMM^+(X)$ and assume that $\kappa:\Gamma_X\times X\to[0,\infty)$ satisfies \ref{property:cocycle}, \ref{property:locallystable} and~\ref{property:finiterange}.  For any  $B\in\BBB_0(X)$ and $t\geq 0$, it holds that \begin{align}\label{eq:intensity1}
           \mathrm{S}_t^{\kappa,\lambda} \xi(B)=\E[Y_t^{\kappa,\lambda,(\xi)}(B)]\leq e^{-t}\xi(B)+A\lambda(B)(1-e^{-t}),
           \end{align}
           and
           \begin{align} \label{eq:intensity2}
               & \mathrm{S}_t^{\kappa,\lambda} \xi(B)^2=\E[Y_t^{\kappa,\lambda,(\xi)}(B)^2]\leq e^{-2t}\xi(B)^2
+C_t^\prime(e^{-t}-e^{-2t})\xi(B)
+ C_t^{\prime\prime},
           \end{align}
           with  $$C_t^\prime= 2A\lambda(B)+1\text{ and }C_t^{\prime\prime}=((A\lambda(B))^2+(A\lambda(B))-(2(A\lambda(B))^2+A\lambda(B))e^{-t}+(A\lambda(B))^2e^{-2t}).$$
    \end{lem}
    
    \begin{proof}
        Let us first assume that $B\in\BBB_0(X)$ is an open set. Then, there exists a monotone increasing sequence of functions $(f_i)_{i\in\N}\subset\CCC_c(X)$ s.t.\   $\lim_{i\to\infty}f_i=\ind_B$. For $i\in\N$, let \begin{align*}
            F_i(\xi):=\int f_i(x)\xi(\mathrm{d}x)\wedge i,\quad \xi\in\Gamma_X,
        \end{align*}
        and note that $F_i\in\mathcal{G}.$ The process \begin{align*}
            M_t^{\kappa,\lambda,(\xi)}(F_i):= F_i(Y_t^{\kappa,\lambda,(\xi)})-F_i(Y_0^{\kappa,\lambda,(\xi)})-\int_0^t \mathfrak{L}_{\kappa,\lambda} F_i(Y_s^{\kappa,\lambda,(\xi)})\mathrm{d}s
        \end{align*}
        is a martingale w.r.t.\ the filtration $\FFF=(\FFF_t)_{t\geq0}$ by Theorem~\ref{thm:SDESemigroup}. For $i>1$, we define the stopping time $$\tau_i:=\inf\{t\geq0:Y_t^{\kappa,\lambda,(\xi)}(B)\geq i-1\}$$ and note that the stopped process $(M_{t\wedge\tau_i}^{\kappa,\lambda,(\xi)}(F_i))_{t\geq0}$ is also a $\FFF$-martingale.
        Hence, 
        \begin{align}\label{eq:1050} \E[F_i(Y_{t\wedge\tau_i}^{\kappa,\lambda,(\xi)})]-F_i(Y_0^{\kappa,\lambda,(\xi)})=\E\left[\int_0^{t\wedge\tau_i} \mathfrak{L}_{\kappa,\lambda} F_i(Y_s^{\kappa,\lambda,(\xi)})\mathrm{d}s\right],
        \end{align}
        for all $t\geq0$.
        If $\eta\in\Gamma_X$ satisfies  $\eta(B)\leq i-1$, we have that 
        \begin{align*}
            F_i(\eta-\delta_x)-F_i(\eta)=-f_i(x)
        \end{align*}
        and 
        \begin{align*}
            F_i(\eta+\delta_x)-F_i(\eta)=f_i(x)\leq \ind_B(x).
        \end{align*}
        Thus, using \ref{property:locallystable}, 
        \begin{align*}
            \mathfrak{L}_{\kappa,\lambda} F_i(Y_t^{\kappa,\lambda,(\xi)})&\leq -\sum_{x\in Y_t^{\kappa,\lambda,(\xi)} }f_i(x)+\int \ind_B(x)\kappa(Y_t^{\kappa,\lambda,(\xi)},x)\lambda(\mathrm{d}x)\leq -\sum_{x\in Y_t^{\kappa,\lambda,(\xi)} }f_i(x)+A\lambda(B).
        \end{align*}
        Inserting the previous inequality into equation~\eqref{eq:1050}, taking the limit $i\to\infty$, and using the monotone convergence theorem yields
        \begin{align*}
            \E[Y_{t}^{\kappa,\lambda,(\xi)}(B)]-Y_0^{\kappa,\lambda,(\xi)}(B)\leq\E\left[\int_0^{t} -Y_s^{\kappa,\lambda,(\xi)}(B)+A\lambda(B)\mathrm{d}s \right].
        \end{align*}
        Finally, by Grönwall's inequality, we conclude that 
        \begin{align*}
            \mathrm{S}_t^{\kappa,\lambda} \xi(B)=\E[Y_t^{\kappa,\lambda,(\xi)}(B)]&\leq e^{-t}Y_0^{\kappa,\lambda,(\xi)}(B)+A\lambda(B)(1-e^{-t})\\
            &= e^{-t}\xi(B)+A\lambda(B)(1-e^{-t}),
        \end{align*}
        finishing the proof of \eqref{eq:intensity1} under the additional assumption that $B\in\BBB_0(X)$ is an open set. Now, let $B\in\BBB_0(X)$ be an arbitrary set. For any open set $O\in\BBB_0(X)$ s.t.\   $B\subseteq O$, it holds that 
        \begin{align*}
          \mathrm{S}_t^{\kappa,\lambda} \xi(B)\leq  \mathrm{S}_t^{\kappa,\lambda} \xi(O)\leq e^{-t}\xi(O)+A\lambda(O)(1-e^{-t}).
        \end{align*}
        As $\xi$ and $\lambda$ are outer regular measures, this implies 
        \begin{align*}
            \mathrm{S}_t^{\kappa,\lambda} \xi(B)\leq \inf\{ e^{-t}\xi(O)+A\lambda(O)(1-e^{-t}):B\subset O \text{ open}\} = e^{-t}\xi(B)+A\lambda(B)(1-e^{-t}),
        \end{align*}
        finishing the proof of the first claim. We continue with the proof of \eqref{eq:intensity2}. By a martingale argument analogous to the one used in the first part of the proof, we obtain that
        \begin{align*}
             \E[Y_{t}^{\kappa,\lambda,(\xi)}(B)^2]-Y_0^{\kappa,\lambda,(\xi)}(B)^2\leq\E\left[\int_0^{t} -2Y_s^{\kappa,\lambda,(\xi)}(B)^2+(2A\lambda(B)+1)Y_s^{\kappa,\lambda,(\xi)}(B)+A\lambda(B)\mathrm{d}s\right].
        \end{align*}
       Inserting inequality \eqref{eq:intensity1} for $\E[Y_s^{\kappa,\lambda,(\xi)}(B)]$ and using again Grönwall's inequality, we conclude similarly as before.
    \end{proof}
    If $X$ is a bounded Polish metric space, then the long-term behavior of the process $(Y_t^{\kappa,\lambda,(\xi)})_{t \ge 0}$ can be analyzed using the following results.

\begin{satz}[{\cite[Theorem 2.1]{Moller},\cite[Proposition 7.1, Theorem 7.1]{Preston}}] \label{thm:convergencesemigroup}
 Let $(Y_t)_{t\ge 0}$ be a  spatial birth-and-death process on a bounded Polish metric space $X$ with total birth rate  $\beta^{\mathrm{tot}}$ and
 total death rate $\Delta^{\mathrm{tot}}$.  For each $i\ge 0$, define 
\begin{align*}
\beta_i^{\mathrm{tot}} := \sup_{\xi\in\Gamma_X:\xi(X)=i}\beta^{\mathrm{tot}}(\xi),
\quad \text{ and }\quad 
\Delta_i^{\mathrm{tot}} := \inf_{\xi\in\Gamma_X:\xi(X)=i}\Delta^{\mathrm{tot}}(\xi),
\end{align*}
with the convention $\Delta_0^{\mathrm{tot}}=0$. Let $Q_t^{(\xi)}$ denote the law of $Y_t$
given $Y_0=\xi$. Assume that $\beta_i^{\mathrm{tot}}>0  \text{ and } \Delta_i^{\mathrm{tot}}>0
 \text{ for all } i\ge 1$,
and that the following conditions hold:
\begin{equation}\label{eq:condition1}
\sum_{i=2}^{\infty}
\frac{\beta_1^{\mathrm{tot}}\beta_2^{\mathrm{tot}}\cdots\beta_{i-1}^{\mathrm{tot}}}
     {\Delta_1^{\mathrm{tot}}\Delta_2^{\mathrm{tot}}\cdots\Delta_i^{\mathrm{tot}}}
<\infty\quad \text{ and } \quad \sum_{i=1}^{\infty}
\frac{\Delta_1^{\mathrm{tot}}\Delta_2^{\mathrm{tot}}\cdots\Delta_{i-1}^{\mathrm{tot}}}
     {\beta_1^{\mathrm{tot}}\beta_2^{\mathrm{tot}}\cdots\beta_i^{\mathrm{tot}}}
=\infty.
\end{equation}

There exists a {unique invariant probability measure}
$\dP\in\PPP_1(\Gamma_X)$ s.t.\ , for every $\xi\in\Gamma_X$ and every
 $K\in\BBB(\Gamma_X)$,
$$\lim_{t\to\infty} Q_t^{(\xi)}(K)=\dP(K).$$

\end{satz}
    \begin{kor}\label{kor:convergencesemigroup}
      Let $X$ be a bounded Polish metric space and $\lambda\in\MMM^+(X).$ Assume that $\kappa:\Gamma_X\times X\to[0,\infty)$ satisfies \ref{property:cocycle}  and  \ref{property:locallystable}. Further, let $\dR = \RRR_{\kappa,\lambda}(\Gamma_{X})$ and  $Q_t^{\kappa,\lambda,(\xi)}$ be the distribution of  $Y_t^{\kappa,\lambda,(\xi)}$. It holds that 
        \begin{align*}
            \lim_{t\to\infty}Q_t^{\kappa,\lambda,(\xi)}(K)=\dR(K),
        \end{align*}
        for all $K\in\BBB(\Gamma_X)$ and $\xi\in\Gamma_X$. Further, $\dR$ is the unique invariant measure of the process $(Y_t^{\kappa,\lambda,(\xi)})_{t\geq 0}$.
    \end{kor}
    \begin{proof}
        The total birth rate of $(Y_t^{\kappa,\lambda,(\xi)})_{t\geq 0}$ is given by $\beta^{\mathrm{tot}}(\xi)=\int\kappa(\xi,x)\lambda(\mathrm{d}x)$ and the total death rate by $\Delta^{\mathrm{tot}}(\xi)=\xi(X)$. Hence, by \ref{property:locallystable}, we obtain that $$\beta_i^{\mathrm{tot}}=\sup_{\xi\in\Gamma_X:\xi(X)=i}\beta^{\mathrm{tot}}(\xi)\leq A\lambda(X) \quad \text{ and }\quad \Delta_i^{\mathrm{tot}}= \inf_{\xi\in\Gamma_X:\xi(X)=i}\Delta^{\mathrm{tot}}(\xi)=i.$$
        Hence, \eqref{eq:condition1} follow by
\begin{align*}
\sum_{i=2}^\infty
\frac{\beta_1^{\mathrm{tot}}\cdots\beta_{i-1}^{\mathrm{tot}}}
{\Delta_1^{\mathrm{tot}}\cdots\Delta_i^{\mathrm{tot}}}
\le
\sum_{i=2}^\infty
\frac{(A\lambda(X))^{i-1}}{i!}
< \infty \quad \text{ and }\quad  \sum_{i=1}^\infty
\frac{\Delta_1^{\mathrm{tot}}\cdots\Delta_{i-1}^{\mathrm{tot}}}
{\beta_1^{\mathrm{tot}}\cdots\beta_{i}^{\mathrm{tot}}}
\ge
\sum_{i=1}^\infty
\frac{(i-1)!}{(A\lambda(X))^i}
= \infty.
\end{align*}

The claim follows directly from Theorem~\ref{thm:convergencesemigroup}.\end{proof}
    We finish this subsection by an integration by parts formula for the generator $\mathfrak{L}_{\kappa,\lambda}$.
     \begin{satz}\label{thm:ibp}
 Let  $\lambda\in\MMM^+(X)$ and assume that $\kappa:\Gamma_X\times X\to[0,\infty)$ satisfies \ref{property:cocycle}, \ref{property:locallystable}, and~\ref{property:finiterange}. Further, let $\dR \in \RRR_{\kappa,\lambda}(\Gamma_{X})$ and assume that one of the following assumptions is satisfied: 
 \begin{enumerate}
     \item[$(i)$] $F\in L^2(\dR)$ and $G\in \mathcal{G}$;
     \item[$(ii)$] $F\in L^1(\dR)$ s.t.\   $\dP=F\dR\in\PPP_1(\Gamma_X)$ and $G\in\mathcal{G}$;
     \item[$(iii)$] $F,G\in\mathcal{H}$ bounded local functions.
 \end{enumerate}
 It holds that 
        \begin{align}\label{eq:ibp}
            \int F(\xi)\mathfrak{L}_{\kappa,\lambda} G(\xi)\dR(\mathrm{d}\xi)=-\int\int D_xF(\xi)D_xG(\xi)\kappa(\xi,x)\dR(\mathrm{d}\xi)\lambda(\mathrm{d}x).
        \end{align}

     Moreover, if $F,G\in\mathcal{G}$ or $F,G\in\mathcal{H}$, then
         \begin{align}\label{eq:ibp2}
            \int F(\xi)\mathfrak{L}_{\kappa,\lambda} G(\xi)\dR(\mathrm{d}\xi)= \int \mathfrak{L}_{\kappa,\lambda} F(\xi) G(\xi)\dR(\mathrm{d}\xi).
        \end{align}
    \end{satz}
    \begin{proof}
    We first verify that under each of the assumptions $(i)$--$(iii)$, it holds that ${F\mathfrak{L}_{\kappa,\lambda}G\in L^1(\dR)}$.
    
    $(i)$ If $G\in \mathcal{G}$, then $\mathfrak{L}_{\kappa,\lambda} G\in L^2(\dR)$ by \cite[Theorem 3.2]{li2010noteequilibriumglauberkawasaki}. Since $F\in L^2(\dR)$, we obtain that $F\mathfrak{L}_{\kappa,\lambda}G\in L^1(\dR)$. 
    
    $(ii)$ For  $G\in\mathcal{G}$, there exists a constant $C<\infty$ and a set $B\in\mathcal{B}_0(X)$ s.t.\  
    \begin{align*}
        |G(\xi+\delta_x)-G(\xi)|\leq C\ind_B(x)\quad \text{ and }\quad |G(\xi-\delta_x)-G(\xi)|\leq C\ind_B(x).
    \end{align*}
     Thus, using \ref{property:locallystable} and $\dP=F\dR\in\PPP_1(\Gamma_X)$, we obtain 
     \begin{align*}
        \int |F(\xi)\mathfrak{L}_{\kappa,\lambda} G(\xi)|\dR(\mathrm{d}\xi)\leq \int (CA\lambda(B)+C\xi(B))  F(\xi)\dR(\mathrm{d}\xi)=CA\lambda(B)+C\int \xi(B)\dP(\mathrm{d}\xi) <\infty.
     \end{align*}
     $(iii)$ If $F,G \in \mathcal{H}$, there exist again a constant $C<\infty$ and a set $B\in\BBB_0(X)$ s.t.\ 
     \begin{align*}
         \int |F(\xi)\mathfrak{L}_{\kappa,\lambda} G(\xi)|\dR(\mathrm{d}\xi)\leq \int C(2CA\lambda(B)+2C\xi(B))\dR(\mathrm{d}\xi)<\infty.
     \end{align*}
     
     In each of the three cases $(i)$--$(iii)$, we obtain by  the GNZ equation~\eqref{eq:papangelou}, 
        \begin{align*}
              &\int F(\xi)\mathfrak{L}_{\kappa,\lambda} G(\xi)\dR(\mathrm{d}\xi)\\&\quad=-\int F(\xi)\int D_xG(\xi-\delta_x)\xi(\mathrm{d}x)\dR(\mathrm{d}\xi)+\int F(\xi)\int D_xG(\xi)\kappa(\xi,x)\dR(\mathrm{d}\xi)\lambda(\mathrm{d}x)\\
              &\quad= -\int\int \kappa(\xi,x) F(\xi+\delta_x) D_xG(\xi)\dR(\mathrm{d}\xi)\lambda(\mathrm{d}x)+\int\int F(\xi) \kappa(\xi,x)D_xG(\xi)\dR(\mathrm{d}\xi)\lambda(\mathrm{d}x)\\
              &\quad=-\int\int D_x F(\xi) D_x G(\xi)\kappa(\xi,x)\dR(\mathrm{d}\xi) \lambda(\mathrm{d}x),
        \end{align*}
        proving \eqref{eq:ibp}. Further, \eqref{eq:ibp2} follows directly from \eqref{eq:ibp}.
    \end{proof}
    \subsection{Dual Semigroup}\label{subsec:dualsemigroup}
        We define the \textit{dual semigroup} of $(\mathrm{S}_t^{\kappa,\lambda})_{t\geq0}$ acting on $\PPP_1(\Gamma_X)$ by
        \begin{align*} \mathrm{S}_t^{\kappa,\lambda*}\dP(K)=\int\mathrm{S}_t^{\kappa,\lambda}\ind_K(\xi)\dP(\mathrm{d}\xi),\quad t\geq0,\; K\in\BBB(\Gamma_X).
        \end{align*}
        \begin{satz}\label{thm:dualsemigroup}
      Let $\lambda\in\MMM^+(X)$ and assume that $\kappa:\Gamma_X\times X\to[0,\infty)$ satisfies \ref{property:cocycle}, \ref{property:locallystable} and~\ref{property:finiterange}. The dual semigroup $(\mathrm{S}_t^{\kappa,\lambda*})_{t\geq0}$ is a $C_0$-semigroup in the sense that for all $\dP\in\PPP_1(\Gamma_X)$ it holds that:
        \begin{itemize}
            \item[$(i)$] $\mathrm{S}_t^{\kappa,\lambda*}\dP\in\PPP_1(\Gamma_X)$ for all $t\geq0$;
            \item [$(ii)$] $\mathrm{S}_0^{\kappa,\lambda*}\dP=\dP$;
            \item [$(iii)$] $\mathrm{S}_{s+t}^{\kappa,\lambda*}\dP=\mathrm{S}_s^{\kappa,\lambda*}\mathrm{S}_t^{\kappa,\lambda*}\dP$ for all $s,t\geq0$;
            \item [$(iv)$] $\lim_{t\to 0}\mathrm{S}_{t}^{\kappa,\lambda*}\dP=\dP$ in $\PPP_1(\Gamma_X)$.
        \end{itemize}
    \end{satz}
      \begin{proof}
      We start with the proof of property $(i)$. Let $\dP\in\PPP_1(\Gamma_X)$ and $t\geq0$.  By the definition of the dual semigroup, we directly deduce that $\mathrm{S}_{t}^{\kappa,\lambda*}\dP(\Gamma_X)=\int\mathrm{S}_{t}^{\kappa,\lambda}1\dP(\mathrm{d}\xi)=1$, which implies that $\mathrm{S}_{t}^{\kappa,\lambda*}\dP\in\PPP(\Gamma_X)$. It remains to show that $ \I_{\mathrm{S}_{t}^{\kappa,\lambda*}\dP}(B)<\infty$ for all $B\in\BBB_0(X)$. To this end, let 
      $B\in\BBB_0(X)$. By Lemma~\ref{lem:intensitymeasuresemigroup}, 
       \begin{align*}
           \I_{\mathrm{S}_{t}^{\kappa,\lambda*}\dP}(B)&= \int \mathrm{S}_{t}^{\kappa,\lambda}\xi(B) \dP(\mathrm{d}\xi)
           \leq e^{-t}\int \xi(B)\dP(\mathrm{d}\xi)+(1-e^{-t})A\lambda(B)=e^{-t}\I_\dP(B)+(1-e^{-t})A\lambda(B)<\infty,
       \end{align*}
       proving that $\mathrm{S}_{t}^{\kappa,\lambda*}\dP\in\PPP_1(\Gamma_X)$.\\
        Since $(\mathrm{S}_t^{\kappa,\lambda})_{t\geq0}$ is  the Markov semigroup associated with the process $(Y_t^{\kappa,\lambda,(\xi)})_{t\geq0}$, properties $(ii)$ and $(iii)$ follow from the definition of $(\mathrm{S}_t^{\kappa,\lambda*})_{t\geq0}$.\\
        It remains to prove  property $(iv)$. Let $\dP\in\PPP_1(\Gamma_X)$. By Proposition~\ref{prop:convergencep1} it suffices to prove that $\lim_{t\to 0}L_{\mathrm{S}_t^{\kappa,\lambda*}\dP}(f)=L_{\dP}(f)$ for all $f\in\CCC_c^+(X)$ and $\lim_{t\to0}\I_{\mathrm{S}_t^{\kappa,\lambda*}\dP}=\I_{\dP}$ vaguely.  Thus, let $f\in\CCC_c^+(X)$ and recall that $(Y_t^{\kappa,\lambda,(\xi)})_{t\geq 0}$ is càdlàg. Hence, by  the dominated convergence theorem
        \begin{align*}
            \lim_{t\to 0} L_{\mathrm{S}_t^{\kappa,\lambda*}\dP}(f)&=\lim_{t\to 0}\int \E[e^{-\int f(x)Y_t^{\kappa,\lambda,(\xi)}(\mathrm{d}x)}]\dP(\mathrm{d}\xi)=\int \E[e^{-\int f(x)\xi(\mathrm{d}x)}]\dP(\mathrm{d}\xi)= L_{\dP}(f).
        \end{align*}
        Further, let $B\in\BBB_0(X)$ be an open set. We can find an increasing sequence of compact sets $K_i\in\BBB_c(X)$ s.t.\  $\bigcup_{i\in\N}K_i=B$ and a monotone increasing sequence of functions $(f_i)_{i\in\N}\subset\CCC_c(X)$ s.t.\   $0\leq f_i\leq \ind_B$ and $f_i(x)=1$ for $x\in K_i$. Let \begin{align*}
            F_i(\xi):=\int f_i(x)\xi(\mathrm{d}x),\quad i\in\N.
        \end{align*}
        Using \eqref{eq:intensity1} of Lemma~\ref{lem:intensitymeasuresemigroup}, we note that 
        \begin{align*}
            \E[ F_i(Y_t^{\kappa,\lambda,(\xi)})]\leq\E[ Y_t^{\kappa,\lambda,(\xi)}(B)]\leq e^{-t}\xi(B)+A\lambda(B)(1-e^{-t})\leq\xi(B)+A\lambda(B).
        \end{align*}
        Since $\dP\in\PPP_1(\Gamma_X)$, the right-hand side is $\dP$-integrable. Hence, using that  $(Y_t^{\kappa,\lambda,(\xi)})_{t\geq 0}$ is càdlàg and the dominated convergence theorem, we obtain
        \begin{align*}
          \lim_{t\to 0} \int F_i(\xi) \mathrm{S}_{t}^{\kappa,\lambda*}\dP(\mathrm{d}\xi)=   \lim_{t\to 0}\int \E[F_i(Y_t^{\kappa,\lambda,(\xi)})]\dP(\mathrm{d}\xi)= \int \E[F_i(Y_0^{\kappa,\lambda,(\xi)})]\dP(\mathrm{d}\xi)=\int F_i(\xi) \dP(\mathrm{d}\xi),
        \end{align*}
        for every $i\in\N$. Thus, 
        \begin{align*}
            \lim_{i\to\infty}   \lim_{t\to 0} \int F_i(\xi)\mathrm{S}_{t}^{\kappa,\lambda*}\dP(\mathrm{d}\xi)=\lim_{i\to\infty} \int F_i(\xi) \dP(\mathrm{d}\xi).
        \end{align*}
        If $\int F_i(\xi)\mathrm{S}_{t}^{\kappa,\lambda*}\dP(\mathrm{d}\xi)$ converges uniformly in $t$ for $i\to\infty$, then the two limits of the left-hand side may be interchanged and we obtain by applying the monotone convergence theorem on both sides that 
        \begin{align*}
             \lim_{t\to 0}\I_{\mathrm{S}_{t}^{\kappa,\lambda*}\dP}(B)=\lim_{t\to 0} \int \xi(B) \mathrm{S}_{t}^{\kappa,\lambda*}\dP(\mathrm{d}\xi)= \int \xi(B) \dP(\mathrm{d}\xi)=\I_{\dP}(B).
        \end{align*}
        It remains to show that $\int F_i(\xi) \mathrm{S}_{t}^{\kappa,\lambda*}\dP(\mathrm{d}\xi)$ converges uniformly in $t$ for $i\to\infty$. To this end, we deduce from the definition of $f_i$ and \eqref{eq:intensity1} of Lemma~\ref{lem:intensitymeasuresemigroup}  that
        \begin{align*}
            \lim_{i\to\infty}\sup_{t\geq0}|\int \xi(B) \mathrm{S}_{t}^{\kappa,\lambda*}\dP(\mathrm{d}\xi)-\int F_i(\xi) \mathrm{S}_{t}^{\kappa,\lambda*}\dP(\mathrm{d}\xi)|
            &= \lim_{i\to\infty}\sup_{t\geq0}\int \E[Y_t^{\kappa,\lambda,(\xi)}(B)-F_i(Y_t^{\kappa,\lambda,(\xi)})] \dP(\mathrm{d}\xi)\\
            &\leq \lim_{i\to\infty}\sup_{t\geq0}\int\E[Y_t^{\kappa,\lambda,(\xi)}(B\setminus K_i)]\dP(\mathrm{d}\xi)\\
          &\leq \lim_{i\to\infty}\sup_{t\geq0} e^{-t}\I_\dP(B\setminus K_i)+A\lambda(B\setminus K_i)(1-e^{-t})\\
          &=0,
            \end{align*}
       finishing the proof for open sets $B\in\BBB_0(X)$. Using the additivity of 
$\I_{\mathrm{S}_{t}^{\kappa,\lambda*}\dP}$ and $\I_{\dP}$ on disjoint sets,
the result extends immediately to all bounded closed sets.
The general case then follows from the regularity of the intensity measures by approximation of bounded Borel sets by open and closed sets.
       \end{proof}
    \begin{lem}\label{lem:solutionsemigroup}
    Let $\lambda\in\MMM^+(X)$ and assume that $\kappa:\Gamma_X\times X\to[0,\infty)$ satisfies \ref{property:cocycle}, \ref{property:locallystable}, and~\ref{property:finiterange}. Further, let $\dR \in \RRR_{\kappa,\lambda}(\Gamma_{X})$,  $\dP=\rho\dR\in\PPP_1(\Gamma_X)$, and $T>0$. Then $$(\mathrm{S}_t^{\kappa,\lambda*}\dP,-\kappa D\mathrm{S}_t^{\kappa,\lambda}\rho(\dR\otimes \lambda))_{t\in[0,T]}\in\CCC\EEE_T^X(\dP,\mathrm{S}_T^{\kappa,\lambda*}\dP).$$
    \end{lem}
    \begin{proof}
        The proof follows the same line of argumentation as the proof of \cite[Proposition 4.3]{schiavo2023wasserstein}. Due to the additional Papangelou intensity, we give a full proof here.  Property $(i)$ of the definition of the continuity equation, Definition~\ref{dfn:ce}, is obviously satisfied.
        
        We continue with the proof of property $(ii)$, the $\PPP_1(\Gamma_X)$-continuity of $(\mathrm{S}_t^{\kappa,\lambda*}\dP)_{t\in[0,T]}$. The right-continuity of the curve $(\mathrm{S}_t^{\kappa,\lambda*}\dP)_{t\in[0,T]}$ follows from Theorem~\ref{thm:dualsemigroup} $(iii)$ and $(iv)$. The proof of the left-continuity is analogous to the proof of Theorem~\ref{thm:dualsemigroup} $(iv)$ using Lemma~\ref{lem:strongcontinuity}. 
        
   For  the proof of property $(iii)$ of Definition~\ref{dfn:ce}, let $\phi\in\CCC_c^\infty((0,T))$ and $f\in\CCC_c^+(X).$ Then  $F(\xi)=e^{-\int f(x)\xi(\mathrm{d}x)}\in\mathcal{G}$. 
   Using integration by parts in $t$, \eqref{eq:derivativegenerator}, the reversibility  of $\dR$ \eqref{eq:Rinvariant}, and finally the integration by parts formula \eqref{eq:ibp}, we obtain
    \begin{align*}
        \int_0^T\dt\phi(t)\int e^{-\int f(x)\xi(\mathrm{d}x)}\mathrm{S}_t^{\kappa,\lambda*}\dP(\mathrm{d}\xi)\mathrm{d}t&=  \int_0^T\dt\phi(t)\int \mathrm{S}_t^{\kappa,\lambda} e^{-\int f(x)\xi(\mathrm{d}x)}\dP(\mathrm{d}\xi)\mathrm{d}t\\
        &= - \int_0^T\phi(t)\int\dt \mathrm{S}_t^{\kappa,\lambda} e^{-\int f(x)\xi(\mathrm{d}x)}\dP(\mathrm{d}\xi)\mathrm{d}t\\
        &\overset{\eqref{eq:derivativegenerator}}{=}  -\int_0^T\phi(t)\int  \mathrm{S}_t^{\kappa,\lambda} \mathfrak{L}_{\kappa,\lambda} e^{-\int f(x)\xi(\mathrm{d}x)}\dP(\mathrm{d}\xi)\mathrm{d}t\\
        & \overset{\eqref{eq:Rinvariant}}{=} -\int_0^T\phi(t)\int   \mathfrak{L}_{\kappa,\lambda} e^{-\int f(x)\xi(\mathrm{d}x)} \mathrm{S}_t^{\kappa,\lambda}\rho(\xi)\dR(\mathrm{d}\xi)\mathrm{d}t\\ 
        &\overset{\eqref{eq:ibp}}{=} \int_0^T\phi(t)\int   D_x  e^{-\int f(x)\xi(\mathrm{d}x)} D_x\mathrm{S}_t^{\kappa,\lambda}\rho(\xi)\kappa(\xi,x)\dR(\mathrm{d}\xi)\lambda(\mathrm{d}x)\mathrm{d}t,
    \end{align*}
    proving property $(iii)$. 
    
    For the proof of property $(iv)$ of Definition~\ref{dfn:ce}, we note that by the GNZ equation~\eqref{eq:papangelou}
        \begin{align*}
            \int_B\int \kappa(\xi,x)\mathrm{S}_t^{\kappa,\lambda}\rho(\xi+\delta_x)\dR(\mathrm{d}\xi)\lambda (\mathrm{d}x)&=\int\int_B \mathrm{S}_t^{\kappa,\lambda}\rho(\xi)\xi(\mathrm{d}x)\dR(\mathrm{d}\xi)=\int \xi(B) \mathrm{S}_t^{\kappa,\lambda}\rho(\xi)\dR(\mathrm{d}\xi)=\I_{\mathrm{S}_t^{\kappa,\lambda*}\dP}(B).
        \end{align*}
        Therefore, we get by \ref{property:locallystable},  Lemma~\ref{lem:intensitymeasuresemigroup}, and the reversibility of $\dR$ \eqref{eq:Rinvariant} that
        \begin{align*}
            \int_0^T\!\int_B\!\int\!| \kappa(\xi,x)D_x\mathrm{S}_t^{\kappa,\lambda}\rho(\xi)| \dR(\mathrm{d}\xi) \lambda(\mathrm{d}x)\mathrm{d}t&\leq  \int_0^T\!\int_B\!\int\! \kappa(\xi,x)\mathrm{S}_t^{\kappa,\lambda}\rho(\xi+\delta_x)+  \kappa(\xi,x)\mathrm{S}_t^{\kappa,\lambda}\rho(\xi) \dR(\mathrm{d}\xi) \lambda(\mathrm{d}x)\mathrm{d}t\\
            &\leq\int_0^T \I_{\mathrm{S}_t^{\kappa,\lambda*}\dP}(B) \mathrm{d}t+TA\int_B\int \rho(\xi)\dR(\mathrm{d}\xi)\lambda(\mathrm{d}x)\\
            &\leq \int_0^T e^{-t}\I_{\dP}(B)+ (1-e^{-t})A\lambda(B) \mathrm{d}t+TA \lambda(B)<\infty,
        \end{align*}
         proving property $(iv)$ and finishing the proof.
    \end{proof}
    \begin{satz}\label{thm:convergencedualsemigroup}
          Let $X$ be a bounded Polish metric space and $\lambda\in\MMM^+(X).$ Assume that \linebreak $\kappa:\Gamma_X\times X\to[0,\infty)$ satisfies \ref{property:cocycle} and~\ref{property:locallystable}. Then $\lim_{t\to\infty}\mathrm{S}_t^{\kappa,\lambda*}\dP=\dR$  in $\PPP_1(\Gamma_X)$  for all $\dP\in\PPP_1(\Gamma_X)$.
    \end{satz}
    \begin{proof}
        Let $\dP\in\PPP_1(\Gamma_X)$.  By Proposition~\ref{prop:convergencep1} it suffices to show that ${\lim_{t\to \infty}L_{\mathrm{S}_t^{\kappa,\lambda*}\dP}(f)=L_{\dR}(f)}$ for all $f\in\CCC_c^+(X)$ and $\lim_{t\to \infty}\I_{\mathrm{S}_t^{\kappa,\lambda*}\dP}=\I_{\dR}$ vaguely. To this end, let $f\in\CCC_c^+(X)$. By  Corollary~\ref{kor:convergencesemigroup} and the dominated convergence theorem, we deduce that 
         \begin{align*}
             \lim_{t\to \infty}L_{\mathrm{S}_t^{\kappa,\lambda*}\dP}(f)= \lim_{t\to \infty}\int \mathrm{S}_t^{\kappa,\lambda} e^{-\int f(x)\xi(\mathrm{d}x)}\dP(\mathrm{d}\xi)= \int  e^{-\int f(x)\xi(\mathrm{d}x)}\dR(\mathrm{d}\xi)=L_\dR(f).
         \end{align*}
         Further, let $B\in\BBB_0(X)$. 
        By an analogous argument to that in the proof of Theorem~\ref{thm:dualsemigroup} $(iv),$ we deduce that $\lim_{t\to \infty}\I_{\mathrm{S}_t^{\kappa,\lambda*}\dP}=\I_{\dR}$ vaguely,
        finishing the proof.
    \end{proof}
     \subsection{Entropy Dissipation}\label{subsec:energyfunctionals}
        In this subsection, we consider the evolution of the entropy along the dual spatial birth-and-death semigroup $(\mathrm{S}_t^{\kappa,\lambda*})_{t\geq0}$ and establish a de Bruijn identity.
 \begin{lem}[De Bruijn identity]\label{lem:energy}
     Let $X$ be a bounded Polish metric space and $\lambda\in\MMM^+(X).$ Assume that $\kappa:\Gamma_X\times X\to[0,\infty)$ satisfies \ref{property:cocycle} and~\ref{property:locallystable}. Further,   let $\dR=\RRR_{\kappa,\lambda}(\Gamma_X)$ and $\dP\in\PPP_1(\Gamma_{X})$ such that  $\Ent(\dP|\dR)<\infty$. It holds that 
        \begin{align*}
            \Ent(\dP|\dR)-\Ent(\mathrm{S}_T^{\kappa,\lambda*}\dP|\dR)=\int_0^T \I(\mathrm{S}_t^{\kappa,\lambda*}\dP|\dR)\mathrm{d}t,\quad T\geq0. 
        \end{align*}
        
    \end{lem}
    \begin{proof}
    The proof follows the same line of argumentation as the proof of the de Bruijn identity in the Poisson setting $\kappa\equiv1$ \cite[Theorem 3.3 (ii)]{schiavo2023wasserstein}. However, due to the slightly different domain of the generator $\mathfrak{L}_{\kappa,\lambda}$ and the additional Papangelou intensity in the definition of the Fisher information, we give an adapted proof here. We first assume that $\Ent(\dP|\dR)<\infty$ and  $\dP=\rho\dR$ with  $\rho \in\FFF_{b}(\Gamma_X)$. Define 
    \begin{align*}
        f(u):=u\log(u)-u+1,\quad u\geq0,
    \end{align*}
    and 
    \begin{align*}
        f_i(u):=\int_1^ui\wedge\log(v)\vee (-i)\mathrm{d}v,\quad i\in\N,\;u\geq0.
    \end{align*}
    As $\rho \in\FFF_{b}(\Gamma_X)$, also $\mathrm{S}_t^{\kappa,\lambda}\rho\in\FFF_{b}(\Gamma_X)$ and $f_i^\prime(\mathrm{S}_t^{\kappa,\lambda}\rho)\in\FFF_{b}(\Gamma_X)$. Therefore, we obtain by the integration by parts formula of Theorem~\ref{thm:ibp} that 
    \begin{align}
        \int f_i(\rho(\xi))\dR(\mathrm{d}\xi)- \int f_i(\mathrm{S}_T^{\kappa,\lambda}\rho(\xi))\dR(\mathrm{d}\xi)
        &=-\int_0^T\int f_i^\prime(\mathrm{S}_t^{\kappa,\lambda}\rho(\xi))\mathfrak{L}_{\kappa,\lambda}\mathrm{S}_t^{\kappa,\lambda}\rho(\xi)\dR(\mathrm{d}\xi)\mathrm{d}t\nonumber\\
        &=\int_0^T\int \int D_x f_i^\prime(\mathrm{S}_t^{\kappa,\lambda}\rho(\xi)) D_x\mathrm{S}_t^{\kappa,\lambda}\rho(\xi)\kappa(\xi,x)\dR(\mathrm{d}\xi)\lambda(\mathrm{d}x)\mathrm{d}t.\label{eq:3001}
    \end{align}
    By the monotone convergence theorem, we obtain for the left-hand side
    \begin{align}\label{eq:3002}
        \lim_{i\to\infty}\int f_i(\rho(\xi))\dR(\mathrm{d}\xi)- \int f_i(\mathrm{S}_T^{\kappa,\lambda}\rho(\xi))\dR(\mathrm{d}\xi)=\Ent(\dP|\dR)- \Ent(\mathrm{S}_T^{\kappa,\lambda*}\dP|\dR).
    \end{align}
   The integrand of the right-hand side $D_x f_i^\prime(\mathrm{S}_t^{\kappa,\lambda}\rho)(\xi) D_x\mathrm{S}_t^{\kappa,\lambda}\rho(\xi)\kappa(\xi,x)$ is non-negative by the convexity of $f_i$,  and it holds that  \begin{align}\label{eq:3004}
       (f^\prime_{i+1}(u)-f^\prime_{i+1}(v))(u-v)\geq (f^\prime_{i}(u)-f^\prime_{i}(v))(u-v),\quad u,v\geq0,
   \end{align}
   since $f_{i+1}-f_i$ is convex (see \cite[p. 20]{schiavo2023wasserstein}). Hence, by the monotone convergence theorem,
   \begin{align}\label{eq:3003}
       \lim_{i\to\infty} \int_0^T\int D_x f_i^\prime(\mathrm{S}_t^{\kappa,\lambda}\rho(\xi)) D_x\mathrm{S}_t^{\kappa,\lambda}\rho(\xi)\kappa(\xi,x)\dR(\mathrm{d}\xi)\lambda(\mathrm{d}x)\mathrm{d}t=\int_0^T \I(\mathrm{S}_t^{\kappa,\lambda*}\dP|\dR)\mathrm{d}t.
   \end{align}
   Combining \eqref{eq:3001}, \eqref{eq:3002}, and \eqref{eq:3003} yields the claim. Now, we assume that $\dP=\rho\dR$ such that  $\rho\in\FFF(\Gamma_X)$ and $\Ent(\dP|\dR)<\infty$. Let $$\dP_i(\mathrm{d}\xi):=\frac{(\rho\wedge i)}{\int \rho(\xi)\wedge i \dR(\mathrm{d}\xi)}\dR(\mathrm{d}\xi).$$
   As $(\rho\wedge i)\in\FFF_b(\Gamma_X)$, we obtain by the previous proof for bounded densities that 
   \begin{align}\label{eq:3005}
        \Ent(\dP_i|\dR)-\Ent(\mathrm{S}_T^{\kappa,\lambda*}\dP_i|\dR)=\int_0^T \I(\mathrm{S}_t^{\kappa,\lambda*}\dP_i|\dR)\mathrm{d}t.
   \end{align}
   Further, by a direct computation we get that
   \begin{align*}
       \I(\mathrm{S}_t^{\kappa,\lambda*}\dP_i|\dR)=\frac{1}{\int \rho(\xi)\wedge i \dR(\mathrm{d}\xi)}\int D_x(\mathrm{S}_t^{\kappa,\lambda}\rho\wedge i)(\xi)D_xf^\prime_i(\mathrm{S}_t^{\kappa,\lambda}\rho)(\xi)\kappa(\xi,x)\dR(\mathrm{d}\xi)\lambda(\mathrm{d}x),
   \end{align*}
   and similarly to \eqref{eq:3004} it holds that 
   \begin{align*}
       (\log(u\wedge(i+1))-\log(v\wedge(i+1)))(u\wedge(i+1)-v\wedge(i+1))\geq (\log(u\wedge i)-\log(v\wedge i))(u\wedge i-v\wedge i).
   \end{align*}
   Hence, by taking the limit in \eqref{eq:3005} and using the monotone convergence theorem on both sides, the claim follows.
    \end{proof}

    \begin{bem}
Note that Lemma~\ref{lem:energy} is restricted to compact base spaces $X$. This assumption is necessary to ensure that the integration by parts formula from Theorem~\ref{thm:ibp} is applicable in \eqref{eq:3001}.
\end{bem}

     The variation of the transport distance $\mathcal{W}_{\kappa,\lambda}$ along the dual semigroup $(\mathrm{S}_t^{\kappa,\lambda*})_{t\geq0}$ can be bounded in terms of the modified Fisher information. 
    \begin{lem}\label{lem:derivativeFisher}
Let $\lambda\in\MMM^+(X)$, assume that $\kappa:\Gamma_X\times X\to[0,\infty)$ satisfies \ref{property:cocycle}, \ref{property:locallystable}, \ref{property:continuity}, and~\ref{property:finiterange} and let $\dR\in\RRR_{\kappa,\lambda}(\Gamma_X)$. If the mean function $\theta$ is given by the logarithmic mean \eqref{eq:logmean} and $\dP_0,\dP_1\in\PPP_1(\Gamma_X)$, then
\begin{align*}
    \mathcal{W}_{\kappa,\lambda}(\mathrm{S}_{t+s}^{\kappa,\lambda*}{\dP}_0,\dP_1)
    -\mathcal{W}_{\kappa,\lambda}(\mathrm{S}_{t}^{\kappa,\lambda*}{\dP}_0,\dP_1)
    \leq
    \int_t^{s+t}\sqrt{\mathrm{I}(\mathrm{S}_{u}^{\kappa,\lambda*}{\dP}_0|\dR)}\mathrm{d}u,
    \quad s,t\geq0.
\end{align*}
\end{lem}

\begin{proof}
        The proof follows the same line of argumentation as the proof of \cite[Proposition 5.22]{schiavo2023wasserstein}. For the sake of completeness, we provide a full proof here. We assume that $\dP_0=\rho\dR$ and $\mathrm{I}(\dP_0|\dR)<\infty$ since otherwise the claim is immediate. Let $\bar{\dV}_t=-\kappa D\mathrm{S}_t^{\kappa,\lambda}\rho(\dR\otimes\lambda),\; t\geq0$, and recall that $(\mathrm{S}_{t}^{\kappa,\lambda*}{\dP}_0,\bar{\dV}_t)_{t\in[0,T]}\in\CCC\EEE^{X}_T(\dP_0,\mathrm{S}_{T}^{\kappa,\lambda*}\dP_0)$ by Lemma~\ref{lem:solutionsemigroup}. 
        By the definition of the modified Fisher information we directly get that
        \begin{align} \mathcal{L}_{\kappa,\lambda}(\mathrm{S}_{t}^{\kappa,\lambda*}{\dP}_0,\bar{\dV}_t)&=\int\int\alpha(\mathrm{S}_t^{\kappa,\lambda}\rho(\xi),\mathrm{S}_t^{\kappa,\lambda}\rho(\xi+\delta_x),-D_x\mathrm{S}_t^{\kappa,\lambda}\rho(\xi))\kappa(\xi,x)\dR(\mathrm{d}\xi)\lambda(\mathrm{d}x)\nonumber\\
        &=\int\frac{(D_x\mathrm{S}_t^{\kappa,\lambda}\rho(\xi))^2}{\frac{D_x\mathrm{S}_t^{\kappa,\lambda}\rho(\xi)}{D_x\log(\mathrm{S}_t^{\kappa,\lambda}\rho)(\xi)}}\kappa(\xi,x)\dR\otimes\lambda(\mathrm{d}\xi,\mathrm{d}x)\nonumber\\
        &=\mathrm{I}(\mathrm{S}_{t}^{\kappa,\lambda*}{\dP}_0|\dR).\label{eq:11137}
        \end{align}
         Thus, by the definition of the metric derivative (see Definition~\ref{dfn:metricderivative}) and Proposition~\ref{prop:metricderivative}, we conclude that 
        \begin{align*} \mathcal{W}_{\kappa,\lambda}(\mathrm{S}_{t+s}^{\kappa,\lambda*}{\dP}_0,\dP_1)-\mathcal{W}_{\kappa,\lambda}(\mathrm{S}_t^{\kappa,\lambda*}{\dP}_0,\dP_1)&\leq \mathcal{W}_{\kappa,\lambda}(\mathrm{S}_{t+s}^{\kappa,\lambda*}{\dP}_0,\mathrm{S}_{t}^{\kappa,\lambda*}{\dP}_0)\\
        &\leq \int_t^{s+t} |\dot{\mathrm{S}_u^{\kappa,\lambda*}{\dP}_0}|\mathrm{d}u\\
        &\leq \int_t^{s+t}  \sqrt{\mathrm{I}(\mathrm{S}_u^{\kappa,\lambda*}{\dP}_0|\dR)}\mathrm{d}u,
        \end{align*}
        proving the claim.
    \end{proof}

    \section{Convergence of the Discretization}\label{sec:convergencediscretization}       
\subsection{Papangelou Measures}\label{subsec:Papangelou}
We first recall relevant notation and prove the convergence of the discrete Papangelou measures. Let $n\in\N$, $\lambda=\Leb_{|\Lambda_n}$, and assume that  $\kappa:\Gamma_{\Lambda_n}\times \Lambda_n\to[0,\infty)$ satisfies \ref{property:cocycle} and~\ref{property:locallystable}. Let $\lambda^M=\Leb_{|\Lambda_n^M}^M$ and
\begin{align*}
    \kappa^M(\xi^M,x^M)=\kappa(\xi^M,x^M)\ind_{\{\xi^M(\{x^M\})=0\}},\quad \xi^M\in\Gamma_{\Lambda_n^M},\ x^M\in\Lambda_n^M.
    \end{align*}
$\tilde{\dR}^M=\RRR_{\kappa^M,\lambda^M}(\Gamma_{\Lambda_n^M})$
admits by Lemma~\ref{lem:finitevolumePapangelou} the explicit representation
      \begin{align*}
          \int F(\xi^M)\tilde\dR^M(\mathrm{d}\xi^M)= \frac{1}{Z_{\kappa^M,\lambda^M}^{\Lambda_n^M}} \sum_{j\geq0}\frac{1}{j!}\int\ldots\int F\left(\sum_{i=1}^j\delta_{x_i^M}\right)\varphi_{\kappa^M}\left(\sum_{i=1}^j\delta_{x_i^M}\right)\lambda^M(\mathrm{d}x_1^M)\ldots \lambda^M(\mathrm{d}x_j^M),
      \end{align*}
    with \begin{align*}
      \varphi_{\kappa^M}\left(\sum_{i=1}^j\delta_{x_i^M}\right)=\prod_{i=1}^j\kappa^M\left(\sum_{l=1}^{i-1}\delta_{x_l^M},x_i^M\right) 
    \end{align*}
    and
    \begin{align*}
        Z_{\kappa^M,\lambda^M}^{\Lambda_n^M}= \sum_{j\geq0}\frac{1}{j!}\int\ldots\int \varphi_{\kappa^M}\left(\sum_{i=1}^j\delta_{x_i^M}\right)\lambda^M(\mathrm{d}x_1^M)\ldots \lambda^M(\mathrm{d}x_j^M).
    \end{align*}

    \begin{satz}\label{thm:convergencePapangeloumeasure}
   Let $\lambda=\Leb_{|\Lambda_n}$ and assume that $\kappa:\Gamma_{\Lambda_n}\times \Lambda_n\to[0,\infty)$ satisfies \ref{property:cocycle}, \ref{property:locallystable}, and~\ref{property:continuity2}. For  $\dR=\RRR_{\kappa,\lambda}(\Gamma_{\Lambda_n})$ and  $\tilde{\dR}^M=\RRR_{\kappa^M,\lambda^M}(\Gamma_{\Lambda_n^M}),\; M\in\N,$ it holds that $\lim_{M\to\infty}\tilde{\dR}^M=\dR$ in $\PPP_1(\Gamma_{\Lambda_n}).$
         \end{satz}
          \begin{proof} By Proposition~\ref{prop:convergencep1} it suffices to show that $\lim_{M\to\infty}L_{\tilde{\dR}^M}(f)=L_{\dR}(f)$ for all functions $f\in\CCC^+_c(\Lambda_n)$ and $\lim_{M\to\infty}\I_{\tilde{\dR}^M}=\I_{\dR}$ vaguely.
          
            We first show that the  partition function $Z_{\kappa^M,\lambda^M}^{\Lambda_n^M}$ converges to $Z_{\kappa,\lambda}^{\Lambda_n}$ for $M\to\infty$. To this end, we note that 
            \begin{align*}
                Z_{\kappa^M,\lambda^M}^{\Lambda_n^M}&= \sum_{j\geq0}\sum_{\substack{(x_1^M,\ldots,x_j^M)\in(\Lambda_n^M)^j\\ x_i^M\not=x_k^M\forall i\not=k}}\frac{1}{j!}\left(\frac{1}{M^d}\right)^j\varphi_{\kappa^M}\left(\sum_{i=1}^j\delta_{x_i^M}\right).
            \end{align*}
            By \ref{property:locallystable}, it holds that 
         \begin{align*}
           \sum_{\substack{(x_1^M,\ldots,x_j^M)\in(\Lambda_n^M)^j\\ x_i^M\not=x_k^M\forall i\not=k}}\frac{1}{j!}\left(\frac{1}{M^d}\right)^j\varphi_{\kappa^M}\left(\sum_{i=1}^j\delta_{x_i^M}\right) \leq  \frac{1}{j!}\left(\frac{1}{M^d}\right)^jA^j(n^dM^d)^j=\frac{(An^d)^j}{j!},
         \end{align*}
         which is summable over $j$. Therefore, by the dominated convergence theorem, the Riemann approximation, and~\ref{property:continuity2}, we conclude 
         \begin{align}\label{eq:4001}
             \lim_{M\to\infty} Z_{\kappa^M,\lambda^M}^{\Lambda_n^M}=\sum_{j\geq0}\int_{\Lambda_n}\ldots\int_{\Lambda_n} \frac{1}{j!}\varphi_{\kappa}\left(\sum_{i=1}^j\delta_{x_i}\right)\lambda(\mathrm{d}x_1)\ldots \lambda(\mathrm{d}x_j)= Z_{\kappa,\lambda}^{\Lambda_n}.
         \end{align}
        Now, let $f\in\CCC_c^+(\Lambda_n)$. 
        Similar to before, we obtain by the Riemann approximation and the dominated convergence theorem that 
        \begin{align*}
            &\lim_{M\to\infty}\int e^{-\int f(x^M)\xi^M(\mathrm{d}x^M)}\tilde{\dR}^M(\mathrm{d}\xi^M)\\
            &\quad=  \lim_{M\to\infty}\frac{1}{ Z_{\kappa^M,\lambda^M}^{\Lambda_n^M}}\sum_{j\geq0}\sum_{\substack{(x_1^M,\ldots,x_j^M)\in(\Lambda_n^M)^j\\ x_i^M\not=x_k^M\forall i\not=k}}\frac{1}{j!}\left(\frac{1}{M^d}\right)^j\varphi_{\kappa^M}\left(\sum_{i=1}^j\delta_{x_i^M}\right) e^{-\int f(y)\sum_{i=1}^j\delta_{x_i^M}(\mathrm{d}y)}\\
            &\quad= \frac{1}{ Z_{\kappa,\lambda}^{\Lambda_n}} \sum_{j\geq0} \int_{\Lambda_n}\ldots\int_{\Lambda_n} \frac{1}{j!}\varphi_{\kappa}\left(\sum_{i=1}^j\delta_{x_i}\right)e^{-\int f(y)\sum_{i=1}^j\delta_{x_i}(\mathrm{d}y)}\lambda(\mathrm{d}x_1)\ldots \lambda(\mathrm{d}x_j)\\
            &\quad= \int e^{-\int f(x)\xi(\mathrm{d}x)}{\dR}(\mathrm{d}\xi).
        \end{align*}
        An analogous argument for the intensity measure finishes the proof.
        \end{proof}

 \subsection{Dual Semigroup}
      In this subsection, we prove Theorem~\ref{thm:convergenceresultsdiscretization} $(i)$. More precisely, we show that if $\lambda=\Leb_{|\Lambda_n}$ and $\kappa:\Gamma_{\Lambda_n}\times\Lambda_n\to[0,\infty)$ satisfies \ref{property:cocycle}, \ref{property:locallystable}, and~\ref{property:continuity2}, then \begin{align*}
          \lim_{M\to\infty}\mathrm{S}_t^{\kappa^M,\lambda^M*}\dP^M=\mathrm{S}_t^{\kappa,\lambda*}\dP \quad \text{ in }\PPP_1(\Gamma_{\Lambda_n}),
      \end{align*}
      for all $t\geq0$ and $\dP\in\PPP_1(\Gamma_{\Lambda_n})$ s.t.\   $\dP\ll\dR$. We begin by establishing the corresponding result for the semigroups $(\mathrm{S}_t^{\kappa^M,\lambda^M})_{t\geq0}$. \begin{satz}\label{thm:discretizedsemigroup}
      Let $\lambda=\Leb_{|\Lambda_n}$ and assume that $\kappa:\Gamma_{\Lambda_n}\times \Lambda_n\to[0,\infty)$ satisfies \ref{property:cocycle}, \ref{property:locallystable}, and~\ref{property:continuity2}. Let $\dR=\RRR_{\kappa,\lambda}(\Gamma_{\Lambda_n})$. It holds that
        \begin{align*}
            \lim_{M\to\infty}\|\mathrm{S}_t^{\kappa^M,\lambda^M}\Psi^M(F)-\mathrm{S}_t^{\kappa,\lambda}F\|_{L^2(\dR)}=0,\quad F\in L^2(\dR),\;  t\geq0,
        \end{align*}
        with $\Psi^M(F)(\xi):=F(\xi^M),\; \xi\in\Gamma_{\Lambda_n}$.
    \end{satz}
    \begin{proof}[Proof of Theorem~\ref{thm:discretizedsemigroup}]
      We first note that  the cylinder functions $\mathcal{G}$ are a core of $\mathfrak{L}_{\kappa,\lambda}$. Using \eqref{eq:derivativegenerator}, we have 
        \begin{align}
        \mathcal{G}\subseteq \mathrm{Dom}(\mathfrak{L}_{\kappa,\lambda}):=\left\{F\in L^2(\dR): \lim_{t\to 0}\frac{\mathrm{S}_t^{\kappa,\lambda}F-F}{t} \text{ exists in } L^2(\dR)\right\}\label{eq:core}
    \end{align} 
      and  $\overline{\mathfrak{L}_{\kappa,\lambda}\big|_{\mathcal{G}}}=\mathfrak{L}_{\kappa,\lambda}$, see e.g.\ \cite[Theorem 4.1]{KoLy05}. Since $(\mathrm{S}_t^{\kappa,\lambda})_{t\geq0}$ and $(\mathrm{S}_t^{\kappa^M,\lambda^M})_{t\geq0}$ are contractive $C_0$-semigroups (Lemma~\ref{lem:strongcontinuity}), the claim follows from the Kato--Trotter approximation, see e.g.\ \cite[Theorem~6.1]{EthierKurtz}, if the corresponding convergence for the generator 
      \begin{align*} \lim_{M\to\infty}\|\mathfrak{L}_{\kappa^M,\lambda^M}\Psi^M(F)-\mathfrak{L}_{\kappa,\lambda} F\|_{L^2(\dR)}=0
       \end{align*}
     holds for all functions $F\in\mathcal{G}$.
      Therefore, let $F\in\mathcal{G}$ and note that
     \begin{align*}
           &\lim_{M\to\infty}\|\mathfrak{L}_{\kappa^M,\lambda^M}\Psi^M(F)-\mathfrak{L}_{\kappa,\lambda} F\|^2_{L^2(\dR)}\\
           &\quad=\lim_{M\to\infty}\int | \int\kappa^M(\xi^M,x^M)(F(\xi^M+\delta_{x^M})-F(\xi^M))-\kappa(\xi,x)(F(\xi+\delta_{x})-F(\xi))\lambda(\mathrm{d}x)\\
           &\qquad +\sum_{x^M\in\xi^M}(F(\xi^M-\delta_{x^M})-F(\xi^M))-\sum_{x\in\xi}(F(\xi-\delta_{x})-F(\xi))|^2\dR(\mathrm{d}\xi).
       \end{align*}
       For the birth term of the generator, we obtain
       \begin{align*}
           &\lim_{M\to\infty}\int |\int \kappa^M(\xi^M,x^M)(F(\xi^M+\delta_{x^M})-F(\xi^M))-\kappa(\xi,x)(F(\xi+\delta_{x})-F(\xi))\lambda(\mathrm{d}x)|^2\dR(\mathrm{d}\xi)\\
           &\quad=\lim_{M\to\infty}\int |\int\kappa^M(\xi^M,x^M)(F(\xi^M+\delta_{x^M})-F(\xi^M)-F(\xi+\delta_{x})+F(\xi))\\
           &\qquad-(\kappa(\xi,x)-\kappa(\xi^M,x^M))(F(\xi+\delta_{x})-F(\xi))\lambda(\mathrm{d}x)|^2\dR(\mathrm{d}\xi)\\
           &\quad \leq\lim_{M\to\infty}2\lambda(\Lambda_n)\int \int |\kappa^M(\xi^M,x^M)|^2|(F(\xi^M+\delta_{x^M})-F(\xi^M)-F(\xi+\delta_{x})+F(\xi))|^2\lambda(\mathrm{d}x)\dR(\mathrm{d}\xi)\\
           &\qquad+ 2\lambda(\Lambda_n)\int \int|(\kappa(\xi,x)-\kappa(\xi^M,x^M))|^2|(F(\xi+\delta_{x})-F(\xi))|^2\lambda(\mathrm{d}x)\dR(\mathrm{d}\xi).
       \end{align*}
       Since $F \in \mathcal{G}$ and $\lim_{M\to\infty}\kappa^M(\xi^M, x^M) = \kappa(\xi, x)$ $\dR\otimes\lambda$-a.e.\ \ref{property:continuity2}, it follows that  
\begin{align*}
\lim_{M \to \infty} 
\bigl|\kappa^M(\xi^M, x^M)\bigr|^2 
\bigl|F(\xi^M + \delta_{x^M}) - F(\xi^M) - F(\xi + \delta_x) + F(\xi)\bigr|^2 
= 0\quad   \dR\otimes\lambda-\text{a.e.}
\end{align*}
Moreover, since $F \in \mathcal{G}$, there exists a constant $c \geq 1$ s.t.\    
$|F(\xi)| \le c^{\xi(\Lambda_n)}$. Hence, by \ref{property:locallystable},  
\begin{align*}
\bigl|\kappa^M(\xi^M, x^M)\bigr|^2 
\bigl|F(\xi^M + \delta_{x^M}) - F(\xi^M) - F(\xi + \delta_x) + F(\xi)\bigr|^2
\le 2A^2\bigl(c^{2\xi(\Lambda_n) + 2} + c^{2\xi(\Lambda_n)}\bigr).
\end{align*}
The upper bound is $ \dR\otimes \lambda$-integrable, since  
$\int c^{\xi(\Lambda_n)} \dR(\mathrm{d}\xi) < \infty$ by \eqref{property:integrability}.  
Therefore, using the dominated convergence theorem, we conclude that
\begin{align*}
           \lim_{M\to\infty}\int \int |\kappa^M(\xi^M,x^M)|^2|(F(\xi^M+\delta_{x^M})-F(\xi^M)-F(\xi+\delta_{x})+F(\xi))|^2\lambda(\mathrm{d}x)\dR(\mathrm{d}\xi)=0.
       \end{align*}
      Similarly, we obtain that 
      \begin{align*}
         \lim_{M\to\infty} \int \int|(\kappa(\xi,x)-\kappa(\xi^M,x^M))|^2|(F(\xi+\delta_{x})-F(\xi))|^2\lambda(\mathrm{d}x)\dR(\mathrm{d}\xi)=0.
      \end{align*}
      This proves the convergence of the birth term.
       Using the GNZ equation~\eqref{eq:papangelou}, we obtain an analogous result for the death term. All in all, for any function $F\in\mathcal{G}$ it holds that \begin{align*}
           \lim_{M\to\infty}\|\mathfrak{L}_{\kappa^M,\lambda^M}\Psi^M(F)-\mathfrak{L}_{\kappa,\lambda} F\|_{L^2(\dR)}=0,
       \end{align*}
       finishing the proof.
      \end{proof}

    Using the previous convergence result for the semigroup, we now turn to the proof of Theorem~\ref{thm:convergenceresultsdiscretization}~$(i)$ and lift the convergence to the dual semigroup.
    \begin{proof}[Proof of Theorem~\ref{thm:convergenceresultsdiscretization} $(i)$]
     Let $\dP=\rho\dR\in\PPP_1(\Gamma_{\Lambda_n})$ and $t\geq0$. By Proposition~\ref{prop:convergencep1}, it suffices to show that $\lim_{M\to\infty}L_{\mathrm{S}_t^{\kappa^M,\lambda^M*}\dP^M}(f)=L_{\mathrm{S}_t^{\kappa,\lambda*}\dP}(f)$ for all $f\in\CCC_c^+(\Lambda_n)$ and $\lim_{M\to\infty}\I_{\mathrm{S}_t^{\kappa^M,\lambda^M*}\dP^M}=\I_{\mathrm{S}_t^{\kappa,\lambda*}\dP}$ vaguely.  Let $f\in\CCC_c^+(\Lambda_n)$. It holds that 
     \begin{align}
         &|\int e^{-\int f(x)\xi^M(\mathrm{d}x)}\mathrm{S}_t^{\kappa^M,\lambda^M*}\dP^M(\mathrm{d}\xi^M)-\int e^{-\int f(x)\xi(\mathrm{d}x)}\mathrm{S}_t^{\kappa,\lambda*}\dP(\mathrm{d}\xi)|\nonumber\\&\quad=|\int (\mathrm{S}_t^{\kappa^M,\lambda^M}e^{-\int f(x)\xi^M(\mathrm{d}x)}-\mathrm{S}_t^{\kappa,\lambda}e^{-\int f(x)\xi(\mathrm{d}x)})\rho(\xi)\dR(\mathrm{d}\xi)|\nonumber\\
         &\quad\leq \int |\mathrm{S}_t^{\kappa^M,\lambda^M}e^{-\int f(x)\xi^M(\mathrm{d}x)}-\mathrm{S}_t^{\kappa,\lambda}e^{-\int f(x)\xi(\mathrm{d}x)}|\rho(\xi)\ind_{\{\rho(\xi)\leq i\}}\dR(\mathrm{d}\xi)+2\int \ind_{\{\rho(\xi)>i\}}\rho(\xi)\dR(\mathrm{d}\xi),\label{eq:65}
     \end{align}
     for all $i\in\N$. By Theorem~\ref{thm:discretizedsemigroup}, we obtain for the first term of the right-hand side of \eqref{eq:65} that \begin{align*}
         \lim_{M\to\infty}\int |\mathrm{S}_t^{\kappa^M,\lambda^M}e^{-\int f(x)\xi^M(\mathrm{d}x)}-\mathrm{S}_t^{\kappa,\lambda}e^{-\int f(x)\xi(\mathrm{d}x)}|\rho(\xi)\ind_{\{\rho(\xi)\leq i\}}\dR(\mathrm{d}\xi)=0.
     \end{align*}
     Since $\rho$ is a probability density, the second term of the right-hand side of \eqref{eq:65} tends to zero for $i\to\infty.$ Hence, $\lim_{M\to\infty}L_{\mathrm{S}_t^{\kappa^M,\lambda^M*}\dP^M}(f)=L_{\mathrm{S}_t^{\kappa,\lambda*}\dP}(f)$. An analogous argument for the intensity measure using additionally Lemma~\ref{lem:intensitymeasuresemigroup} finishes the proof.
    \end{proof}

\subsection{Entropy}
In this subsection, we prove Theorem~\ref{thm:convergenceresultsdiscretization} $(ii)$. More precisely,  let $\lambda=\Leb_{|\Lambda_n}$ and  assume that $\kappa:\Gamma_{\Lambda_n}\times\Lambda_n\to[0,\infty)$ satisfies \ref{property:cocycle}, \ref{property:locallystable}, \ref{property:lowerbound}, \ref{property:continuity2}, and~\ref{property:repulsive}. We show that if $\dR=\RRR_{\kappa,\lambda}(\Gamma_{\Lambda_n})$ and $\tilde\dR^M=\RRR_{\kappa^M,\lambda^M}(\Gamma_{\Lambda_n^M}),\; M\in\N$, then $$\Gamma\text{-}\lim_{M\to\infty}\Ent(\cdot|\tilde{\dR}^M)
    =\Ent(\cdot|\dR)$$ on $\PPP_1(\Gamma_{\Lambda_n})$. A recovery sequence for $\dP\in\PPP_1(\Gamma_{\Lambda_n})$ is $(\dP^M)_{M\in\N}$.
\begin{proof}[Proof of Theorem~\ref{thm:convergenceresultsdiscretization} $(ii)$]
Let $\dP_M\in\PPP_1(\Gamma_{\Lambda_n^M}),M\in\N,$ and $\dP\in\PPP_1(\Gamma_{\Lambda_n}),$ be s.t.\  $\lim_{M\to\infty}\dP_M=\dP$ in $\PPP_1(\Gamma_{\Lambda_n})$. Since $\lim_{M\to\infty}\tilde\dR^M=\dR$ in $\PPP_1(\Gamma_{\Lambda_n})$  by Theorem~\ref{thm:convergencePapangeloumeasure},  we obtain by the lower semicontinuity of the relative entropy (Lemma~\ref{lem:propertiesent} $(i)$) that
    \begin{align*}
    \Ent(\dP|\dR)\leq\liminf_{M\to\infty}\Ent(\dP_M|\tilde\dR^M).
    \end{align*}
     It remains to show that  for $\dP\in\PPP_1(\Gamma_{\Lambda_n})$ and $\dP^M=(T^M_\Gamma)_\#\dP,\;M\in\N$, it holds that
         \begin{align*}
            \limsup_{M\to\infty} \Ent(\dP^M|\tilde{\dR}^M)\leq  \Ent(\dP|\dR). 
        \end{align*}
        If $\Ent(\dP|\dR)=\infty$, the claim is true. Therefore, let $\Ent(\dP|\dR)<\infty$ and $\dP=\rho\dR$.
        By \ref{property:locallystable} and~\ref{property:lowerbound}, there exists a  set
$E\in\BBB(\Gamma_{\Lambda_n})$ with $\dR(E)=\dR^M(E)=1$ and  for all $\xi+\delta_x\in E$ it holds that $\tfrac{1}{a } \leq \kappa(\xi,x) \leq A$.  
This implies that $\tilde{\dR}^M(\{\xi^M\}) > 0$ for all  $\xi^M \in E$. Hence, $\dR^M \ll \tilde{\dR}^M$, with
\begin{align*}
   f^M(\xi^M):= \frac{\mathrm{d}{\dR}^M}{\mathrm{d}\tilde\dR^M}(\xi^M)=\frac{{\dR}^M(\{\xi^M\})}{\tilde\dR^M(\{\xi^M\})}.
\end{align*}
Further, by disintegration, we find measures $\dR^{*,M}(\mathrm{d}\eta|\xi^M),\; \xi^M\in\Gamma_{\Lambda_n^M}$, on $\Gamma_{\Lambda_n}$ s.t.\   
           \begin{align*}
               \int F(\xi)\dR(\mathrm{d}\xi)=\int\int F(\eta)\dR^{*,M}(\mathrm{d}\eta|\xi^M)\dR^M(\mathrm{d}\xi^M),
           \end{align*}
           for all $F\in\FFF^+(\Gamma_{\Lambda_n})$. Let $K\in\BBB(\Gamma_{\Lambda_n^M})$. It holds that 
           \begin{align*}
               \int \ind_K(\xi^M)\dP^M(\mathrm{d}\xi^M)&= \int \ind_K(\xi^M)\rho(\xi)\dR(\mathrm{d}\xi)= \int \ind_K(\xi^M)\int\rho(\eta)\dR^{*,M}(\mathrm{d}\eta|\xi^M)\dR^M(\mathrm{d}\xi^M),
           \end{align*}
           implying that $\dP^M\ll\dR^M$ with
           \begin{align*}
            g^M(\xi^M):=   \frac{\mathrm{d}\dP^M}{\mathrm{d}\dR^M}(\xi^M)=\int \rho(\eta)\dR^{*,M}(\mathrm{d}\eta|\xi^M).
           \end{align*}
     Using this observation, we have that
     \begin{align}
            \Ent(\dP^M|\tilde{\dR}^M)&=\int  g^M(\xi^M)f^M(\xi^M)\log\left( g^M(\xi^M)f^M(\xi^M)\right)\tilde{\dR}^M(\mathrm{d}\xi^M)\nonumber\\
            &= \Ent(\dP^M|\dR^M) +  \int  \log\left( \frac{{\dR}^M(\{\xi^M\})}{\tilde\dR^M(\{\xi^M\})}\right)\rho(\xi){\dR}(\mathrm{d}\xi).\label{eq:3000}
        \end{align}
       Further, by the contraction property of the relative entropy (Lemma~\ref{lem:propertiesent} $(ii)$), 
        \begin{align}
            \Ent(\dP^M|\dR^M)=\Ent(({T_\Gamma^M})_\#\dP|({T_\Gamma^M})_\#\dR)\leq \Ent(\dP|\dR).\label{eq:4010}
        \end{align}
        Combining equation~\eqref{eq:3000} and \eqref{eq:4010}, we obtain that 
        \begin{align*}
            \limsup_{M\to\infty} \Ent(\dP^M|\tilde{\dR}^M)\leq \Ent(\dP|\dR)+\limsup_{M\to\infty} \int  \log\left( \frac{{\dR}^M(\{\xi^M\})}{\tilde\dR^M(\{\xi^M\})}\right)\rho(\xi){\dR}(\mathrm{d}\xi).
        \end{align*}
        Hence, it remains to show that $$\lim_{M\to\infty}\int  \log\left( \frac{{\dR}^M(\{\xi^M\})}{\tilde\dR^M(\{\xi^M\})}\right)\rho(\xi){\dR}(\mathrm{d}\xi)=0.$$
        By the definition of $\tilde\dR^M$ and Lemma~\ref{lem:finitevolumePapangelou}, it holds that 
        \begin{align}\label{eq:954}
            \tilde\dR^M(\{\xi^M\})=\frac{1}{Z_{\kappa^M,\lambda^M}^{\Lambda_n^M}}\left(\frac{1}{M^d}\right)^{\xi^M(\Lambda_n^M)}\varphi_{\kappa^M}(\xi^M).
        \end{align}
        For finite $\xi\in\Gamma_{\Lambda_n}$, we denote by $j_\xi:=\xi(\Lambda_n)$ the total number of points in $\Lambda_n$ and by $\xi(i)$ the $i$-th point of $\xi,$ where the order of the points is chosen arbitrarily. Further, we introduce the function
         \begin{align}\label{eq:definitionnumberofpointsinabox}
             n_{\xl}^M(\xi)=\begin{cases}
                 \xi(\xl+\Lambda_{\frac{1}{M}}), &\text{if } \xi(\xl+\Lambda_{\frac{1}{M}})>0,\\
                 \infty, &\text{otherwise, } 
             \end{cases}
             \quad   \xl\in\Lambda_n^M,\; \xi\in\Gamma_{\Lambda_n}.
         \end{align} 
        Using the definition of $n_{\xl}^M$ and the  GNZ equation~\eqref{eq:papangelou}, we obtain for $\xi\in\Gamma_{\Lambda_n}$ that 
    \begin{align*}
         {\dR}^M(\{\xi^M\})&=\int\prod_{i=1}^{j_\xi}\ind_{\{\eta(\Lambda_{\frac{1}{M}}+{\xi(i)}^M)\geq1\}}\prod_{\substack{z^M\in\Lambda_n^M\setminus\xi^M}}\ind_{\{\eta(\Lambda_{\frac{1}{M}}+z^M)=0\}}\dR(\mathrm{d}\eta)\\
         &= \int\int_{\Lambda_{\frac{1}{M}}+{\xi(1)}^M}\frac{1}{n_{{\xi(1)}^M}^M(\eta)}\ind_{\{\eta(\Lambda_{\frac{1}{M}}+{\xi(1)}^M)\geq1\}}\eta(\mathrm{d}y_1)\prod_{i=2}^{j_\xi}\ind_{\{\eta(\Lambda_{\frac{1}{M}}+{\xi(i)}^M)\geq1\}}\\
         &\quad\times\prod_{\substack{z^M\in\Lambda_n^M\setminus\xi^M}}\ind_{\{\eta(\Lambda_{\frac{1}{M}}+z^M)=0\}}\dR(\mathrm{d}\eta)\\
         &= \int\left(\int_{\Lambda_{\frac{1}{M}}+{\xi(1)}^M}\frac{1}{n_{{\xi(1)}^M}^M(\eta)+1}\kappa(\eta,y_1)\ind_{\{\eta(\Lambda_{\frac{1}{M}}+\xi(1)^M)\geq0\}}\lambda(\mathrm{d}y_1)\right)\\
         &\quad \times \prod_{i=2}^{j_\xi}\ind_{\{\eta(\Lambda_{\frac{1}{M}}+{\xi(i)}^M)\geq1\}}\prod_{\substack{z^M\in\Lambda_n^M\setminus\xi^M}}\ind_{\{\eta(\Lambda_{\frac{1}{M}}+z^M)=0\}}\dR(\mathrm{d}\eta)\\
         &\leq \int\left(\int_{\Lambda_{\frac{1}{M}}+{\xi(1)}^M}\kappa(\eta,y_1)\lambda(\mathrm{d}y_1)\right)\prod_{i=2}^{j_\xi}\ind_{\{\eta(\Lambda_{\frac{1}{M}}+{\xi(i)}^M)\geq1\}}\prod_{\substack{z^M\in\Lambda_n^M\setminus\xi^M}}\ind_{\{\eta(\Lambda_{\frac{1}{M}}+z^M)=0\}}\dR(\mathrm{d}\eta).
         \end{align*}
         Repeating this step $j_\xi-1$ times and using that $\kappa$ is repulsive \ref{property:repulsive}, we conclude that  
         \begin{align}
             {\dR}^M(\{\xi^M\}) \nonumber
             &\leq \int\!\!\int_{\Lambda_{\frac{1}{M}}\!+\xi(1)^M}\!\!\!\!\!\!\ldots\int_{\Lambda_{\frac{1}{M}}\!+\xi({j_\xi})^M}\prod_{i=1}^{j_\xi}\kappa\left(\!\eta\!+\!\sum_{l=1}^{i-1}\delta_{y_l},y_i\!\right)\!\lambda(\mathrm{d}y_1)\ldots \lambda(\mathrm{d}y_{j_\xi})\!\!\!\!\!\prod_{\substack{z^M\in\Lambda_n^M\setminus\xi^M}}\!\!\!\!\ind_{\{\eta(\Lambda_{\frac{1}{M}}\!+z^M)=0\}}\dR(\mathrm{d}\eta)\nonumber\\
         &\leq  \int_{\Lambda_{\frac{1}{M}}\!+\xi(1)^M}\!\!\!\!\!\!\ldots\int_{\Lambda_{\frac{1}{M}}\!+\xi({j_\xi})^M}\prod_{i=1}^{j_\xi}\kappa\left(\sum_{l=1}^{i-1}\delta_{y_l},y_i\!\right)\lambda(\mathrm{d}y_1)\ldots \lambda(\mathrm{d}y_{j_\xi})\!\int \!\!\!\prod_{\substack{z^M\in\Lambda_n^M\setminus\xi^M}}\ind_{\{\eta(\Lambda_{\frac{1}{M}}\!+z^M)=0\}}\dR(\mathrm{d}\eta)\nonumber\\
         &= \int_{\Lambda_{\frac{1}{M}}\!+\xi(1)^M}\ldots\int_{\Lambda_{\frac{1}{M}}\!+\xi({j_\xi})^M}\varphi_\kappa\left(\sum_{i=1}^{j_\xi}\delta_{y_i}\right)\lambda(\mathrm{d}y_1)\ldots \lambda(\mathrm{d}y_{j_\xi}) \int \!\!\!\prod_{\substack{z^M\in\Lambda_n^M\setminus\xi^M}}\ind_{\{\eta(\Lambda_{\frac{1}{M}}\!+z^M)=0\}}\dR(\mathrm{d}\eta).\label{eq:4002}
    \end{align}
    Therefore, by \eqref{eq:954} and \eqref{eq:4002}, we have that 
    \begin{align}
       &\int \log\left(\frac{{\dR}^M(\{\xi^M\})}{\tilde\dR^M(\{\xi^M\})}\right)\rho(\xi)\dR(\mathrm{d}\xi)\nonumber\\
       &\quad\leq \int\log\left(\frac{\int_{\Lambda_{\frac{1}{M}}+\xi(1)^M}\ldots\int_{\Lambda_{\frac{1}{M}}+\xi({j_\xi})^M}\varphi_\kappa(\sum_{i=1}^{j_\xi}\delta_{y_i})\lambda(\mathrm{d}y_1)\ldots \lambda(\mathrm{d}y_{j_\xi}) }{\frac{1}{(M^d)^{{j_\xi}}}\varphi_{\kappa^M}(\xi^M)}\right)\rho(\xi)\dR(\mathrm{d}\xi)\nonumber\\
        &\qquad+\int\log\left(\frac{ \int \prod_{\substack{z^M\in\Lambda_n^M\setminus\xi^M}}\ind_{\{\eta(\Lambda_{\frac{1}{M}}+z^M)=0\}}\dR(\mathrm{d}\eta)}{\frac{1}{Z_{\kappa^M,\lambda^M}^{\Lambda_n^M}}}\right)\rho(\xi)\dR(\mathrm{d}\xi).\label{eq:7000}
    \end{align}
    We continue by proving that both terms of the right-hand side converge to zero. Since $\kappa$ is  $\dR\otimes\lambda$-a.e.\ continuous \ref{property:continuity2}, we obtain by the Lebesgue differentiation theorem  that 
    \begin{align*}
        \limsup_{M\to\infty}\log\left({\frac{\int_{\Lambda_{\frac{1}{M}}+\xi(1)^M}\ldots\int_{\Lambda_{\frac{1}{M}}+\xi({j_\xi})^M}\varphi_\kappa\left(\sum_{i=1}^{j_\xi}\delta_{y_i}\right)\lambda(\mathrm{d}y_1)\ldots \lambda(\mathrm{d}y_{j_\xi})}{(\frac{1}{M^d})^{j_\xi}\varphi_{\kappa^M}(\xi^M)}}\right)&=\log\left(\frac{\varphi_\kappa\left(\sum_{i=1}^{j_\xi}\delta_{{\xi(i)}}\right)}{\varphi_\kappa\left(\sum_{i=1}^{j_\xi}\delta_{{\xi(i)}}\right) }\right)\\
        &=0\ \dR\text{-a.s.}
    \end{align*}
    By \ref{property:locallystable} and~\ref{property:lowerbound},
    \begin{align*}
        \xi(\Lambda_n)\log\left(\frac{1}{Aa }\right)&\leq \log\left(\frac{\int_{\Lambda_{\frac{1}{M}}+\xi(1)^M}\ldots\int_{\Lambda_{\frac{1}{M}}+\xi({j_\xi})^M}\varphi_\kappa(\sum_{i=1}^{j_\xi}\delta_{y_i})\lambda(\mathrm{d}y_1)\ldots \lambda(\mathrm{d}y_{j_\xi}) }{\frac{1}{(M^d)^{{j_\xi}}}\varphi_{\kappa^M}(\xi^M)}\right)\leq \xi(\Lambda_n)\log\left(Aa \right) 
    \end{align*}
 $ \dR\text{-a.s.}$ Since $\dP\in\PPP_1(\Gamma_{\Lambda_n})$, we have $\int \xi(\Lambda_n)\dP(\mathrm{d}\xi)<\infty$. Thus, the dominated convergence theorem yields for the first term of \eqref{eq:7000} that 
   \begin{align}\limsup_{M\to\infty}\int\log\left(\frac{\int_{\Lambda_{\frac{1}{M}}+\xi(1)^M}\ldots\int_{\Lambda_{\frac{1}{M}}+\xi({j_\xi})^M}\varphi_\kappa(\sum_{i=1}^{j_\xi}\delta_{y_i})\lambda(\mathrm{d}y_1)\ldots \lambda(\mathrm{d}y_{j_\xi}) }{\frac{1}{(M^d)^{{j_\xi}}}\varphi_{\kappa^M}(\xi^M)}\right)\rho(\xi)\dR(\mathrm{d}\xi)=0.\label{eq:7001}
   \end{align}
    For the second term of \eqref{eq:7000}, we deduce by the dominated convergence theorem and Lemma~\ref{lem:finitevolumePapangelou} that \begin{align*}
      \limsup_{M\to\infty}  \int \prod_{\substack{z^M\in\Lambda_n^M\setminus\xi^M}}\ind_{\{\eta(\Lambda_{\frac{1}{M}}+z^M)=0\}}\dR(\mathrm{d}\eta)&=  \int \ind_{\{\eta(\Lambda_{n})=0\}}\dR(\mathrm{d}\eta)\nonumber= \frac{1}{Z_{\kappa,\lambda}^{\Lambda_n}}.
    \end{align*}
    Moreover, we recall that by \eqref{eq:4001} it holds that
    \begin{align*}
        \lim_{M\to\infty}Z_{\kappa^M,\lambda^M}^{\Lambda_n^M}&=\sum_{{j_\xi}\geq0}\int_{\Lambda_n}\ldots\int_{\Lambda_n} \frac{1}{{j_\xi}!}\varphi_{\kappa}\left(\sum_{i=1}^{j_\xi}\delta_{x_{i}}\right)\lambda(\mathrm{d}x_{1})\ldots \lambda(\mathrm{d}x_{{j_\xi}})={Z_{\kappa,\lambda}^{\Lambda_n}}.
    \end{align*}
    Hence, by the dominated convergence theorem, \begin{align}\label{eq:4030}
        \lim_{M\to\infty} \int\log\left(\frac{ \int \prod_{\substack{z^M\in\Lambda_n^M\setminus\xi^M}}\ind_{\{\eta(\Lambda_{\frac{1}{M}}\!+z^M)=0\}}\dR(\mathrm{d}\eta)}{\frac{1}{Z_{\kappa^M,\lambda^M}^{\Lambda_n^M}}}\right)\rho(\xi)\dR(\mathrm{d}\xi)=0.
    \end{align}
    Combining \eqref{eq:7000}, \eqref{eq:7001}, and \eqref{eq:4030}, we conclude for the error term that
    \begin{align*}
       & \limsup_{M\to\infty} \int  \log\left( \frac{{\dR}^M(\{\xi^M\})}{\tilde\dR^M(\{\xi^M\})}\right)\rho(\xi){\dR}(\mathrm{d}\xi)\leq0,
        \end{align*}
        finishing the proof.
        \end{proof}
\subsection{Transport Distance}\label{subsec:discretizationdistance}

In this subsection, we prove Theorem~\ref{thm:convergenceresultsdiscretization}  $(iii)$. More precisely,  let $\theta$ be the logarithmic mean \eqref{eq:logmean}, $\lambda=\Leb_{|\Lambda_n}$, and assume that $\kappa:\Gamma_{\Lambda_n}\times\Lambda_n\to[0,\infty)$ satisfies \ref{property:cocycle}, \hyperref[property:locallystable]{$(\mathbf{S}(1))$}, \ref{property:lowerbound}, \ref{property:continuity2}, and~\ref{property:repulsive}. We show that  $$\Gamma\text{-}\lim_{M\to\infty}\mathcal{W}_{\kappa^M,\lambda^M}(\cdot,\cdot)
    =\mathcal{W}_{\kappa,\lambda}(\cdot,\cdot)$$ 
   on $\PPP_1(\Gamma_{\Lambda_n})\times\{\dP\in\PPP_1(\Gamma_{\Lambda_n}):\dP\ll\dR\}$ (or with the two spaces interchanged). A recovery sequence for $\PPP_1(\Gamma_{\Lambda_n})\times\{\dP\in\PPP_1(\Gamma_{\Lambda_n}):\dP\ll\dR\}$ is  $(\dP_0^M,\dP_1^M)_{M\in\N}$.\\

We prove the $\Gamma$-convergence in two steps. The limsup inequality is established in Lemma~\ref{lem:upperbounddistance}, and the liminf inequality in Lemma~\ref{lem:lowersemicontinuitydiscretization}.
Combining these two lemmas yields the proof of Theorem~\ref{thm:convergenceresultsdiscretization}~$(iii)$.

        \subsubsection{Limsup Inequality}  
       \begin{lem}\label{lem:upperbounddistance}
        Let $\lambda=\Leb_{|\Lambda_n}$ and assume that $\kappa:\Gamma_{\Lambda_n}\times\Lambda_n\to[0,\infty)$ satisfies \ref{property:cocycle},  \hyperref[property:locallystable]{$(\mathbf{S}(1))$}, \ref{property:lowerbound},  \ref{property:continuity2}, and~\ref{property:repulsive}. Further, let $\dR=\RRR_{\kappa,\lambda}(\Gamma_{\Lambda_n})$,  $\dP_0,\dP_1\in\PPP_1(\Gamma_{\Lambda_n})$ s.t.\   $\dP_0\ll\dR$ or $\dP_1\ll\dR$. It holds that \begin{align*}
              \limsup_{M\to\infty} \mathcal{W}_{\kappa^M,\lambda^M}(\dP_0^M,\dP_1^M)\leq\mathcal{W}_{\kappa,\lambda}(\dP_0,\dP_1).
           \end{align*}
       \end{lem}
        The general proof strategy for Lemma \ref{lem:upperbounddistance} is as follows. We select a solution to the continuity equation minimizing the optimization problem in the definition of $\mathcal{W}_{\kappa,\lambda}(\dP_0,\dP_1)$. We then discretize this minimizer to obtain a competitor for 
$\mathcal{W}_{\kappa^M,\lambda^M}(\dP_0^M,\dP_1^M)$. Finally, we pass to the limit and obtain the desired upper bound.

        For this strategy to work, it is essential  that discretizations of solutions to the continuity equation on $\Lambda_n$ again satisfy the continuity equation on $\Lambda_n^M$.

        \begin{lem}\label{lem:discretizedCE}
       Let $\dP_0,\dP_1\in\PPP_1(\Gamma_{\Lambda_n})$ and 
       $(\Bar{\dP},\Bar{\dV})\in\CCC\EEE^{\Lambda_n}(\dP_0,\dP_1)$. For $\bar{\dV}_t\in\MMM_{b,0}(\Gamma_{\Lambda_n}\times\Lambda_n)$ define \begin{align*}
           \Bar{\dV}_t^M:=(T^M_\Gamma,T^M)_\#\Bar{\dV}_t.
       \end{align*} If $x^M+\Lambda_{\frac{1}{M}}$ are $\dP_0,\dP_1$-continuity sets for all $x^M\in\Lambda_n^M$, then the pair $(\Bar{\dP}^M,\ind_{\{\xi^M(\{x^M\})=0\}}\Bar{\dV}^M)$ satisfies the continuity equation~\eqref{eq:continuityequation} and there exists a continuous representative $\Tilde{\dP}^{M}$ of $\Bar{\dP}^{M}$ s.t.\   $$(\Tilde{\dP}^{M},\ind_{\{\xi^M(\{x^M\})=0\}}\Bar{\dV}^M)\in\CCC\EEE^{\Lambda_n^M}(\dP_0^M,\dP_1^M).$$ 
    \end{lem}
    \begin{proof}
          We start the proof by showing that $(\Bar{\dP}^M,\ind_{\{\xi^M(\{x^M\})=0\}}\Bar{\dV}^M)$ satisfies the continuity equation~\eqref{eq:continuityequation}. Let  $f\in\CCC_c^+(\Lambda_n^M)$ and $\phi\in \mathcal{C}_c^\infty((0,1)).$ Since $({\Bar{\dP}},{\Bar{\dV}})\in\CCC\EEE^{\Lambda_n}(\dP_0, \dP_1)$ and $D_xe^{-\int f(y)T_{\Gamma}^M(\xi)(\mathrm{d}y)}\in\FFF_{b,0}(\Gamma\times\R^d)$, we obtain by  Theorem~\ref{thm:propertiesce} $(ii)$ that
             \begin{align*}
                 &\int_0^1 \dt \phi(t)\int e^{-\int f(y^M)\xil(\mathrm{d}y^M)}\Bar{\dP}_t^{M}(\mathrm{d}\xil)\mathrm{d}t \\
                 &\quad= \int_0^1 \dt \phi(t)\int e^{-\int f(y^M)T_{\Gamma}^M(\xi)(\mathrm{d}y^M)}\Bar{\dP}_t(\mathrm{d}\xi)\mathrm{d}t \\
                 &\quad= -  \int_0^1  \phi(t)\int D_x e^{-\int f(y^M)T_{\Gamma}^M(\xi)(\mathrm{d}y^M)}\Bar{\dV}_t(\mathrm{d}\xi,\mathrm{d}x)\mathrm{d}t \\
                 &\quad= -  \int_0^1  \phi(t)\int \ind_{\{T_{\Gamma}^M(\xi)(T^M(x))=0\}}D_{T^M(x)} e^{-\int f(y^M)T_{\Gamma}^M(\xi)(\mathrm{d}y^M)}\Bar{\dV}_t(\mathrm{d}\xi,\mathrm{d}x)\mathrm{d}t\\
                 &\quad=  -\int_0^1  \phi(t)\int D_{\xl} e^{-\int f(y^M)\xil(\mathrm{d}y^M)}\ind_{\{\xil(\{\xl\})=0\}}\Bar{\dV}_t^{M}(\mathrm{d}\xil,\mathrm{d}\xl)\mathrm{d}t,
             \end{align*}
    proving the first claim. 
    
     By  Theorem~\ref{thm:propertiesce} $(v)$, we can find a unique continuous representative $\Tilde{\dP}^M$ of ${\Bar{\dP}}^{M}$.  Property $(iv)$ of Definition \ref{dfn:ce} is satisfied by the definition of $\ind_{\{\xi^M(\{x^M\})=0\}}\Bar{\dV}^M$. It remains to prove that $\dP_0^M=\Tilde{\dP}_0^M$ and $\dP_1^M=\Tilde{\dP}_1^M$. Let $(t_m)_{m\in\N}$ s.t.\   $\lim_{m\to\infty}t_m=0$ and   $\Tilde{\dP}_{t_m}^M=\Bar{\dP}_{t_m}^M$ for all $m\in\N$. As $x^M+\Lambda_{\frac{1}{M}}$ is a $\dP_0$-continuity set for any $x^M\in\Lambda_n^M$ and $\Bar{\dP}$ is $\PPP_1(\Gamma_{\Lambda_n})$-continuous, we obtain by Proposition \ref{prop:convergencep1} that 
     \begin{align*}
        &\int  \ind_{\{\xi^M(x^M_1)=0,\ldots,\xi^M(x^M_j)=0,\xi^M(x^M_{j+1})=1,\ldots,\xi^M(x^M_{j+k})=1\}}\tilde\dP^M_{0}(\mathrm{d}\xi^M)\\ 
        &\quad=\lim_{m\to\infty}\int  \ind_{\{\xi^M(x^M_1)=0,\ldots,\xi^M(x^M_j)=0,\xi^M(x^M_{j+1})=1,\ldots,\xi^M(x^M_{j+k})=1\}}\tilde{\dP}^M_{t_m}(\mathrm{d}\xi^M)\\
        &\quad=\lim_{m\to\infty}\int  \ind_{\{\xi^M(x^M_1)=0,\ldots,\xi^M(x^M_j)=0,\xi^M(x^M_{j+1})=1,\ldots,\xi^M(x^M_{j+k})=1\}}\Bar\dP^M_{t_m}(\mathrm{d}\xi^M)\\
         &\quad= \lim_{m\to\infty}\int  \ind_{\{\xi^M(x^M_1)=0,\ldots,\xi^M(x^M_j)=0,\xi^M(x^M_{j+1})\geq1,\ldots,\xi^M(x^M_{j+k})\geq1\}}\Bar{\dP}_{t_m}(\mathrm{d}\xi)\\
         &\quad=\lim_{m\to\infty}\int  \ind_{\{\xi(x^M_1+\Lambda_{\frac{1}{M}})=0,\ldots,\xi(x^M_j+\Lambda_{\frac{1}{M}})=0,\xi(x^M_{j+1}+\Lambda_{\frac{1}{M}})\geq1,\ldots,\xi(x^M_{j+k}+\Lambda_{\frac{1}{M}})\geq1\}}\Bar\dP_{t_m}(\mathrm{d}\xi)\\
         &\quad= \int  \ind_{\{\xi(x^M_1+\Lambda_{\frac{1}{M}})=0,\ldots,\xi(x^M_j+\Lambda_{\frac{1}{M}})=0,\xi(x^M_{j+1}+\Lambda_{\frac{1}{M}})\geq1,\ldots,\xi(x^M_{j+k}+\Lambda_{\frac{1}{M}})\geq1\}}\dP_0(\mathrm{d}\xi)\\
         &\quad= \int  \ind_{\{\xi^M(x^M_1)=0,\ldots,\xi^M(x^M_j)=0,\xi^M(x^M_{j+1})=1,\ldots,\xi^M(x^M_{j+k})=1\}}\dP^M_0(\mathrm{d}\xi^M).
     \end{align*}
     A similar argument for the intensity measures proves that $\I_{\dP_0^M}=\I_{\tilde\dP_0^M}$.  
       Hence, by Proposition \ref{prop:convergencep1}   we have $\dP_0^M=\Tilde{\dP}_0^M$. Analogously, it follows that  $\dP_1^M=\Tilde{\dP}_1^M$, finishing the proof.
    \end{proof}
        With these preparations in place, we now prove Lemma~\ref{lem:upperbounddistance}.
       \begin{proof}[Proof of Lemma~\ref{lem:upperbounddistance}]
       If $\mathcal{W}_{\kappa,\lambda}(\dP_0,\dP_1)=\infty$, the claim is immediate. Hence, let $\mathcal{W}_{\kappa,\lambda}(\dP_0,\dP_1)<\infty$.

           By Theorem~\ref{thm:minimizer}, there exists  a solution $(\bar\dP,\bar{\dV})\in\CCC\EEE^{\Lambda_n}(\dP_0,\dP_1)$ s.t.\  
           \begin{align*} \mathcal{W}_{\kappa,\lambda}(\dP_0,\dP_1)^2=\int_0^1\mathcal{L}_{\kappa,\lambda}(\bar{\dP}_t,\bar{\dV}_t)\mathrm{d}t.
           \end{align*}
         By Lemma~\ref{lem:discretizedCE}, the pair
$(\bar\dP^M,\ind_{\{\xi^M(\{x^M\})=0\}}\bar{\dV}^M)$ solves the continuity equation \eqref{eq:continuityequation} and, up to choosing a continuous representative of $\bar\dP^M$, 
\(
(\bar\dP^M,\ind_{\{\xi^M(\{x^M\})=0\}}\bar{\dV}^M)\in\CCC\EEE^{\Lambda_n^M}(\dP_0^M,\dP_1^M).
\) Thus, 
          \begin{align}
               &\limsup_{M\to\infty} \mathcal{W}_{\kappa^M,\lambda^M}^2(\dP_0^M,\dP_1^M)\leq \limsup_{M\to\infty}\int_0^1\mathcal{L}_{\kappa^M,\lambda^M}(\bar{\dP}_t^M,\ind_{\{\xi^M(\{x^M\})=0\}}\bar{\dV}_t^M)\mathrm{d}t.\label{eq:1113}
          \end{align}
            Hence, it remains to prove that 
          \begin{align*} \limsup_{M\to\infty}\int_0^1\mathcal{L}_{\kappa^M,\lambda^M}(\bar{\dP}_t^M,\ind_{\{\xi^M(\{x^M\})=0\}}\bar{\dV}_t^M)\mathrm{d}t\leq \int_0^1\mathcal{L}_{\kappa,\lambda}(\bar{\dP}_t,\bar{\dV}_t)\mathrm{d}t.
          \end{align*} 
           Let $$\sigma_t=\kappa(\bar{\dP}_t\otimes\lambda)+C_{\bar{\dP}_t}+|\bar{\dV}_t|$$ and denote  $\kappa(\bar{\dP}_t\otimes\lambda)=\psi_t\sigma_t,\ C_{\bar{\dP}_t}=\phi_t\sigma_t, \text{ and } \bar{\dV}_t=w_t\sigma_t.$
            Since $\bar{\dP}_0\ll\dR$ or $\bar{\dP}_1\ll\dR$ also $\bar{\dP}_t\ll\dR$  (see Lemma \ref{lem:preservationabsolutecontinuity}),   $C_{\bar{\dP}_t}\ll\kappa(\dR\otimes{\lambda})$ (Corollary~\ref{kor:CampbellPapangelou}), and $\bar{\dV}_t\ll\kappa(\dR\otimes\lambda)$ (Lemma~\ref{lem:absolutecontinuityvelocity}). Hence, we conclude that $\sigma_t\ll\kappa(\dR\otimes{\lambda}).$ 
           By disintegration, we find measures $\sigma^{*,M}_t(\cdot|\xi^M,x^M),\; \xi^M\in\Gamma_{\Lambda_n^M},\; x^M\in \Lambda_n^M$, on $\Gamma_{\Lambda_n}\times \Lambda_n$ s.t.\   
           \begin{align*}
               \int F(\xi,x)\sigma_t(\mathrm{d}\xi,\mathrm{d}x)=\int\int F(\eta,y)\sigma^{*,M}_t(\mathrm{d}\eta,\mathrm{d}y|\xi^M,x^M)\sigma^M_t(\mathrm{d}\xi^M,\mathrm{d}x^M)
           \end{align*}
           for all $F\in\FFF^+(\Gamma_{\Lambda_n}\times \Lambda_n)$. Since $\sigma_t\ll \kappa(\dR\otimes{\lambda})$ and $\frac 1 a\leq\kappa(\xi,x)$  $\kappa(\dR\otimes{\lambda})$-a.e.\ by \ref{property:lowerbound} and Remark~\ref{bem:lowerboundcampbell}, for any $K\in\BBB(\Gamma_{\Lambda_n^M})$ and $B\in\BBB(\Lambda_n^M)$, it holds that 
         \begin{align*}
             &\int\int\ind_{K\times B}(\xi^M,x^M)\kappa^M(\xi^M,x^M)\bar\dP^M_t(\mathrm{d}\xi^M)\lambda^M(\mathrm{d}x^M)\\
             &\quad= \int\int\ind_{K\times B}(\xi^M,x^M) \kappa^M(\xi^M,x^M)\bar\dP_t(\mathrm{d}\xi)\lambda(\mathrm{d}x)\\
             &\quad=  \int\ind_{K\times B}(\xi^M,x^M) \frac{\kappa^M(\xi^M,x^M)}{\kappa(\xi,x)}\psi_t(\xi,x)\sigma_t(\mathrm{d}\xi,\mathrm{d}x)\\
             &\quad= \int\ind_{K\times B}(\xi^M,x^M) \kappa^M(\xi^M,x^M)\int \frac{1}{\kappa(\eta,y)}\psi_t(\eta,y)\sigma^{*,M}_t(\mathrm{d}\eta,\mathrm{d}y|\xi^M,x^M)\sigma^M_t(\mathrm{d}\xi^M,\mathrm{d}x^M).
         \end{align*}
         This implies that \begin{align*}
             \kappa^M(\xi^M,x^M)(\bar\dP^M_t\otimes\lambda^M)(\mathrm{d}\xi^M,\mathrm{d}x^M)&=\kappa^M(\xi^M,x^M)\int \frac{1}{\kappa(\eta,y)}\psi_t(\eta,y)\sigma^{*,M}_t(\mathrm{d}\eta,\mathrm{d}y|\xi^M,x^M)\sigma^M_t(\mathrm{d}\xi^M,\mathrm{d}x^M)\\
         &=:\psi^{M}_t(\xi^M,x^M)\sigma^M_t(\mathrm{d}\xi^M,\mathrm{d}x^M).\end{align*} 
         Similarly, we find that 
         \begin{align*}
            \ind_{\{\xi^M(\{x^M\})=0\}}\bar{\dV}^M_t(\mathrm{d}\xi^M,\mathrm{d}x^M)&=\ind_{\{\xi^M(\{x^M\})=0\}}\int w_t(\eta,y)\sigma^{*,M}_t(\mathrm{d}\eta,\mathrm{d}y|\xi^M,x^M)\sigma^M_t(\mathrm{d}\xi^M,\mathrm{d}x^M)\\
            &=:w^{M}_t(\xi^M,x^M)\sigma^M_t(\mathrm{d}\xi^M,\mathrm{d}x^M).
         \end{align*}
         To compute the density of the reduced Campbell measure $C_{\bar\dP^M_t}$, we recall the definition of the function
         \begin{align*}
             n_{\xl}^M(\xi)=\begin{cases}
                 \xi(\xl+\Lambda_{\frac{1}{M}}), &\text{if } \xi(\xl+\Lambda_{\frac{1}{M}})>0,\\
                 \infty, &\text{otherwise, } 
             \end{cases}
             \quad \xl\in\Lambda_n^M,\;\xi\in\Gamma_{\Lambda_n},
         \end{align*}
         see \eqref{eq:definitionnumberofpointsinabox}.
         For any $K\in\BBB(\Gamma_{\Lambda_n^M})$ and $B\in\BBB(\Lambda_n^M)$, it holds that 
         \begin{align*}
             \int \ind_{K\times B}C_{\bar\dP^M_t}(\mathrm{d}\xi^M,\mathrm{d}x^M)&=\int\int \ind_{K\times B}(\xi^M-\delta_{x^M},x^M)\xi^M(\mathrm{d}x^M)\bar\dP_t(\mathrm{d}\xi)\\
             &= \int\int \ind_{K\times B}(\xi^M-\delta_{x^M},x^M)\frac{1}{n_{x^M}^M(\xi)}\xi(\mathrm{d}x)\bar\dP_t(\mathrm{d}\xi)\\
             &= \int \ind_{K\times B}((\xi+\delta_x)^M-\delta_{x^M},x^M)\frac{1}{n_{x^M}^M(\xi+\delta_x)}C_{\bar\dP_t}(\mathrm{d}\xi,\mathrm{d}x)\\
             &\geq \int \ind_{K\times B}(\xi^M,x^M)\ind_{\{\xi^M(\{x^M\})=0\}}C_{\bar\dP_t}(\mathrm{d}\xi,\mathrm{d}x)\\
             &= \int \ind_{K\times B}(\xi^M,x^M)\ind_{\{\xi^M(\{x^M\})=0\}}\int\phi_t(\eta,y)\sigma^{*,M}_t(\mathrm{d}\eta,\mathrm{d}y|\xi^M,x^M)\sigma^M_t(\mathrm{d}\xi^M,\mathrm{d}x^M),
         \end{align*}
         which implies that 
         \begin{align*}
             C_{\bar\dP^M_t}(\mathrm{d}\xi^M,\mathrm{d}x^M)&\geq\ind_{\{\xi^M(\{x^M\})=0\}}\int\phi_t(\eta,y)\sigma^{*,M}_t(\mathrm{d}\eta,\mathrm{d}y|\xi^M,x^M)\sigma^M_t(\mathrm{d}\xi^M,\mathrm{d}x^M)\\
             &=:\phi^M_t(\xi^M,x^M) \sigma^M_t(\mathrm{d}\xi^M,\mathrm{d}x^M).
         \end{align*}
         Using that $\alpha$ is monotone decreasing in the first two components and  Jensen's inequality together with the convexity and the homogeneity of $\alpha,$ we obtain that 
         \begin{align}
             \mathcal{L}_{\kappa^M,\lambda^M}(\bar\dP^M_t,\ind_{\{\xi^M(\{x^M\})=0\}}\bar\dV^M_t)\nonumber&\leq\int\alpha(\psi^M_t(\xi^M,x^M),\phi^M_t(\xi^M,x^M),w^M_t(\xi^M,x^M))\sigma^M_t(\mathrm{d}\xi^M,\mathrm{d}x^M)\nonumber\\
             &\leq\!\int \!\alpha\left(\!\frac{\kappa^M(\xi^M,x^M)}{\kappa(\xi,x)}\psi_t(\xi,x), \phi_t(\xi,x), w_t(\xi,x)\right)\ind_{\{\xi^M(\{x^M\})=0\}}\sigma_t(\mathrm{d}\xi,\mathrm{d}x). \label{eq:2030}
         \end{align}
         By \ref{property:continuity2}, $\lim_{M\to\infty} \kappa^M(\xi^M,x^M)=\kappa(\xi,x)$ $\dR\otimes\lambda$-a.e.  Using the monotonicity of $\alpha$, \ref{property:locallystable}, and~\ref{property:lowerbound}, we have that 
        \begin{align*}
            \alpha\left(\frac{\kappa^M(\xi^M,x^M)}{\kappa(\xi,x)}\psi_t(\xi,x), \phi_t(\xi,x), w_t(\xi,x)\right)\ind_{\{\xi^M(\{x^M\})=0\}}&\leq \alpha\left(\frac{ 1}{Aa }\psi_t(\xi,x), \phi_t(\xi,x), w_t(\xi,x)\right)\\
            &\leq Aa  \alpha\left(\psi_t(\xi,x), \phi_t(\xi,x), w_t(\xi,x)\right)
        \end{align*}
        $\kappa(\dR\otimes\lambda)$-a.e.
        Since \begin{align*}
            \int_0^1\int Aa \alpha\left(\psi_t(\xi,x), \phi_t(\xi,x), w_t(\xi,x)\right)\sigma_t(\mathrm{d}\xi,\mathrm{d}x)\mathrm{d}t&=Aa  \int_0^1\mathcal{L}_{\kappa,\lambda}(\bar{\dP}_t,\bar{\dV}_t)\mathrm{d}t\leq Aa  \int_0^1\mathcal{L}_{\kappa,\lambda}(\bar{\dP}_t,\bar{\dV}_t)\mathrm{d}t<\infty,
        \end{align*}
       and $\sigma_t\ll\kappa(\dR\otimes\lambda)$, we obtain by the dominated convergence theorem and equation~\eqref{eq:2030} 
        \begin{align}
           & \limsup_{M\to\infty}\int_0^1\mathcal{L}_{\kappa^M,\lambda^M}(\bar{\dP}_t^M,\ind_{\{\xi^M(\{x^M\})=0\}}\bar{\dV}_t^M)\mathrm{d}t\nonumber\\
            &\quad\leq \int_0^1\int \limsup_{M\to\infty} \alpha\left(\frac{\kappa^M(\xi^M,x^M)}{\kappa(\xi,x)}\psi_t(\xi,x), \phi_t(\xi,x), w_t(\xi,x)\right)\ind_{\{\xi^M(\{x^M\})=0\}}\sigma_t(\mathrm{d}\xi,\mathrm{d}x)\mathrm{d}t\nonumber\\ 
            &\quad\leq \int_0^1\int \alpha\left(\psi_t(\xi,x), \phi_t(\xi,x), w_t(\xi,x)\right)\sigma_t(\mathrm{d}\xi,\mathrm{d}x)\mathrm{d}t \nonumber\\
            &\quad=\int_0^1\mathcal{L}_{\kappa,\lambda}(\bar{\dP}_t,\bar{\dV}_t)\mathrm{d}t.\label{eq:1112}
        \end{align}
        All in all, combining \eqref{eq:1113} and \eqref{eq:1112} and using that $(\bar{\dP},\bar{\dV})$ is a minimizer of the optimization problem of  $\mathcal{W}_{\kappa,\lambda}(\dP_0,\dP_1)$, we obtain that 
        \begin{align*}
            \limsup_{M\to\infty} \mathcal{W}_{\kappa^M,\lambda^M}^2(\dP_0^M,\dP_1^M)&\leq \limsup_{M\to\infty}\int_0^1\mathcal{L}_{\kappa^M,\lambda^M}(\bar{\dP}_t^M,\ind_{\{\xi^M(\{x^M\})=0\}}\bar{\dV}_t^M)\mathrm{d}t\\
            &\leq \int_0^1\mathcal{L}_{\kappa,\lambda}(\bar{\dP}_t,\bar{\dV}_t)\mathrm{d}t\\
            &= \mathcal{W}_{\kappa,\lambda}^2(\dP_0,\dP_1),
        \end{align*}
        finishing the proof.
       \end{proof}

       \subsubsection{Liminf Inequality}
We continue with the proof of the liminf inequality for the $\Gamma$-convergence of $(\mathcal{W}_{\kappa^M,\lambda^M})_{M\in\N}$. We first prove the corresponding result for the Lagrange functional in Lemma~\ref{lem:lscwithvariablekappa}. Thereafter, we establish a compactness result for the action functional in Lemma~\ref{lem:compactnessdiscretization}, which allows us to lift the liminf inequality to the transport distance.
         \begin{lem}\label{lem:lscwithvariablekappa}  Let $\lambda=\Leb_{|\Lambda_n}$ and assume that $\kappa:\Gamma_{\Lambda_n}\times\Lambda_n\to[0,\infty)$ satisfies \ref{property:locallystable} and~\ref{property:continuity}.
            Let $\dP_M\in\PPP_1(\Gamma_{\Lambda_n^M})$, $M\in\N$, and $\dP\in\PPP_1(\Gamma_{\Lambda_n})$ s.t.\   $\lim_{M\to\infty}\dP_M=\dP$ in $\PPP_1(\Gamma_{\Lambda_n})$. Further, assume that $\dV_M\in\MMM_{b,0}(\Gamma_{\Lambda_n^M}\times\Lambda_n^M)$, $M\in\N$, and $\dV\in\MMM_{b,0}(\Gamma_{\Lambda_n}\times\Lambda_n)$ s.t.\   $\lim_{M\to\infty}\dV_M=\dV$ in $\MMM_{b,0}(\Gamma_{\Lambda_n}\times\Lambda_n).$ It holds that 
             \begin{align*}
                 \mathcal{L}_{\kappa,\lambda}(\dP,\dV)\leq\liminf_{M\to\infty}\mathcal{L}_{\kappa^M,\lambda^M}(\dP_M,\dV_M). 
             \end{align*}
         \end{lem}
         \begin{proof}
        First, we assume that $\kappa$ not only satisfies \ref{property:continuity} but is even continuous.  If $F\in\CCC_{b,0}(\Gamma_{\Lambda_n}\times \Lambda_n)$, then $\kappa F\in\CCC_{b,0}(\Gamma_{\Lambda_n}\times\Lambda_n)$ as $\kappa$ is continuous and bounded above by \ref{property:locallystable}.  Using additionally that $\dP$ is simple, we deduce that \begin{align*}
            &\lim_{M\to\infty}\int F(\xi^M,x^M)\kappa^M(\xi^M,x^M)\dP_M\otimes\lambda^M(\mathrm{d}\xi^M,\mathrm{d}x^M)\\
            &\quad=  \lim_{M\to\infty}\int F(\xi^M,x^M)\kappa(\xi^M,x^M)\ind_{\{\xi^M(\{x^M\})=0\}}\dP_M\otimes\lambda^M(\mathrm{d}\xi^M,\mathrm{d}x^M)\\
            &\quad= \int F(\xi,x)\kappa(\xi,x)\dP\otimes\lambda(\mathrm{d}\xi,\mathrm{d}x).
        \end{align*} 
         Thus, $\lim_{M\to\infty}\kappa^M(\dP_M\otimes\lambda^M)=\kappa(\dP\otimes\lambda)$ in $\MMM_{b,0}(\Gamma_{\Lambda_n}\times\Lambda_n)$. As $\lim_{M\to\infty}C_{\dP_M}= C_\dP$ in $\MMM_{b,0}(\Gamma_{\Lambda_n}\times\Lambda_n)$  by  Proposition~\ref{prop:convergencep1},  Theorem~\ref{athm:lowersemicontinuity} implies that
        \begin{align*}
            \mathcal{L}_{\kappa,\lambda}(\dP,\dV)&=\int\alpha\left(\frac{\mathrm{d}\kappa(\dP\otimes\lambda)}{\mathrm{d}\sigma}(\xi,x),\frac{\mathrm{d}C_{\dP}}{\mathrm{d}\sigma}(\xi,x),\frac{\mathrm{d}\dV}{\mathrm{d}\sigma}(\xi,x)\right)\sigma(\mathrm{d}\xi,\mathrm{d}x)\\
            &\leq\liminf_{M\to\infty}\int\alpha\left(\frac{\mathrm{d}\kappa^M(\dP_M\otimes\lambda^M)}{\mathrm{d}\sigma}(\xi,x),\frac{\mathrm{d}C_{\dP_M}}{\mathrm{d}\sigma}(\xi,x),\frac{\mathrm{d}\dV_M}{\mathrm{d}\sigma}(\xi,x)\right)\sigma(\mathrm{d}\xi,\mathrm{d}x)\\
            &=\liminf_{M\to\infty}\mathcal{L}_{\kappa^M,\lambda^M}(\dP_M,\dV_M),
        \end{align*}
        finishing the proof under the additional assumption that $\kappa$ is continuous. Now, let $\kappa$ be upper semicontinuous \ref{property:continuity}. There exists a sequence of functions $(\kappa_i)_{i\in\N}\subset\CCC_b(\Gamma_{\Lambda_n}\times \Lambda_n)$ s.t.\   $\kappa_i\downarrow\kappa$ for $i\to\infty.$ By the monotone convergence theorem, the monotonicity of $\alpha$, and the previous result for the continuous functions $\kappa_i$, we obtain 
       \begin{align*}
           \mathcal{L}_{\kappa,\lambda}(\dP,\dV)&=\int \alpha\left(\frac{\mathrm{d}\kappa(\dP\otimes\lambda)}{\mathrm{d}\sigma}(\xi,x),\frac{\mathrm{d}C_\dP}{\mathrm{d}\sigma}(\xi,x), \frac{\mathrm{d}\dV}{\mathrm{d}\sigma}(\xi,x)\right)\sigma(\mathrm{d}\xi,\mathrm{d}x)\\
           &=\sup_{i\in\N}\int \alpha\left(\frac{\mathrm{d}\kappa_i(\dP\otimes\lambda)}{\mathrm{d}\sigma}(\xi,x),\frac{\mathrm{d}C_\dP}{\mathrm{d}\sigma}(\xi,x), \frac{\mathrm{d}\dV}{\mathrm{d}\sigma}(\xi,x)\right)\sigma(\mathrm{d}\xi,\mathrm{d}x)\\
           &\leq \sup_{i\in\N}\liminf_{M\to\infty}\int \alpha\left(\frac{\mathrm{d}\kappa_i^M(\dP_M\otimes\lambda^M)}{\mathrm{d}\sigma}(\xi,x),\frac{\mathrm{d}C_{\dP_M}}{\mathrm{d}\sigma}(\xi,x), \frac{\mathrm{d}\dV_M}{\mathrm{d}\sigma}(\xi,x)\right)\sigma(\mathrm{d}\xi,\mathrm{d}x)\\
           &\leq \liminf_{M\to\infty}\sup_{i\in\N}\int \alpha\left(\frac{\mathrm{d}\kappa_i^M(\dP_M\otimes\lambda^M)}{\mathrm{d}\sigma}(\xi,x),\frac{\mathrm{d}C_{\dP_M}}{\mathrm{d}\sigma}(\xi,x), \frac{\mathrm{d}\dV_M}{\mathrm{d}\sigma}(\xi,x)\right)\sigma(\mathrm{d}\xi,\mathrm{d}x)\\
           &= \liminf_{M\to\infty} \mathcal{L}_{\kappa^M,\lambda^M}(\dP_M,\dV_M),
       \end{align*}
       finishing the proof.
        \end{proof}
\begin{lem}\label{lem:compactnessdiscretization}
   Let $\lambda=\Leb_{|\Lambda_n}$ and assume that $\kappa:\Gamma_{\Lambda_n}\times\Lambda_n\to[0,\infty)$ satisfies \ref{property:locallystable} and~\ref{property:continuity}. Further, let   $\dP_{0,M},\dP_{1,M}\in\PPP_1(\Gamma_{\Lambda_n^M}),\;M\in\N,$ and $\dP_0,\dP_1\in\PPP_1(\Gamma_{\Lambda_n})$ s.t.\  $\lim_{M\to\infty}\dP_{0,M}=\dP_0$ and $\lim_{M\to\infty}\dP_{1,M}=\dP_1$ in $\PPP_1(\Gamma_{\Lambda_n})$. If $(\bar{\dP}_M,\bar{\dV}_M)\in\CCC\EEE^{\Lambda_n^M}(\dP_{0,M},\dP_{1,M}),$ $M\in\N,$ s.t.\ 
\begin{align*}
    \sup_{M\in\N}\mathcal{A}_{\kappa^M,\lambda^M}(\bar{\dP}_M,\bar{\dV}_M)\leq C<\infty,
\end{align*}
        then there exists a subsequence $(\bar{\dP}_{M_N},\bar{\dV}_{M_N})_{N\in\N}$ and a solution $(\bar{\dP},\bar{\dV})\in\CCC\EEE^{\Lambda_n}(\dP_0,\dP_1)$ s.t.\   
        \begin{align*}
            \lim_{N\to\infty}\bar{\dP}_{M_N,t}=\bar{\dP}_t\quad &\text{ in }\PPP_1(\Gamma_{\Lambda_n}) \text{ for all } t\in [0,1],\\
            \lim_{N\to\infty}\bar{\dV}_{M_N}=\bar{\dV}\quad &\text{ in }\MMM_{b,0}(\Gamma_{\Lambda_n}\times {\Lambda_n}\times[0,1]).
             \end{align*}
    \end{lem}
     \begin{proof}
         The proof is analogous to the proof of Lemma~\ref{lem:compactnessofsolutions} $(iii)$ by additionally using that 
         \begin{align*}
             \lambda^{M}(B)\leq\lambda^{out}(B):=\sum_{z\in\Z^d:B\cap\Lambda_1+z\not=\emptyset}\lambda(\Lambda_1+z).
         \end{align*}
     \end{proof}   
\begin{lem}\label{lem:lowersemicontinuitydiscretization}
Let $\lambda=\Leb_{|\Lambda_n}$ and assume that $\kappa:\Gamma_{\Lambda_n}\times\Lambda_n\to[0,\infty)$ satisfies \ref{property:cocycle}, \ref{property:locallystable}, and~\ref{property:continuity}. Let  $\dP_{0,M},\dP_{1,M}\in\PPP_1(\Gamma_{\Lambda_n^M}),\;M\in\N,$ and $\dP_0,\dP_1\in\PPP_1(\Gamma_{\Lambda_n})$ such that   $\lim_{M\to\infty}\dP_{0,M}=\dP_0$ and  $\lim_{M\to\infty}\dP_{1,M}=\dP_1$ in $\PPP_1(\Gamma_{\Lambda_n})$. It holds that \begin{align*}
                \mathcal{W}_{\kappa,\lambda}(\dP_0,\dP_1)\leq \liminf_{M\to\infty}\mathcal{W}_{\kappa^M,\lambda^M}(\dP_{0,M},\dP_{1,M}).
            \end{align*}
            \end{lem}
            \begin{proof}
            If $\liminf_{M\to\infty}\mathcal{W}_{\kappa^M,\lambda^M}(\dP_{0,M},\dP_{1,M})=\infty$, the claim is trivial. Hence, let 
            \begin{align*}
                \liminf_{M\to\infty}\mathcal{W}_{\kappa^M,\lambda^M}(\dP_{0,M},\dP_{1,M})\leq C<\infty.
            \end{align*}
            Consider an increasing sequence of natural numbers $(M_N)_{N\in\N}$ s.t.\   
                \begin{align*}
                    \liminf_{M\to\infty}\mathcal{W}_{\kappa^M,\lambda^M}(\dP_{0,M},\dP_{1,M})=\lim_{N\to\infty} \mathcal{W}_{\kappa^{M_N},\lambda^{M_N}}(\dP_{0,M_N},\dP_{1,M_N}).
                \end{align*}
          Theorem~\ref{thm:minimizer} implies that there exist solutions $(\bar{\dP}_{M_N},\bar{\dV}_{M_N})\in\CCC\EEE^{\Lambda_n^{M_N}}(\dP_{0,M_N},\dP_{1,M_N})$ s.t.\   
                \begin{align*} {\mathcal{W}_{\kappa^{M_N},\lambda^{M_N}}^2}(\dP_{0,M_N},\dP_{1,M_N})=\int_0^1\mathcal{L}_{\kappa^{M_N},\lambda^{M_N}}(\bar{\dP}_{t,M_N},\bar{\dV}_{t,M_N})\mathrm{d}t\leq C^2+1.
                \end{align*}
                
            Since additionally $\lim_{M\to\infty}\dP_{0,M}=\dP_0$ and $\lim_{M\to\infty}\dP_{1,M}=\dP_1$ in $\PPP_1(\Gamma_{\Lambda_n})$, by Lemma~\ref{lem:compactnessdiscretization} there exists a subsequence $(M_{N_O})_{O\in\N}$ of $(M_N)_{N\in\N}$ and a solution $(\bar\dP,\bar\dV)\in\CCC\EEE^{\Lambda_n}(\dP_0,\dP_1)$ s.t.\   
               \begin{align*} \lim_{O\to\infty}\bar{\dP}_{t,M_{N_O}}=\bar{\dP}_t\quad &\text{ in }\PPP_1(\Gamma_{\Lambda_n})  \text{ for all } t\in [0,1],\\
            \lim_{O\to\infty}\bar{\dV}_{M_{N_O}}=\bar{\dV}\quad &\text{ in }\MMM_{b,0}(\Gamma_{\Lambda_n}\times {\Lambda_n}\times[0,1]).
             \end{align*}
                Using  Lemma~\ref{lem:lscwithvariablekappa} and Fatou's lemma, we obtain that 
            \begin{align*}
                {\mathcal{W}_{\kappa,\lambda}^2}(\dP_0,\dP_1)&\leq \liminf_{O\to\infty}\int_0^1\mathcal{L}_{\kappa^{M_{N_O}},\lambda^{M_{N_O}}}(\bar\dP_{t,M_{N_O}},\bar\dV_{t,M_{N_O}})\mathrm{d}t\\ &=\liminf_{O\to\infty}{\mathcal{W}^2_{\kappa^{M_{N_O}},\lambda^{M_{N_O}}}}(\dP_{0,M_{N_O}},\dP_{1,M_{N_O}})\\
                &=\lim_{N\to\infty}{\mathcal{W}^2_{\kappa^{M_{N}},\lambda^{M_{N}}}}(\dP_{0,M_{N}},\dP_{1,M_{N}})\\
                &=\liminf_{M\to\infty}{\mathcal{W}^2_{\kappa^M,\lambda^M}}(\dP_{0,M},\dP_{1,M}),
            \end{align*}
             which finishes the proof.
       \end{proof}
          \section{Convergence of the Approximation}\label{sec:truncation}
        Let $\lambda=\Leb_{|\Lambda_n}$ and assume that $\kappa:\Gamma_{\Lambda_n}\times \Lambda_n\to[0,\infty)$ satisfies \ref{property:ALB} with  approximating sequence $(\kappa_a)_{a\in\N}$. In this section, we prove Theorem~\ref{thm:convergenceresultstruncation} and show that   the dual spatial birth-and-death semigroup, the relative entropy, and the transport distance converge along this approximation. 
    \subsection{Papangelou Measures}
    We first recall relevant notation and prove the convergence of the Papangelou measures along the approximation.
     \begin{lem}\label{lem:convergencetruncation}
          Let $\lambda=\Leb_{|\Lambda_n}$ and assume that $\kappa:\Gamma_{\Lambda_n}\times\Lambda_n\to[0,\infty)$ satisfies \ref{property:cocycle}, \ref{property:locallystable}, and~\ref{property:continuity2}. Moreover, assume that $\kappa$ satisfies \ref{property:ALB} with  approximating sequence $(\kappa_a)_{a\in\N}$. Let $\dR=\RRR_{\kappa,\lambda}(\Gamma_{\Lambda_n})$ and $\dR_a=\RRR_{\kappa_a ,\lambda}(\Gamma_{\Lambda_n})$, $a \in\N$. It holds that $\lim_{a \to\infty}\dR_a =\dR\ \text{ in } \PPP_1(\Gamma_{\Lambda_n}).$
            \end{lem}
       \begin{proof}
            By Proposition~\ref{prop:convergencep1} it suffices to prove that $\lim_{a \to\infty}L_{\dR_a }(f)=L_{\dR}(f)$ for all $f\in\CCC^+_c(\Lambda_n)$ and $\lim_{a \to\infty}\I_{\dR_a }=\I_{\dR}$ vaguely. Recall that by Lemma~\ref{lem:finitevolumePapangelou},
            \begin{align*}
              L_{\dR_a }(f)&=\int e^{-\int f(y)\xi(\mathrm{d}y)}\dR_a (\mathrm{d}\xi)= \frac{1}{Z_{\kappa_a ,\lambda}^{\Lambda_n}}\sum_{j\geq0}\frac{1}{j!}\int\!\ldots\!\int e^{-\int f(y)\sum_{i=1}^j\delta_{x_i}(\mathrm{d}y)}\varphi_{\kappa_a }\left(\sum_{i=1}^j\delta_{x_i}\!\right)\lambda(\mathrm{d}x_1)\ldots\lambda(\mathrm{d}x_j),
            \end{align*}
            with $\varphi_{\kappa_a }(\sum_{i=1}^j\delta_{x_i})=\prod_{i=1}^j\kappa_a (\sum_{l=1}^{i-1}\delta_{x_l},x_i)$.
            We first prove that the partition function $Z_{\kappa_a ,\lambda}^{\Lambda_n}$ converges. 
       By definition, we obtain that $\lim_{a \to\infty}\varphi_{\kappa_{a }}(\sum_{i=1}^j\delta_{x_i})=\varphi_{\kappa}(\sum_{i=1}^j\delta_{x_i})$ $\lambda^{j}$-a.s. Further, since $ \kappa\leq \kappa_{a }\leq A$ by \ref{property:ALB}, we have 
       \begin{align*}
          \frac{1}{j!}\int\ldots\int \varphi_{\kappa_a }\left(\sum_{i=1}^j\delta_{x_i}\right)\lambda(\mathrm{d}x_1)\ldots\lambda(\mathrm{d}x_j)\leq \frac{1}{j!}(n^d)^j A^j,
       \end{align*}
        which is summable in $j$. Combining the previous two observations, we obtain by the dominated convergence theorem that 
       \begin{align*}
           \lim_{a \to\infty} Z_{\kappa_a ,\lambda}^{\Lambda_n}&=\lim_{a \to\infty} \sum_{j\geq0}\frac{1}{j!}\int\ldots\int \varphi_{\kappa_a }\left(\sum_{i=1}^j\delta_{x_i}\right)\lambda(\mathrm{d}x_1)\ldots\lambda(\mathrm{d}x_j)
           =Z_{\kappa,\lambda}^{\Lambda_n}.
       \end{align*}
       Analogously, we obtain by the dominated convergence theorem that 
       \begin{align*}
           \lim_{a \to\infty}\sum_{j\geq0}\frac{1}{j!}\int\ldots\int e^{-\int f(y)\sum_{i=1}^j\delta_{x_i}(\mathrm{d}y)}\varphi_{\kappa_a }\left(\sum_{i=1}^j\delta_{x_i}\right)\lambda(\mathrm{d}x_1)\ldots\lambda(\mathrm{d}x_j)\\
           \quad=\sum_{j\geq0}\frac{1}{j!}\int\ldots\int e^{-\int f(y)\sum_{i=1}^j\delta_{x_i}(\mathrm{d}y)}\varphi_{\kappa}\left(\sum_{i=1}^j\delta_{x_i}\right)\lambda(\mathrm{d}x_1)\ldots\lambda(\mathrm{d}x_j).
       \end{align*}
       All in all, we conclude that 
       \begin{align*}
       \lim_{a \to\infty}L_{\dR_a }(f)&=\lim_{a \to\infty}\frac{1}{Z_{\kappa_a ,\lambda}^{\Lambda_n}}\sum_{j\geq0}\frac{1}{j!}\int\ldots\int e^{-\int f(y)\sum_{i=1}^j\delta_{x_i}(\mathrm{d}y)}\varphi_{\kappa_a }\left(\sum_{i=1}^j\delta_{x_i}\right)\lambda(\mathrm{d}x_1)\ldots\lambda(\mathrm{d}x_j)\\
       &=\frac{1}{Z_{\kappa,\lambda}^{\Lambda_n}}\sum_{j\geq0}\frac{1}{j!}\int\ldots\int e^{-\int f(y)\sum_{i=1}^j\delta_{x_i}(\mathrm{d}y)}\varphi_{\kappa}\left(\sum_{i=1}^j\delta_{x_i}\right)\lambda(\mathrm{d}x_1)\ldots\lambda(\mathrm{d}x_j)\\
       &=L_{\dR}(f).
       \end{align*}
    An analogous argument for the intensity measure finishes the proof.
       \end{proof}
\subsection{Dual Semigroup}
           In this subsection, we prove Theorem~\ref{thm:convergenceresultstruncation} $(i)$. More precisely, we show that if $\lambda=\Leb_{|\Lambda_n}$ and $\kappa:\Gamma_{\Lambda_n}\times\Lambda_n\to[0,\infty)$ satisfies \ref{property:cocycle}, \ref{property:locallystable}, and~\ref{property:ALB} with approximating sequence $(\kappa_a)_{a\in\N}$, then \begin{align*}
          \lim_{a\to\infty}\mathrm{S}_t^{\kappa_a,\lambda*}\dP=\mathrm{S}_t^{\kappa,\lambda*}\dP \quad \text{ in }\PPP_1(\Gamma_{\Lambda_n}),
      \end{align*}
      for all $t\geq0$ and $\dP\in\PPP_1(\Gamma_{\Lambda_n})$ s.t.\   $\dP\ll\dR$. We begin by establishing the corresponding result for the semigroup $(\mathrm{S}_t^{\kappa_a,\lambda})_{t\geq0}$. 
    \begin{satz}\label{thm:convergencetruncatedsemigroup}
     Let $\lambda=\Leb_{|\Lambda_n}$  and assume that $\kappa:\Gamma_{\Lambda_n}\times \Lambda_n \to [0,\infty)$  satisfies \ref{property:cocycle}, \ref{property:locallystable}, and~\ref{property:ALB} with  approximating sequence $(\kappa_a)_{a\in\N}$. Further, let $\dR=\RRR_{\kappa,\lambda}(\Gamma_{\Lambda_n})$ and $F\in L^2(\dR)$. It holds that 
        \begin{align*}
            \lim_{a \to\infty}\|\mathrm{S}_t^{\kappa_a ,\lambda}F-\mathrm{S}_t^{\kappa,\lambda}F\|_{L^2(\dR)}=0,\quad t\geq0.
        \end{align*}
    \end{satz}
    \begin{proof}
    The proof follows the same line of argumentation as the proof of Theorem~\ref{thm:discretizedsemigroup}. We recall that $(\mathrm{S}_t^{\kappa,\lambda})_{t\geq0}$ and $(\mathrm{S}_t^{\kappa_a,\lambda})_{t\geq0}$ are contractive $C_0$-semigroups and  that  $\mathcal{G}$ is a core of $\mathfrak{L}_{\kappa,\lambda}$ and $\mathfrak{L}_{\kappa_a,\lambda}$, see \eqref{eq:core}. Hence, using the Kato--Trotter approximation, see e.g.\ \cite[Theorem 6.1]{EthierKurtz}, the claim follows if 
    \begin{align*}
        \lim_{a \to\infty}\|\mathfrak{L}_{\kappa_a ,\lambda}F-\mathfrak{L}_{\kappa,\lambda} F\|^2_{L^2(\dR)}
          =0
    \end{align*}
    for all $F\in\mathcal{G}.$ Using Jensen's inequality,  we obtain for $F\in\mathcal{G}$ that \begin{align*}
           \lim_{a \to\infty}\|\mathfrak{L}_{\kappa_a ,\lambda}F-\mathfrak{L}_{\kappa,\lambda} F\|^2_{L^2(\dR)}
           &=\lim_{a \to\infty}\int | \int\kappa_a (\xi,x)(F(\xi+\delta_{x})-F(\xi))-\kappa(\xi,x)(F(\xi+\delta_{x})-F(\xi))\lambda(\mathrm{d}x)\\
           &\quad +\sum_{x\in\xi}F(\xi-\delta_{x})-F(\xi)-\sum_{x\in\xi}F(\xi-\delta_{x})-F(\xi)|^2\dR(\mathrm{d}\xi)\\
           &=\lim_{a \to\infty}\int | \int(\kappa_a (\xi,x)-\kappa(\xi,x))(F(\xi+\delta_{x})-F(\xi))\lambda(\mathrm{d}x)|^2\dR(\mathrm{d}\xi)\\
           &\leq \lim_{a \to\infty}\lambda(\Lambda_n)\int  \int(\kappa_a (\xi,x)-\kappa(\xi,x))^2(F(\xi+\delta_{x})-F(\xi))^2\lambda(\mathrm{d}x)\dR(\mathrm{d}\xi).
       \end{align*}
       By \ref{property:ALB}, we have that $0\leq \kappa_a (\xi,x)-\kappa(\xi,x)\leq g(a)$ with $\lim_{a\to\infty}g(a)=0$. Moreover, since $F \in \mathcal{G}$, there exists a constant $c \geq 1$ s.t.\    
$|F(\xi)| \le c^{\xi(\Lambda_n)}$. Hence, we obtain  
\begin{align*}
(F(\xi+\delta_{x})-F(\xi))^2
\leq 2c^{2(\xi(\Lambda_n) + 1)} + 2c^{2\xi(\Lambda_n)}.
\end{align*}
The upper bound is $  \dR$-integrable since  $\kappa$ satisfies \ref{property:locallystable} and hence 
$\int c^{2\xi(\Lambda_n)} \dR(\mathrm{d}\xi) < \infty$, see  \eqref{property:integrability}.  
Therefore, using the dominated convergence theorem, we conclude that
\begin{align*}
\lim_{a \to\infty}\|\mathfrak{L}_{\kappa_a ,\lambda}F-\mathfrak{L}_{\kappa,\lambda} F\|^2_{L^2(\dR)}&\leq \lambda(\Lambda_n)\lim_{a \to\infty} \int\!  \int(\kappa_a (\xi,x)-\kappa(\xi,x))^2(F(\xi+\delta_{x})-F(\xi))^2\lambda(\mathrm{d}x)\dR(\mathrm{d}\xi) =0,
       \end{align*}
finishing the proof.
    \end{proof}
    \begin{proof}[Proof of Theorem~\ref{thm:convergenceresultstruncation} $(i)$]
      The proof follows the same line of argument as the proof of Theorem~\ref{thm:convergenceresultsdiscretization} $(i)$, with Theorem~\ref{thm:discretizedsemigroup} replaced by Theorem~\ref{thm:convergencetruncatedsemigroup}. 
    \end{proof}

  \subsection{Entropy}
          In this subsection, we prove Theorem~\ref{thm:convergenceresultstruncation} $(ii)$. More precisely, let $\lambda=\Leb_{|\Lambda_n}$ and assume that $\kappa:\Gamma_{\Lambda_n}\times\Lambda_n\to[0,\infty)$ satisfies \ref{property:cocycle}, \ref{property:locallystable}, \ref{property:continuity2}, and~\ref{property:ALB} with approximating sequence $(\kappa_a)_{a\in\N}$. We show that if $\dR=\RRR_{\kappa,\lambda}(\Gamma_{\Lambda_n})$ and $\dR_a=\RRR_{\kappa_a,\lambda}(\Gamma_{\Lambda_n}),\;a\in\N$, then $$\Gamma\text{-}\lim_{a\to\infty}\Ent(\cdot|\dR_a)=\Ent(\cdot|\dR)$$ on $\PPP_1(\Gamma_{\Lambda_n})$. A recovery sequence for $\dP\in\PPP_1(\Gamma_{\Lambda_n})$ is given by the constant sequence $(\dP)_{a\in\N}$.
\begin{proof}[Proof of Theorem~\ref{thm:convergenceresultstruncation} $(ii)$]
    Let $\lambda=\Leb_{|\Lambda_n}$ and assume that $\kappa:\Gamma_{\Lambda_n}\times\Lambda_n\to[0,\infty)$ satisfies \ref{property:cocycle}, \ref{property:locallystable}, \ref{property:continuity2}, and~\ref{property:ALB}. 
    Let $(\dP_a)_{a\in\N}\subset\PPP_1(\Gamma_{\Lambda_n})$ and  $\dP\in\PPP_1(\Gamma_{\Lambda_n})$
    s.t.\  $\lim_{a\to\infty}\dP_a=\dP$ in $\PPP_1(\Gamma_{\Lambda_n})$. For    $\dR_a=\RRR_{\kappa_a,\lambda}(\Gamma_{\Lambda_n}),\; a\in\N$, and $\dR=\RRR_{\kappa,\lambda}(\Gamma_{\Lambda_n})$ it holds that  $\lim_{a\to\infty}\dR_a=\dR$ in $\PPP_1(\Gamma_{\Lambda_n})$, see Lemma~\ref{lem:convergencetruncation}. Hence,  we obtain by the lower semicontinuity of the relative entropy (Lemma~\ref{lem:propertiesent} $(i)$) that $$\Ent(\dP|\dR)\leq\liminf_{a\to\infty}\Ent(\dP_a|\dR_a).$$
   It remains to prove that $ \limsup_{a \to\infty}\Ent(\dP|\dR_a )\leq\Ent(\dP|\dR) .$ If $\Ent(\dP|\dR)=\infty$, the claim is immediate. Hence, let $\Ent(\dP|\dR)<\infty$ and $\dP=\rho\dR$.  Since $\varphi_{\kappa}\leq\varphi_{\kappa_a}$, Lemma~\ref{lem:finitevolumePapangelou} yields, for all $K\in\mathcal{B}(\Gamma_{\Lambda_n})$, that
            \begin{align*}
             \int \ind_K(\xi)\dP(\mathrm{d}\xi)&=\int \ind_K(\xi)\frac{Z_{\kappa_a ,\lambda}^{\Lambda_n}}{Z_{\kappa,\lambda}^{\Lambda_n}}\frac{\varphi_{\kappa}(\xi)}{\varphi_{\kappa_a }(\xi)}\rho(\xi) \dR_a (\mathrm{d}\xi).
            \end{align*}
            Hence, $\dP=\rho_a \dR_a $ with 
            \begin{align*}
                \rho_a (\xi)=\frac{Z_{\kappa_a ,\lambda}^{\Lambda_n}}{Z_{\kappa,\lambda}^{\Lambda_n}}\frac{\varphi_{\kappa}(\xi)}{\varphi_{\kappa_a }(\xi)}\rho(\xi).
            \end{align*}
            Inserting this into the definition of the relative entropy yields 
            \begin{align*}
                \limsup_{a \to\infty}\Ent(\dP|\dR_a )&=\limsup_{a \to\infty}\int\log\left(\frac{Z_{\kappa_a ,\lambda}^{\Lambda_n}}{Z_{\kappa,\lambda}^{\Lambda_n}}\frac{\varphi_{\kappa}(\xi)}{\varphi_{\kappa_a }(\xi)}\rho(\xi)\right)\dP(\mathrm{d}\xi)\\
                &\leq \limsup_{a \to\infty}\log\left(\frac{Z_{\kappa_a ,\lambda}^{\Lambda_n}}{Z_{\kappa,\lambda}^{\Lambda_n}}\right)+ \limsup_{a \to\infty}  \int\log\left(\frac{\varphi_{\kappa}(\xi)}{\varphi_{\kappa_a }(\xi)}\rho(\xi)\right)\dP(\mathrm{d}\xi).
            \end{align*}
             Analogously to the  proof of Lemma~\ref{lem:convergencetruncation}, we get by  dominated convergence  $\lim_{a \to\infty} Z_{\kappa_a ,\lambda}^{\Lambda_n}=Z_{\kappa,\lambda}^{\Lambda_n}$. Since additionally $\varphi_{\kappa}\leq\varphi_{\kappa_a }$, we conclude
            \begin{align*}
                \limsup_{a \to\infty}\Ent(\dP|\dR_a )&\leq \limsup_{a \to\infty}\log\left(\frac{Z_{\kappa_a ,\lambda}^{\Lambda_n}}{Z_{\kappa,\lambda}^{\Lambda_n}}\right)+   \int\log\left(\frac{\varphi_{\kappa}(\xi)}{\varphi_{\kappa}(\xi)}\rho(\xi)\right)\dP(\mathrm{d}\xi)\\
                &= 0+\int \log(\rho(\xi))\rho(\xi)\dR(\mathrm{d}\xi)\\
                &=\Ent(\dP|\dR),
            \end{align*}
            finishing the proof.
        \end{proof}

       \subsection{Transport Distance}

       In this subsection, we prove Theorem~\ref{thm:convergenceresultstruncation}  $(iii)$. Let $\lambda=\Leb_{|\Lambda_n}$ and assume that $\kappa:\Gamma_{\Lambda_n}\times\Lambda_n\to[0,\infty)$ satisfies \ref{property:cocycle}, \ref{property:locallystable}, \ref{property:continuity}, and~\ref{property:ALB} with approximating sequence $(\kappa_a)_{a\in\N}$. We show that $$\Gamma\text{-}\lim_{a\to\infty}\mathcal{W}_{\kappa_a,\lambda}(\cdot,\cdot)
    =\mathcal{W}_{\kappa,\lambda}(\cdot,\cdot)$$ on $\PPP_1(\Gamma_{\Lambda_n})\times \PPP_1(\Gamma_{\Lambda_n})$. A recovery sequence for $(\dP_0,\dP_1)\in\PPP_1(\Gamma_{\Lambda_n})\times \PPP_1(\Gamma_{\Lambda_n})$ is given by the constant sequence $(\dP_0,\dP_1)_{a\in\N}$.\\

     The proof is divided into two parts. We establish the limsup inequality in Lemma~\ref{lem:upperboundtruncation} and the liminf inequality in Lemma~\ref{lem:lowerboundtruncation}. Combining both lemmas yields the proof of Theorem~\ref{thm:convergenceresultstruncation} $(iii)$.
       \begin{lem}\label{lem:upperboundtruncation} Let $\lambda=\Leb_{|\Lambda_n}$ and assume that $\kappa:\Gamma_{\Lambda_n}\times\Lambda_n\to[0,\infty)$ satisfies \ref{property:cocycle}, \ref{property:locallystable}, and~\ref{property:ALB} with approximating sequence $(\kappa_a)_{a\in\N}$. Further, let $\dP_0,\dP_1\in\PPP_1(\Gamma_{\Lambda_n})$. It holds that 
           \begin{align*}
              \limsup_{a \to\infty} \mathcal{W}_{\kappa_a ,\lambda}(\dP_0,\dP_1)\leq \mathcal{W}_{\kappa,\lambda}(\dP_0,\dP_1).
           \end{align*}
       \end{lem}
       \begin{proof}
       If   $\mathcal{W}_{\kappa,\lambda}(\dP_0,\dP_1)=\infty$, the claim is immediate. Therefore, let $\mathcal{W}_{\kappa,\lambda}(\dP_0,\dP_1)<\infty.$ By Theorem~\ref{thm:minimizer} exists a solution $(\bar{\dP},\bar{\dV})\in\CCC\EEE^{\Lambda_n}(\dP_0,\dP_1)$ s.t.\   
           \begin{align*}
               \mathcal{W}_{\kappa,\lambda}^2(\dP_0,\dP_1)=\int_0^1\mathcal{L}_{\kappa,\lambda}(\bar{\dP}_t,\bar{\dV}_t)\mathrm{d}t.
           \end{align*}
           Using that $\kappa \leq\kappa_a$ together with the monotonicity of $\alpha$, we obtain that 
           \begin{align*}
               \limsup_{a \to\infty}\mathcal{W}_{\kappa_a ,\lambda}^2(\dP_0,\dP_1)&\leq \limsup_{a \to\infty} \int_0^1\mathcal{L}_{\kappa_a ,\lambda}(\bar{\dP}_t,\bar{\dV}_t)\mathrm{d}t\leq \int_0^1\mathcal{L}_{\kappa,\lambda}(\bar{\dP}_t,\bar{\dV}_t)\mathrm{d}t=\mathcal{W}_{\kappa,\lambda}^2(\dP_0,\dP_1),
           \end{align*}
           proving the claim.
       \end{proof}
       We continue with the proof of the liminf inequality for the $\Gamma$-convergence of $(\mathcal{W}_{\kappa_a,\lambda})_{a\in\N}$. We first prove the corresponding result for the Lagrange functional in Lemma~\ref{lem:lscLagrangetruncatedkappa}. Thereafter, we establish a compactness result for the action functional (Lemma~\ref{lem:compactnesstruncation}) which allows us to lift the liminf inequality to the transport distance.
         \begin{lem}\label{lem:lscLagrangetruncatedkappa}
          Let $\lambda=\Leb_{|\Lambda_n}$ and assume that $\kappa:\Gamma_{\Lambda_n}\times\Lambda_n\to[0,\infty)$ satisfies \ref{property:cocycle}, \ref{property:locallystable}, \ref{property:continuity}, and~\ref{property:ALB} with approximating sequence $(\kappa_a)_{a\in\N}$.  If $(\dP_a)_{a\in\N}\subset\PPP_1(\Gamma_{\Lambda_n})$ and $\dP\in\PPP_1(\Gamma_{\Lambda_n})$ such that   $\lim_{a\to\infty}\dP_a=\dP$ in $\PPP_1(\Gamma_{\Lambda_n})$ and $(\dV_a)_{a\in\N}\subset\MMM_{b,0}(\Gamma_{\Lambda_n}\times {\Lambda_n})$ and $\dV\in \MMM_{b,0}(\Gamma_{\Lambda_n}\times {\Lambda_n})$ such that   $\lim_{a\to\infty}\dV_a=\dV$ in $\MMM_{b,0}(\Gamma_{\Lambda_n}\times \Lambda_n)$, then     
             \begin{align*}
                 \mathcal{L}_{\kappa,\lambda}(\dP,\dV)\leq\liminf_{a \to\infty}\mathcal{L}_{\kappa_a ,\lambda}(\dP_a ,\dV_a ). 
             \end{align*}
         \end{lem}
         \begin{proof}
          Assume that $\kappa$ is not only upper semicontinuous \ref{property:continuity} but continuous. Let $F\in\CCC_{b,0}(\Gamma_{\Lambda_n}\times \Lambda_n)$ s.t.\   $|F|\leq C$. As $\kappa$ is continuous and bounded above  \ref{property:locallystable}, we note that $\kappa F\in\CCC_{b,0}(\Gamma_{\Lambda_n}\times\Lambda_n)$. 
        Hence, for $\varepsilon>0$, there exists some $a _1\in\N$ s.t.\   \begin{align*}
            |\int\int F(\xi,x)\kappa(\xi,x)\dP_a (\mathrm{d}\xi)\lambda(\mathrm{d}x)-\int \int F(\xi,x)\kappa(\xi,x)\dP(\mathrm{d}\xi)\lambda(\mathrm{d}x)|<\frac{\varepsilon}{2}
        \end{align*}
        for all $a \geq a _1$.
        Further, by \ref{property:ALB} we have that \begin{align*}
            &|\int\int F(\xi,x)\kappa_a (\xi,x)\dP_a (\mathrm{d}\xi)\lambda(\mathrm{d}x)-\int \int F(\xi,x)\kappa(\xi,x)\dP_a (\mathrm{d}\xi)\lambda(\mathrm{d}x)|\\
            &\quad\leq\int \int |F(\xi,x)||\kappa_a (\xi,x)-\kappa(\xi,x)|\dP_a (\mathrm{d}\xi)\lambda(\mathrm{d}x)\\
            &\quad\leq C\lambda(\Lambda_n)g(a)\\&\quad<\frac{\varepsilon}{2},
        \end{align*}
        for all $a \geq a _2$. Using the triangle inequality, we conclude that 
        \begin{align*}
            &|\int\int F(\xi,x)\kappa_a (\xi,x)\dP_a (\mathrm{d}\xi)\lambda(\mathrm{d}x)-\int\int F(\xi,x)\kappa(\xi,x)\dP(\mathrm{d}\xi)\lambda(\mathrm{d}x)|\\
            &\quad\leq |\int\int F(\xi,x)\kappa_a (\xi,x)\dP_a (\mathrm{d}\xi)\lambda(\mathrm{d}x)-\int\int F(\xi,x)\kappa(\xi,x)\dP_a (\mathrm{d}\xi)\lambda(\mathrm{d}x)| \\
            &\qquad+|\int\int F(\xi,x)\kappa(\xi,x)\dP_a (\mathrm{d}\xi)\lambda(\mathrm{d}x)-\int\int F(\xi,x)\kappa(\xi,x)\dP(\mathrm{d}\xi)\lambda(\mathrm{d}x)|\\
            &\quad \leq \frac{\varepsilon}{2}+\frac{\varepsilon}{2},
        \end{align*}
        for all $a \geq a _1\vee a _2$.
         Thus, $\lim_{a \to\infty}\kappa_a (\dP_a \otimes\lambda)=\kappa(\dP\otimes\lambda)$ in $\MMM_{b,0}(\Gamma_{\Lambda_n}\times\Lambda_n)$. Since further $\lim_{a \to\infty}C_{\dP_a }= C_\dP$ in $\MMM_{b,0}(\Gamma_{\Lambda_n}\times\Lambda_n)$  by  Proposition~\ref{prop:convergencep1},  Theorem~\ref{athm:lowersemicontinuity} and the monotonicity of $\alpha$ imply 
        \begin{align*}
            \mathcal{L}_{\kappa,\lambda}(\dP,\dV)&=\int\alpha\left(\frac{\mathrm{d}\kappa(\dP\otimes\lambda)}{\mathrm{d}\sigma}(\xi,x),\frac{\mathrm{d}C_{\dP}}{\mathrm{d}\sigma}(\xi,x),\frac{\mathrm{d}\dV}{\mathrm{d}\sigma}(\xi,x)\right)\sigma(\mathrm{d}\xi,\mathrm{d}x)\\
            &\leq\liminf_{a \to\infty}\int\alpha\left(\frac{\mathrm{d}\kappa_a (\dP_a \otimes\lambda)}{\mathrm{d}\sigma}(\xi,x),\frac{\mathrm{d}C_{\dP_a }}{\mathrm{d}\sigma}(\xi,x),\frac{\mathrm{d}\dV_a }{\mathrm{d}\sigma}(\xi,x)\right)\sigma(\mathrm{d}\xi,\mathrm{d}x)\\
            &=\liminf_{a \to\infty}\mathcal{L}_{\kappa_a ,\lambda}(\dP_a ,\dV_a ),
        \end{align*}
        finishing the proof for $\kappa$ being continuous. Now, assume that $\kappa$ is upper semicontinuous \ref{property:continuity}. There exists a sequence of functions $(\kappa_i)_{i\in\N}\subset\CCC_b(\Gamma_{\Lambda_n}\times \Lambda_n)$ s.t.\  $\kappa_i\downarrow\kappa$ for $i\to\infty.$ Since $0\leq \kappa_a-\kappa\leq g(a)$ with $g(a)\to0$, we have
\[
    \kappa_a\leq \kappa_i+g(a)
    \quad \text{for all } a,i\in\N .
\]
For fixed $i\in\N$, the functions $\kappa_i+g(a)\in\CCC_b(\Gamma_{\Lambda_n}\times \Lambda_n)$ converge uniformly to $\kappa_i$ as $a\to\infty$. Hence,  by the previous result for continuous functions and analogously to the proof of Lemma \ref{lem:lscwithvariablekappa}, we obtain that
\begin{align*}
           \mathcal{L}_{\kappa,\lambda}(\dP,\dV)
           &\leq \sup_{i\in\N}\liminf_{a \to\infty}\mathcal{L}_{\kappa_i+g(a),\lambda}(\dP_a,\dV_a)\leq \liminf_{a \to\infty}\sup_{i\in\N}\mathcal{L}_{\kappa_i+g(a),\lambda}(\dP_a,\dV_a)\leq \liminf_{a \to\infty} \mathcal{L}_{\kappa_a,\lambda}(\dP_a,\dV_a),
       \end{align*}
       finishing the proof.
        \end{proof}

         \begin{lem}\label{lem:compactnesstruncation}
        Let $\lambda=\Leb_{|\Lambda_n}$ and assume that $\kappa:\Gamma_{\Lambda_n}\times\Lambda_n\to[0,\infty)$ satisfies \ref{property:cocycle}, \ref{property:locallystable}, \ref{property:continuity}, and~\ref{property:ALB} with approximating sequence $(\kappa_a)_{a\in\N}$. Let $(\dP_{0,a })_{a \in\N},(\dP_{1,a })_{a \in\N}\subset\PPP_1(\Gamma_{\Lambda_n})$ and $\dP_0,\dP_1\in\PPP_1(\Gamma_{\Lambda_n})$ s.t.\   $\lim_{a \to\infty}\dP_{0,a }=\dP_0$ and 
        $\lim_{a \to\infty}\dP_{1,a }=\dP_1$ in $\PPP_1(\Gamma_{\Lambda_n})$. 
        If $(\bar{\dP}_a ,\bar{\dV}_a )_{a \in\N}\subset\CCC\EEE^{\Lambda_n}(\dP_{0,a },\dP_{1,a })$  s.t.\   
        \begin{align*}
            \sup_{a \in\N}\mathcal{A}_{\kappa_a ,\lambda}(\bar{\dP}_a ,\bar{\dV}_a )\leq C<\infty,
        \end{align*}
        then there exists a subsequence $(\bar{\dP}_{a_b} ,\bar{\dV}_{a_b} )_{b \in\N}$ and a solution $(\bar{\dP},\bar{\dV})\in\CCC\EEE^{\Lambda_n}(\dP_0,\dP_1)$ s.t.\   
        \begin{align*}
            \lim_{b\to\infty}\bar{\dP}_{t,a _b}=\bar{\dP}_t\quad &\text{ in }\PPP_1(\Gamma_{\Lambda_n}) 
            \text{ for all } t\in [0,1],\\
            \lim_{b\to\infty}\bar{\dV}_{a _b}=\bar{\dV}\quad &\text{ in }\MMM_{b,0}(\Gamma_{\Lambda_n}\times \Lambda_n\times[0,1]).
             \end{align*}
    \end{lem}
        \begin{proof}
            The proof is identical to the proof of Lemma~\ref{lem:compactnessofsolutions} $(iii)$.
        \end{proof} \begin{lem}\label{lem:lowerboundtruncation}
       Let $\lambda=\Leb_{|\Lambda_n}$ and assume that $\kappa:\Gamma_{\Lambda_n}\times\Lambda_n\to[0,\infty)$ satisfies \ref{property:cocycle}, \ref{property:locallystable}, \ref{property:continuity}, and~\ref{property:ALB} with approximating sequence $(\kappa_a)_{a\in\N}$.  Let $(\dP_{0,a })_{a \in\N},(\dP_{1,a })_{a \in\N}\subset\PPP_1(\Gamma_{\Lambda_n})$ and $\dP_0,\dP_1\in\PPP_1(\Gamma_{\Lambda_n})$ s.t.\   $\lim_{a \to\infty}\dP_{0,a }=\dP_0$ and 
        $\lim_{a \to\infty}\dP_{1,a }=\dP_1$ in $\PPP_1(\Gamma_{\Lambda_n})$. It holds that 
            \begin{align*}
                \mathcal{W}_{\kappa,\lambda}(\dP_0,\dP_1)\leq\liminf_{a \to\infty}\mathcal{W}_{\kappa_a ,\lambda}(\dP_{0,a},\dP_{1,a}).
            \end{align*}
        \end{lem}
        \begin{proof}
            The proof is analogous to the proof of Lemma~\ref{lem:lowersemicontinuitydiscretization} by replacing Lemma~\ref{lem:lscwithvariablekappa} by Lemma~\ref{lem:lscLagrangetruncatedkappa} and Lemma~\ref{lem:compactnessdiscretization} by Lemma~\ref{lem:compactnesstruncation}.
        \end{proof}
         \section{Transport Geometry for Stationary Measures}
             \subsection{Gradient Flows, Convexity, and Functional Inequalities}
              In this section, we extend the results of Section \ref{sec:gradientflows} to  {stationary} measures. More precisely, we identify the spatial birth-and-death process w.r.t.\ a stationary Papangelou measure $\dR\in\RRR_{\kappa,s}(\Gamma)$ (see Subsection \ref{subsec:semigroup}) as the EVI gradient flow of the specific relative entropy $\mathcal{E}(\cdot|\dR)$ in a geometry induced by a suitable transport distance.  

 Let $\kappa:\Gamma\times\R^d\to[0,\infty)$ satisfy \ref{property:locallystable}, \ref{property:continuity}, and \ref{property:shiftinvariance}. Following the approach of \cite{huesmann2025nonlocalwassersteingeometrygradient}, we define a \textit{specific} non-local transport distance  $\mathcal{W}_{\kappa,s}:\PPP_s(\Gamma)\times\PPP_s(\Gamma)\to[0,\infty]$ by 
\begin{align*}
\mathcal{W}_{\kappa,s}^2(\dP_0, \dP_1) &:= \inf \left\{ \mathcal{A}_{\kappa,s}(\Bar{\dP}, \Bar{\dV}) : (\Bar{\dP}, \Bar{\dV}) \in \CCC\EEE^s(\dP_0, \dP_1) \right\},
\end{align*}
where $\CCC\EEE^s(\dP_0,\dP_1)$ is the set of stationary solutions to the continuity equation connecting $\dP_0$ and $\dP_1$ (see Definition \ref{dfn:stationarysolutions}) and $\mathcal{A}_{\kappa,s}$ is the normalized large-volume limit of the action functional introduced in Section \ref{sec:stationarydistance} defined by
    \begin{align*}
         \mathcal{A}_{\kappa,s}({\Bar{\dP}},{\Bar{\dV}}):=\sup_{n\in\N}\frac{1}{\Leb(\Lambda_n)}\int_0^1\mathcal{L}_{\kappa_n,\Leb_{|\Lambda_n}}(\Bar{\dP}_{t|\Lambda_n},\Bar{\dV}_{t|\Lambda_n})\mathrm{d}t,
     \end{align*}
     with 
     \begin{align*}
         \kappa_n(\xi,x):=\kappa(\xi_{|\Lambda_n},x)\ind_{\Lambda_n}(x),\quad \xi\in\Gamma,x\in\R^d
     \end{align*}
     and the mean function $\theta$ in the definition of the Lagrange functional is given by the logarithmic mean \ref{eq:logmean}.
 We refer to Section \ref{sec:stationarydistance} for a more detailed discussion.
\begin{bem}
    By definition, $\mathcal{L}_{\kappa_n,\Leb_{|\Lambda_n}}({\dP},{\dV})=\mathcal{L}_{\kappa,\Leb_{|\Lambda_n}}({\dP},{\dV})$. For notational convenience, we briefly write $\mathcal{L}_{\kappa_n}$ and $\mathcal{W}_{\kappa_n}$ instead of $\mathcal{L}_{\kappa_n,\Leb_{|\Lambda_n}}$ and $\mathcal{W}_{\kappa_n,\Leb_{|\Lambda_n}},$ respectively. Analogously, we write $(\mathrm{S}_t^{\kappa_n*})_{t\geq0}$ instead of $(\mathrm{S}_t^{\kappa_n,\Leb_{|\Lambda_n} *})_{t\geq0}$.
\end{bem}
The main result of this section is the following EVI.

 \begin{satz}\label{thm:evistationary}
Assume that $\kappa:\Gamma\times\R^d\to[0,\infty)$ satisfies \ref{property:cocycle}, \hyperref[property:locallystable]{$(\mathbf{S}(1))$}, \ref{property:continuity2}, \ref{property:repulsive}, \ref{property:finiterange}, and~\ref{property:shiftinvariance}. Moreover, assume that $\kappa_n$ satisfies \ref{property:ALB} and \hyperref[property:dobrushin2]{$(\mathbf{D}^\prime(\zeta(\kappa_n)))$} for all $n\in\N$. Let $\dR\in\RRR_{\kappa,s}(\Gamma)$.  Then $(\mathrm{S}_t^{\kappa*})_{t\geq0}$ satisfies the EVI$(\zeta(\kappa))$ for $\mathcal{E}(\cdot|\dR)$, i.e.\ for all $\dP_0,\dP_1\in\PPP_s(\Gamma)$ s.t.\  $ \mathcal{E}(\dP_1|\dR)<\infty$ and $\mathcal{W}_{\kappa,s}(\dP_0,\dP_1)<\infty$, it holds
\begin{align*}
            \frac{1}{2}\frac{\mathrm{d}}{\mathrm{d}t}
\mathcal{W}_{\kappa,s}^2(\mathrm{S}_t^{\kappa*}\dP_0,\dP_1)
+
\frac{\zeta(\kappa)}{2}
\mathcal{W}_{\kappa,s}^2(\mathrm{S}_t^{\kappa*}\dP_0,\dP_1)
\leq
\mathcal{E}(\dP_1 | \dR)
-
\mathcal{E}(\mathrm{S}_t^{\kappa*}\dP_0 | \dR),
\quad
 t>0,
\end{align*}
or equivalently in integrated form
\begin{align*}
        \frac{e^{\zeta(\kappa)t}}{2}\mathcal{W}_{\kappa,s}^2(\mathrm{S}_t^{\kappa*}\dP_0,\dP_1)-\frac{1}{2}\mathcal{W}_{\kappa,s}^2(\dP_0,\dP_1)\leq E_{\zeta(\kappa)}(t)(\mathcal{E}(\dP_1 | {\dR})
-
\mathcal{E}(\mathrm{S}_t^{\kappa*}\dP_0 | {\dR})),\quad t>0,
\end{align*}
with $E_{\zeta(\kappa)}(t)=\int_0^te^{\zeta(\kappa)s}\mathrm{d}s$.
\end{satz}

Before proving Theorem~\ref{thm:evistationary}, we discuss some of its consequences. The following theorem may be regarded as the counterpart of McCann's displacement convexity in our non-local stationary setting.

\begin{satz}\label{thm:geodconvexitystationary}
Assume that $\kappa:\Gamma\times\R^d\to[0,\infty)$ satisfies \ref{property:cocycle}, \hyperref[property:locallystable]{$(\mathbf{S}(1))$}, \ref{property:continuity2}, \ref{property:repulsive}, \ref{property:finiterange}, and~\ref{property:shiftinvariance}. Moreover, assume that $\kappa_n$ satisfies \ref{property:ALB} and \hyperref[property:dobrushin2]{$(\mathbf{D}^\prime(\zeta(\kappa_n)))$} for all $n\in\N$. Let $\dR\in\RRR_{\kappa,s}(\Gamma)$. The specific relative entropy $\mathcal{E}(\cdot|\dR)$ is $\zeta(\kappa)$-geodesically convex w.r.t.\ $\mathcal{W}_{\kappa,s}$. More precisely, if $\dP_0,\dP_1\in\PPP_s(\Gamma)$ satisfy $\mathcal{W}_{\kappa,s}(\dP_0,\dP_1)<\infty$, then for every $\mathcal{W}_{\kappa,s}$-geodesic $\bar{\dP}$ connecting $\dP_0$ and $\dP_1$, it holds 
\begin{align*}
        \mathcal{E}(\bar{\dP}_t|\dR)\leq (1-t) \mathcal{E}({\dP}_0|\dR)+t\mathcal{E}({\dP}_1|\dR)-\frac{\zeta(\kappa)t(1-t)}{2}\mathcal{W}_{\kappa,s}^2(\dP_0,\dP_1), \quad 0\leq t\leq 1.
\end{align*}
\end{satz}
\begin{proof}
The proof follows directly from \cite[Proposition 3.2]{Daneri_2008} and Theorem~\ref{thm:evistationary}.
\end{proof}

Furthermore, we derive several new functional inequalities for stationary Papangelou measures. These inequalities are the stationary counterparts of the finite-volume functional inequalities of Theorem~\ref{thm:FIbounded}. 
\begin{satz}\label{thm:FIstationary}
Assume that $\kappa:\Gamma\times\R^d\to[0,\infty)$ satisfies \ref{property:cocycle}, \hyperref[property:locallystable]{$(\mathbf{S}(1))$}, \ref{property:continuity2}, \ref{property:repulsive}, \ref{property:finiterange}, and~\ref{property:shiftinvariance}. Moreover, assume that $\kappa_n$ satisfies \ref{property:ALB} and \hyperref[property:dobrushin2]{$(\mathbf{D}^\prime(\zeta(\kappa_n)))$} for all $n\in\N$. For $\dR\in\RRR_{\kappa,s}(\Gamma)$ and $\dP_0,\dP_1\in\PPP_s(\Gamma)$ s.t.\  $\mathcal{E}(\dP_0|\dR),\mathcal{E}(\dP_1|\dR)<\infty$, the following functional inequalities are satisfied:
\begin{enumerate}
    \item[${(i)}$]
    The dual semigroup $(\mathrm{S}_t^{\kappa*})_{t\geq0}$ is $\zeta(\kappa)$-exponentially contractive with respect to 
    $\mathcal{W}_{\kappa,s}$: 
    \[
        \mathcal{W}_{\kappa,s}
        \bigl(\mathrm{S}_t^{\kappa*}\dP_0,\mathrm{S}_t^{\kappa*}\dP_1\bigr)
        \leq
        e^{-\zeta(\kappa)t}
        \mathcal{W}_{\kappa,s}(\dP_0,\dP_1),
        \quad t\geq0.
    \]
     Further, $(\mathrm{S}_t^{\kappa*})_{t\geq0}$ is $2\zeta(\kappa)$-exponentially contractive  w.r.t.\ $\mathcal{E}(\cdot|\dR)$:
\begin{align*}
    \mathcal{E}
    (
        \mathrm{S}_t^{\kappa*}\dP_0|\dR
    )
    &\leq
    e^{-2\zeta(\kappa)t}
    \mathcal{E}(\dP_0|\dR),\quad t\geq0.
    \end{align*}

    \item[${(ii)}$]
    The specific HWI inequality with parameter $\zeta(\kappa)$ is satisfied:
    \[
        \mathcal{E}(\dP_0|\dR)
        \leq
        \mathcal{W}_{\kappa,s}(\dP_0,\dR)
        \sqrt{\mathcal{I}(\dP_0|\dR)}
        -
        \frac{\zeta(\kappa)}{2}
        \mathcal{W}_{\kappa,s}^2(\dP_0,\dR).
    \]

    \item[${(iii)}$]
    The specific modified Log-Sobolev inequality with parameter 
    $\zeta(\kappa)$ is satisfied: 
    \[
        \mathcal{E}(\dP_0|\dR)
        \leq
        \frac{1}{2\zeta(\kappa)}
        \mathcal{I}(\dP_0|\dR).
    \]

    \item[${(iv)}$]
    The specific Talagrand inequality with parameter 
    $\zeta(\kappa)$ is satisfied: 
    \[
        \mathcal{W}_{\kappa,s}(\dP_0,\dR)
        \leq
        \sqrt{
        \frac{2}{\zeta(\kappa)}
        \mathcal{E}(\dP_0|\dR)
        }.
    \]
\end{enumerate}
\end{satz}
\begin{bem}
  Since we assume in Theorem~\ref{thm:FIstationary} that $\kappa_n$ satisfies \hyperref[property:dobrushin2]{$(\mathbf{D}^\prime(\zeta(\kappa_n)))$} for all $n\in\N$, Theorem~\ref{thm:convergenceresultsstationary}~$(v)$ implies that $\zeta(\kappa)>0$. Hence, all functional inequalities in Theorem~\ref{thm:FIstationary} are non-trivial.
\end{bem}
 \begin{bem}
Extending the exponential contraction of the specific relative entropy as in Section \ref{sec:phientropies} to more general $\Phi$-entropies requires a suitable
normalization, since these functionals need not grow linearly
with the volume. For instance, the power entropies grow
exponentially with the volume even for homogeneous Poisson
processes with different intensities. Normalization by
$\Leb(\Lambda_n)$ therefore generally yields an infinite limit.
\end{bem}
\subsubsection{Proof of Theorem~\ref{thm:evistationary} and Theorem~\ref{thm:FIstationary}}\label{subsec:proofevistationary}
The general proof strategy of Theorem~\ref{thm:evistationary} and Theorem~\ref{thm:FIstationary} is to lift the respective finite-volume results of Theorem~\ref{thm:evi} and Theorem~\ref{thm:FIbounded} to the stationary setting by taking the normalized large-volume limit. 
The proof relies on the following convergence results.

\begin{satz}\label{thm:convergenceresultsstationary}
Assume that $\kappa:\Gamma\times\R^d\to[0,\infty)$ satisfies \ref{property:cocycle}, \ref{property:locallystable}, \ref{property:continuity}, \ref{property:repulsive}, \ref{property:shiftinvariance}, and~\ref{property:finiterange}. Let $\dR\in\RRR_{\kappa,s}(\Gamma)$, $\dR^{\Lambda_n}=\RRR_{\kappa_n,\Leb_{|\Lambda_n}}(\Gamma_{\Lambda_n})$, $n\in\N$, and $\dP_0,\dP_1\in\PPP_s(\Gamma)$. Then 
\begin{enumerate}
    \item[$(i)$]
       $ \lim_{n\to\infty}\mathrm{S}_t^{\kappa_n*}\dP_{|\Lambda_n}
        =
        \mathrm{S}_t^{\kappa*}\dP$ in $\PPP_1(\Gamma)$ for all $\dP\in\PPP_s(\Gamma)$ and $t\geq0$;

    \item[$(ii)$] $\Gamma\text{-}\lim_{n\to\infty}\frac{1}{\Leb(\Lambda_n)}
        \Ent(\cdot
            |
            \dR^{\Lambda_n}
        )=\mathcal{E}(\cdot|\dR)$ on $\PPP_s(\Gamma)\subset\PPP_1(\Gamma)$.  A recovery sequence for $\dP\in\PPP_s(\Gamma)$ is $(\dP_{|\Lambda_n})_{n\in\N}$;

    \item[$(iii)$] $\Gamma\text{-}\lim_{n\to\infty}\frac{1}{\Leb(\Lambda_n)}
        \mathcal{W}_{\kappa_n}^2
        (
            \cdot, \cdot
        )=\mathcal{W}_{\kappa,s}^2(\cdot,\cdot)$ on $\PPP_s(\Gamma)\times\PPP_s(\Gamma)\subset\PPP_1(\Gamma)\times\PPP_1(\Gamma)$.  A recovery sequence for $(\dP_0,\dP_1)\in\PPP_s(\Gamma)\times\PPP_s(\Gamma)$ is  $(\dP_{0|\Lambda_n},\dP_{1|\Lambda_n})_{n\in\N}$;
    \item[$(iv)$] $\lim_{n\to\infty}
        \frac{1}{\Leb(\Lambda_n)}
        \I(\dP_{|\Lambda_n}|\dR^{\Lambda_n})=\mathcal{I}(\dP|\dR)$ for all $\dP\in\PPP_s(\Gamma)$.
\end{enumerate}
Moreover, if $\kappa_n$ satisfies \hyperref[property:dobrushin]{$(\mathbf{D}(\zeta(\kappa_n)))$} for all $n\in\N$, then $\kappa$ satisfies \ref{property:dobrushin} and
\begin{enumerate}
    \item[$(v)$] 
    $0<\lim_{n\to\infty}\zeta(\kappa_n)=\zeta(\kappa)$.
    
\end{enumerate}
\end{satz}
\begin{bem}
  Theorem~\ref{thm:convergenceresultsstationary}~$(iii)$ implies that, for all $\dP_0,\dP_1\in\PPP_s(\Gamma)$, it holds that
\begin{align*}
     \mathcal{W}_{\kappa,s}^2(\dP_0,\dP_1)
     &=\inf\left\{\lim_{n\to\infty}\frac{1}{\Leb(\Lambda_n)}
     \mathcal{A}_{\kappa_n}({\bar{\dP}_{|\Lambda_n}},{\bar{\dV}_{|\Lambda_n}})
     :({\bar{\dP}},{\bar{\dV}})\in\CCC\EEE^s(\dP_0,\dP_1)\right\}\\
     &=\lim_{n\to\infty}\frac{1}{\Leb(\Lambda_n)}
     \mathcal W_{\kappa_n}^2(\dP_{0|\Lambda_n},\dP_{1|\Lambda_n}).
\end{align*}
Hence, $\mathcal{W}_{\kappa,s}^2$ can be represented  as the normalized large-volume limit of the finite-volume transport distances $\mathcal W_{\kappa_n}^2$, justifying the name \emph{specific non-local transport distance}.
\end{bem}
\begin{proof}[Proof of Theorem~\ref{thm:evistationary}]
For every $n\in\N$, the finite-volume Papangelou intensity
$\kappa_n$ satisfies \ref{property:cocycle},
\hyperref[property:locallystable]{$(\mathbf{S}(1))$},
\ref{property:continuity2}, \ref{property:repulsive},
\ref{property:ALB}, and
\hyperref[property:dobrushin2]{$(\mathbf{D}^\prime(\zeta(\kappa_n)))$}. By Lemma~\ref{lem:boundcands},
\[
\frac{1}{\Leb(\Lambda_n)}
\mathcal{W}_{\kappa_n}^2
(\dP_{0|\Lambda_n},\dP_{1|\Lambda_n})
\leq
\mathcal{W}_{\kappa,s}^2(\dP_0,\dP_1)<\infty.
\]
Moreover, Lemma~\ref{lem:limsupentropiestat} yields
$\Ent(\dP_{1|\Lambda_n}|\dR^{\Lambda_n})<\infty$
for all sufficiently large $n$. Hence, we obatin by Theorem~\ref{thm:evi} that
\begin{align*}
    &\liminf_{n\to\infty}\frac{e^{\zeta(\kappa_n)t}}{2}\frac{1}{\Leb(\Lambda_n)}\mathcal{W}_{\kappa_n}^2(\mathrm{S}_t^{\kappa_n*}\dP_{0|\Lambda_n},\dP_{1|\Lambda_n})-\frac{1}{\Leb(\Lambda_n)}\frac{1}{2}\mathcal{W}_{\kappa_n}^2(\dP_{0|\Lambda_n},\dP_{1|\Lambda_n})\\
    &\quad\leq \liminf_{n\to\infty} E_{\zeta(\kappa_n)}(t)\left(\frac{1}{\Leb(\Lambda_n)}\Ent(\dP_{1|\Lambda_n} | \dR^{\Lambda_n})-
\frac{1}{\Leb(\Lambda_n)}\Ent(\mathrm{S}_t^{\kappa_n*}\dP_{0|\Lambda_n} | \dR^{\Lambda_n})\right),\quad t>0.
\end{align*}
Applying the convergence results of Theorem~\ref{thm:convergenceresultsstationary} finishes the proof.
\end{proof}

\begin{proof}[Proof of Theorem~\ref{thm:FIstationary}]
The proof is analogous to that of Theorem~\ref{thm:evistationary} and follows directly from the finite-volume functional inequalities in Theorem~\ref{thm:FIbounded} and the convergence results in Theorem~\ref{thm:convergenceresultsstationary}.
\end{proof}
\subsection{Transport Distance for Stationary Measures}\label{sec:stationarydistance}
In this section, we discuss the central properties of $$\mathcal{W}_{\kappa,s}^2(\dP_0,\dP_1)= \inf \left\{ \mathcal{A}_{\kappa,s}(\Bar{\dP}, \Bar{\dV}) : (\Bar{\dP}, \Bar{\dV}) \in \CCC\EEE^s(\dP_0, \dP_1) \right\},\quad \dP_0,\dP_1\in\PPP_s(\Gamma).$$ 
Since the action functional
\begin{align*}
\mathcal{A}_{\kappa,s}({\Bar{\dP}},{\Bar{\dV}})
=\sup_{n\in\N}\frac{1}{\Leb(\Lambda_n)}\int_0^1
\mathcal{L}_{\kappa_n}(\Bar{\dP}_{t|\Lambda_n},\Bar{\dV}_{t|\Lambda_n})\mathrm{d}t
\end{align*}
is defined analogously to the Poisson setting $\kappa \equiv 1$ \cite{huesmann2025nonlocalwassersteingeometrygradient}, they share the same structural properties. Therefore, we present these properties in a concise form and either give short proofs or refer to the corresponding proofs in \cite{huesmann2025nonlocalwassersteingeometrygradient}. The result of this Section hold for general mean function $\theta$ satisfying the assumptions (A1) -- (A7).

\begin{lem}\label{lem:propertiesactionstationary}
Assume that $\kappa:\Gamma\times\R^d\to[0,\infty)$ satisfies \ref{property:locallystable}, \ref{property:continuity}, and~\ref{property:shiftinvariance}. 
\begin{enumerate}
    \item [$(i)$]  The specific action functional $\mathcal{A}_{\kappa,s}$ is jointly  convex, i.e.\ for all $({\bar{\dP}}^1,{\bar{\dV}}^1),({\bar{\dP}}^2,{\bar{\dV}}^2)$  $\in\CCC\EEE^s(\dP_0,\dP_1)$ it holds that 
             \begin{align*}
                 \mathcal{A}_{\kappa,s}(c{\bar{\dP}}^1+(1-c){\bar{\dP}}^2,c{\bar{\dV}}^1+(1-c){\bar{\dV}}^2)\leq c\mathcal{A}_{\kappa,s} ({\bar{\dP}}^1,{\bar{\dV}}^1)+(1-c) \mathcal{A}_{\kappa,s}({\bar{\dP}}^2,{\bar{\dV}}^2),\quad c\in[0,1].
             \end{align*}
    \item [$(ii)$] The specific action functional $\mathcal{A}_{\kappa,s}$ is lower semicontinuous, i.e.\ let $\dP_0,\dP_1\in\PPP_s(\Gamma)$ and $(\dP_0^i)_{i\in\N},$ $(\dP_1^i)_{i\in\N}\subset\PPP_s(\Gamma)$ s.t.\   $\lim_{i\to\infty}\dP_0^i=\dP_0$ and $\lim_{i\to\infty}\dP_1^i=\dP_1$ in $\PPP_s(\Gamma)$. Furthermore,  let $(\bar{\dP}^i,\bar{\dV}^i)\in\CCC\EEE^s(\dP_0^i,\dP_1^i), \;i\in\N$, and $(\bar{\dP},\bar{\dV})\in\CCC\EEE^s(\dP_0,\dP_1)$  s.t.\   $\lim_{i\to\infty}\bar\dP_t^i=\bar\dP_t$ in $\PPP_s(\Gamma)$ for all $t\in[0,1]$ and $\lim_{i\to\infty}\bar\dV^i=\bar\dV$ in $\MMM_{s}(\Gamma\times \R^d\times[0,1])$, then
    \begin{align*}
        \mathcal{A}_{\kappa,s}(\bar{\dP},\bar{\dV})\leq \liminf_{i\to\infty}\mathcal{A}_{\kappa,s}(\bar{\dP}^i,\bar{\dV}^i).
    \end{align*}
    \item[$(iii)$]  Let $(\dP_0^i)_{i\in\N},(\dP_1^i)_{i\in\N}\subset\PPP_s(\Gamma)$ and $\dP_0,\dP_1\in\PPP_s(\Gamma)$ s.t.\   $\lim_{i\to\infty}\dP_0^i=\dP_0$ and $\lim_{i\to\infty}\dP_1^i=\dP_1$ in $\PPP_s(\Gamma)$. Further, let $(\bar{\dP}^i,\bar{\dV}^i)\in\CCC\EEE^s(\dP_0^i,\dP_1^i),$ $i\in\N,$ s.t.\   
        \begin{align*}
            \sup_{i\in\N}\mathcal{A}_{\kappa,s}(\bar{\dP}^i,\bar{\dV}^i)\leq C<\infty.
        \end{align*}
        There exists a subsequence $(\bar{\dP}^{i_j},\bar{\dV}^{i_j})_{j\in\N}$ and a solution $(\bar{\dP},\bar{\dV})\in\CCC\EEE^s(\dP_0,\dP_1)$ s.t.\   
        \begin{align*}
            \lim_{j\to\infty}\bar{\dP}_t^{i_j}=\bar{\dP}_t\quad &\text{ in }\PPP_s(\Gamma) \text{ for all } t\in [0,1],\\
            \lim_{j\to\infty}\bar{\dV}^{i_j}=\bar{\dV}\quad &\text{ in }\MMM_{s}(\Gamma\times \R^d\times[0,1]).
             \end{align*}
\end{enumerate}
\end{lem}
\begin{proof}
      The proof is identical to the proof of \cite[Lemma 2.7]{huesmann2025nonlocalwassersteingeometrygradient}.
\end{proof}

The following theorem establishes the existence of minimizers of the optimization problem in the definition of $\mathcal{W}_{\kappa,s}.$
  \begin{satz} \label{thm:minimizerstationary}
    Assume that $\kappa:\Gamma\times\R^d\to[0,\infty)$ satisfies \ref{property:locallystable}, \ref{property:continuity}, and~\ref{property:shiftinvariance}.   For all $\dP_0,\dP_1\in\PPP_s(\Gamma)$   there exists a stationary solution $({\bar{\dP}},{\bar{\dV}})\in\CCC\EEE^s(\dP_0,\dP_1)$ s.t.\   
    \begin{align*}
      \mathcal{W}_{\kappa,s}^2(\dP_0,\dP_1)=\mathcal{A}_{\kappa,s}({\bar{\dP}},{\bar{\dV}}).
    \end{align*}
    \end{satz}
    \begin{proof}
        The claim is a consequence of the direct method of calculus of variations and Lemma~\ref{lem:propertiesactionstationary}.
    \end{proof}
 By a   standard reparameterization argument, see e.g.\ \cite[Theorem 5.4]{Dolbeault_2008}, \cite[Lemma 1.1.4]{Ambrosio2008GradientFI}, $\mathcal{W}_{\kappa,s}$ can be represented as follows.
      \begin{lem}\label{lem:reparametrization}
      Assume that $\kappa:\Gamma\times\R^d\to[0,\infty)$ satisfies \ref{property:locallystable}, \ref{property:continuity}, and~\ref{property:shiftinvariance}.   Let $\dP_0,\dP_1\in\PPP_s(\Gamma).$ It holds that 
        \begin{align*}
            \mathcal{W}_{\kappa,s}(\dP_0,\dP_1)=\inf\left\{\sup_{n\in\N}\frac{1}{\Leb(\Lambda_n)^{1/2}}\int_0^T\mathcal{L}^{1/2}_{\kappa_n}(\bar{\dP}_{t|\Lambda_n},{\bar{\dV}_{t|\Lambda_n}})\mathrm{d}t:({\bar{\dP}},{\bar{\dV}})\in\CCC\EEE^s_T(\dP_0,\dP_1)\right\},
        \end{align*}
        for all $T>0.$
        \end{lem}
        
        The following theorem shows that $\mathcal{W}_{\kappa,s}$ is indeed a distance and collects several fundamental properties.
        \begin{satz}\label{thm:propertiesdistancestationary}
        Assume that $\kappa:\Gamma\times\R^d\to[0,\infty)$ satisfies \ref{property:locallystable}, \ref{property:continuity}, and~\ref{property:shiftinvariance}.
             \begin{enumerate}
                 \item[$(i)$] The map $\mathcal{W}_{\kappa,s}:\PPP_s(\Gamma)\times\PPP_s(\Gamma)\to[0,\infty]$ is an extended distance on $\PPP_s(\Gamma).$
             
                 \item[${(ii)}$]  The topology induced by $\mathcal{W}_{\kappa,s}$ on $\PPP_s(\Gamma)$ is stronger than the usual $\CCC_1$-weak topology  on $\PPP_s(\Gamma)$.   
                 \item[${(iii)}$] The specific non-local transport distance $\mathcal{W}_{\kappa,s}$ is lower semicontinuous on $\PPP_s(\Gamma)\times \PPP_s(\Gamma)$. 
                 \item[${(iv)}$] 
                 Bounded sets w.r.t.\ $\mathcal{W}_{\kappa,s}$ are $\PPP_s$-relatively compact.
                 \item[${(v)}$]  For every $\dP \in \PPP_s(\Gamma)$, the accessible component $(\{\mathcal{W}_{\kappa,s}(\dP, \cdot) < \infty\}, \mathcal{W}_{\kappa,s}) $ is a complete geodesic space. The geodesics are given by the minimizing curves $\bar{\dP}$ of Theorem~\ref{thm:minimizerstationary}.
            \end{enumerate}
         \end{satz}
         \begin{proof} The proof is analogous to the proof of  \cite[Theorem 5.4]{huesmann2025nonlocalwassersteingeometrygradient} and  \cite[Lemma 5.5]{huesmann2025nonlocalwassersteingeometrygradient}.
         \end{proof}
        We finish this section by showing that the supremum in the definition of the stationary action functional is realized as a limit and thus exists. 
   \begin{lem}\label{lem:changelimsup}
        Assume that $\kappa:\Gamma\times\R^d\to[0,\infty)$ satisfies \ref{property:locallystable}, \ref{property:continuity}, \ref{property:repulsive}, and~\ref{property:shiftinvariance}.  For $\dP_0,\dP_1\in\PPP_s(\Gamma),$ it holds that 
         \begin{align*}
             \mathcal{W}^2_{\kappa,s} (\dP_0,\dP_1)&=\inf\left\{\sup_{n\in\N}\frac{1}{\Leb(\Lambda_n)}\int_0^1\mathcal{L}_{\kappa_n}({\bar{\dP}_{t|\Lambda_n}},{\bar{\dV}_{t|\Lambda_n}})\mathrm{d}t: ({\bar{\dP}},{\bar{\dV}})\in\CCC\EEE^s( \dP_0,\dP_1)\right\}\\
             & =\inf\left\{\lim_{n\to\infty}\frac{1}{\Leb(\Lambda_n)}\int_0^1\mathcal{L}_{\kappa_n}({\bar{\dP}_{t|\Lambda_n}},{\bar{\dV}_{t|\Lambda_n}})\mathrm{d}t: ({\bar{\dP}},{\bar{\dV}})\in\CCC\EEE^s( \dP_0,\dP_1)\right\}.
         \end{align*}
         Moreover,  the limit
         \begin{align*}
             \lim_{n\to\infty}\frac{1}{\Leb(\Lambda_n)}\int_0^1\mathcal{L}_{\kappa_n}({\bar{\dP}_{t|\Lambda_n}},\bar{\dV}_{t|\Lambda_n})\mathrm{d}t
         \end{align*}
         exists in $[0,\infty]$ for every $({\bar{\dP}},{\bar{\dV}})\in\CCC\EEE^s(\dP_0,\dP_1).$
     \end{lem}
     
         The proof is analogous to \cite[Lemma 5.6]{huesmann2025nonlocalwassersteingeometrygradient} and follows directly from Fekete's Lemma, see e.g., \cite[p. 74]{fekete} and the following shift invariance and super-additivity of the Lagrange functional.
         \begin{lem}\label{lem:shiftinvariance}
  Assume that $\kappa:\Gamma\times\R^d\to[0,\infty)$ satisfies \ref{property:shiftinvariance}. Let $z\in\R^d$ and $B\in\BBB(\R^d)$. Consider $\dP_1\in\PPP_1(\Gamma_{B}),\;\dP_2\in\PPP_1(\Gamma_{B+z}),\;{\dV}_1\in\MMM_{b,0}(\Gamma_{B}\times B)$, and ${\dV}_2\in\MMM_{b,0}(\Gamma_{B+z}\times (B+z)),$ and assume that $\dP_1=({\Theta_{z}})_\#\dP_2$ and ${\dV}_1=({\Theta_{z}^{\mathrm{prod}}})_\#{\dV}_2.$ Then $$\mathcal{L}_{\kappa,\Leb_{|B}}(\dP_1,{\dV}_1)=\mathcal{L}_{\kappa,\Leb_{|B}}(({\Theta_{z}})_\#\dP_2,({\Theta_{z}^{\mathrm{prod}}})_\#{\dV}_2)=\mathcal{L}_{\kappa,\Leb_{|B+z}}(\dP_2,{\dV}_2).$$
    \end{lem}
    \begin{proof}
        The proof follows directly by choosing $\sigma_1=({\Theta_z^{\mathrm{prod}}})_\#\sigma_2$ in the definition of the Lagrange functional and using the shift-invariance of the Lebesgue measure as well as \ref{property:shiftinvariance}.
    \end{proof}
     \begin{lem}\label{lem:superadditivity}
      Assume that $\kappa:\Gamma\times\R^d\to[0,\infty)$ satisfies \ref{property:repulsive}.  Let $B_1,{B_2}\in\BBB(\R^d)$ be disjoint sets. For $\dP\in\PPP_1(\Gamma_{B_1\cup {B_2}})$ and $\dV\in\MMM_{b,0}(\Gamma_{B_1\cup {B_2}}\times (B_1\cup {B_2}))$, it holds that 
        \begin{align}
               \mathcal{L}_{\kappa,\Leb_{|B_1} }(\dP_{|B_1},\dV_{|B_1})+\mathcal{L}_{\kappa,\Leb_{|{B_2}} }(\dP_{|{B_2}},\dV_{|{B_2}})\leq \mathcal{L}_{\kappa, \Leb_{|B_1\cup {B_2}} }(\dP,\dV).
           \end{align}
    \end{lem}
    \begin{proof}
       Let $\sigma=\dP\otimes\Leb_{|B_1\cup {B_2}}+C_{\dP}+|\dV|$ and denote  $\dP\otimes\Leb_{|B_1\cup {B_2}}=\psi\sigma,\; C_\dP=\phi\sigma,$ and $ \dV=w\sigma.$
       By disintegration, we find probability measures $\sigma^{*,1}(\cdot|\xi,x),\; \xi\in\Gamma_{B_2},\;x\in B_1\cup {B_2},$ on $\Gamma_{B_1}$ s.t.\   
           \begin{align*}
               \int_{\Gamma_{B_1\cup {B_2}}\times (B_1\cup {B_2})} F(\xi,x)\sigma(\mathrm{d}\xi,\mathrm{d}x)=\int_{\Gamma_{B_2}\times (B_1\cup {B_2})}\int _{\Gamma_{B_1}}F(\xi_{|{B_2}}+\xi_1,x)\sigma^{*,1}(\mathrm{d}\xi_1|\xi_{|{B_2}},x)\sigma(\mathrm{d}\xi,\mathrm{d}x),
           \end{align*}
           for all $F\in\FFF^+(\Gamma_{B_1\cup {B_2}}\times (B_1\cup {B_2}))$. Analogously, we define $\sigma^{*,{2}}(\cdot|\xi,x),$ $\xi\in\Gamma_{B_1},\;x\in {B_1}\cup {B_2},$ on $\Gamma_{B_2}.$
           Using these measures, the Lagrange functional can be rewritten by
           \begin{align}
               &\mathcal{L}_{\kappa,\Leb_{| {B_1}\cup {B_2}}}(\dP,\dV)\\
               &\quad=\int_{\Gamma_{{B_1}\cup {B_2}}\times ({B_1}\cup {B_2})}\alpha (\kappa(\xi,x)\psi(\xi,x),\phi(\xi,x),w(\xi,x))\sigma(\mathrm{d}\xi,\mathrm{d}x)\nonumber\\
               &\quad= \int_{\Gamma_{{B_1}\cup {B_2}}\times ({B_1}\cup {B_2})}\ind_{\{\Gamma_{{B_1}\cup {B_2}}\times {B_1}\}}(\xi,x)\alpha (\kappa(\xi,x)\psi(\xi,x),\phi(\xi,x),w(\xi,x))\sigma(\mathrm{d}\xi,\mathrm{d}x)\nonumber\\
               &\qquad+\int_{\Gamma_{{B_1}\cup {B_2}}\times ({B_1}\cup {B_2})}\ind_{\{\Gamma_{{B_1}\cup {B_2}}\times {B_2}\}}(\xi,x)\alpha (\kappa(\xi,x)\psi(\xi,x),\phi(\xi,x),w(\xi,x))\sigma(\mathrm{d}\xi,\mathrm{d}x)\nonumber\\
               &\quad= \int_{\Gamma_{B_1}\times {B_1}} \int_{\Gamma_{B_2}} \alpha (\kappa(\xi_1+\xi_2,x)\psi(\xi_1+\xi_2,x),\phi(\xi_1+\xi_2,x),w(\xi_1+\xi_2,x))\sigma^{*,{2}}(\mathrm{d}\xi_2|\xi_{1},x)\sigma_{|{B_1}}(\mathrm{d}\xi_1,\mathrm{d}x)\nonumber\\
               &\qquad+ \int_{\Gamma_{B_2}\times {B_2}} \int_{\Gamma_{B_1}}  \alpha (\kappa(\xi_1+\xi_2,x)\psi(\xi_1+\xi_2,x),\phi(\xi_1+\xi_2,x),w(\xi_1+\xi_2,x))\sigma^{*,{1}}(\mathrm{d}\xi_1|\xi_{2},x)\sigma_{|{B_2}}(\mathrm{d}\xi_2,\mathrm{d}x) .\label{eq:1022}
           \end{align}
           For any $K\in\BBB(\Gamma_{B_1})$ and $B\in\BBB({B_1})$, it holds that 
         \begin{align*}
             \int\ind_{K\times B}(\xi,x)\dP_{|{B_1}}\otimes\Leb_{|{B_1}}(\mathrm{d}\xi,\mathrm{d}x)&= \int\ind_{K\times B}(\xi_{|{B_1}},x)\ind_{{B_1}}(x) \dP\otimes\Leb_{|{B_1}\cup {B_2}}(\mathrm{d}\xi,\mathrm{d}x)\\
             &= \int\ind_{K\times B}(\xi_{|{B_1}},x)\ind_{{B_1}}(x)\psi(\xi,x)\sigma(\mathrm{d}\xi,\mathrm{d}x)\\
             &= \int\ind_{K\times B}(\xi_{1},x)\int\psi(\xi_{1}+\xi_2,x)\sigma^{*,{2}}(\mathrm{d}\xi_2|\xi_{1},x)\sigma_{|{B_1}}(\mathrm{d}\xi_1,\mathrm{d}x).
         \end{align*}
         This computation implies that $$\dP_{|{B_1}}\otimes\Leb_{|{B_1}}(\mathrm{d}\xi_1,\mathrm{d}x)=\int\psi(\xi_{1}+\xi_2,x)\sigma^{*,{2}}(\mathrm{d}\xi_2|\xi_{1},x)\sigma_{|{B_1}}(\mathrm{d}\xi_1,\mathrm{d}x)=:\psi^{{1}}(\xi_1,x)\sigma_{|{B_1}}(\mathrm{d}\xi_1,\mathrm{d}x).$$ Similarly, we obtain  
         \begin{align*}
             C_{\dP_{|{B_1}}}(\mathrm{d}\xi_1,\mathrm{d}x)&=\int\phi(\xi_{1}+\xi_2,x)\sigma^{*,{2}}(\mathrm{d}\xi_2|\xi_{1},x)\sigma_{|{B_1}}(\mathrm{d}\xi_1,\mathrm{d}x)=:\phi^{{1}}(\xi_1,x)\sigma_{|{B_1}}(\mathrm{d}\xi_1,\mathrm{d}x),\\
            \dV_{|{B_1}}(\mathrm{d}\xi_1,\mathrm{d}x)&=\int w(\xi_{1}+\xi_2,x)\sigma^{*,{2}}(\mathrm{d}\xi_2|\xi_{1},x)\sigma_{|{B_1}}(\mathrm{d}\xi_1,\mathrm{d}x)=:w^{{1}}(\xi_1,x)\sigma_{|{B_1}}(\mathrm{d}\xi_1,\mathrm{d}x).
         \end{align*}
         Since $\sigma^{*,{2}}$ is a probability measure, we may apply Jensen's inequality below. By first using  $\kappa(\xi_1+\xi_2,x)\leq \kappa(\xi_1,x)$ \ref{property:repulsive} together with the monotonicity of $\alpha$ in its first two components and then applying the convexity and homogeneity of $\alpha$ together with Jensen's inequality, we obtain for the first summand of \eqref{eq:1022} that
         \begin{align}
             &\int_{\Gamma_{B_1}\times {B_1}} \int_{\Gamma_{B_2}} \alpha (\kappa(\xi_1+\xi_2,x)\psi(\xi_1+\xi_2,x),\phi(\xi_1+\xi_2,x),w(\xi_1+\xi_2,x))\sigma^{*,{2}}(\mathrm{d}\xi_2|\xi_{1},x)\sigma_{|{B_1}}(\mathrm{d}\xi_1,\mathrm{d}x)\nonumber\\
             &\quad\geq \int_{\Gamma_{B_1}\times {B_1}} \int_{\Gamma_{B_2}} \alpha (\kappa(\xi_1,x)\psi(\xi_1+\xi_2,x),\phi(\xi_1+\xi_2,x),w(\xi_1+\xi_2,x))\sigma^{*,{2}}(\mathrm{d}\xi_2|\xi_{1},x)\sigma_{|{B_1}}(\mathrm{d}\xi_1,\mathrm{d}x)\nonumber\\
             &\quad\geq \int_{\Gamma_{B_1}\times {B_1}}\alpha(\kappa(\xi_1,x)\psi^{1}(\xi_1,x),\phi^{1}(\xi_1,x),w^{1}(\xi_1,x))\sigma_{|{B_1}}(\mathrm{d}\xi_1,\mathrm{d}x)\nonumber\\
             &\quad=\mathcal{L}_{\kappa,\Leb_{|{B_1}}}(\dP_{|{B_1}},\dV_{|{B_1}}).\label{eq:1023}
         \end{align}
         Analogously, we obtain for the second summand of \eqref{eq:1022} that 
         \begin{align}
             &\int_{\Gamma_{B_2}\times {B_2}} \int_{\Gamma_{B_1}} \alpha (\kappa(\xi_1+\xi_2,x)\psi(\xi_1+\xi_2,x),\phi(\xi_1+\xi_2,x),w(\xi_1+\xi_2,x))\sigma^{*,{1}}(\mathrm{d}\xi_1|\xi_{2},x)\sigma_{|{B_2}}(\mathrm{d}\xi_2,\mathrm{d}x)\nonumber\\
             &\quad\geq \mathcal{L}_{\kappa,\Leb_{|{B_2}}}(\dP_{|{B_2}},\dV_{|{B_2}}).\label{eq:1024}
         \end{align}
         Inserting \eqref{eq:1023} and \eqref{eq:1024} into \eqref{eq:1022} finishes the proof.
    \end{proof}
              \section{Proof of Theorem~\ref{thm:convergenceresultsstationary}}\label{sec:Proofconvergenceresultsstationary}
   Of particular importance in the proof of Theorem~\ref{thm:convergenceresultsstationary} is the finite-range property~\ref{property:finiterange}. In the following, let $0\leq r<\infty$ always denote the range of the interaction. For $n\in\N$ with $n>2r$, we define and recall
\[
    \Lambda_n^{r-}=[-n/2+r,n/2-r)^d
    \quad \text{and} \quad
    \Lambda_n^{r+}:=[-n/2-r,n/2+r)^d .
\]
Further, let
\(
    \kappa_n^{r\pm}(\xi,x):=\kappa(\xi_{|\Lambda_n^{r\pm}},x)\ind_{\Lambda_n^{r\pm}}(x),\;\xi\in\Gamma,\; x\in\R^d .
\)
Consistent with the notation introduced before, we write $\mathcal{L}_{\kappa_n^{r\pm}}$ and $\mathcal{W}_{\kappa_n^{r\pm}}$ instead of $\mathcal{L}_{\kappa_n^{r\pm},\Leb_{|\Lambda_n^{r\pm}}}$ and $\mathcal{W}_{\kappa_n^{r\pm},\Leb_{|\Lambda_n^{r\pm}}}$, respectively.\\

\subsection{Papangelou Measures}\label{subsec:papangeloumeasuresconvergence}
For the proof of Theorem~\ref{thm:convergenceresultsstationary}, we need to
compare the finite-volume Papangelou measure
$\dR^{\Lambda_n}=\RRR_{\kappa_n,\Leb_{|\Lambda_n}}(\Gamma_{\Lambda_n})$
with the restriction $\dR_{|\Lambda_n}$ of an infinite-volume Papangelou measure
$\dR\in\RRR_{\kappa,s}(\Gamma)$. The following lemma shows that $\dR_{|\Lambda_n}$ is itself a finite-volume Papangelou measure with Papangelou intensity obtained by averaging the infinite-volume intensity $\kappa$ over all possible configurations on $\Lambda_n^C$.
\begin{lem}\label{lem:restrictionpapangelou}
Assume that $\kappa:\Gamma\times\R^d\to[0,\infty)$ satisfies
\ref{property:cocycle}, \ref{property:locallystable},  and~\ref{property:finiterange}, and let $\dR\in\RRR_{\kappa,s}(\Gamma)$. The restriction $\dR_{|\Lambda_n}$ is a Papangelou measure with reference measure $\Leb_{|\Lambda_n}$ and Papangelou intensity
\begin{align*}
\kappa_{|\Lambda_n}(\xi_1,x)
:=
\int
\kappa(\xi_1+\xi_2,x)
\dR^{\Lambda_n^C}(\mathrm d\xi_2| \xi_1),
\quad
\xi_1\in\Gamma_{\Lambda_n},\; x\in\Lambda_n,
\end{align*}
where $\dR^{\Lambda_n^C}(\cdot| \xi_1)$ is defined by 
\begin{align*}
\int F(\xi)\dR(\mathrm d\xi)
=
\int\int
F(\xi_1+\xi_2)
\dR^{\Lambda_n^C}(\mathrm d\xi_2| \xi_1)
\dR_{|\Lambda_n}(\mathrm d\xi_1),
\quad F\in\FFF^+(\Gamma).
\end{align*}
\end{lem}
\begin{proof}
Let $F\in\FFF^+(\Gamma_{\Lambda_n}\times\Lambda_n)$. By the definition of $\kappa_{|\Lambda_n}$ and the GNZ equation~\eqref{eq:papangelou}, we obtain
\begin{align*}
\int\int
F(\xi+\delta_x,x)\kappa_{|\Lambda_n}(\xi,x)
\Leb_{|\Lambda_n}(\mathrm dx)\dR_{|\Lambda_n}(\mathrm d\xi)& =
\int\int\ind_{\Lambda_n}(x)
F(\xi_{|\Lambda_n}+\delta_x,x)\kappa(\xi,x)
\Leb(\mathrm dx)\dR(\mathrm d\xi) \\
&=
\int\int
F(\xi,x)\xi(\mathrm{d}x)\dR_{|\Lambda_n}(\mathrm d\xi),
\end{align*}
finishing the proof.
\end{proof}
We proceed with a discussion of the relation between the Papangelou measures $\dR^{\Lambda_n}$ and $\dR_{|\Lambda_n}$.
\begin{lem}\label{lem:propertiespapangelou}
Assume that $\kappa:\Gamma\times\R^d\to[0,\infty)$ satisfies
\ref{property:cocycle}, \ref{property:locallystable}, \ref{property:repulsive},
\ref{property:finiterange}, and~\ref{property:shiftinvariance}. Let
$\dR\in\RRR_{\kappa,s}(\Gamma)$ and
$\dR^{\Lambda_n}=\RRR_{\kappa_n,\Leb_{|\Lambda_n}}(\Gamma_{\Lambda_n})$.
Then $\dR_{|\Lambda_n}\ll\dR^{\Lambda_n}$ and
\begin{align*}
  f_n(\xi)
  := 
    \frac{\mathrm d\dR_{|\Lambda_n}}
    {\mathrm d\dR^{\Lambda_n}}(\xi)=
  \int
  \frac{Z_{\kappa_n,\Leb_{|\Lambda_n}}^{\Lambda_n}}
  {Z_{\kappa_n,\Leb_{|\Lambda_n}}^{\Lambda_n}
  (\eta_{|\Lambda_n^C})}
  \frac{\varphi^{\eta_{|\Lambda_n^C}}_{\kappa_n}(\xi)}
  {\varphi_{\kappa_n}(\xi)}
  \dR(\mathrm{d}\eta),
  \quad \xi\in\Gamma_{\Lambda_n}.
\end{align*}
The density $f_n$ is uniformly bounded
\begin{align}
    e^{-A\Leb(\Lambda_n^{r+}\setminus\Lambda_n)}
    \leq f_n(\xi)
    \leq e^{\Leb(\Lambda_n\setminus\Lambda_n^{r-})},
    \quad \xi\in\Gamma_{\Lambda_n},
    \label{eq:bounddensity}
\end{align}
and it holds that
\begin{align}
    \kappa(\xi,x)f_n(\xi+\delta_x)
    =
    \kappa_{|\Lambda_n}(\xi,x)f_n(\xi)
    \quad
    \dR^{\Lambda_n}\otimes\Leb_{|\Lambda_n}\text{-a.e.}
    \label{eq:57}
\end{align}
\end{lem}
\begin{proof}
 Let $K\in\BBB(\Gamma_{\Lambda_n})$. Using the DLR equation
 (Theorem~\ref{thm:DLR}) and Lemma~\ref{lem:finitevolumePapangelou},
we obtain
\begin{align*}
    \int \ind_K(\xi)\dR_{|\Lambda_n}(\mathrm d\xi)
     =
    \int
    \ind_K(\xi)
    \int
    \frac{Z_{\kappa_n,\Leb_{|\Lambda_n}}^{\Lambda_n}}
    {Z_{\kappa_n,\Leb_{|\Lambda_n}}^{\Lambda_n}
    (\eta_{|\Lambda_n^C})}
    \frac{
    \varphi^{\eta_{|\Lambda_n^C}}_{\kappa_n}(\xi)}
    {
    \varphi_{\kappa_n}(\xi)}
    \dR(\mathrm d\eta)
    \dR^{\Lambda_n}(\mathrm d\xi).
\end{align*}
Note that all expressions are well-defined by
\ref{property:repulsive}. Hence, $\dR_{|\Lambda_n}\ll\dR^{\Lambda_n}$
and
\begin{align*}
    f_n(\xi)
    &=
    \frac{\mathrm d\dR_{|\Lambda_n}}
    {\mathrm d\dR^{\Lambda_n}}(\xi)=
    \int
    \frac{Z_{\kappa_n,\Leb_{|\Lambda_n}}^{\Lambda_n}}
    {Z_{\kappa_n,\Leb_{|\Lambda_n}}^{\Lambda_n}
    (\eta_{|\Lambda_n^C})}
    \frac{\varphi^{\eta_{|\Lambda_n^C}}_{\kappa_n}(\xi)}
    {\varphi_{\kappa_n}(\xi)}
    \dR(\mathrm d\eta).
\end{align*}
We continue with the proof of \eqref{eq:bounddensity} and note that by \ref{property:finiterange},
\begin{align}
    Z_{\kappa_n,\Leb_{|\Lambda_n}}^{\Lambda_n}(\eta_{|\Lambda_n^C})&=\sum_{j\geq0}\frac{1}{j!}\int_{\Lambda_n}\ldots\int_{\Lambda_n}\varphi^{\eta_{|\Lambda_n^C}}_{\kappa_n}\left(\sum_{i=1}^j\delta_{x_i}\right)\Leb(\mathrm dx_1)\ldots\Leb(\mathrm dx_j)\nonumber\\
    &\geq \sum_{j\geq0}\frac{1}{j!}\int_{\Lambda_n^{r-}}\ldots\int_{\Lambda_n^{r-}}\varphi^{\eta_{|\Lambda_n^C}}_{\kappa_n}\left(\sum_{i=1}^j\delta_{x_i}\right)\Leb(\mathrm dx_1)\ldots\Leb(\mathrm dx_j)\nonumber\\
    &=\sum_{j\geq0}\frac{1}{j!}\int_{\Lambda_n^{r-}}\ldots\int_{\Lambda_n^{r-}}\varphi_{\kappa_n}\left(\sum_{i=1}^j\delta_{x_i}\right)\Leb(\mathrm dx_1)\ldots\Leb(\mathrm dx_j)\nonumber\\
    &= Z_{\kappa_n,\Leb_{|\Lambda_n^{r-}}}^{\Lambda_n^{r-}},\label{eq:37}
\end{align}
and 
\begin{align}
Z_{\kappa_n,\Leb_{|\Lambda_n}}^{\Lambda_n}&=\sum_{j\geq0}\frac{1}{j!}\int_{\Lambda_n}\ldots\int_{\Lambda_n}\varphi_{\kappa_n}\left(\sum_{i=1}^j\delta_{x_i}\right)\Leb(\mathrm dx_1)\ldots\Leb(\mathrm dx_j)\nonumber\\
    &=e^{\Leb(\Lambda_n)}\int \varphi_{\kappa_n}(\xi)\Poi^{\Leb_{|\Lambda_n}}(\mathrm d\xi)\nonumber\\
    &=e^{\Leb(\Lambda_n\setminus\Lambda_n^{r-})}e^{\Leb(\Lambda_n^{r-})}\int\int\varphi_{\kappa_n}(\xi+\eta)\Poi^{\Leb_{|\Lambda_n^{r-}}}(\mathrm d\xi)\Poi^{\Leb_{|\Lambda_n\setminus\Lambda_n^{r-}}}(\mathrm d\eta)\nonumber\\
    &\leq e^{\Leb(\Lambda_n\setminus\Lambda_n^{r-})}e^{\Leb(\Lambda_n^{r-})}\int\int\varphi_{\kappa_n}(\xi)\Poi^{\Leb_{|\Lambda_n^{r-}}}(\mathrm d\xi)\Poi^{\Leb_{|\Lambda_n\setminus\Lambda_n^{r-}}}(\mathrm d\eta)\nonumber\\
    &=e^{\Leb(\Lambda_n\setminus\Lambda_n^{r-})} Z_{\kappa_n,\Leb_{|\Lambda_n^{r-}}}^{\Lambda_n^{r-}}.\label{eq:38}
\end{align}
Combining \eqref{eq:37} and \eqref{eq:38} and using \ref{property:repulsive}, we deduce that 
\begin{align}
f_n(\xi)&=\int\frac{Z_{\kappa_n,\Leb_{|\Lambda_n}}^{\Lambda_n}}{Z_{\kappa_n,\Leb_{|\Lambda_n}}^{\Lambda_n}(\eta_{|\Lambda_n^C})}\frac{\varphi^{\eta_{|\Lambda_n^C}}_{\kappa_n}(\xi)}{\varphi_{\kappa_n}(\xi)}\dR(\mathrm d\eta)\leq
    \int\frac{e^{\Leb(\Lambda_n\setminus\Lambda_n^{r-})} Z_{\kappa_n,\Leb_{|\Lambda_n^{r-}}}^{\Lambda_n^{r-}}}{ Z_{\kappa_n,\Leb_{|\Lambda_n^{r-}}}^{\Lambda_n^{r-}}}\frac{\varphi_{\kappa_n}(\xi)}{\varphi_{\kappa_n}(\xi)}\dR(\mathrm d\eta)=e^{\Leb(\Lambda_n\setminus\Lambda_n^{r-})}.\label{eq:39}
\end{align}
On the other hand, \ref{property:repulsive} yields \begin{align}\label{eq:41}
    Z_{\kappa_n,\Leb_{|\Lambda_n}}^{\Lambda_n}\geq Z_{\kappa_n,\Leb_{|\Lambda_n}}^{\Lambda_n}(\eta_{|\Lambda_n^C}),
\end{align}
and by the DLR equation (Theorem~\ref{thm:DLR}) and~\ref{property:locallystable}, we obtain \begin{align}
    \int\frac{\varphi^{\eta_{|\Lambda_n^C}}_{\kappa_n}(\xi)}{\varphi_{\kappa_n}(\xi)}\dR(\mathrm d\eta)&\geq \int\frac{\varphi^{\eta_{|\Lambda_n^C}}_{\kappa_n}(\xi)}{\varphi_{\kappa_n}(\xi)}\ind_{\{\eta(\Lambda_n^{r+}\setminus\Lambda_n)=0\}}\dR(\mathrm d\eta)\nonumber\\
    &=\int \ind_{\{\eta(\Lambda_n^{r+}\setminus\Lambda_n)=0\}}\dR(\mathrm d\eta)\nonumber\\
    &=\int\int \ind_{\{\eta(\Lambda_n^{r+}\setminus\Lambda_n)=0\}}\dR^{\Lambda_n^{r+}\setminus\Lambda_n}(\mathrm d\eta|\xi_{|(\Lambda_n^{r+}\setminus\Lambda_n)^C})\dR(\mathrm d\xi)\nonumber\\
    &=\int \frac{1}{Z_{\kappa,\Leb_{|\Lambda_n^{r+}\setminus\Lambda_n}}^{\Lambda_n^{r+}\setminus\Lambda_n}(\xi_{|(\Lambda_n^{r+}\setminus\Lambda_n)^C})}\dR(\mathrm d\xi)\nonumber\\
    &\geq e^{-A\Leb(\Lambda_n^{r+}\setminus\Lambda_n)}.\label{eq:40}
\end{align}
Combining \eqref{eq:41} and \eqref{eq:40} yields 
\begin{align} f_n(\xi)&=\int\frac{Z_{\kappa_n,\Leb_{|\Lambda_n}}^{\Lambda_n}}{Z_{\kappa_n,\Leb_{|\Lambda_n}}^{\Lambda_n}(\eta_{|\Lambda_n^C})}\frac{\varphi^{\eta_{|\Lambda_n^C}}_{\kappa_n}(\xi)}{\varphi_{\kappa_n}(\xi)}\dR(\mathrm d\eta)\geq e^{- A\Leb(\Lambda_n^{r+}\setminus\Lambda_n)},\label{eq:42}
\end{align}
finishing the proof of \eqref{eq:bounddensity}. It remains to prove \eqref{eq:57}. Let $F\in\FFF^+(\Gamma_{\Lambda_n}\times\Lambda_n)$. By the GNZ equation~\eqref{eq:papangelou}, we obtain
\begin{align*}
    \int\int
    F(\xi,x)
    \kappa(\xi,x)
    f_n(\xi+\delta_x)
    \dR^{\Lambda_n}(\mathrm d\xi)
    \Leb_{|\Lambda_n}(\mathrm dx)
    & =
    \int\int
    F(\xi-\delta_x,x)
    f_n(\xi)
    \xi(\mathrm dx)
    \dR^{\Lambda_n}(\mathrm d\xi)\\
    & =
    \int\int
    F(\xi-\delta_x,x)
    \xi(\mathrm dx)
    \dR_{|\Lambda_n}(\mathrm d\xi)\\
    &=
    \int\int
    F(\xi,x)
    \kappa_{|\Lambda_n}(\xi,x)
    \dR_{|\Lambda_n}(\mathrm d\xi)
    \Leb_{|\Lambda_n}(\mathrm dx)\\
    &=
    \int\int
    F(\xi,x)
    \kappa_{|\Lambda_n}(\xi,x)
    f_n(\xi)
    \dR^{\Lambda_n}(\mathrm d\xi)
    \Leb_{|\Lambda_n}(\mathrm dx), 
\end{align*}
proving \eqref{eq:57} and finishing the proof.
\end{proof}

\subsection{Tiled and Stationarized Measures}\label{subsec:tiledstationarized}
In this subsection, we introduce tiled and stationarized measures and discuss their main properties.
\begin{dfn}[Tiled and stationarized measures]\label{dfn:tiledstationarized}
Let $n\in\N$, $\dP\in\PPP_1(\Gamma_{\Lambda_n})$ and
$\dV\in\MMM_{b,0}(\Gamma_{\Lambda_n}\times \Lambda_n)$. The \textit{tiled measure} of $\dP$ is defined by
\begin{align*}
    \dP^{\mathrm{til}}
    :=
    \bigoplus_{z\in\Z^d}(\Theta_{nz})_{\#}\dP,
\end{align*}
and its \textit{stationarized measure} is defined by
\begin{align*}
    \dP^{\mathrm{stat}}
    :=
    \fint_{\Lambda_n}(\Theta_u)_\#\dP^{\mathrm{til}}\Leb(\mathrm{d}u),
\end{align*}
where $\fint_{\Lambda_n}:=\frac{1}{\Leb(\Lambda_n)}\int_{\Lambda_n}$. Further, the \textit{tiled velocity field} of $\dV$ w.r.t.\ $\dP$ is defined by
\begin{align*}
    \dV^{\mathrm{til}}
    :=
    \bigoplus_{z\in\Z^d}^{\mathrm{prod}}
    \left((\Theta_{nz}^{\mathrm{prod}})_\#\dV,(\Theta_{nz})_\#\dP\right),
\end{align*}
see \eqref{eq:superpositionvelocity} for the definition of the superposition on the product space. The corresponding \textit{stationarized velocity field} is defined by
\begin{align*}
    \dV^{\mathrm{stat}}
    :=
    \fint_{\Lambda_n}(\Theta_u^{\mathrm{prod}})_\#\dV^{\mathrm{til}}\Leb(\mathrm{d}u).
\end{align*}
\end{dfn}
\begin{bem}\label{rem:limstat}
Every stationary measure can be approximated by the stationarizations of its
finite-volume restrictions, i.e.\ for every $\dP\in\PPP_s(\Gamma)$ it holds that $\lim_{n\to\infty}(\dP_{|\Lambda_n})^{\mathrm{stat}}=\dP$ in $\PPP_s(\Gamma),$ see Remark \cite[Remark 6.4 $(ii)$]{huesmann2025nonlocalwassersteingeometrygradient}.
\end{bem}
The following lemma shows that the tiled and stationarized versions of a solution to the continuity equation are again solutions to the continuity equation.
\begin{lem}[{\cite[Lemma 6.6]{huesmann2025nonlocalwassersteingeometrygradient}}]\label{lem:stationarysolutions}
    Let $\dP_0, \dP_1\in\PPP_1(\Gamma_{\Lambda_n})$. If $({\bar{\dP}},{\bar{\dV}})\in\CCC\EEE^{\Lambda_n}(\dP_0, \dP_1),$ then \begin{align*}
        ({\bar{\dP}}^{\mathrm{til}},{\bar{\dV}}^{\mathrm{til}})\in\CCC\EEE(\dP_0^{\mathrm{til}}, \dP_1^{\mathrm{til}})\quad \text{ and } \quad ({\bar{\dP}}^{\mathrm{stat}},{\bar{\dV}}^{\mathrm{stat}})\in\CCC\EEE^s(\dP_0^{\mathrm{stat}}, \dP_1^{\mathrm{stat}}).
    \end{align*}
\end{lem}     
       \subsection{Dual Semigroup}
In this subsection, we prove  Theorem~\ref{thm:convergenceresultsstationary} $(i)$ and show that if $\kappa:\Gamma\times\R^d\to[0,\infty)$ satisfies \ref{property:cocycle}, \ref{property:locallystable}, \ref{property:finiterange}, and~\ref{property:shiftinvariance}, then  \begin{align*}
    \lim_{n\to\infty}\mathrm{S}_t^{\kappa_n*}\dP_{|\Lambda_n}=\mathrm{S}_t^{\kappa*}\dP\quad \text{ in } \PPP_1(\Gamma),
\end{align*} 
for all  $\dP\in\PPP_s(\Gamma)$ and $ t\geq0$. The proof is a consequence of the following finite speed of propagation result. 
\begin{lem}[{\cite[Lemma 4.4, Lemma 4.5]{steenbeck2025reversiblebirthanddeathdynamicscontinuum}, \cite[Lemma 6.5]{jahnel2026reversiblebirthanddeathdynamicscontinuum}}]\label{lem:finitespeedofpropagation}
  Assume that ${\kappa:\Gamma\times\R^d\to[0,\infty)}$ satisfies \ref{property:cocycle}, \ref{property:locallystable}, \ref{property:finiterange}, and~\ref{property:shiftinvariance}. For any bounded local function $F\in\mathcal{H}$, it holds that \begin{align*}
       \lim_{n\to\infty} \int F(\xi)\mathrm{S}_t^{\kappa_n*}\dP_{|\Lambda_n}(\mathrm{d}\xi)=\int F(\xi)\mathrm{S}_t^{\kappa*}\dP(\mathrm{d}\xi),
    \end{align*}
    for all $\dP\in\PPP_s(\Gamma)$ and $t\geq0$.
\end{lem}
\begin{proof}
Note that the lemma is stated under stronger assumptions in \cite{steenbeck2025reversiblebirthanddeathdynamicscontinuum} and \cite{jahnel2026reversiblebirthanddeathdynamicscontinuum}. However, its proof only uses \ref{property:cocycle}, \ref{property:locallystable}, \ref{property:finiterange}, and~\ref{property:shiftinvariance}.
\end{proof}
\begin{proof}[Proof of Theorem~\ref{thm:convergenceresultsstationary} $(i)$]
Let $t\geq0$ and $\dP\in\PPP_s(\Gamma)$.
By Proposition~\ref{prop:convergencep1} it suffices to show that \linebreak$\lim_{n\to\infty}L_{\mathrm{S}_t^{\kappa_n*}\dP_{|\Lambda_n}}(f)=L_{\mathrm{S}_t^{\kappa*}\dP}(f)$ for all $f\in\CCC^+_c(\R^d)$ and $\lim_{n\to\infty}\I_{\mathrm{S}_t^{\kappa_n*}\dP_{|\Lambda_n}}=\I_{\mathrm{S}_t^{\kappa*}\dP}$ vaguely.  Lemma~\ref{lem:finitespeedofpropagation} directly implies that $\lim_{n\to\infty}L_{\mathrm{S}_t^{\kappa_n*}\dP_{|\Lambda_n}}(f)=L_{\mathrm{S}_t^{\kappa*}\dP}(f)$ for all $f\in\CCC^+_c(\R^d)$. Hence, it remains to show that the intensity measures converge vaguely. To this end, let $B\in\BBB_0(\R^d)$. It holds that 
\begin{align}
    |\int \xi(B)\mathrm{S}_t^{\kappa_n*}\dP_{|\Lambda_n}(\mathrm{d}\xi)-\int \xi(B)\mathrm{S}_t^{\kappa*}\dP(\mathrm{d}\xi)|\nonumber&\leq|\int \xi(B)\mathrm{S}_t^{\kappa_n*}\dP_{|\Lambda_n}(\mathrm{d}\xi)-\int (\xi(B)\wedge N)\mathrm{S}_t^{\kappa_n*}\dP_{|\Lambda_n}(\mathrm{d}\xi)|\nonumber\\
    &\quad+|\int (\xi(B)\wedge N)\mathrm{S}_t^{\kappa_n*}\dP_{|\Lambda_n}(\mathrm{d}\xi)-\int (\xi(B)\wedge N)\mathrm{S}_t^{\kappa*}\dP(\mathrm{d}\xi)|\nonumber\\
    &\quad+|\int (\xi(B)\wedge N)\mathrm{S}_t^{\kappa*}\dP(\mathrm{d}\xi)-\int \xi(B)\mathrm{S}_t^{\kappa*}\dP(\mathrm{d}\xi)|.\label{eq:413}
\end{align}
   Since $\kappa$ satisfies \ref{property:locallystable}, Theorem~\ref{thm:SDESemigroup} directly yields that there exists some random variable $Z_t^B$ that is Poisson distributed with parameter $A(1-e^{-t})\Leb(B)$ s.t.\   
   \begin{align*}
       \mathrm{S}_t^{\kappa_n}(\xi(B)-N)_+=\E[(Y_t^{\kappa_n,\Leb_{|\Lambda_n},(\xi)}(B)-N)_+]\leq\E[(\xi(B)+Z_t^B-N)_+]. 
   \end{align*}
   Together with $\dP\in\PPP_s(\Gamma)\subset\PPP_1(\Gamma)$, we obtain for the first term of \eqref{eq:413} that
    \begin{align*}
        \lim_{N\to\infty}\sup_{n\in\N}|\int \xi(B)\mathrm{S}_t^{\kappa_n*}\dP_{|\Lambda_n}(\mathrm{d}\xi)-\int (\xi(B)\wedge N)\mathrm{S}_t^{\kappa_n*}\dP_{|\Lambda_n}(\mathrm{d}\xi)|
       &=\lim_{N\to\infty}\sup_{n\in\N}\int\mathrm{S}_t^{\kappa_n}(\xi(B)-N)_+\dP_{|\Lambda_n}(\mathrm{d}\xi)\\
        &\leq \lim_{N\to\infty} \int \E[(\xi(B)+Z_t^B-N)_+] \dP(\mathrm{d}\xi)\\
        &=0.
    \end{align*}
    Similarly, for the last term in \eqref{eq:413} it holds that
    \begin{align*}
       \lim_{N\to\infty} |\int (\xi(B)\wedge N)\mathrm{S}_t^{\kappa*}\dP(\mathrm{d}\xi)-\int \xi(B)\mathrm{S}_t^{\kappa*}\dP(\mathrm{d}\xi)|=0.
    \end{align*}
    Finally, for the second term of \eqref{eq:413}, we obtain by Lemma~\ref{lem:finitespeedofpropagation} 
    \begin{align*}
       \lim_{n\to\infty} |\int (\xi(B)\wedge N)\mathrm{S}_t^{\kappa_n*}\dP_{|\Lambda_n}(\mathrm{d}\xi)-\int (\xi(B)\wedge N)\mathrm{S}_t^{\kappa*}\dP(\mathrm{d}\xi)|=0,
    \end{align*}
    finishing the proof.
\end{proof}
\subsection{Specific Relative Entropy}
In this subsection, we prove Theorem~\ref{thm:convergenceresultsstationary} $(ii)$. More precisely, let $\kappa:\Gamma\times\R^d\to[0,\infty)$ satisfy \ref{property:cocycle}, \ref{property:locallystable}, \ref{property:repulsive}, \ref{property:shiftinvariance}, and~\ref{property:finiterange}, and consider $\dR\in\RRR_{\kappa,s}(\Gamma)$, $\dR^{\Lambda_n}=\RRR_{\kappa_n,\Leb_{|\Lambda_n}}(\Gamma_{\Lambda_n}),\;n\in\N$. We show that \begin{align*}
    \Gamma\text{-}\lim_{n\to\infty}\frac{1}{\Leb(\Lambda_n)}\Ent(\cdot|\dR^{\Lambda_n})=\mathcal{E}(\cdot|\dR)
\end{align*}
 on $\PPP_s(\Gamma)\subset\PPP_1(\Gamma)$. A recovery sequence for $\dP\in\PPP_s(\Gamma)$ is  $(\dP_{|\Lambda_n})_{n\in\N}$.\\

We divide the proof into two parts and prove the limsup inequality for the $\Gamma$-convergence in Lemma~\ref{lem:limsupentropiestat} and the liminf inequality in Lemma~\ref{lem:8}. Combining these two lemmas directly yields the proof of Theorem~\ref{thm:convergenceresultsstationary} $(ii)$.

\begin{lem}\label{lem:limsupentropiestat}
   Assume that $\kappa:\Gamma\times\R^d\to[0,\infty)$ satisfies \ref{property:cocycle}, \ref{property:locallystable}, \ref{property:repulsive}, \ref{property:shiftinvariance}, and~\ref{property:finiterange}. For $\dR\in\RRR_{\kappa,s}(\Gamma)$, $\dR^{\Lambda_n}=\RRR_{\kappa_n,\Leb_{|\Lambda_n}}(\Gamma_{\Lambda_n}),\;n\in\N$, and $\dP\in\PPP_s(\Gamma)$, it holds that
   \begin{align*}
       \limsup_{n\to\infty} \frac{1}{\Leb(\Lambda_n)}\Ent(\dP_{|\Lambda_n}|\dR^{\Lambda_n})\leq \mathcal{E}(\dP|\dR).
   \end{align*}
\end{lem}

\begin{proof}
If $\mathcal{E}(\dP|\dR)=\infty$, the claim is immediate. Hence, let $\mathcal{E}(\dP|\dR)<\infty$. Then, by the definition of the specific relative entropy, it follows that $\dP_{|\Lambda_n}\ll\dR_{|\Lambda_n}$ for $n\in\N$ large enough. Furthermore, by Lemma~\ref{lem:propertiespapangelou}, we have $\dR_{|\Lambda_n}\ll\dR^{\Lambda_n}$ and \begin{align}
      f_n(\xi)=  \frac{\mathrm{d}\dR_{|\Lambda_n}}{\mathrm{d}\dR^{\Lambda_n}}(\xi)=\int\frac{Z_{\kappa_n,\Leb_{|\Lambda_n}}^{\Lambda_n}}{Z_{\kappa_n,\Leb_{|\Lambda_n}}^{\Lambda_n}(\eta_{|\Lambda_n^C})}\frac{\varphi^{\eta_{|\Lambda_n^C}}_{\kappa_n}(\xi)}{\varphi_{\kappa_n}(\xi)}\dR(\mathrm{d}\eta).\label{eq:1605}
    \end{align}
Therefore, we obtain
\begin{align*}
\Ent(\dP_{|\Lambda_n}|\dR^{\Lambda_n})&=\int\frac{\mathrm{d}\dP_{|\Lambda_n}}{\mathrm{d}\dR_{|\Lambda_n}}(\xi)f_n(\xi)\log\left(\frac{\mathrm{d}\dP_{|\Lambda_n}}{\mathrm{d}\dR_{|\Lambda_n}}(\xi)f_n(\xi)\right)\dR^{\Lambda_n}(\mathrm{d}\xi)\\
    &=\Ent(\dP_{|\Lambda_n}|\dR_{|\Lambda_n})+\int \log (f_n(\xi))\dP_{|\Lambda_n}(\mathrm{d}\xi).
\end{align*}
Taking the normalized large-volume limit and using \eqref{eq:bounddensity} of Lemma~\ref{lem:propertiespapangelou} yields
\begin{align*}
    \limsup_{n\to\infty} \frac{1}{\Leb(\Lambda_n)}\Ent(\dP_{|\Lambda_n}|\dR^{\Lambda_n})&=\mathcal{E}(\dP|\dR)+\limsup_{n\to\infty} \frac{1}{\Leb(\Lambda_n)}\int\log (f_n(\xi))\dP_{|\Lambda_n}(\mathrm{d}\xi)\\
    &\leq  \mathcal{E}(\dP|\dR)+\limsup_{n\to\infty}\frac{\Leb(\Lambda_n\setminus\Lambda_n^{r-})}{\Leb(\Lambda_n)}\\
    &= \mathcal{E}(\dP|\dR),
\end{align*}
finishing the proof.
\end{proof}

\begin{lem}\label{lem:8}
    Assume that $\kappa:\Gamma\times\R^d\to[0,\infty)$ satisfies \ref{property:cocycle}, \ref{property:locallystable}, \ref{property:repulsive}, \ref{property:shiftinvariance}, and~\ref{property:finiterange}. Let $\dR\in\RRR_{\kappa,s}(\Gamma)$. For $\dP_n\in\PPP_1(\Gamma_{\Lambda_n}),\;n\in\N,$ and $\dP\in\PPP_s(\Gamma)$ s.t.\  $\lim_{n\to\infty}\dP_n=\dP$ in $\PPP_1(\Gamma)$, it holds that \begin{align*}
    \mathcal{E}(\dP|\dR)\leq\liminf_{n\to\infty}\frac{1}{\Leb(\Lambda_n)}\Ent(\dP_n|\dR^{\Lambda_n}).
\end{align*}
\end{lem}

\begin{proof}
    Assume that $\liminf_{n\to\infty}\frac{1}{\Leb(\Lambda_n)}\Ent(\dP_{n}|\dR^{\Lambda_n})<\infty$ since otherwise the claim is immediate.
    Recall that $0\leq r<\infty$ is the range of the interaction of $\kappa$ and consider the tiled and stationarized measures\begin{align*}
        (\dP_{n|\Lambda_n^{r-}})^\mathrm{til}=\bigoplus_{z\in\Z^d}(\Theta_{nz})_\#\dP_{n|\Lambda_n^{r-}} \quad\text{ and }\quad (\dP_{n|\Lambda_n^{r-}})^\mathrm{stat}=\fint_{\Lambda_n}(\Theta_u)_\#(\dP_{n|\Lambda_n^{r-}})^\mathrm{til}\Leb(\mathrm{d}u).
    \end{align*}
    Note that we intentionally leave a strip around each tile empty. This prevents interactions between neighboring tiles. Analogously to the proof of Remark~\ref{rem:limstat}, we  obtain that $$\lim_{n\to\infty}(\dP_{n|\Lambda_n^{r-}})^\mathrm{stat}=\dP \quad\text{in } \PPP_1(\Gamma).$$ Thus, by the lower semicontinuity of the specific relative entropy (Lemma~\ref{lem:entropy} $(ii)$), it follows that \begin{align*}
        \mathcal{E}(\dP|\dR)\leq\liminf_{n\to\infty}\mathcal{E}((\dP_{n|\Lambda_n^{r-}})^\mathrm{stat}|\dR).
    \end{align*}
    By Jensen's inequality, the contraction property of the relative entropy (Lemma~\ref{lem:propertiesent} $(ii)$), and the stationarity of $\dR$, we obtain
    \begin{align}
        \mathcal{E}((\dP_{n|\Lambda_n^{r-}})^\mathrm{stat}|\dR)\nonumber&=\lim_{m\to\infty}\frac{1}{\Leb(\Lambda_{(2m-1)n})}\Ent\left(\left(\fint_{\Lambda_n}(\Theta_u)_\#(\dP_{n|\Lambda_n^{r-}})^\mathrm{til}\Leb(\mathrm{d}u)\right)_{|\Lambda_{(2m-1)n}}|\dR_{|\Lambda_{(2m-1)n}}\right)\nonumber\\
        & =\lim_{m\to\infty}\frac{1}{\Leb(\Lambda_{(2m-1)n})}\Ent\left(\fint_{\Lambda_n}((\Theta_u)_\#(\dP_{n|\Lambda_n^{r-}})^\mathrm{til})_{|\Lambda_{(2m-1)n}}\Leb(\mathrm{d}u)|\dR_{|\Lambda_{(2m-1)n}}\right)\nonumber\\
        &\leq \lim_{m\to\infty}\frac{1}{\Leb(\Lambda_{(2m-1)n})}\fint_{\Lambda_n}\Ent(((\Theta_u)_\#(\dP_{n|\Lambda_n^{r-}})^\mathrm{til})_{|\Lambda_{(2m-1)n}}|\dR_{|\Lambda_{(2m-1)n}})\Leb(\mathrm{d}u)\nonumber\\
        &\leq \lim_{m\to\infty}\frac{1}{\Leb(\Lambda_{(2m-1)n})}\fint_{\Lambda_n}\Ent(((\Theta_u)_\#(\dP_{n|\Lambda_n^{r-}})^\mathrm{til})_{|\Lambda_{(2m+1)n}-u}|\dR_{|\Lambda_{(2m+1)n}-u})\Leb(\mathrm{d}u)\nonumber\\
        &= \lim_{m\to\infty}\frac{1}{\Leb(\Lambda_{(2m-1)n})}\Ent((\dP_{n|\Lambda_n^{r-}})^\mathrm{til}_{|\Lambda_{(2m+1)n}}|\dR_{|\Lambda_{(2m+1)n}}).\label{eq:137}
    \end{align}
Furthermore, we have
\begin{align}
&\Ent((\dP_{n|\Lambda_n^{r-}})^\mathrm{til}_{|\Lambda_{(2m+1)n}}|\dR_{|\Lambda_{(2m+1)n}})\nonumber\\
    &= \Ent((\dP_{n|\Lambda_n^{r-}})^\mathrm{til}_{|\Lambda_{(2m+1)n}}|(\dR^{\Lambda_n^{r-}})^\mathrm{til}_{|\Lambda_{(2m+1)n}})+\int\log\left(\frac{\mathrm{d}(\dR^{\Lambda_n^{r-}})^\mathrm{til}_{|\Lambda_{(2m+1)n}}}{\mathrm{d}\dR_{|\Lambda_{(2m+1)n}}}(\xi)\right)(\dP_{n|\Lambda_n^{r-}})^\mathrm{til}_{|\Lambda_{(2m+1)n}}(\mathrm{d}\xi).\label{eq:138}
\end{align}
By the construction of the tiled process, it holds that
\begin{align}
    \Ent((\dP_{n|\Lambda_n^{r-}})^\mathrm{til}_{|\Lambda_{(2m+1)n}}|(\dR^{\Lambda_n^{r-}})^\mathrm{til}_{|\Lambda_{(2m+1)n}})=(2m+1)^d\Ent(\dP_{n|\Lambda_n^{r-}}|\dR^{\Lambda_n^{r-}}).\label{eq:140}
\end{align}
For the second term, we note that by \ref{property:finiterange} $(\dR^{\Lambda_n^{r-}})^\mathrm{til}_{|\Lambda_{(2m+1)n}}=\dR^{\Lambda_{(2m+1)n}}(\cdot|E)$ with \begin{align*}
    E:=\left\{\xi\in\Gamma_{\Lambda_{(2m+1)n}}:\xi\bigg(\Lambda_{(2m+1)n}\setminus\big(\bigcup_{z\in\Z^d:\|z\|_\infty\leq m}\Lambda_n^{r-}+nz\big)\bigg)=0\right\},
\end{align*}
where $\|\cdot\|_\infty$ is the infinity norm on $\Z^d$. Since $(\dP_{n|\Lambda_n^{r-}})^\mathrm{til}_{|\Lambda_{(2m+1)n}}(E)=1$, we obtain
\begin{align}
  &\int\log\left(\frac{\mathrm{d}(\dR^{\Lambda_n^{r-}})^\mathrm{til}_{|\Lambda_{(2m+1)n}}}{\mathrm{d}\dR_{|\Lambda_{(2m+1)n}}}(\xi)\right)(\dP_{n|\Lambda_n^{r-}})^\mathrm{til}_{|\Lambda_{(2m+1)n}}(\mathrm{d}\xi)\nonumber\\&\quad=  \int\log\left(\frac{\mathrm{d}\dR^{\Lambda_{(2m+1)n}}(\cdot|E)}{\mathrm{d}\dR^{\Lambda_{(2m+1)n}}}(\xi)\frac{\mathrm{d}\dR^{\Lambda_{(2m+1)n}}}{\dR_{|\Lambda_{(2m+1)n}}}(\xi)\right)(\dP_{n|\Lambda_n^{r-}})^\mathrm{til}_{|\Lambda_{(2m+1)n}}(\mathrm{d}\xi)\nonumber\\&\quad=\int\left(\log\left(\frac{\ind_E(\xi)}{\dR^{\Lambda_{(2m+1)n}}(E)}\right)+\log\left(\frac{\mathrm{d}\dR^{\Lambda_{(2m+1)n}}}{\mathrm{d}\dR_{|\Lambda_{(2m+1)n}}}(\xi)\right)\right)(\dP_{n|\Lambda_n^{r-}})^\mathrm{til}_{|\Lambda_{(2m+1)n}}(\mathrm{d}\xi)\nonumber\\
    &\quad=-\log(\dR^{\Lambda_{(2m+1)n}}(E))+\int\log\left(\frac{\mathrm{d}\dR^{\Lambda_{(2m+1)n}}}{\mathrm{d}\dR_{|\Lambda_{(2m+1)n}}}(\xi)\right)(\dP_{n|\Lambda_n^{r-}})^\mathrm{til}_{|\Lambda_{(2m+1)n}}(\mathrm{d}\xi).\label{eq:141}
\end{align}
Using \eqref{eq:bounddensity} from Lemma~\ref{lem:propertiespapangelou}, we deduce that
\begin{align}\label{eq:149}
\int\log\left(\frac{\mathrm{d}\dR^{\Lambda_{(2m+1)n}}}{\mathrm{d}\dR_{|\Lambda_{(2m+1)n}}}(\xi)\right)(\dP_{n|\Lambda_n^{r-}})^\mathrm{til}_{|\Lambda_{(2m+1)n}}(\mathrm{d}\xi)\leq A\Leb(\Lambda_{(2m+1)n}^{r+}\setminus\Lambda_{(2m+1)n}).
\end{align}
Moreover, analogously to the computation in \eqref{eq:40}, we obtain
\begin{align}\label{eq:142}
    -\log(\dR^{\Lambda_{(2m+1)n}}(E))\leq -\log\left(e^{-A\Leb(\Lambda_{(2m+1)n}\setminus\bigcup_{z\in\Z^d:\|z\|_\infty\leq m}(\Lambda_n^{r-}+nz)}\right)=A(2m+1)^d\Leb(\Lambda_n\setminus\Lambda_n^{r-}).
\end{align}
Combining \eqref{eq:137}, \eqref{eq:138}, \eqref{eq:140}, \eqref{eq:149}, and \eqref{eq:142}, we conclude that
\begin{align*}
    \mathcal{E}(\dP|\dR)
    &\leq\liminf_{n\to\infty}\lim_{m\to\infty}\frac{(2m+1)^d}{(2m-1)^d}\left(\frac{1}{\Leb(\Lambda_n)}\Ent(\dP_{n|\Lambda_n^{r-}}|\dR^{\Lambda_n^{r-}})+\frac{A\Leb(\Lambda_n\setminus\Lambda_n^{r-})}{\Leb(\Lambda_n)}\right)\\
    &\quad+\frac{A\Leb(\Lambda_{(2m+1)n}^{r+}\setminus\Lambda_{(2m+1)n})}{\Leb(\Lambda_{(2m-1)n})}\nonumber\\
    &=\liminf_{n\to\infty}\frac{1}{\Leb(\Lambda_n)}\Ent(\dP_{n|\Lambda_n^{r-}}|\dR^{\Lambda_n^{r-}}).
\end{align*}
Analogously to the proof of Lemma~\ref{lem:limsupentropiestat}, we obtain that
\begin{align*}
    \liminf_{n\to\infty}\frac{1}{\Leb(\Lambda_n)}\Ent(\dP_{n|\Lambda_n^{r-}}|\dR^{\Lambda_n^{r-}})\leq\liminf_{n\to\infty}\frac{1}{\Leb(\Lambda_n)}\Ent(\dP_{n|\Lambda_n^{r-}}|\dR_{|\Lambda_n^{r-}}).
\end{align*}
By the contraction property of the relative entropy (Lemma~\ref{lem:propertiesent} $(ii)$), $\Ent(\dP_{n|\Lambda_n^{r-}}|\dR_{|\Lambda_n^{r-}})\leq \Ent(\dP_{n}|\dR_{|\Lambda_n})$. Hence, following the same line of argumentation as in Lemma~\ref{lem:limsupentropiestat}, we finally obtain from \eqref{eq:bounddensity} that
\begin{align*}
    \mathcal{E}(\dP|\dR)
    &\leq \liminf_{n\to\infty}\frac{1}{\Leb(\Lambda_n)}\Ent(\dP_{n}|\dR_{|\Lambda_n})\\
    &=\liminf_{n\to\infty}\frac{1}{\Leb(\Lambda_n)}\Ent(\dP_{n}|\dR^{\Lambda_n})+ \liminf_{n\to\infty}\frac{1}{\Leb(\Lambda_n)}\int\log\left(\frac{\mathrm{d}\dR^{\Lambda_n}}{\mathrm{d}\dR_{|\Lambda_n}}(\xi)\right)\dP_{n}(\mathrm{d}\xi)\\
    &\leq \liminf_{n\to\infty}\frac{1}{\Leb(\Lambda_n)}\Ent(\dP_{n}|\dR^{\Lambda_n})+ \liminf_{n\to\infty}\frac{A\Leb(\Lambda_n^{r+}\setminus\Lambda_n)}{\Leb(\Lambda_n)}\\
    &=\liminf_{n\to\infty}\frac{1}{\Leb(\Lambda_n)}\Ent(\dP_{n}|\dR^{\Lambda_n}),
\end{align*}
finishing the proof.
\end{proof}

     \subsection{Specific Transport Distance}\label{subsec:convergencestationary}
     In this subsection, we prove Theorem~\ref{thm:convergenceresultsstationary} $(iii)$. More precisely, assume that $\kappa:\Gamma\times\R^d\to[0,\infty)$ satisfies \ref{property:cocycle}, \ref{property:locallystable}, \ref{property:continuity},\ref{property:repulsive}, \ref{property:shiftinvariance}, and~\ref{property:finiterange}. We show that 
     \begin{align*}
          \Gamma\text{-}\lim_{n\to\infty}\frac{1}{\Leb(\Lambda_n)}
        \mathcal{W}_{\kappa_n}^2
        (
            \cdot, \cdot
        )=\mathcal{W}_{\kappa,s}^2(\cdot,\cdot)
     \end{align*}
        on $\PPP_s(\Gamma)\times\PPP_s(\Gamma)\subset\PPP_1(\Gamma)\times\PPP_1(\Gamma)$.  A recovery sequence for $(\dP_0,\dP_1)\in\PPP_s(\Gamma)\times\PPP_s(\Gamma)$ is $(\dP_{0|\Lambda_n},\dP_{1|\Lambda_n})_{n\in\N}$.\\
        
        We divide the proof into two parts. The limsup inequality for the $\Gamma$-convergence is proven in Lemma~\ref{lem:boundcands} and the liminf inequality in Lemma~\ref{lem:boundstat}. Combining these two lemmas directly yields the proof of Theorem~\ref{thm:convergenceresultsstationary} $(iii)$. 
     \begin{lem}\label{lem:boundcands}
          Assume that $\kappa:\Gamma\times\R^d\to[0,\infty)$ satisfies \ref{property:locallystable} and \ref{property:shiftinvariance}. Let $\dP_0,\dP_1\in\PPP_s(\Gamma).$ Then
            \begin{align*}
       \frac{1}{\Leb(\Lambda_n)}\mathcal{W}_{\kappa_n}^2({\dP_0}_{|\Lambda_n},{\dP_1}_{|\Lambda_n}) \leq \mathcal{W}_{\kappa,s}^2(\dP_0,\dP_1),
            \end{align*}
            for all $n\in\N$.
        \end{lem}
        \begin{proof}
        Let $n\in\N.$ As the restrictions of solutions to the continuity equation are again solutions on the restricted domain (Theorem~\ref{thm:transformationofce} $(iii)$), we directly obtain
        \begin{align*}
            \mathcal{W}_{\kappa,s}^2(\dP_0,\dP_1)&=\inf\left\{\sup_{m\in\N}\frac{1}{\Leb(\Lambda_m)}\int_0^1\mathcal{L}_{\kappa_m}({\bar{\dP}_{t|\Lambda_m}},{\bar{\dV}_{t|\Lambda_m}})\mathrm{d}t:({\bar{\dP}},{\bar{\dV}})\in\CCC\EEE^s(\dP_0,\dP_1)\right\}\\
            &\geq \inf\left\{\frac{1}{\Leb(\Lambda_n)}\int_0^1\mathcal{L}_{\kappa_n}({\bar{\dP}_t^n},{\bar{\dV}_t^n})\mathrm{d}t:({\bar{\dP}}^n,{\bar{\dV}}^n)\in\CCC\EEE^{\Lambda_n}({\dP_0}_{|\Lambda_n},{\dP_1}_{|\Lambda_n})\right\}\\
            &= \frac{1}{\Leb(\Lambda_n)}\mathcal{W}_{\kappa_n}^2({\dP_0}_{|\Lambda_n},{\dP_1}_{|\Lambda_n}),
        \end{align*}
        finishing the proof.
        \end{proof}
       In order to prove the liminf inequality for the $\Gamma$-convergence, we select a minimizer of the optimization problem defining $\mathcal{W}_{\kappa_n}$ and construct by tiling and stationarizing a competitor for the optimization problem in the definition of $\mathcal{W}_{\kappa,s}$. In a next step, we then aim to  exploit the particular structure of tiled and stationarized measures.  The main technical result is the following lemma, which shows that the Lagrange functional is additive for  disjoint base spaces that are separated by an additional buffer zone. 
    \begin{lem}\label{lem:additivityLagrangefunctional}
Assume that $\kappa:\Gamma\times\R^d\to[0,\infty)$ satisfies
\ref{property:locallystable}, \ref{property:shiftinvariance},
and~\ref{property:finiterange}.
Let $m\in\N$, $I:=\{-m,\ldots,m\}^d$, and set $z_i:=ni$ for $i\in I$. 
Let $\dP^0\in\PPP_1(\Gamma_{\Lambda_n^{r-}})$ and
$\dV^0\in\MMM_{b,0}(\Gamma_{\Lambda_n^{r-}}\times \Lambda_n^{r-})$
and define $\dP^{i}:=(\Theta_{z_i})_\#\dP^{0}$ and
$\dV^{i}:=(\Theta_{z_i}^{\mathrm{prod}})_\#\dV^{0}$ for all
$i\in I\setminus\{0\}$. Finally, set
${\dP}:=\bigoplus_{i\in I}{\dP^i}$ and
${\dV}:=\bigoplus^{\mathrm{prod}}_{i\in I}(\dV^i,\dP^i)$, see
\eqref{eq:superpositionvelocity} for the definition of the superposition.
It holds that
\begin{align*}
    \mathcal{L}_{\kappa_{(2m+1)n}}({\dP}, {\dV})
    &= (2m+1)^d\mathcal{L}_{\kappa_{n}^{r-}}({\dP}^{0},{\dV}^{0}).
\end{align*}
\end{lem}

\begin{proof}
For $i\in I$, let $\sigma_i=\dP^i\otimes\Leb_{|(\Lambda_n-z_i)}+C_{\dP^i}+|\dV^i|$ and denote
$\dP^i\otimes\Leb_{|(\Lambda_n-z_i)}=\rho_i\sigma_i$,
$C_{\dP^i}=\phi_i\sigma_i$, and $\dV^i=w_i\sigma_i$. Let
$K\in\BBB(\Gamma_{\Lambda_{(2m+1)n}})$ and
$B\in\BBB(\Lambda_{(2m+1)n})$. It holds that
\begin{align*}
   & \int \int\ind_{K\times B}(\xi,x)\kappa(\xi,x)\dP(\mathrm{d}\xi)\Leb_{|\Lambda_{(2m+1)n}}(\mathrm{d}x)\\
   &\quad=\sum_{i\in I} \int \int
   \ind_{K\times B}\left(\sum_{j\in I}\xi_j,x\right)
   \kappa\left(\sum_{j\in I}\xi_j,x\right)
   \bigotimes_{j\in I}\dP^j(\mathrm{d}\xi_j)
   \Leb_{|(\Lambda_{n}-z_i)}(\mathrm{d}x)\\
   &\quad= \sum_{i\in I} \int \int
   \ind_{K\times B}\left(\sum_{j\in I}\xi_j,x\right)
   \kappa\left(\sum_{j\in I}\xi_j,x\right)\rho_i(\xi_i,x)
   \sigma_i(\mathrm{d}\xi_i,\mathrm{d}x)
   \bigotimes_{\substack{j\in I\\j\not=i}}\dP^j(\mathrm{d}\xi_j).
\end{align*}
Hence, we can deduce that
\[
\kappa(\xi,x)\dP\otimes\Leb_{|\Lambda_{(2m+1)n}}(\mathrm{d}\xi,\mathrm{d}x)
=\kappa(\xi,x)\sum_{i\in I}\rho_i(\xi_{|(\Lambda_n-z_i)},x)
\mathrm{sp}^\mathrm{prod}_\#(\sigma_i\otimes(\oplus_{j\in I\setminus\{i\}}\dP^j))
(\mathrm{d}\xi,\mathrm{d}x).
\]
Further,
\begin{align*}
   \int\ind_{K\times B}(\xi,x)C_\dP(\mathrm{d}\xi,\mathrm{d}x)
   &=\int\int
   \ind_{K\times B}\left(\sum_{j\in I}\xi_j-\delta_x,x\right)
   \sum_{j\in I}\xi_j(\mathrm{d}x)
   \bigotimes_{j\in I}\dP^j(\mathrm{d}\xi_j)\\
   &=\sum_{i\in I}\int
   \ind_{K\times B}\left(\sum_{j\in I}\xi_j,x\right)
   C_{\dP^i}(\mathrm{d}\xi_i,\mathrm{d}x)
   \bigotimes_{\substack{j\in I\\ j\not= i}}\dP^j(\mathrm{d}\xi_j)\\
   &=\sum_{i\in I}\int
   \ind_{K\times B}\left(\sum_{j\in I}\xi_j,x\right)
   \phi_i(\xi_i,x)\sigma_i(\mathrm{d}\xi_i,\mathrm{d}x)
   \bigotimes_{\substack{j\in I\\ j\not= i}}\dP^j(\mathrm{d}\xi_j).
\end{align*}
This implies that
\[
C_\dP(\mathrm{d}\xi,\mathrm{d}x)
=\sum_{i\in I}\phi_i(\xi_{|(\Lambda_n-z_i)},x)
\mathrm{sp}^\mathrm{prod}_\#(\sigma_i\otimes(\oplus_{j\in I\setminus\{i\}}\dP^j))
(\mathrm{d}\xi,\mathrm{d}x).
\]
By an analogous computation, we obtain that
\[
\dV(\mathrm{d}\xi,\mathrm{d}x)
=\sum_{i\in I}w_i(\xi_{|(\Lambda_n-z_i)},x)
\mathrm{sp}^\mathrm{prod}_\#(\sigma_i\otimes(\oplus_{j\in I\setminus\{i\}}\dP^j))
(\mathrm{d}\xi,\mathrm{d}x).
\]
Inserting these densities into the definition of the Lagrange functional yields
\begin{align}
   \mathcal{L}_{\kappa_{(2m+1)n}}(\dP,\dV)\nonumber&=\int \alpha\left(
   \frac{\mathrm{d}\kappa(\dP\otimes\Leb_{|\Lambda_{(2m+1)n}})}{\mathrm{d}\sigma}(\xi,x),
   \frac{\mathrm{d}C_\dP}{\mathrm{d}\sigma}(\xi,x),
   \frac{\mathrm{d}\dV}{\mathrm{d}\sigma}(\xi,x)
   \right)\sigma(\mathrm{d}\xi,\mathrm{d}x)\nonumber\\
   &=\sum_{i\in I}\int \alpha\left(
   \frac{\mathrm{d}\kappa(\dP\otimes\Leb_{|\Lambda_{(2m+1)n}})}{\mathrm{d}\sigma}(\xi,x),
   \frac{\mathrm{d}C_\dP}{\mathrm{d}\sigma}(\xi,x),
   \frac{\mathrm{d}\dV}{\mathrm{d}\sigma}(\xi,x)
   \right)\ind_{(\Lambda_n-z_i)}(x)\sigma(\mathrm{d}\xi,\mathrm{d}x)\nonumber\\
   &=\sum_{i\in I}\int\ldots\int \int
   \alpha\left(
   \kappa\left(\sum_{j\in I}\xi_j,x\right)\rho_i(\xi_i,x),
   \phi_i(\xi_i,x),
   w_i(\xi_i,x)
   \right)
   \sigma_i(\mathrm{d}\xi_i,\mathrm{d}x)
   \bigotimes_{\substack{j\in I\\j\not=i}}\dP^j(\mathrm{d}\xi_j).
   \label{eq:5200}
\end{align}
Since $\kappa$ satisfies the finite-range property
\ref{property:finiterange}, we obtain
\begin{align}
   &\int\ldots\int \int
   \alpha\left(
   \kappa\left(\sum_{j\in I}\xi_j,x\right)\rho_i(\xi_i,x),
   \phi_i(\xi_i,x),
   w_i(\xi_i,x)
   \right)
   \sigma_i(\mathrm{d}\xi_i,\mathrm{d}x)
   \bigotimes_{\substack{j\in I\\j\not=i}}\dP^j(\mathrm{d}\xi_j)\nonumber\\
   &\quad=\int\ldots\int \int
   \alpha\left(
   \kappa\left(\xi_i,x\right)\rho_i(\xi_i,x),
   \phi_i(\xi_i,x),
   w_i(\xi_i,x)
   \right)
   \sigma_i(\mathrm{d}\xi_i,\mathrm{d}x)
   \bigotimes_{\substack{j\in I\\j\not=i}}\dP^j(\mathrm{d}\xi_j)\nonumber\\
   &\quad=\mathcal{L}_{\kappa,\Leb_{|(\Lambda_{n}-z_i)}}(\dP^i,\dV^i).
   \label{eq:5201}
\end{align}
Combining \eqref{eq:5200} and \eqref{eq:5201} and using the shift-invariance
of the Lagrange functional (Lemma~\ref{lem:shiftinvariance}) yields
\[
   \mathcal{L}_{\kappa_{(2m+1)n}}(\dP,\dV)
   =(2m+1)^d\mathcal{L}_{\kappa_{n}}({\dP}^{0},{\dV}^{0}).
\]
Finally, by the definition of the Lagrange functional, it holds that
\[
\mathcal{L}_{\kappa_{n}}({\dP}^{0},{\dV}^{0})
=\mathcal{L}_{\kappa_{n}^{r-}}({\dP}^{0},{\dV}^{0}),
\]
finishing the proof.
\end{proof}
    \begin{bem}
    For $\kappa\equiv 1$, the interaction radius is zero. Hence, Lemma~\ref{lem:additivityLagrangefunctional} yields additivity of $\mathcal{L}_{1,\lambda}$ for disjoint domains. See also \cite[Lemma 4.4]{huesmann2025nonlocalwassersteingeometrygradient}.
\end{bem}
\begin{lem}\label{lem:boundstat}
   Assume that $\kappa:\Gamma\times\R^d\to[0,\infty)$ satisfies \ref{property:cocycle}, \ref{property:locallystable}, 
    \ref{property:continuity}, \ref{property:repulsive}, \ref{property:shiftinvariance}, and~\ref{property:finiterange}. For $\dP_0^n,\dP_1^n\in\PPP_1(\Gamma_{\Lambda_n}
    ),\;n\in\N,$ and $\dP_0,\dP_1\in\PPP_s(\Gamma)$ s.t.\  $\lim_{n\to\infty}\dP_0^n=\dP_0$ and $\lim_{n\to\infty}\dP_1^n=\dP_1$ in $\PPP_1(\Gamma)$, it holds that 
     \begin{align*}
        \mathcal{W}_{\kappa,s}^2(\dP_{0},\dP_1)\leq \liminf_{n\to\infty}\frac{1}{\Leb(\Lambda_n)}\mathcal{W}_{\kappa_n}^2(\dP_{0}^n,\dP_{1}^n).
     \end{align*}
\end{lem}
\begin{proof}
If 
$\liminf_{n\to\infty}\frac{1}{\Leb(\Lambda_n)}\mathcal{W}_{\kappa_n}^2(\dP_{0}^n,\dP_{1}^n)=\infty$,
the claim is immediate. Therefore, we assume that 
$\liminf_{n\to\infty}\frac{1}{\Leb(\Lambda_n)}\mathcal{W}_{\kappa_n}^2(\dP_{0}^n,\dP_{1}^n)<\infty$. Consider the tiled and stationarized measures
\begin{align*}
        (\dP_{i|\Lambda_n^{r-}}^n)^\mathrm{til}=\bigoplus_{z\in\Z^d}(\Theta_{nz})_\#\dP_{i|\Lambda_n^{r-}}^n \quad\text{ and }\quad (\dP_{i|\Lambda_n^{r-}}^n)^\mathrm{stat}=\fint_{\Lambda_n}(\Theta_u)_\#(\dP_{i|\Lambda_n^{r-}}^n)^\mathrm{til}\Leb(\mathrm{d}u),\quad i=0,1.
\end{align*}
Note that, similar to the proof of Lemma~\ref{lem:8}, we intentionally leave an additional strip around each tile empty to prevent interactions between neighboring tiles.
Analogously to the proof of Remark~\ref{rem:limstat}, we obtain that 
\begin{align*}
    \lim_{n\to\infty} (\dP_{i|\Lambda_n^{r-}}^n)^\mathrm{stat}=\dP_i\quad\text{ in }\PPP_1(\Gamma),\; i=0,1.
\end{align*}
Hence, by the lower semicontinuity of $\mathcal{W}_{\kappa,s}$ (Theorem~\ref{thm:propertiesdistancestationary} $(iii)$), it holds that 
\begin{align*}
    \mathcal{W}_{\kappa,s}^2(\dP_0,\dP_1)&\leq \liminf_{n\to\infty}\mathcal{W}_{\kappa,s}^2((\dP_{0|\Lambda_n^{r-}}^n)^\mathrm{stat},(\dP_{1|\Lambda_n^{r-}}^n)^\mathrm{stat}).
\end{align*}
By Theorem~\ref{thm:minimizer}, there exist solutions $({\bar{\dP}}^{n},{\bar{\dV}}^{n})\in\CCC\EEE^{\Lambda_n}(\dP_{0|\Lambda_n}^n,\dP_{1|\Lambda_n}^n)$ s.t.\   
\begin{align*}
    \mathcal{W}_{\kappa_n}^2(\dP_{0|\Lambda_n}^n,\dP_{1|\Lambda_n}^n)
    =\int_0^1\mathcal{L}_{\kappa_n}({\bar{\dP}_t^{n}}, {\bar{\dV}_t^{n}})\mathrm{d}t.
\end{align*}
Recall the definition of tiled and stationarized velocity fields from
Definition~\ref{dfn:tiledstationarized}. 
Lemma~\ref{lem:stationarysolutions} entails that $(({\bar{\dP}}^{n}_{|\Lambda_n^{r-}})^\mathrm{stat},({\bar{\dV}}^{n}_{|\Lambda_n^{r-}})^\mathrm{stat})\in\CCC\EEE^s((\dP_{0|\Lambda_n^{r-}}^n)^\mathrm{stat},(\dP_{1|\Lambda_n^{r-}}^n)^\mathrm{stat})$.
Combining these observations with Lemma~\ref{lem:changelimsup} yields
\begin{align}
\mathcal{W}_{\kappa,s}^2((\dP_{0|\Lambda_n^{r-}}^n)^\mathrm{stat},(\dP_{1|\Lambda_n^{r-}}^n)^\mathrm{stat})
&\leq  \lim_{m\to\infty}\frac{1}{\Leb(\Lambda_{(2m-1)n})}\int_0^1\!\!\mathcal{L}_{\kappa_{(2m-1)n}}((\bar{\dP}^{n}_{t|\Lambda_n^{r-}})^\mathrm{stat}_{|\Lambda_{(2m-1)n}}, (\bar{\dV}^{n}_{t|\Lambda_n^{r-}})^\mathrm{stat}_{|\Lambda_{(2m-1)n}})\mathrm{d}t.
\label{eq:204}
\end{align}
Using the convexity of the Lagrange functional  (Lemma~\ref{lem:propertiesLagrange} $(i)$) and Jensen's inequality, we obtain
\begin{align}
&\mathcal{L}_{\kappa_{(2m-1)n}}((\bar{\dP}^{n}_{t|\Lambda_n^{r-}})^\mathrm{stat}_{|\Lambda_{(2m-1)n}}, (\bar{\dV}^{n}_{t|\Lambda_n^{r-}})^\mathrm{stat}_{|\Lambda_{(2m-1)n}})\nonumber\\
&\quad\leq \fint_{\Lambda_n}\mathcal{L}_{\kappa_{(2m-1)n}}( {(({\Theta_u})_\#(\bar{\dP}^{n}_{t|\Lambda_n^{r-}})^\mathrm{til}})_{|\Lambda_{(2m-1)n}},(({\Theta_u^\mathrm{prod}})_\#(\bar{\dV}^{n}_{t|\Lambda_n^{r-}})^\mathrm{til})_{|\Lambda_{(2m-1)n}} )\Leb(\mathrm{d}u).
\label{eq:200}
\end{align}
In the following, we will make use of the specific form of $(({\Theta_u})_\#(\bar{\dP}_{t|\Lambda_n^{r-}}^{n})^\mathrm{til}, ({\Theta_u^\mathrm{prod}})_\#(\bar{\dV}_{t|\Lambda_n^{r-}}^{n})^\mathrm{til})$. Therefore, let $u\in\Lambda_n.$ Since $\Lambda_{(2m-1)n}\subset\Lambda_{(2m+1)n}-u$, we obtain by the super-additivity of the Lagrange functional (Lemma~\ref{lem:superadditivity}) that
\begin{align*}
&\mathcal{L}_{\kappa_{(2m-1)n}}( {(({\Theta_u})_\# (\bar{\dP}^{n}_{t|\Lambda_n^{r-}})^\mathrm{til}})_{|\Lambda_{(2m-1)n}},(({\Theta_u^\mathrm{prod}})_\#(\bar{\dV}^{n}_{t|\Lambda_n^{r-}})^\mathrm{til})_{|\Lambda_{(2m-1)n}} )\\
&\leq \mathcal{L}_{\kappa,\Leb_{|\Lambda_{(2m+1)n}-u}}( {({(\Theta_u)}_\#(\bar{\dP}^{n}_{t|\Lambda_n^{r-}})^\mathrm{til}})_{|\Lambda_{(2m+1)n}-u},({(\Theta_u^\mathrm{prod})}_\#(\bar{\dV}^{n}_{t|\Lambda_n^{r-}})^\mathrm{til})_{|\Lambda_{(2m+1)n}-u} )\\
&=\mathcal{L}_{\kappa_{(2m+1)n}}( ({\Theta_{-u}})_\#({(({\Theta_u})_\#(\bar{\dP}^{n}_{t|\Lambda_n^{r-}})^\mathrm{til}})_{|\Lambda_{(2m+1)n}-u}),{(\Theta_{-u}^\mathrm{prod})}_\#(({(\Theta_u^\mathrm{prod})}_\#(\bar{\dV}^{n}_{t|\Lambda_n^{r-}})^\mathrm{til})_{|\Lambda_{(2m+1)n}-u} ))\nonumber\\
&= \mathcal{L}_{\kappa_{(2m+1)n}}( {((\bar{\dP}^{n}_{t|\Lambda_n^{r-}})^\mathrm{til}})_{|\Lambda_{(2m+1)n}},((\bar{\dV}^{n}_{t|\Lambda_n^{r-}})^\mathrm{til})_{|\Lambda_{(2m+1)n}} ).
\end{align*}
Inserting this into equation~\eqref{eq:200} yields
\begin{align}
&\fint_{\Lambda_n}\mathcal{L}_{\kappa_{(2m-1)n}}( {(({\Theta_u})_\#(\bar{\dP}^{n}_{t|\Lambda_n^{r-}})^\mathrm{til}})_{|\Lambda_{(2m-1)n}},(({\Theta_u^\mathrm{prod}})_\#(\bar{\dV}^{n}_{t|\Lambda_n^{r-}})^\mathrm{til})_{|\Lambda_{(2m-1)n}} )\Leb(\mathrm{d}u)\nonumber\\
&\quad\leq\fint_{\Lambda_n} \mathcal{L}_{\kappa_{(2m+1)n}}({((\bar{\dP}^{n}_{t|\Lambda_n^{r-}})^\mathrm{til}})_{|\Lambda_{(2m+1)n}},((\bar{\dV}^{n}_{t|\Lambda_n^{r-}})^\mathrm{til})_{|\Lambda_{(2m+1)n}}) \Leb(\mathrm{d}u)\nonumber\\
&\quad=\mathcal{L}_{\kappa_{(2m+1)n}}({((\bar{\dP}^{n}_{t|\Lambda_n^{r-}})^\mathrm{til}})_{|\Lambda_{(2m+1)n}},((\bar{\dV}^{n}_{t|\Lambda_n^{r-}})^\mathrm{til})_{|\Lambda_{(2m+1)n}}).
\label{eq:5301}
\end{align}
By Lemma~\ref{lem:additivityLagrangefunctional} and the super-additivity of the Lagrange functional (Lemma~\ref{lem:superadditivity}), we have 
\begin{align}
\mathcal{L}_{\kappa_{(2m+1)n}}({((\bar{\dP}^{n}_{t|\Lambda_n^{r-}})^\mathrm{til}})_{|\Lambda_{(2m+1)n}},((\bar{\dV}^{n}_{t|\Lambda_n^{r-}})^\mathrm{til})_{|\Lambda_{(2m+1)n}})
&= (2m+1)^d \mathcal{L}_{\kappa_{n}}({\bar{\dP}^{n}_{t|\Lambda_n^{r-}}},\bar{\dV}^{n}_{t|\Lambda_n^{r-}})\nonumber\\
&\leq (2m+1)^d \mathcal{L}_{\kappa_{n}}({\bar{\dP}^{n}_{t},\bar{\dV}^{n}_{t}}).
\label{eq:5300}
\end{align}
Combining \eqref{eq:204}, \eqref{eq:200}, \eqref{eq:5301}, and \eqref{eq:5300}, we obtain  
\begin{align}
&\liminf_{n\to\infty}\mathcal{W}_{\kappa,s}^2((\dP_{0|\Lambda_n^{r-}}^n)^\mathrm{stat},(\dP_{1|\Lambda_n^{r-}}^n)^\mathrm{stat})\nonumber\\
&\quad\leq \liminf_{n\to\infty}\lim_{m\to\infty}\frac{1}{\Leb(\Lambda_{(2m-1)n})}\int_0^1\mathcal{L}_{\kappa_{(2m+1)n}}(((\bar{\dP}^{n}_{t|\Lambda_n^{r-}})^\mathrm{til})_{|\Lambda_{(2m+1)n}},((\bar{\dV}^{n}_{t|\Lambda_n^{r-}})^\mathrm{til})_{|\Lambda_{(2m+1)n}}) \mathrm{d}t\nonumber\\
&\quad\leq \liminf_{n\to\infty}\lim_{m\to\infty}\frac{(2m+1)^d}{(2m-1)^d\Leb(\Lambda_n)}\int_0^1 \mathcal{L}_{\kappa_{n}}({\bar{\dP}^{n}_t},\bar{\dV}^{n}_t)\mathrm{d}t\nonumber\\
&\quad\leq \liminf_{n\to\infty}\frac{1}{\Leb(\Lambda_n)} \mathcal{W}_{\kappa_n}^2(\dP_{0}^n,\dP_{1}^n),
\label{eq:44}
\end{align}
finishing the proof. 
\end{proof}

\subsection{Specific Modified Fisher Information}
In this subsection, we prove Theorem~\ref{thm:convergenceresultsstationary}~$(iv)$ and show that if $\kappa:\Gamma\times\R^d\to[0,\infty)$ satisfies \ref{property:cocycle}, \ref{property:locallystable}, \ref{property:repulsive}, \ref{property:finiterange}, and~\ref{property:shiftinvariance}, then for $\dR\in\RRR_{\kappa,s}(\Gamma),\;\dR^{\Lambda_n}=\RRR_{\kappa_n,\Leb_{|\Lambda_n}}(\Gamma_{\Lambda_n}),\;n\in\N,$ and $\dP\in\PPP_s(\Gamma)$, it holds that
\begin{align*}
    \lim_{n\to\infty}\frac{1}{\Leb(\Lambda_n)}
    \I(\dP_{|\Lambda_n}|\dR^{\Lambda_n})
    =
    \mathcal{I}(\dP|\dR).
\end{align*}

Before turning to the proof of Theorem~\ref{thm:convergenceresultsstationary}~$(iv)$, we first show that finiteness of the specific modified Fisher information $\mathcal{I}(\dP|\dR)<\infty$ implies that the restrictions $\dP_{|\Lambda_n}$ are themselves Papangelou measures. For this purpose, we recall that the restrictions $\dR_{|\Lambda_n}$ of infinite-volume Papangelou measures $\dR\in\RRR_{\kappa,s}(\Gamma)$ are themselves Papangelou measures with Papangelou intensity
\[
    \kappa_{|\Lambda_n}(\xi_1,x)
    :=
    \int
    \kappa(\xi_1+\xi_2,x)
    \dR^{\Lambda_n^C}(\mathrm d\xi_2| \xi_1),
    \quad
    \xi_1\in\Gamma_{\Lambda_n},\; x\in\Lambda_n,
\]
see Lemma~\ref{lem:restrictionpapangelou}.
\begin{lem}\label{lem:regularisationFisher}
    Assume that $\kappa:\Gamma\times\R^d\to[0,\infty)$ satisfies \ref{property:cocycle}, \ref{property:locallystable}, \ref{property:repulsive}, \ref{property:finiterange}, and~\ref{property:shiftinvariance}. Let $\dR\in\RRR_{\kappa,s}(\Gamma)$ and $\dP\in\PPP_s(\Gamma)$ s.t.\  $\mathcal{I}(\dP|\dR)<\infty$ and $\dP_{|\Lambda_n}=\rho_n\dR_{|\Lambda_n}$, $n\in\N$. Then $\dP_{|\Lambda_n}$ is a Papangelou measure with reference measure $\Leb_{|\Lambda_n}$ and Papangelou intensity
    \begin{align*}
        \kappa^\dP_n(\xi,x)
        :=
        \frac{\rho_n(\xi+\delta_x)}{\rho_n(\xi)}
        \kappa_{|\Lambda_n}(\xi,x),  \quad
    \xi\in\Gamma_{\Lambda_n},\; x\in\Lambda_n.
    \end{align*}
\end{lem}

\begin{proof}
    Let $F\in\FFF^+(\Gamma_{\Lambda_n}\times\Lambda_n)$. By Lemma~\ref{lem:propertiespapangelou}, we obtain
    \begin{align*}
        \int\int
        F(\xi-\delta_x,x)
        \xi(\mathrm{d}x)
        \dP_{|\Lambda_n}(\mathrm{d}\xi)
        &=
        \int\int
        F(\xi-\delta_x,x)
        \rho_n(\xi)
        \xi(\mathrm{d}x)
        \dR_{|\Lambda_n}(\mathrm{d}\xi)\\
        &=
        \int\int
        F(\xi,x)
        \rho_n(\xi+\delta_x)
        \kappa_{|\Lambda_n}(\xi,x)
        \dR_{|\Lambda_n}(\mathrm{d}\xi)
        \Leb_{|\Lambda_n}(\mathrm{d}x)\\
        &=
        \int\int
        F(\xi,x)
        \frac{\rho_n(\xi+\delta_x)}
        {\rho_n(\xi)}
        \kappa_{|\Lambda_n}(\xi,x)
        \dP_{|\Lambda_n}(\mathrm{d}\xi)
        \Leb_{|\Lambda_n}(\mathrm{d}x).
    \end{align*}
    Thus, $\dP_{|\Lambda_n}$ satisfies the GNZ equation~\eqref{eq:papangelou} with reference measure $\Leb_{|\Lambda_n}$ and Papangelou intensity
    \[
        \kappa^\dP_n(\xi,x)
        =
        \frac{\rho_n(\xi+\delta_x)}
        {\rho_n(\xi)}
        \kappa_{|\Lambda_n}(\xi,x).
    \]
    It remains to show that $\kappa^\dP_n(\xi,x)$ is well-defined. Since
    \begin{align*}
        &\int\int
        \left(\rho_n(\xi+\delta_x)-\rho_n(\xi)\right)
        \left(\log(\rho_n(\xi+\delta_x))-\log(\rho_n(\xi))\right)
        \kappa_{|\Lambda_n}(\xi,x)
        \dR_{|\Lambda_n}(\mathrm{d}\xi)
        \Leb_{|\Lambda_n}(\mathrm{d}x)
        <\infty,
    \end{align*}
    we obtain
    \begin{align*}
        &\kappa_{|\Lambda_n}(\xi,x)
        \dR_{|\Lambda_n}(\mathrm{d}\xi)
        \Leb_{|\Lambda_n}(\mathrm{d}x)
        \big(
        \{(\xi,x):\rho_n(\xi)=0,\ \rho_n(\xi+\delta_x)>0\}
        \big)
        =
        0,
    \end{align*}
    which finishes the proof.
\end{proof}

\begin{proof}[Proof of Theorem~\ref{thm:convergenceresultsstationary}~$(iv)$]
    Assume that $\mathcal{I}(\dP|\dR)<\infty$ and let $\dP_{|\Lambda_n}=\rho_n\dR_{|\Lambda_n},\;n\in\N$. By Lemma~\ref{lem:propertiespapangelou}, $\dR_{|\Lambda_n}\ll\dR^{\Lambda_n}$ and
    \begin{align*}
        f_n(\xi)
        =
        \frac{\mathrm{d}\dR_{|\Lambda_n}}{\mathrm{d}\dR^{\Lambda_n}}(\xi)
        =
        \int
        \frac{Z_{\kappa_n,\Leb_{|\Lambda_n}}^{\Lambda_n}}
        {Z_{\kappa_n,\Leb_{|\Lambda_n}}^{\Lambda_n}(\eta_{|\Lambda_n^C})}
        \frac{\varphi^{\eta_{|\Lambda_n^C}}_{\kappa_n}(\xi)}
        {\varphi_{\kappa_n}(\xi)}
        \dR(\mathrm{d}\eta),
        \quad \xi\in\Gamma_{\Lambda_n}.
    \end{align*}
    Hence, $\dP_{|\Lambda_n}=\rho_nf_n\dR^{\Lambda_n}$ and by the definition of the Fisher information it holds that
    \begin{align}
        &\left|
        \I(\dP_{|\Lambda_n}|\dR_{|\Lambda_n})
        -
        \I(\dP_{|\Lambda_n}|\dR^{\Lambda_n})
        \right| \nonumber\\
        &\quad=
        \left|
        \int\int
        D_x\rho_n(\xi)D_x\log(\rho_n)(\xi)
        \kappa_{|\Lambda_n}(\xi,x)
        \dR_{|\Lambda_n}(\mathrm{d}\xi)
        \Leb_{|\Lambda_n}(\mathrm{d}x) \right. \nonumber\\
        &\qquad \left.
        -
        \int\int
        D_x(\rho_nf_n)(\xi)D_x\log(\rho_nf_n)(\xi)
        \kappa(\xi,x)
        \dR^{\Lambda_n}(\mathrm{d}\xi)
        \Leb_{|\Lambda_n}(\mathrm{d}x)
        \right| \nonumber\\
        &\quad=
        \left|
        \int\int
        \left(
        \frac{\rho_n(\xi+\delta_x)}{\rho_n(\xi)}-1
        \right)
        \log\left(
        \frac{\rho_n(\xi+\delta_x)}{\rho_n(\xi)}
        \right)
        \kappa_{|\Lambda_n}(\xi,x)
        \dP_{|\Lambda_n}(\mathrm{d}\xi)
        \Leb_{|\Lambda_n}(\mathrm{d}x) \right. \nonumber\\
        &\qquad \left.
        -
        \int\int
        \left(
        \frac{\rho_n(\xi+\delta_x)f_n(\xi+\delta_x)}
        {\rho_n(\xi)f_n(\xi)}-1
        \right)
        \log\left(
        \frac{\rho_n(\xi+\delta_x)f_n(\xi+\delta_x)}
        {\rho_n(\xi)f_n(\xi)}
        \right)
        \kappa(\xi,x)
        \dP_{|\Lambda_n}(\mathrm{d}\xi)
        \Leb_{|\Lambda_n}(\mathrm{d}x)
        \right|.
        \label{eq:5678}
    \end{align}

    Since $\kappa$ satisfies \ref{property:finiterange}, it holds that
    \begin{align}
        f_n(\xi+\delta_x)
        &=
        \int
        \frac{Z_{\kappa_n,\Leb_{|\Lambda_n}}^{\Lambda_n}}
        {Z_{\kappa_n,\Leb_{|\Lambda_n}}^{\Lambda_n}(\eta_{|\Lambda_n^C})}
        \frac{\varphi^{\eta_{|\Lambda_n^C}}_{\kappa_n}(\xi+\delta_x)}
        {\varphi_{\kappa_n}(\xi+\delta_x)}
        \dR(\mathrm{d}\eta) =
        \int
        \frac{Z_{\kappa_n,\Leb_{|\Lambda_n}}^{\Lambda_n}}
        {Z_{\kappa_n,\Leb_{|\Lambda_n}}^{\Lambda_n}(\eta_{|\Lambda_n^C})}
        \frac{\varphi^{\eta_{|\Lambda_n^C}}_{\kappa_n}(\xi)}
        {\varphi_{\kappa_n}(\xi)}
        \dR(\mathrm{d}\eta) =
        f_n(\xi)
        \label{eq:ffiniterange}
    \end{align}
    and $\kappa_{|\Lambda_n}(\xi,x)=\kappa(\xi,x)$ for all $\xi\in\Gamma_{\Lambda_n}$ and $x\in\Lambda_n^{r-}$. Inserting these observations into~\eqref{eq:5678} yields
    \begin{align*}
        &\left|
        \I(\dP_{|\Lambda_n}|\dR_{|\Lambda_n})
        -
        \I(\dP_{|\Lambda_n}|\dR^{\Lambda_n})
        \right|\\
        &\quad=
        \left|
        \int_{\Lambda_n\setminus\Lambda_n^{r-}}\int
        \left(
        \frac{\rho_n(\xi+\delta_x)}{\rho_n(\xi)}-1
        \right)
        \log\left(
        \frac{\rho_n(\xi+\delta_x)}{\rho_n(\xi)}
        \right)
        \kappa_{|\Lambda_n}(\xi,x)
        \dP_{|\Lambda_n}(\mathrm{d}\xi)
        \Leb_{|\Lambda_n}(\mathrm{d}x) \right.\\
        & \qquad\left.
        -
        \int_{\Lambda_n\setminus\Lambda_n^{r-}}\int
        \left(
        \frac{\rho_n(\xi+\delta_x)f_n(\xi+\delta_x)}
        {\rho_n(\xi)f_n(\xi)}-1
        \right)
        \log\left(
        \frac{\rho_n(\xi+\delta_x)f_n(\xi+\delta_x)}
        {\rho_n(\xi)f_n(\xi)}
        \right)
        \kappa(\xi,x)
        \dP_{|\Lambda_n}(\mathrm{d}\xi)
        \Leb_{|\Lambda_n}(\mathrm{d}x)
        \right|.
    \end{align*}

    Analogously to the proof of~\eqref{eq:bounddensity} in Lemma~\ref{lem:propertiespapangelou}, we obtain by \ref{property:locallystable} that
    \begin{align*}
        \kappa_{|\Lambda_n}(\xi,x)
        &\geq
        \int
        \kappa(\xi,x)
        \ind_{\{\eta(B_r(x))=0\}}
        \dR^{\Lambda_n^C}(\mathrm{d}\eta|\xi)\geq
        \kappa(\xi,x)e^{-A\Leb(B_r(0))}.
    \end{align*}
    On the other hand,  $\kappa_{|\Lambda_n}(\xi,x)\leq\kappa(\xi,x)$ by repulsiveness~\ref{property:repulsive}. Together, the two estimates imply
    \begin{align}
        e^{-A\Leb(B_r(0))}
        \leq
        \frac{\kappa_{|\Lambda_n}(\xi,x)}{\kappa(\xi,x)}
        \leq
        1.
        \label{eq:58}
    \end{align}
    Applying these bounds and using Lemma~\ref{lem:propertiespapangelou} yields
    \begin{align}
        &\left|
        \I(\dP_{|\Lambda_n}|\dR_{|\Lambda_n})
        -
        \I(\dP_{|\Lambda_n}|\dR^{\Lambda_n})
        \right| \nonumber\\
        &\lesssim
        \left|
        \int_{\Lambda_n\setminus\Lambda_n^{r-}}\int
        \bigg[
        \left(
        \frac{\rho_n(\xi+\delta_x)}{\rho_n(\xi)}-1
        \right)
        \log\left(
        \frac{\rho_n(\xi+\delta_x)}{\rho_n(\xi)}
        \right)\right. \nonumber\\
        &\quad \left.
        -
        \left(
        \frac{\rho_n(\xi+\delta_x)f_n(\xi+\delta_x)}
        {\rho_n(\xi)f_n(\xi)}-1
        \right)
        \log\left(
        \frac{\rho_n(\xi+\delta_x)f_n(\xi+\delta_x)}
        {\rho_n(\xi)f_n(\xi)}
        \right)
        \bigg]
        \kappa_{|\Lambda_n}(\xi,x)
        \dP_{|\Lambda_n}(\mathrm{d}\xi)
        \Leb_{|\Lambda_n}(\mathrm{d}x)
        \right| \nonumber\\
        &=
        \left|
        \int_{\Lambda_n\setminus\Lambda_n^{r-}}\int
        \bigg[
        \left(
        \frac{\rho_n(\xi+\delta_x)}{\rho_n(\xi)}-1
        \right)
        \log\left(
        \frac{\rho_n(\xi+\delta_x)}{\rho_n(\xi)}
        \right)\right. \nonumber\\
        &\quad \left.
        -
        \left(
        \frac{\rho_n(\xi+\delta_x)\kappa_{|\Lambda_n}(\xi,x)}
        {\rho_n(\xi)\kappa(\xi,x)}-1
        \right)
        \log\left(
        \frac{\rho_n(\xi+\delta_x)\kappa_{|\Lambda_n}(\xi,x)}
        {\rho_n(\xi)\kappa(\xi,x)}
        \right)
        \bigg]
        \kappa_{|\Lambda_n}(\xi,x)
        \dP_{|\Lambda_n}(\mathrm{d}\xi)
        \Leb_{|\Lambda_n}(\mathrm{d}x)
        \right|,
        \label{eq:6789}
    \end{align}
    where $\lesssim$ denotes an inequality up to a finite constant. Define
    \begin{align*}
       u=u(\xi,x):=
       \frac{\rho_n(\xi+\delta_x)}{\rho_n(\xi)},
       \quad
       v=v(\xi,x):=
       \frac{\kappa_{|\Lambda_n}(\xi,x)}{\kappa(\xi,x)},
       \quad\text{ and }\quad
       h(w):=(w-1)\log(w).
    \end{align*}
    Since $0<v\leq1$ by~\eqref{eq:58}, we have $uv\leq u$. By Taylor's formula with integral remainder around $uv$ and the fact that $h''(w)=w^{-1}+w^{-2}$ is monotone decreasing, we obtain
    \begin{align*}
        \left|h(u)-h(uv)\right|
        &=
        \left|
        u(1-v)h'(uv)
        +
        \int_{uv}^{u}(u-s)h''(s)\mathrm{d}s
        \right|\\
        &\leq
        u(1-v)\left|h'(uv)\right|
        +
        \int_{uv}^{u}(u-s)h''(s)\mathrm{d}s\\
        &\leq
        u(1-v)\left|h'(uv)\right|
        +
        h''(uv)\int_{uv}^{u}(u-s)\mathrm{d}s\\
        &=
        u(1-v)
        \left|1-\frac{1}{uv}+\log(uv)\right|
        +
        \frac{1}{2}(1-v)^2
        \left(\frac{u}{v}+\frac{1}{v^2}\right).
    \end{align*}
    Using~\eqref{eq:58} once more, we may simplify this inequality to
    \begin{align*}
        \left|h(uv)-h(u)\right|
        \lesssim
        u+1+u\left|\log(uv)\right|.
    \end{align*}

    Inserting this estimate back into~\eqref{eq:6789} yields
    \begin{align}
        &\left|
        \I(\dP_{|\Lambda_n}|\dR_{|\Lambda_n})
        -
        \I(\dP_{|\Lambda_n}|\dR^{\Lambda_n})
        \right| \nonumber\\
        &\quad\lesssim
        \int_{\Lambda_n\setminus\Lambda_n^{r-}}\int
        \bigg(
        \frac{\rho_n(\xi+\delta_x)}{\rho_n(\xi)}
        +
        1
        +
        \frac{\rho_n(\xi+\delta_x)}{\rho_n(\xi)}
        \left|
        \log\left(
        \frac{\rho_n(\xi+\delta_x)f_n(\xi+\delta_x)}
        {\rho_n(\xi)f_n(\xi)}
        \right)
        \right|
        \bigg)
        \kappa_{|\Lambda_n}(\xi,x)
        \dP_{|\Lambda_n}(\mathrm{d}\xi)
        \Leb(\mathrm{d}x).
        \label{eq:92637}
    \end{align}
    It remains to show that the three terms on the right-hand side vanish in the large-volume limit. For the first term of~\eqref{eq:92637}, Corollary~\ref{kor:CampbellPapangelou} gives
    \begin{align*}
        &\lim_{n\to\infty}
        \frac{1}{\Leb(\Lambda_n)}
        \int_{\Lambda_n\setminus\Lambda_n^{r-}}\int
        \frac{\rho_n(\xi+\delta_x)}{\rho_n(\xi)}
        \kappa_{|\Lambda_n}(\xi,x)
        \dP_{|\Lambda_n}(\mathrm{d}\xi)
        \Leb_{|\Lambda_n}(\mathrm{d}x)\\
        &\quad =
        \lim_{n\to\infty}
        \frac{1}{\Leb(\Lambda_n)}
        \int_{\Lambda_n\setminus\Lambda_n^{r-}}\int
        \rho_n(\xi+\delta_x)
        \kappa_{|\Lambda_n}(\xi,x)
        \dR_{|\Lambda_n}(\mathrm{d}\xi)
        \Leb_{|\Lambda_n}(\mathrm{d}x)\\
        &\quad =
        \lim_{n\to\infty}
        \frac{1}{\Leb(\Lambda_n)}
        \I_{\dP}(\Lambda_n\setminus\Lambda_n^{r-})\\
        &\quad =
        0.
    \end{align*}
    For the second term of~\eqref{eq:92637}, \ref{property:locallystable} and stationarity imply
    \begin{align*}
        &\lim_{n\to\infty}
        \frac{1}{\Leb(\Lambda_n)}
        \int_{\Lambda_n\setminus\Lambda_n^{r-}}\int
        \kappa_{|\Lambda_n}(\xi,x)
        \dP_{|\Lambda_n}(\mathrm{d}\xi)
        \Leb_{|\Lambda_n}(\mathrm{d}x) \lesssim
        \lim_{n\to\infty}
        \frac{1}{\Leb(\Lambda_n)}
        \Leb(\Lambda_n\setminus\Lambda_n^{r-})\ =
        0.
    \end{align*}
    For the third term of~\eqref{eq:92637} we note that
    \begin{align}
        &\int_{\Lambda_n\setminus\Lambda_n^{r-}}\int
        \frac{\rho_n(\xi+\delta_x)}{\rho_n(\xi)}
        \left|
        \log\left(
        \frac{\rho_n(\xi+\delta_x)}{\rho_n(\xi)}
        \frac{f_n(\xi+\delta_x)}{f_n(\xi)}
        \right)
        \right|
        \kappa_{|\Lambda_n}(\xi,x)
        \dP_{|\Lambda_n}(\mathrm{d}\xi)
        \Leb_{|\Lambda_n}(\mathrm{d}x) \nonumber\\
        &\quad\leq
        \int_{\Lambda_n\setminus\Lambda_n^{r-}}\int
        \frac{\rho_n(\xi+\delta_x)}{\rho_n(\xi)}
        \kappa_{|\Lambda_n}(\xi,x)
        \left|
        \log\left(
        \frac{\rho_n(\xi+\delta_x)}{\rho_n(\xi)}
        \kappa_{|\Lambda_n}(\xi,x)
        \right)
        \right|
        \dP_{|\Lambda_n}(\mathrm{d}\xi)
        \Leb_{|\Lambda_n}(\mathrm{d}x) \nonumber\\
        &\qquad +
        \int_{\Lambda_n\setminus\Lambda_n^{r-}}\int
        \frac{\rho_n(\xi+\delta_x)}{\rho_n(\xi)}
        \kappa_{|\Lambda_n}(\xi,x)
        \left|
        \log\left(\frac{1}{\kappa(\xi,x)}\right)
        \right|
        \dP_{|\Lambda_n}(\mathrm{d}\xi)
        \Leb_{|\Lambda_n}(\mathrm{d}x).
        \label{eq:7890}
    \end{align}
    We first consider the first term on the right-hand side of~\eqref{eq:7890}. By Lemma~\ref{lem:regularisationFisher}, $\dP_{|\Lambda_n}$ is a Papangelou measure with Papangelou intensity
\begin{align*}
    \kappa^\dP_n(\xi,x)
    =
    \frac{\rho_n(\xi+\delta_x)}{\rho_n(\xi)}
    \kappa_{|\Lambda_n}(\xi,x),
    \quad \xi\in\Gamma_{\Lambda_n},\; x\in\Lambda_n.
\end{align*}
Moreover, analogously to Lemma~\ref{lem:restrictionpapangelou}, $\dP_{|\Lambda_{n-1}}$ is again a Papangelou measure with Papangelou intensity
\begin{align*}
    \kappa^\dP_{n-1}(\xi_1,x)
    =
    \int
    \kappa^\dP_n(\xi_1+\xi_2,x)
    \dP^{\Lambda_n\setminus\Lambda_{n-1}}_{|\Lambda_n}(\mathrm{d}\xi_2|\xi_1),
    \quad \xi_1\in\Gamma_{\Lambda_{n-1}},\; x\in\Lambda_{n-1},
\end{align*}
where $\dP^{\Lambda_n\setminus\Lambda_{n-1}}_{|\Lambda_n}(\cdot|\xi_1)$ is defined by
\begin{align*}
    \int F(\xi_1+\xi_2)
    \dP_{|\Lambda_n}(\mathrm{d}(\xi_1+\xi_2))
    =
    \int
    \int F(\xi_1+\xi_2)
    \dP^{\Lambda_n\setminus\Lambda_{n-1}}_{|\Lambda_n}(\mathrm{d}\xi_2|\xi_1)
    \dP_{|\Lambda_{n-1}}(\mathrm{d}\xi_1),\quad F\in\FFF^+(\Gamma_{\Lambda_n}).
\end{align*}
By Jensen's inequality, we obtain
\begin{align*}
    g_{\Lambda_{n-1}}(x)
    &:=
    \int
    \kappa^\dP_{n-1}(\xi,x)
    \log\big(\kappa^\dP_{n-1}(\xi,x)\big)
    \dP_{|\Lambda_{n-1}}(\mathrm{d}\xi)\leq
    \int
    \kappa^\dP_n(\xi,x)
    \log\big(\kappa^\dP_n(\xi,x)\big)
    \dP_{|\Lambda_n}(\mathrm{d}\xi)=
    g_{\Lambda_n}(x).
\end{align*}
Thus, $g_{\Lambda_n}(x)$ is monotone increasing in $n$ and the limit $g_\infty(x):=\lim_{n\to\infty}g_{\Lambda_n}(x)$ exists. Further, we note that by stationarity $g_{\Lambda_n}(x)=g_{\Lambda_n+y}(x+y)$. Thus, if $m\in\N$ is s.t.\  $\Lambda_n+y\subset\Lambda_m$ and $l\in\N$ is s.t.\  $\Lambda_m-y\subset\Lambda_l$, then
    \begin{align*}
        g_{\Lambda_n}(x)
        =
        g_{\Lambda_n+y}(x+y)
        \leq
        g_{\Lambda_m}(x+y)
        =
        g_{\Lambda_m-y}(x)
        \leq
        g_{\Lambda_l}(x).
    \end{align*}
    Letting $n\to\infty$ shows that $g_\infty(x)=g_\infty(x+y)$, which implies that $g_\infty(x)$ is constant. We assume for now that $g_\infty(x)\equiv g_\infty<\infty$ and justify this later. For the first term of~\eqref{eq:7890}, we  obtain
    \begin{align*}
        &\lim_{n\to\infty}
        \frac{1}{\Leb(\Lambda_n)}
        \int_{\Lambda_n\setminus\Lambda_n^{r-}}\int
        \frac{\rho_n(\xi+\delta_x)}{\rho_n(\xi)}
        \kappa_{|\Lambda_n}(\xi,x)
        \left|
        \log\left(
        \frac{\rho_n(\xi+\delta_x)}{\rho_n(\xi)}
        \kappa_{|\Lambda_n}(\xi,x)
        \right)
        \right|
        \dP_{|\Lambda_n}(\mathrm{d}\xi)
        \Leb_{|\Lambda_n}(\mathrm{d}x)\\
        &\quad\leq
        \lim_{n\to\infty}
        \frac{1}{\Leb(\Lambda_n)}
        \int_{\Lambda_n\setminus\Lambda_n^{r-}}\int
        \frac{\rho_n(\xi+\delta_x)}{\rho_n(\xi)}
        \kappa_{|\Lambda_n}(\xi,x)
        \log\left(
        \frac{\rho_n(\xi+\delta_x)}{\rho_n(\xi)}
        \kappa_{|\Lambda_n}(\xi,x)
        \right)
        +2e^{-1}
        \dP_{|\Lambda_n}(\mathrm{d}\xi)
        \Leb_{|\Lambda_n}(\mathrm{d}x)\\
        &\quad\leq
        \lim_{n\to\infty}
        \frac{1}{\Leb(\Lambda_n)}
        \int_{\Lambda_n\setminus\Lambda_n^{r-}}
        \big(g_\infty(x)+2e^{-1}\big)
        \Leb_{|\Lambda_n}(\mathrm{d}x)\\
        &\quad\lesssim
        \lim_{n\to\infty}
        \frac{\Leb(\Lambda_n\setminus\Lambda_n^{r-})}
        {\Leb(\Lambda_n)}\\
        &\quad=
        0.
    \end{align*}
    We next consider the second term of~\eqref{eq:7890}. By \ref{property:finiterange}, $\kappa(\xi,x)=\kappa(\xi_{|\Lambda_n},x)$ for all $n>2r+\|x\|$. Hence, for every $k\geq0$ and all such $n$,
    \begin{align*}
        \widetilde g_{\Lambda_{n+k}}(x)
        &:=
        \int
        \kappa^\dP_{n+k}(\xi,x)
        \left|
        \log\left(\frac{1}{\kappa(\xi,x)}\right)
        \right|
        \dP_{|\Lambda_{n+k}}(\mathrm{d}\xi)=
        \int
        \kappa^\dP_n(\xi,x)
        \left|
        \log\left(\frac{1}{\kappa(\xi,x)}\right)
        \right|
        \dP_{|\Lambda_n}(\mathrm{d}\xi)=
        \widetilde g_{\Lambda_n}(x).
    \end{align*}
    By shift-invariance and stationarity, it holds that $\widetilde g_{\Lambda_n}(x)=\widetilde g_{\Lambda_n+y}(x+y)$. The same sandwich argument as for $g_{\Lambda_n}$ shows that
    \[
        \widetilde g_\infty(x)
        :=
        \lim_{n\to\infty}\widetilde g_{\Lambda_n}(x)
    \]
    exists and is independent of $x$. We again assume that $\widetilde g_\infty(x)\equiv\widetilde g_{\infty}<\infty$ and postpone the proof to later. We obtain for the second term of~\eqref{eq:7890} that
    \begin{align*}
        &\lim_{n\to\infty}
        \frac{1}{\Leb(\Lambda_n)}
        \int_{\Lambda_n\setminus\Lambda_n^{r-}}\int
        \frac{\rho_n(\xi+\delta_x)}{\rho_n(\xi)}
        \kappa_{|\Lambda_n}(\xi,x)
        \left|
        \log\left(\frac{1}{\kappa(\xi,x)}\right)
        \right|
        \dP_{|\Lambda_n}(\mathrm{d}\xi)
        \Leb_{|\Lambda_n}(\mathrm{d}x)\\
        &\quad =
        \lim_{n\to\infty}
        \frac{1}{\Leb(\Lambda_n)}
        \int_{\Lambda_n\setminus\Lambda_n^{r-}}
        \widetilde g_{\Lambda_n}(x)
        \Leb_{|\Lambda_n}(\mathrm{d}x)\\
        &\quad \lesssim
        \lim_{n\to\infty}
        \frac{\Leb(\Lambda_n\setminus\Lambda_n^{r-})}
        {\Leb(\Lambda_n)}\\
        &\quad =
        0.
    \end{align*}
    All in all this proves that the right-hand side of~\eqref{eq:92637} vanishes in the normalized large-volume limit. It remains to prove that $g_\infty$ and $\widetilde g_{\infty}$ are finite. Define
    \[
        \Phi(p,q):=(p-q)\log\left(\frac{p}{q}\right),
    \]
    and
    \[
        j_{\Lambda_n}(x)
        :=
        \int
        \Phi\big(
            \kappa^\dP_n(\xi,x),
            \kappa_{|\Lambda_n}(\xi,x)
        \big)
        \dP_{|\Lambda_n}(\mathrm{d}\xi).
    \]
    By the definition of $\kappa^\dP_n$, it holds that
    \begin{align*}
        \int_{\Lambda_n}j_{\Lambda_n}(x)\Leb(\mathrm{d}x)
        &=
        \int_{\Lambda_n}\int
        (\rho_n(\xi+\delta_x)-\rho_n(\xi))
        \log\left(
        \frac{\rho_n(\xi+\delta_x)}{\rho_n(\xi)}
        \right)
        \kappa_{|\Lambda_n}(\xi,x)
        \dR_{|\Lambda_n}(\mathrm{d}\xi)
        \Leb(\mathrm{d}x)\\
        &=
        \I(\dP_{|\Lambda_n}|\dR_{|\Lambda_n}).
    \end{align*}
   By disintegration, we obtain for $m<n$ that
\begin{align*}
    \rho_m(\xi_1)
    =
    \int
    \rho_n(\xi_1+\xi_2)
    \dR_{|\Lambda_n}^{\Lambda_n\setminus\Lambda_m}
    (\mathrm{d}\xi_2|\xi_1),
    \quad \xi_1\in\Gamma_{\Lambda_m},
\end{align*}
where $\dR_{|\Lambda_n}^{\Lambda_n\setminus\Lambda_m}(\cdot|\xi_1)$ is defined by
\begin{align*}
    \int F(\xi_1+\xi_2)
    \dR_{|\Lambda_n}(\mathrm{d}(\xi_1+\xi_2))
    =
    \int
    \int F(\xi_1+\xi_2)
    \dR_{|\Lambda_n}^{\Lambda_n\setminus\Lambda_m}(\mathrm{d}\xi_2|\xi_1)
    \dR_{|\Lambda_m}(\mathrm{d}\xi_1),\quad F\in\FFF^+(\Gamma_{\Lambda_n}).
\end{align*}
 By \ref{property:finiterange}, it holds that
    \begin{align*}
        \rho_m(\xi_1+\delta_x)
        =
        \int
        \rho_n(\xi_1+\xi_2+\delta_x)
        \dR_{|\Lambda_n}^{\Lambda_n\setminus\Lambda_m}
        (\mathrm{d}\xi_2|\xi_1),
    \end{align*}
    for all $x\in\Lambda_m^{r-}$.
    Hence, by the joint convexity of $\Phi$, Jensen's inequality, and~\ref{property:finiterange}, we get 
    \begin{align*}
        j_{\Lambda_m}(x)
        &=
        \int
        (\rho_m(\xi+\delta_x)-\rho_m(\xi))
        \log\left(
        \frac{\rho_m(\xi+\delta_x)}{\rho_m(\xi)}
        \right)
        \kappa_{|\Lambda_m}(\xi,x)
        \dR_{|\Lambda_m}(\mathrm{d}\xi)\\
        &\leq
        \int
        (\rho_n(\xi+\delta_x)-\rho_n(\xi))
        \log\left(
        \frac{\rho_n(\xi+\delta_x)}{\rho_n(\xi)}
        \right)
        \kappa_{|\Lambda_n}(\xi,x)
        \dR_{|\Lambda_n}(\mathrm{d}\xi)\\
        &=
        j_{\Lambda_n}(x).
    \end{align*}
    The same stationarity and sandwich argument as for $g_{\Lambda_n}$ shows that $j_\infty\equiv\lim_{n\to\infty}j_{\Lambda_n}(x)$ exists and is constant. Moreover, for $m\leq n$ large enough and all $x\in\Lambda_n$ s.t.\  $x+\Lambda_m\subset\Lambda_n$, monotonicity and stationarity imply
    \[
        j_{\Lambda_n}(x)
        \geq
        j_{x+\Lambda_m}(x)
        =
        j_{\Lambda_m}(0).
    \]
    Since $j_{\Lambda_n}(x)$ is non-negative, this implies
    \begin{align}
        \frac{1}{\Leb(\Lambda_n)}
        \int_{\Lambda_n}j_{\Lambda_n}(x)\Leb(\mathrm{d}x)
        \geq
        \frac{\Leb(\{x\in\Lambda_n:x+\Lambda_m\subset\Lambda_n\})}
        {\Leb(\Lambda_n)}
        j_{\Lambda_m}(0).
        \label{eq:01234}
    \end{align}
    Consequently, we deduce that
    \[
        j_\infty
        =
        j_\infty(0)
        \leq
        \liminf_{n\to\infty}
        \frac{1}{\Leb(\Lambda_n)}
        \int_{\Lambda_n}j_{\Lambda_n}(x)\Leb(\mathrm{d}x)
        =
        \mathcal I(\dP|\dR)
        <\infty.
    \]

    We proceed with the proof of the finiteness of $g_\infty$. To this end, we first show that for all \(p\geq 0\) and \(0\leq q\leq A\), there exists a constant \(C_A<\infty\) s.t.\ 
\begin{equation}
    p\log(p)
    \leq
    C_A\bigl(1+\Phi(p,q)\bigr).
    \label{eq:10000}
\end{equation}
If \(p\leq 1\), then \(p\log(p)\leq 0\). If \(q=0\), the estimate is trivial since \(\Phi(p,0)=\infty\). Thus, let \(q>0\). 
For \(p\geq \max\{2A,A^2,1\}\), we obtain
\[
    \Phi(p,q)
    =
    (p-q)\log\left(\frac{p}{q}\right)
    \geq
    \frac{p}{2}\log\left(\frac{p}{A}\right)
    \geq
    \frac{1}{4}p\log(p).
\]
It remains to consider \(p\in[0,\max\{2A,A^2,1\})\). On this interval, \(p\log(p)\) is bounded above by a constant depending only on \(A\), which proves~\eqref{eq:10000}. Since $0\leq\kappa_{|\Lambda_n}\leq A$ by \ref{property:locallystable}, we may apply~\eqref{eq:10000} with $p=\kappa^\dP_n(\xi,x)$ and $q=\kappa_{|\Lambda_n}(\xi,x)$ and obtain
    \[
        -e^{-1}
        \leq
        g_\infty
        =
        \lim_{n\to\infty}g_{\Lambda_n}(x)
        \leq
        \lim_{n\to\infty}
        C_A(1+j_{\Lambda_n}(x))
        =
        C_A\big(1+j_\infty\big)
        <\infty.
    \]
	   It remains to prove that $\widetilde g_{\infty}$ is finite. We first show that
	for all $p\geq0$, $0\leq q\leq A$, and $w>0$ with
	$c\leq q/w\leq1$, there exists a constant $C_{A,c}<\infty$ s.t.\ 
\begin{align}
    p
    \left|
    \log\left(\frac{1}{w}\right)
    \right|
    \leq
    C_{A,c}\bigl(1+\Phi(p,q)\bigr).
    \label{eq:20000}
\end{align}
If $q=0$, then $\Phi(p,0)=\infty$. Hence, let
	$q>0$. Since $c\leq q/w\leq1$ it holds that
\[
    \left|
    \log\left(\frac{1}{w}\right)
    \right|
    \leq
    C_c+
    \left|
    \log\left(\frac{1}{q}\right)
    \right|.
\]
Let $\ell:=p/q$. If $\ell\leq2$, then $p=q\ell\leq2q$, and thus
\[
    p
    \left|
    \log\left(\frac{1}{w}\right)
    \right|
    \leq
    C_c q + 2q\left|\log\left(\frac{1}{q}\right)\right|
    \leq C_{A,c}.
\]
The last inequality is true since $q|\log\left(\frac{1}{q}\right)|$ is bounded on $[0,A]$. Now, let $\ell>2$. Then
\[
    p
    \left|
    \log\left(\frac{1}{w}\right)
    \right|
    \leq
    C_{A,c}\left(1+q\ell\log(\ell)\right).
\]
Moreover, since $\ell>2$, we have $\ell\leq2(\ell-1)$, and hence
\[
    q\ell\log(\ell)
    \leq
    2q(\ell-1)\log(\ell)
    =
    2\Phi(p,q).
\]
Consequently,
\[
    p
    \left|
    \log\left(\frac{1}{w}\right)
    \right|
    \leq
    C_{A,c}\bigl(1+\Phi(p,q)\bigr),
\]
which proves \eqref{eq:20000}.
    By~\eqref{eq:58} and~\ref{property:locallystable}, we may apply the inequality of~\eqref{eq:20000} with
    \(
        p=\kappa^\dP_n(\xi,x),\;
        q=\kappa_{|\Lambda_n}(\xi,x),
    \)
    and $w=\kappa(\xi,x)$ and obtain
    \[
        \widetilde g_\infty
        =
        \lim_{n\to\infty}\widetilde g_{\Lambda_n}(x)
        \leq
        \lim_{n\to\infty}
        C_{A,c}\big(1+j_{\Lambda_n}(x)\big)
        =
        C_{A,c}\big(1+j_\infty\big)
        <\infty,
    \]
    which completes the proof for $\mathcal I(\dP|\dR)<\infty$.

    It remains to consider the case $\mathcal I(\dP|\dR)=\infty$. We recall that by \ref{eq:ffiniterange},  $\kappa_{|\Lambda_n}(\xi,x)=\kappa(\xi,x)$ and $f_n(\xi+\delta_x)=f_n(\xi)$, for $\xi\in\Gamma_{\Lambda_n}$ and  $x\in\Lambda_n^{r-}$, see~\eqref{eq:ffiniterange}. Hence, the integrand in $\I(\dP_{|\Lambda_n}|\dR^{\Lambda_n})$ coincides on $\Lambda_n^{r-}$ with the integrand defining $j_{\Lambda_n}(x)$ and we obtain
    \[
        \I(\dP_{|\Lambda_n}|\dR^{\Lambda_n})
        \geq
        \int_{\Lambda_n^{r-}} j_{\Lambda_n}(x)\Leb(\mathrm{d}x).
    \]
    Let $m\in\N$. Analogously to~\eqref{eq:01234}, we obtain
    \[
        \liminf_{n\to\infty}
        \frac{1}{\Leb(\Lambda_n)}
        \I(\dP_{|\Lambda_n}|\dR^{\Lambda_n})
        \geq
        \liminf_{n\to\infty}
        \frac{
            \Leb(\{x\in\Lambda_n^{r-}:x+\Lambda_m\subset\Lambda_n^{r-}\})
        }
        {\Leb(\Lambda_n)}
        j_{\Lambda_m}(0).
    \]
    Finally, $\mathcal I(\dP|\dR)=\infty$ implies that $\lim_{m\to\infty}j_{\Lambda_m}(0)=j_\infty=\infty$, finishing the proof.
\end{proof}

\subsection{Convergence of the Curvature constant}
We continue with the proof of  Theorem~\ref{thm:convergenceresultsstationary} $(v)$ and show that if $\kappa$ satisfies \ref{property:cocycle}, \ref{property:locallystable}, \ref{property:repulsive}, \ref{property:finiterange}, and~\ref{property:shiftinvariance}, and $\kappa_n$ satisfies \hyperref[property:dobrushin]{$(\mathbf{D}(\zeta(\kappa_n)))$} for all $n\in\N$, then $\kappa$ satisfies \ref{property:dobrushin} and it holds that 
\begin{align*}
   0< \lim_{n\to\infty}\zeta(\kappa_n)=\zeta(\kappa).
\end{align*}
   \begin{proof}[Proof of Theorem~\ref{thm:convergenceresultsstationary} $(v)$]
  Since $\kappa$ satisfies \ref{property:repulsive}, it holds that $\kappa(\xi,y)\geq\kappa(\xi+\delta_{x},y)$ for all  $\xi\in\Gamma,\; x,y\in\R^d$. Further, using \ref{property:finiterange} and~\ref{property:shiftinvariance}, we get for $n>4r$ that
\begin{align*}
    \int
    \big(
        \kappa(\xi,y)-\kappa(\xi+\delta_x,y)
    \big)
    \Leb(\mathrm{d}y)& =
    \int_{B_r(x)}
    \big(
        \kappa(\xi_{|B_{2r}(x)},y)
        -
        \kappa(\xi_{|B_{2r}(x)}+\delta_x,y)
    \big)
    \Leb(\mathrm{d}y)\\
    & =
    \int_{B_r(0)}
    \big(
        \kappa(\Theta_x\xi_{|B_{2r}(x)},y)
        -
        \kappa(\Theta_x\xi_{|B_{2r}(x)}+\delta_0,y)
    \big)
    \Leb(\mathrm{d}y)\\
    & \leq
    \sup_{\substack{\eta\in\Gamma_{\Lambda_n}\\ z\in\Lambda_n}}
    \int_{\Lambda_n}
    \big(
        \kappa(\eta,y)-\kappa(\eta+\delta_z,y)
    \big)
    \Leb(\mathrm{d}y),
\end{align*}
for all $\xi\in\Gamma,\;x\in\R^d$. Hence, $\zeta(\kappa_n)\leq \zeta(\kappa)$ for all $n>4r$. Conversely, by \ref{property:repulsive}, we have
\begin{align*}
    &\sup_{\substack{\xi\in\Gamma_{\Lambda_n}\\ x\in\Lambda_n}}
    \int_{\Lambda_n}
    \big(
        \kappa(\xi,y)-\kappa(\xi+\delta_x,y)
    \big)
    \Leb(\mathrm{d}y) \leq
    \sup_{\substack{\xi\in\Gamma\\ x\in\R^d}}
    \int_{\R^d}
    \big(
        \kappa(\xi,y)-\kappa(\xi+\delta_x,y)
    \big)
    \Leb(\mathrm{d}y),
\end{align*}
implying that $\zeta(\kappa)\leq \zeta(\kappa_n)$.  All in all, we obtain that $0<\zeta(\kappa_n)=\zeta(\kappa)$ for all $n>4r$, finishing the proof.
\end{proof} 
\section*{Acknowledgments}
ME, MH, and HS are supported by the Deutsche Forschungsgemeinschaft (DFG, German Research Foundation) through the SPP 2265 { “Random Geometric Systems"}. MH and HS have been funded and supported by the Deutsche Forschungsgemeinschaft (DFG, German Research Foundation) under Germany's Excellence Strategy EXC 2044/2 - 390685587 “Mathematics M\"unster: Dynamics--Geometry--Structure" and through the Research Training Group 3027 - 524444762 “Rigorous Analysis of Complex Random Systems"
\section*{AI Statement}
The mathematical arguments were developed by the authors without any assistance
by AI. Large language models were used for proof and grammar checking, the text
was written solely by the authors.
\addcontentsline{toc}{section}{References}
\bibliography{references}

@article{huesmann2025nonlocalwassersteingeometrygradient,
  title         = {{Non-Local {Wasserstein} Geometry, Gradient Flows, and Functional Inequalities for Stationary Point Processes}},
  author        = {Huesmann, Martin and Stange, Hanna},
  journal       = {{arXiv} preprint {arXiv}:2504.12047},
  year          = {2025},
  doi           = {10.48550/arXiv.2504.12047},
  url           = {https://arxiv.org/abs/2504.12047}
}

@book{daley2003introduction,
  title     = {{An Introduction to the Theory of Point Processes. Volume I: Elementary Theory and Methods}},
  author    = {Daley, Daryl J. and Vere-Jones, David},
  edition   = {2nd},
  series    = {Probability and its Applications},
  publisher = {Springer},
  address   = {New York},
  year      = {2003},
  isbn      = {9780387955413}
}

@article{Che85,
  author  = {Chen, Louis H. Y.},
  title   = {Poincar{\'e}-type inequalities via stochastic integrals},
  journal = {Z. Wahrsch. Verw. Gebiete},
  volume  = {69},
  number  = {2},
  pages   = {251--277},
  year    = {1985},
  doi     = {10.1007/BF02450284}
}

@article{HPA95,
  author  = {Houdr{\'e}, Christian and P{\'e}rez-Abreu, V{\'\i}ctor},
  title   = {Covariance identities and inequalities for functionals on {W}iener and {P}oisson spaces},
  journal = {Ann. Probab.},
  volume  = {23},
  number  = {1},
  pages   = {400--419},
  year    = {1995},
  doi     = {10.1214/aop/1176988392}
}

@article{LP11,
  author  = {Last, G{\"u}nter and Penrose, Mathew D.},
  title   = {Poisson process {F}ock space representation, chaos expansion and covariance inequalities},
  journal = {Probab. Theory Related Fields},
  volume  = {150},
  number  = {3-4},
  pages   = {663--690},
  year    = {2011},
  doi     = {10.1007/s00440-010-0288-5}
}

@article{BL98,
  author  = {Bobkov, Sergey G. and Ledoux, Michel},
  title   = {On modified logarithmic {S}obolev inequalities for {B}ernoulli and {P}oisson measures},
  journal = {J. Funct. Anal.},
  volume  = {156},
  number  = {2},
  pages   = {347--365},
  year    = {1998},
  doi     = {10.1006/jfan.1997.3187}
}

@article{BT06,
  author  = {Bobkov, Sergey G. and Tetali, Prasad},
  title   = {Modified logarithmic {S}obolev inequalities in discrete settings},
  journal = {J. Theoret. Probab.},
  volume  = {19},
  number  = {2},
  pages   = {289--336},
  year    = {2006},
  doi     = {10.1007/s10959-006-0016-3}
}

@article{Cha04,
  author  = {Chafa{\"\i}, Djalil},
  title   = {Entropies, convexity, and functional inequalities: on {$\Phi$}-entropies and {$\Phi$}-{S}obolev inequalities},
  journal = {J. Math. Kyoto Univ.},
  volume  = {44},
  number  = {2},
  pages   = {325--363},
  year    = {2004},
  doi     = {10.1215/kjm/1250283556}
}

@article{Cha06,
  author  = {Chafa{\"\i}, Djalil},
  title   = {Binomial--{P}oisson entropic inequalities and the {M/M/$\infty$} queue},
  journal = {ESAIM Probab. Stat.},
  volume  = {10},
  pages   = {317--339},
  year    = {2006},
  doi     = {10.1051/ps:2006013}
}

@article{CJ13,
  author  = {Chafa{\"\i}, Djalil and Joulin, Ald{\'e}ric},
  title   = {Intertwining and commutation relations for birth--death processes},
  journal = {Bernoulli},
  volume  = {19},
  number  = {5A},
  pages   = {1855--1879},
  year    = {2013},
  doi     = {10.3150/12-BEJ433}
}

@misc{ALS24,
  author        = {Aryan, Shrey and L{\'o}pez-Rivera, Pablo and Shenfeld, Yair},
  title         = {Improvement of {W}u's logarithmic {S}obolev inequality via the {P}oisson--{F}{\"o}llmer process},
  year          = {2024},
  eprint        = {2410.06117},
  archiveprefix = {arXiv},
  primaryclass  = {math.PR},
  howpublished  = {arXiv:2410.06117}
}

@article{BCC02,
  author  = {Bertini, Lorenzo and Cancrini, Nicoletta and Cesi, Filippo},
  title   = {The spectral gap for a {G}lauber-type dynamics in a continuous gas},
  journal = {Ann. Inst. H. Poincar{\'e} Probab. Statist.},
  volume  = {38},
  number  = {1},
  pages   = {91--108},
  year    = {2002},
  doi     = {10.1016/S0246-0203(01)01090-3}
}

@article{SS14,
  author  = {Schuhmacher, Dominic and Stucki, Kaspar},
  title   = {Gibbs point process approximation: total variation bounds using {S}tein's method},
  journal = {Ann. Probab.},
  volume  = {42},
  number  = {5},
  pages   = {1911--1951},
  year    = {2014},
  doi     = {10.1214/13-AOP895}
}

@article{Tor17,
  author  = {Torrisi, Giovanni Luca},
  title   = {Probability approximation of point processes with {P}apangelou conditional intensity},
  journal = {Bernoulli},
  volume  = {23},
  number  = {4A},
  pages   = {2210--2256},
  year    = {2017},
  doi     = {10.3150/16-BEJ808}
}

@misc{DV18,
  author        = {Decreusefond, Laurent and Vasseur, Aur{\'e}lien},
  title         = {Stein's method and {P}apangelou intensity for {P}oisson or {C}ox process approximation},
  year          = {2018},
  eprint        = {1807.02453},
  archiveprefix = {arXiv},
  primaryclass  = {math.PR},
  howpublished  = {arXiv:1807.02453}
}

@article{KKO13,
  author  = {Kondratiev, Yuri and Kuna, Tobias and Ohlerich, Nataliya},
  title   = {Spectral gap for {G}lauber type dynamics for a special class of potentials},
  journal = {Electron. J. Probab.},
  volume  = {18},
  pages   = {paper no. 42, 18 pp.},
  year    = {2013},
  doi     = {10.1214/EJP.v18-2260}
}

@article{CDPP09,
  author  = {Caputo, Pietro and Dai Pra, Paolo and Posta, Gustavo},
  title   = {Convex entropy decay via the {B}ochner--{B}akry--{E}mery approach},
  journal = {Ann. Inst. Henri Poincar{\'e} Probab. Stat.},
  volume  = {45},
  number  = {3},
  pages   = {734--753},
  year    = {2009},
  doi     = {10.1214/08-AIHP183}
}

@article{DPP13,
  author  = {Dai Pra, Paolo and Posta, Gustavo},
  title   = {Entropy decay for interacting systems via the {B}ochner--{B}akry--{\'E}mery approach},
  journal = {Electron. J. Probab.},
  volume  = {18},
  pages   = {paper no. 52, 21 pp.},
  year    = {2013},
  doi     = {10.1214/EJP.v18-2041}
}

@article{Con22,
  author  = {Conforti, Giovanni},
  title   = {A probabilistic approach to convex {$(\phi)$}-entropy decay for {M}arkov chains},
  journal = {Ann. Appl. Probab.},
  volume  = {32},
  number  = {2},
  pages   = {932--973},
  year    = {2022},
  doi     = {10.1214/21-AAP1700}
}

@article{Rei13,
  author  = {Reitzner, Matthias},
  title   = {Poisson point processes: large deviation inequalities for the convex distance},
  journal = {Electron. Commun. Probab.},
  volume  = {18},
  pages   = {paper no. 96, 7 pp.},
  year    = {2013},
  doi     = {10.1214/ECP.v18-2851}
}

@article{BP16,
  author  = {Bachmann, Sascha and Peccati, Giovanni},
  title   = {Concentration bounds for geometric {P}oisson functionals: logarithmic {S}obolev inequalities revisited},
  journal = {Electron. J. Probab.},
  volume  = {21},
  pages   = {paper no. 6, 44 pp.},
  year    = {2016},
  doi     = {10.1214/16-EJP4235}
}

@article{GHP21,
  author  = {Gozlan, Nathael and Herry, Ronan and Peccati, Giovanni},
  title   = {Transport inequalities for random point measures},
  journal = {J. Funct. Anal.},
  volume  = {281},
  number  = {9},
  pages   = {paper no. 109141, 45 pp.},
  year    = {2021},
  doi     = {10.1016/j.jfa.2021.109141}
}

@article{Dec08,
  author  = {Decreusefond, Laurent},
  title   = {Wasserstein distance on configuration space},
  journal = {Potential Anal.},
  volume  = {28},
  number  = {3},
  pages   = {283--300},
  year    = {2008},
  doi     = {10.1007/s11118-008-9077-5}
}

@article{DJS10,
  author  = {Decreusefond, Laurent and Joulin, Ald{\'e}ric and Savy, Nicolas},
  title   = {Upper bounds on {R}ubinstein distances on configuration spaces and applications},
  journal = {Commun. Stoch. Anal.},
  volume  = {4},
  number  = {3},
  pages   = {377--399},
  year    = {2010},
  note    = {CHECK volume/pages}
}

@article{MSWW11,
  author  = {Ma, Yutao and Shen, Shi and Wang, Xinyu and Wu, Liming},
  title   = {Transportation inequalities: from {P}oisson to {G}ibbs measures},
  journal = {Bernoulli},
  volume  = {17},
  number  = {1},
  pages   = {155--169},
  year    = {2011},
  doi     = {10.3150/00-BEJ268}
}

@article{MWW11,
  author  = {Ma, Yutao and Wang, Ran and Wu, Liming},
  title   = {Transportation-information inequalities for continuum {G}ibbs measures},
  journal = {Electron. Commun. Probab.},
  volume  = {16},
  pages   = {600--613},
  year    = {2011},
  doi     = {10.1214/ECP.v16-1670}
}

@book{daley2007introduction,
  title     = {{An Introduction to the Theory of Point Processes. Volume II: General Theory and Structure}},
  author    = {Daley, Daryl J. and Vere-Jones, David},
  edition   = {2nd},
  series    = {Probability and its Applications},
  publisher = {Springer},
  address   = {New York},
  year      = {2007},
  isbn      = {9780387213377}
}

@book{last_penrose_2017,
  title     = {{Lectures on the {Poisson} Process}},
  author    = {Last, G{\"u}nter and Penrose, Mathew D.},
  series    = {Institute of Mathematical Statistics Textbooks},
  publisher = {Cambridge University Press},
  year      = {2017},
  doi       = {10.1017/9781316104476}
}

@article{Preston,
  author  = {Preston, Chris},
  title   = {{Spatial birth-and-death processes}},
  journal = {Bulletin of the International Statistical Institute},
  volume  = {46},
  number  = {2},
  pages   = {371--391},
  year    = {1977}
}

@article{schiavo2023wasserstein,
  author    = {Dello Schiavo, Lorenzo and Herry, Ronan and Suzuki, Kohei},
  title     = {{{Wasserstein} geometry and {Ricci} curvature bounds for {Poisson} spaces}},
  journal   = {Journal de l'{\'E}cole polytechnique --- Math{\'e}matiques},
  volume    = {11},
  pages     = {957--1010},
  year      = {2024},
  publisher = {\'Ecole polytechnique},
  doi       = {10.5802/jep.270},
  url       = {https://jep.centre-mersenne.org/articles/10.5802/jep.270/}
}

@article{erbar2012gradientflowsentropyjump,
  author  = {Erbar, Matthias},
  title   = {{Gradient flows of the entropy for jump processes}},
  journal = {Annales de l'Institut Henri Poincar{\'e} Probabilit{\'e}s et Statistiques},
  volume  = {50},
  number  = {3},
  pages   = {920--945},
  year    = {2014},
  doi     = {10.1214/12-AIHP537},
  url     = {https://doi.org/10.1214/12-AIHP537}
}

@book{fekete,
  author    = {Rassoul-Agha, Firas and Sepp{\"a}l{\"a}inen, Timo},
  title     = {{A Course on Large Deviations with an Introduction to {Gibbs} Measures}},
  series    = {Graduate Studies in Mathematics},
  volume    = {162},
  publisher = {American Mathematical Society},
  address   = {Providence, RI},
  year      = {2015},
  isbn      = {978-0-8218-7578-0}
}

@article{Dolbeault_2008,
  author  = {Dolbeault, Jean and Nazaret, Bruno and Savar{\'e}, Giuseppe},
  title   = {{A new class of transport distances between measures}},
  journal = {Calculus of Variations and Partial Differential Equations},
  volume  = {34},
  number  = {2},
  pages   = {193--231},
  year    = {2008},
  doi     = {10.1007/s00526-008-0182-5},
  url     = {https://doi.org/10.1007/s00526-008-0182-5}
}

@book{Ambrosio2008GradientFI,
  author    = {Ambrosio, Luigi and Gigli, Nicola and Savar{\'e}, Giuseppe},
  title     = {{Gradient Flows: In Metric Spaces and in the Space of Probability Measures}},
  edition   = {2nd},
  series    = {Lectures in Mathematics ETH Z{\"u}rich},
  publisher = {Birkh{\"a}user},
  address   = {Basel},
  year      = {2008},
  doi       = {10.1007/978-3-7643-8722-8}
}

@book{georgii2011gibbs,
  author    = {Georgii, Hans-Otto},
  title     = {{{Gibbs} Measures and Phase Transitions}},
  edition   = {2nd},
  series    = {De Gruyter Studies in Mathematics},
  volume    = {9},
  publisher = {De Gruyter},
  address   = {Berlin},
  year      = {2011},
  doi       = {10.1515/9783110250329},
  isbn      = {978-3-11-025029-9}
}

@article{KoLy05,
  author  = {Kondratiev, Yuri G. and Lytvynov, Eugene},
  title   = {{{Glauber} dynamics of continuous particle systems}},
  journal = {Annales de l'Institut Henri Poincar{\'e} Probabilit{\'e}s et Statistiques},
  volume  = {41},
  number  = {4},
  pages   = {685--702},
  year    = {2005},
  issn    = {0246-0203},
  doi     = {10.1016/j.anihpb.2004.05.002},
  url     = {https://doi.org/10.1016/j.anihpb.2004.05.002}
}

@article{Ruelle,
  author  = {Ruelle, David},
  title   = {{Superstable interactions in classical statistical mechanics}},
  journal = {Communications in Mathematical Physics},
  volume  = {18},
  number  = {2},
  pages   = {127--159},
  year    = {1970},
  doi     = {10.1007/BF01646087},
  url     = {https://doi.org/10.1007/BF01646087}
}

@article{Erbar_2012,
  author  = {Erbar, Matthias and Maas, Jan},
  title   = {{{Ricci} Curvature of Finite {Markov} Chains via Convexity of the Entropy}},
  journal = {Archive for Rational Mechanics and Analysis},
  volume  = {206},
  number  = {3},
  pages   = {997--1038},
  year    = {2012},
  doi     = {10.1007/s00205-012-0554-z},
  url     = {https://doi.org/10.1007/s00205-012-0554-z}
}

@article{maas2011gradientflowsentropyfinite,
  author  = {Maas, Jan},
  title   = {{Gradient flows of the entropy for finite {Markov} chains}},
  journal = {Journal of Functional Analysis},
  volume  = {261},
  number  = {8},
  pages   = {2250--2292},
  year    = {2011},
  issn    = {0022-1236},
  doi     = {10.1016/j.jfa.2011.06.009},
  url     = {https://doi.org/10.1016/j.jfa.2011.06.009}
}

@article{Daneri_2008,
  author  = {Daneri, Sara and Savar{\'e}, Giuseppe},
  title   = {{Eulerian Calculus for the Displacement Convexity in the {Wasserstein} Distance}},
  journal = {SIAM Journal on Mathematical Analysis},
  volume  = {40},
  number  = {3},
  pages   = {1104--1122},
  year    = {2008},
  doi     = {10.1137/08071346X},
  url     = {https://doi.org/10.1137/08071346X}
}

@book{buttazzo1989semicontinuity,
  author    = {Buttazzo, Giuseppe},
  title     = {{Semicontinuity, Relaxation and Integral Representation in the Calculus of Variations}},
  series    = {Pitman Research Notes in Mathematics Series},
  volume    = {207},
  publisher = {Longman Scientific \& Technical},
  address   = {Harlow},
  year      = {1989}
}

@article{Houdebert_Zass_2022,
  author  = {Houdebert, Pierre and Zass, Alexander},
  title   = {{An explicit {Dobrushin} uniqueness region for {Gibbs} point processes with repulsive interactions}},
  journal = {Journal of Applied Probability},
  volume  = {59},
  number  = {2},
  pages   = {541--555},
  year    = {2022},
  doi     = {10.1017/jpr.2021.70},
  url     = {https://doi.org/10.1017/jpr.2021.70}
}

@article{li2010noteequilibriumglauberkawasaki,
  author  = {Li, Guanhua and Lytvynov, Eugene},
  title   = {{A note on equilibrium {Glauber} and {Kawasaki} dynamics for permanental point processes}},
  journal = {Methods of Functional Analysis and Topology},
  volume  = {17},
  number  = {1},
  pages   = {29--46},
  year    = {2011},
  issn    = {1029-3531},
  url     = {https://mfat.imath.kiev.ua/article/?id=563}
}

@article{steenbeck2025reversiblebirthanddeathdynamicscontinuum,
  author        = {Jahnel, Benedikt and K{\"o}ppl, Jonas and Steenbeck, Yannic and Zass, Alexander},
  title         = {{Reversible birth-and-death dynamics in continuum: free-energy dissipation and attractor properties}},
  journal       = {{arXiv} preprint {arXiv}:2508.21196},
  year          = {2025},
  eprint        = {2508.21196},
  archivePrefix = {arXiv},
  primaryClass  = {math.PR},
  doi           = {10.48550/arXiv.2508.21196},
  url           = {https://arxiv.org/abs/2508.21196}
}

@article{garcia2006spatialbirthdeathprocesses,
  author  = {Garcia, Nancy L. and Kurtz, Thomas G.},
  title   = {{Spatial birth and death processes as solutions of stochastic equations}},
  journal = {Electronic Journal of Probability},
  volume  = {12},
  pages   = {1169--1221},
  year    = {2007},
  doi     = {10.1214/EJP.v12-430},
  url     = {https://doi.org/10.1214/EJP.v12-430}
}

@article{Papangelou,
  author  = {Nehring, Benjamin and Zessin, Hans},
  title   = {{The {Papangelou} Process: A Concept for {Gibbs}, {Fermi} and {Bose} Processes}},
  journal = {Journal of Contemporary Mathematical Analysis},
  volume  = {46},
  number  = {6},
  pages   = {327--336},
  year    = {2011},
  doi     = {10.3103/S1068362311060069},
  url     = {https://doi.org/10.3103/S1068362311060069}
}

@book{EthierKurtz,
  author    = {Ethier, Stewart N. and Kurtz, Thomas G.},
  title     = {{{Markov} Processes: Characterization and Convergence}},
  series    = {Wiley Series in Probability and Mathematical Statistics},
  publisher = {Wiley},
  address   = {New York},
  year      = {1986},
  isbn      = {9780471081868}
}

@article{Lytvynov-Ohlerich07,
  author  = {Lytvynov, Eugene and Ohlerich, Nataliya},
  title   = {{A note on equilibrium {Glauber} and {Kawasaki} dynamics for fermion point processes}},
  journal = {Methods of Functional Analysis and Topology},
  volume  = {13},
  number  = {4},
  pages   = {362--372},
  year    = {2007}
}

@article{Moller,
  author  = {M{\o}ller, Jesper},
  title   = {{On the rate of convergence of spatial birth-and-death processes}},
  journal = {Annals of the Institute of Statistical Mathematics},
  volume  = {41},
  number  = {3},
  pages   = {565--581},
  year    = {1989},
  issn    = {0020-3157},
  doi     = {10.1007/BF00050669},
  url     = {https://doi.org/10.1007/BF00050669}
}

@article{EHMT17,
  author  = {Erbar, Matthias and Henderson, Christopher and Menz, Georg and Tetali, Prasad},
  title   = {{{Ricci} curvature bounds for weakly interacting {Markov} chains}},
  journal = {Electronic Journal of Probability},
  volume  = {22},
  pages   = {40},
  year    = {2017},
  issn    = {1083-6489},
  doi     = {10.1214/17-EJP49},
  url     = {https://doi.org/10.1214/17-EJP49}
}

@article{FM16,
  author  = {Fathi, Max and Maas, Jan},
  title   = {{Entropic {Ricci} curvature bounds for discrete interacting systems}},
  journal = {The Annals of Applied Probability},
  volume  = {26},
  number  = {3},
  pages   = {1774--1806},
  year    = {2016},
  issn    = {1050-5164},
  doi     = {10.1214/15-AAP1133},
  url     = {https://doi.org/10.1214/15-AAP1133}
}

@book{Willard,
  author    = {Willard, Stephen},
  title     = {{General Topology}},
  series    = {Addison-Wesley Series in Mathematics},
  publisher = {Addison-Wesley},
  address   = {Reading, MA},
  year      = {1970},
  isbn      = {0201087073}
}

@article{Geogii76,
  author  = {Georgii, Hans-Otto},
  title   = {{Canonical and grand canonical {Gibbs} states for continuum systems}},
  journal = {Communications in Mathematical Physics},
  volume  = {48},
  number  = {1},
  pages   = {31--51},
  year    = {1976},
  doi     = {10.1007/BF01609410},
  url     = {https://doi.org/10.1007/BF01609410}
}

@article{NguyenZessin77,
  author  = {Xanh, Nguyen X. and Zessin, Hans},
  title   = {{Integral and differential characterizations of the {Gibbs} process}},
  journal = {Advances in Applied Probability},
  volume  = {9},
  number  = {3},
  pages   = {446--447},
  year    = {1977},
  doi     = {10.2307/1426106},
  url     = {https://doi.org/10.2307/1426106}
}

@article{Betsch23,
  author  = {Betsch, Steffen},
  title   = {{Structural Properties of {Gibbsian} Point Processes in Abstract Spaces}},
  journal = {Journal of Theoretical Probability},
  volume  = {36},
  number  = {4},
  pages   = {2501--2563},
  year    = {2023},
  doi     = {10.1007/s10959-023-01262-9},
  url     = {https://doi.org/10.1007/s10959-023-01262-9}
}

@article{Mecke1967,
  author  = {Mecke, Joseph},
  title   = {{Station{\"a}re zuf{\"a}llige Ma{\ss}e auf lokalkompakten abelschen Gruppen}},
  journal = {Zeitschrift f{\"u}r Wahrscheinlichkeitstheorie und verwandte Gebiete},
  volume  = {9},
  number  = {1},
  pages   = {36--58},
  year    = {1967},
  doi     = {10.1007/BF00531950},
  url     = {https://doi.org/10.1007/BF00531950}
}

@article{papangelou1974conditional,
  author  = {Papangelou, Fosco},
  title   = {{The conditional intensity of general point processes and an application to line processes}},
  journal = {Zeitschrift f{\"u}r Wahrscheinlichkeitstheorie und verwandte Gebiete},
  volume  = {28},
  number  = {3},
  pages   = {207--226},
  year    = {1974},
  doi     = {10.1007/BF00533244},
  url     = {https://doi.org/10.1007/BF00533244}
}

@article{Dobruschin1968TheDO,
  author  = {Dobrushin, Roland L.},
  title   = {{The Description of a Random Field by Means of Conditional Probabilities and Conditions of Its Regularity}},
  journal = {Theory of Probability and Its Applications},
  volume  = {13},
  number  = {2},
  pages   = {197--224},
  year    = {1968},
  doi     = {10.1137/1113026},
  url     = {https://doi.org/10.1137/1113026}
}

@article{NguyenZessin1979,
  author  = {Xanh, Nguyen X. and Zessin, Hans},
  title   = {{Ergodic theorems for spatial processes}},
  journal = {Zeitschrift f{\"u}r Wahrscheinlichkeitstheorie und verwandte Gebiete},
  volume  = {48},
  number  = {2},
  pages   = {133--158},
  year    = {1979},
  doi     = {10.1007/BF01886869},
  url     = {https://doi.org/10.1007/BF01886869}
}

@article{jahnel2026reversiblebirthanddeathdynamicscontinuum,
  author  = {Jahnel, Benedikt and K{\"o}ppl, Jonas and Steenbeck, Yannic and Zass, Alexander},
  title   = {{Reversible birth-and-death dynamics in continuum: a de {Bruijn}-type identity for free-energy dissipation}},
  journal = {{arXiv} preprint {arXiv}:2602.13474},
  year    = {2026},
  doi     = {10.48550/arXiv.2602.13474},
  url     = {https://arxiv.org/abs/2602.13474}
}

@article{Macchi_1975,
  author  = {Macchi, Odile},
  title   = {{The coincidence approach to stochastic point processes}},
  journal = {Advances in Applied Probability},
  volume  = {7},
  number  = {1},
  pages   = {83--122},
  year    = {1975},
  doi     = {10.2307/1425855},
  url     = {https://doi.org/10.2307/1425855}
}

@article{Soshnikov2000Determinantal,
  author  = {Soshnikov, Alexander},
  title   = {{Determinantal Random Point Fields}},
  journal = {Russian Mathematical Surveys},
  volume  = {55},
  number  = {5},
  pages   = {923--975},
  year    = {2000},
  doi     = {10.1070/RM2000v055n05ABEH000321},
  url     = {https://doi.org/10.1070/RM2000v055n05ABEH000321}
}

@article{GeorgiiYoo2005,
  author  = {Georgii, Hans-Otto and Yoo, Hyun J.},
  title   = {{Conditional Intensity and {Gibbsianness} of Determinantal Point Processes}},
  journal = {Journal of Statistical Physics},
  volume  = {118},
  number  = {1--2},
  pages   = {55--84},
  year    = {2005},
  doi     = {10.1007/s10955-004-8811-9},
  url     = {https://doi.org/10.1007/s10955-004-8811-9},
  eprint  = {math/0401402}
}

@article{ErbarMaas2014GradientFlow,
  author  = {Erbar, Matthias and Maas, Jan},
  title   = {{Gradient Flow Structures for Discrete Porous Medium Equations}},
  journal = {Discrete and Continuous Dynamical Systems - Series A},
  volume  = {34},
  number  = {4},
  pages   = {1355--1374},
  year    = {2014},
  doi     = {10.3934/dcds.2014.34.1355}
}

@incollection{Dereudre2019,
  author    = {Dereudre, David},
  title     = {{Introduction to the Theory of {Gibbs} Point Processes}},
  booktitle = {{Stochastic Geometry: Modern Research Frontiers}},
  publisher = {Springer},
  address   = {Cham},
  pages     = {181--229},
  year      = {2019},
  doi       = {10.1007/978-3-030-13547-8_5},
  url       = {https://doi.org/10.1007/978-3-030-13547-8_5}
}

@article{leble2015uniqueness,
  author  = {Lebl{\'e}, Thomas},
  title   = {{A Uniqueness Result for Minimizers of the {1D} Log-gas Renormalized Energy}},
  journal = {Journal of Functional Analysis},
  volume  = {268},
  number  = {7},
  pages   = {1649--1677},
  year    = {2015},
  doi     = {10.1016/j.jfa.2014.12.014}
}

@article{Mielke2011,
  author  = {Mielke, Alexander},
  title   = {{A Gradient Structure for Reaction--Diffusion Systems and for Energy-Drift-Diffusion Systems}},
  journal = {Nonlinearity},
  volume  = {24},
  number  = {4},
  pages   = {1329--1346},
  year    = {2011}
}

@article{BeBr00,
  author  = {Benamou, Jean-David and Brenier, Yann},
  title   = {{A computational fluid mechanics solution to the {Monge--Kantorovich} mass transfer problem}},
  journal = {Numerische Mathematik},
  volume  = {84},
  number  = {3},
  pages   = {375--393},
  year    = {2000},
  issn    = {0029-599X},
  doi     = {10.1007/s002110050002},
  url     = {https://doi.org/10.1007/s002110050002}
}

@article{St06,
  author  = {Sturm, Karl-Theodor},
  title   = {{On the geometry of metric measure spaces. {II}}},
  journal = {Acta Mathematica},
  volume  = {196},
  number  = {1},
  pages   = {133--177},
  year    = {2006},
  doi     = {10.1007/s11511-006-0003-7},
  url     = {https://doi.org/10.1007/s11511-006-0003-7}
}

@article{LoVi09,
  author  = {Lott, John and Villani, C{\'e}dric},
  title   = {{{Ricci} curvature for metric-measure spaces via optimal transport}},
  journal = {Annals of Mathematics},
  volume  = {169},
  number  = {3},
  pages   = {903--991},
  year    = {2009},
  doi     = {10.4007/annals.2009.169.903},
  url     = {https://doi.org/10.4007/annals.2009.169.903}
}

@article{OtVi00,
  author  = {Otto, Felix and Villani, C{\'e}dric},
  title   = {{Generalization of an inequality by {Talagrand} and links with the logarithmic {Sobolev} inequality}},
  journal = {Journal of Functional Analysis},
  volume  = {173},
  number  = {2},
  pages   = {361--400},
  year    = {2000},
  issn    = {0022-1236},
  doi     = {10.1006/jfan.1999.3557}
}

@article{Wu00,
  author  = {Wu, Li-Ming},
  title   = {{A new modified logarithmic {Sobolev} inequality for {Poisson} point processes and several applications}},
  journal = {Probability Theory and Related Fields},
  volume  = {118},
  number  = {3},
  pages   = {427--438},
  year    = {2000},
  doi     = {10.1007/PL00008749},
  url     = {https://doi.org/10.1007/PL00008749}
}

@article{GrossRothaus1998,
  author  = {Gross, Leonard and Rothaus, Oscar},
  title   = {{Herbst Inequalities for Supercontractive Semigroups}},
  journal = {Journal of Mathematics of Kyoto University},
  volume   = {38},
  number   = {2},
  pages    = {295--318},
  year     = {1998},
  doi      = {10.1215/kjm/1250518120}
}

@article{steenbeck2026modifiedlogsobolevinequalitiesconcentration,
      title={{Modified {log-Sobolev} inequalities, concentration bounds and uniqueness of {Gibbs} measures}},
      author={Steenbeck, Yannic},
      journal = {{arXiv} preprint {arXiv}:2603.25479},
      year={2026},
      archivePrefix={arXiv},
      primaryClass={math.PR},
      url={https://arxiv.org/abs/2603.25479}
}

@article{ErbarFathiSchlichting2020,
  author  = {Erbar, Matthias and Fathi, Max and Schlichting, Andr{\'e}},
  title   = {{Entropic curvature and convergence to equilibrium for mean-field dynamics on discrete spaces}},
  journal = {ALEA, Latin American Journal of Probability and Mathematical Statistics},
  volume  = {17},
  number  = {1},
  pages   = {445--471},
  year    = {2020}
}

@article{ErHu15,
  author  = {Erbar, Matthias and Huesmann, Martin},
  title   = {{Curvature bounds for configuration spaces}},
  journal = {Calculus of Variations and Partial Differential Equations},
  volume  = {54},
  number  = {1},
  pages   = {397--430},
  year    = {2015},
  issn    = {0944-2669},
  doi     = {10.1007/s00526-014-0790-1},
  url     = {https://doi.org/10.1007/s00526-014-0790-1}
}

@article{DeSu22,
  author  = {Dello Schiavo, Lorenzo and Suzuki, Kohei},
  title   = {{Configuration {Spaces} over {Singular} {Spaces} -- {II}. {Curvature}}},
  journal = {{arXiv} preprint {arXiv}:2205.01379},
  year    = {2022},
  doi     = {10.48550/arXiv.2205.01379},
  url     = {https://arxiv.org/abs/2205.01379}
}

@article{DeSu21,
  author  = {Dello Schiavo, Lorenzo and Suzuki, Kohei},
  title   = {{Configuration {Spaces} over {Singular} {Spaces} -- {I}. {Dirichlet} form and metric measure geometry}},
  journal = {{arXiv} preprint {arXiv}:2109.03192},
  year    = {2021},
  doi     = {10.48550/arXiv.2109.03192},
  url     = {https://arxiv.org/abs/2109.03192}
}

@article{HuSt13,
  author  = {Huesmann, Martin and Sturm, Karl-Theodor},
  title   = {{Optimal transport from {Lebesgue} to {Poisson}}},
  journal = {The Annals of Probability},
  volume  = {41},
  number  = {4},
  pages   = {2426--2478},
  year    = {2013},
  issn    = {0091-1798},
  doi     = {10.1214/12-AOP814},
  url     = {https://doi.org/10.1214/12-AOP814}
}

@article{Hu16,
  author  = {Huesmann, Martin},
  title   = {{Optimal transport between random measures}},
  journal = {Annales de l'Institut Henri Poincar{\'e} Probabilit{\'e}s et Statistiques},
  volume  = {52},
  number  = {1},
  pages   = {196--232},
  year    = {2016},
  issn    = {0246-0203},
  doi     = {10.1214/14-AIHP634},
  url     = {https://doi.org/10.1214/14-AIHP634}
}

@article{erbar2018onedimensionalloggasfreeenergy,
  author  = {Erbar, Matthias and Huesmann, Martin and Lebl{\'e}, Thomas},
  title   = {{The One-Dimensional Log-Gas Free Energy Has a Unique Minimizer}},
  journal = {Communications on Pure and Applied Mathematics},
  volume  = {74},
  number  = {3},
  pages   = {615--675},
  year    = {2021},
  doi     = {10.1002/cpa.21977},
  url     = {https://doi.org/10.1002/cpa.21977}
}

@article{EHJM25,
  author  = {Erbar, Matthias and Huesmann, Martin and Jalowy, Jonas and M{\"u}ller, Bastian},
  title   = {{Optimal transport of stationary point processes: Metric structure, gradient flow and convexity of the specific entropy}},
  journal = {Journal of Functional Analysis},
  volume  = {289},
  number  = {4},
  pages   = {110974},
  year    = {2025},
  doi     = {10.1016/j.jfa.2025.110974},
  url     = {https://doi.org/10.1016/j.jfa.2025.110974}
}

@article{HuMu25,
  author        = {Huesmann, Martin and M{\"u}ller, Bastian},
  title         = {{Gap metrics for stationary point processes and quantitative convexity of the free energy}},
  journal       = {{arXiv} preprint {arXiv}:2509.08659},
  year          = {2025},
  archivePrefix = {arXiv},
  eprint        = {2509.08659},
  primaryClass  = {math.PR},
  doi           = {10.48550/arXiv.2509.08659}
}

@article{dereudre2024nonhyperuniformityperturbedlattices,
  title={{{(Non)}-hyperuniformity of perturbed lattices}},
  author={Dereudre, David and Flimmel, Daniela and Huesmann, Martin and Lebl{\'e}, Thomas},
  journal={{arXiv} preprint {arXiv}:2405.19881},
  year={2024}
}

@article{DerDig,
  title={{Rigidity of one-dimensional point processes via optimal transport}},
  author={Dereudre, David and Digneaux, Rafa{\"e}l},
  journal={{arXiv} preprint {arXiv}:2510.17257},
  year={2025}
}

@article{HuLe24,
  author  = {Huesmann, Martin and Lebl{\'e}, Thomas},
  title   = {{The link between hyperuniformity, {Coulomb} energy, and {Wasserstein} distance to {Lebesgue} for two-dimensional point processes}},
  journal = {Probability and Mathematical Physics},
  volume  = {7},
  pages   = {123--173},
  year    = {2026},
  doi     = {10.2140/pmp.2026.7.123},
  url     = {https://doi.org/10.2140/pmp.2026.7.123}
}

@article{JuengelYue2017,
  author  = {J{\"u}ngel, Ansgar and Yue, Wen},
  title   = {Discrete {Beckner} inequalities via the
             {Bochner--Bakry--Emery} approach for {Markov} chains},
  journal = {The Annals of Applied Probability},
  year    = {2017},
  volume  = {27},
  number  = {4},
  pages   = {2238--2269},
  doi     = {10.1214/16-AAP1258}
}

@article{MuSa20,
  author  = {Muratori, Matteo and Savar{\'e}, Giuseppe},
  title   = {Gradient flows and {Evolution Variational Inequalities}
             in metric spaces. {I}: Structural properties},
  journal = {Journal of Functional Analysis},
  year    = {2020},
  volume  = {278},
  number  = {4},
  pages   = {108347},
  doi     = {10.1016/j.jfa.2019.108347}
}
\bibliographystyle{alpha}
\end{document}